%% file: main.tex
\documentclass[a4paper,twoside,10pt]{book}

\input{header}

	\title{Dirac Eigenvalues of higher Multiplicity}
	\author{Nikolai Nowaczyk}
	
	\makeindex
	\usepackage[totoc]{idxlayout}
	
	\makenomenclature

	\nomenclature[N]{$\N$}{the natural numbers $\N=\{0,1,2, \ldots \}$}
	\nomenclature[TN]{$\mathcal{T}(M)$}{smooth vector fields on $M$}

\begin{document}
	
		\input{titlepages}

		\newpage
		\tableofcontents
		\thispagestyle{empty}

		  \chapter{Overview of the Thesis}
			\input{intro.overview}

		 \chapter{A Short Review of Spin Geometry} \label{ChapterReviewSpinGeometry}
			\input{review.abstract}

			\input{review.clifford}
			\input{review.spinstructs}
			\input{review.spinorbundles}
			\input{review.realcomplex}
			\input{review.dahl}

			\input{review.spinmorphisms}

			\input{review.category}

		 \chapter{Continuity of Dirac Spectra} \label{ChapterContDiracSpec}
			 \input{evpaper.abstract}
			 \input{evpaper.intro}

			 \input{evpaper.famdiscrete}
			 \input{evpaper.prfmain}
			 \input{evpaper.specfl}
		 \chapter{The Universal Spinor Field Bundle} \label{ChapterUSFB}
			 \input{usb.abstract}
			 \input{usb.intro}

			 \input{usb.usb}

			 \input{usb.univdirac}
			 \input{usb.specdec}

		 \addtocontents{toc}{\protect\newpage}
		 \chapter{Higher Multiplicities} \label{ChapterHigher}
			 \input{higher.abstract}
			 \input{higher.intro} 
			 \input{higher.orientors}

			 \input{higher.spindiffeos}

			 \input{higher.sphere}
			 \input{higher.surgery}

			 \input{higher.prfmainthm}

		\begin{appendices}
		\chapter{Supplementary Material}
			\input{review.appendix}
			\input{usb.appendix}

			\input{higher.appendix}

			\input{evpaper.appendix}
		\end{appendices}

		\newpage
		\setlength{\nomitemsep}{-0.8\parsep}   
		\printnomenclature[8.5em]
		\newpage

		\printindex
		
		\phantomsection
		\addcontentsline{toc}{chapter}{\listfigurename}
		\listoffigures
	
		\printbibliography[title=References]
		\addcontentsline{toc}{chapter}{References}

\end{document}

%% file: header.tex
	\usepackage{amsmath}				
	\usepackage{amscd}					
	\usepackage{amssymb}				
	\usepackage[titletoc]{appendix}		
	\usepackage[
		style=alphabetic,				
		sorting=nyt,					
		maxnames=5,						
		minnames=3,						
		backend=bibtex8,				
		firstinits=true,				
		doi=false,						
	]{biblatex}
	\usepackage{bm}						
	\usepackage{color}					
	\usepackage{enumitem}				
	\usepackage{etoolbox}				
	\usepackage{fancyhdr}				
	\usepackage{float}					
	\usepackage[T1]{fontenc}			
	\usepackage{gensymb}				
	\usepackage[a4paper,
		outer=3cm,
		inner=4cm,
		top=3.5cm,
		bottom=3cm]{geometry}
	\usepackage[utf8]{inputenc}			
	\usepackage{ifmtarg}				
	\usepackage{mathrsfs}				
	\usepackage{makeidx}				
	\usepackage{multicol}				
	\usepackage{parskip}				
	\usepackage[refpage,intoc]{nomencl}		
	\usepackage[thmmarks,
		hyperref,
		amsmath]{ntheorem}
	\usepackage{slashed}				
	\usepackage{stmaryrd}				
	\usepackage{tensor}					
	\usepackage{textcomp}				
	\usepackage{tikz}
	\usepackage{tocloft}				
	\usepackage{url}					
	\usepackage[all,cmtip]{xy}			
	
	\usepackage{csquotes}				
	
	\usetikzlibrary{decorations.markings}
	
	\usepackage{hyperref}		
	\usepackage[capitalize]{cleveref}	

	\definecolor{linkcolor}{RGB}{00,10,138} 
	\hypersetup{
		colorlinks=true,
		breaklinks=true,
		linkcolor=linkcolor,
		menucolor=linkcolor,
		citecolor=linkcolor,
		filecolor=linkcolor,
		urlcolor=linkcolor,
		frenchlinks=false,
		pageanchor=true
	}

	\bibliography{refs} 
	\DeclareFieldFormat{postnote}{#1}
	\DeclareFieldFormat{multipostnote}{#1}
	
		\numberwithin{equation}{section}
		\crefformat{equation}{(#2#1#3)} 
		
	\fancypagestyle{plain}{%
		\fancyhf{}%
		\fancyfoot[C]{}%
	}

		\renewcommand{\nomname}{List of Symbols}
		\patchcmd{\thenomenclature}
		  {\chapter*{\nomname}}
		  {\phantomsection\chapter*{\nomname}\label{SectNomenclature}}
		  {}
		  {}

		\newcommand{\textname}{chapter }

	\makeatletter
	\def\theorem@checkbold{}
	\makeatother

	\makeatletter
	\newtheoremstyle{myprfstlye}%
		{\item[\theorem@headerfont\hskip\labelsep ##1\theorem@separator]}%
		{\item[\theorem@headerfont\hskip\labelsep ##1 of ##3\theorem@separator]}
	\newtheoremstyle{mythmstyle}%
		{\item[\theorem@headerfont\hskip\labelsep ##1\ ##2\theorem@separator]}%
		{\item[\theorem@headerfont\hskip\labelsep ##1\ ##2\ \normalfont{(##3)}\textbf{\theorem@separator}]}	
	\newtheoremstyle{myabstractstyle}%
		{\item[\theorem@headerfont\hskip\labelsep ##1.]}%
		{\item[\theorem@headerfont\hskip\labelsep ##1.]}			
		
	\makeatother
	
	\theoremstyle{mythmstyle}
	\theoremheaderfont{\normalfont\bfseries}
	\theorembodyfont{\normalfont}
	\theoremseparator{.}
	\theoremsymbol{$\lozenge$}
	\newtheorem{Def}{Definition}[section]

	\newtheorem{Exm}[Def]{Example}
	
	\newtheorem{Cor}[Def]{Corollary}
	
	\newtheorem{Lem}[Def]{Lemma}
	
	\newtheorem{Rem}[Def]{Remark}

	\newtheorem{Cnj}[Def]{Conjecture}
	
	\newtheorem{MainThm}{Main Theorem}
	\newtheorem{Prb}{Problem}
	
	\newtheorem{Thm}[Def]{Theorem}
	\newtheorem{Not}[Def]{Notation}
	
	\theoremstyle{myprfstlye}
	\theoremsymbol{$\square$}
	\newtheorem{Prf}{Proof}
	
	\theoremstyle{empty}
	\theoremsymbol{$\lozenge$}
	\newtheorem{recall}{}
	
	\newenvironment{subabstract}{\quote \textbf{Abstract}.}{\endquote}
	
	\crefname{Def}{Definition}{Definitions}
	\crefname{Exm}{Example}{Examples}
	\crefname{Con}{Convention}{Conventions}
	\crefname{Cor}{Corollary}{Corollaries}
	\crefname{Fct}{Fact}{Facts}
	\crefname{Lem}{Lemma}{Lemmas}
	\crefname{NumText}{Paragraph}{Paragraphs}
	\crefname{Rem}{Remark}{Remarks}
	\crefname{Cnj}{Conjecture}{Conjectures}
	\crefname{Exc}{Exercise}{Exercises}
	\crefname{MainThm}{Main Theorem}{Main Theorems}
	\crefname{Prb}{Problem}{Problems}
	\crefname{Qst}{Question}{Questions}
	\crefname{Thm}{Theorem}{Theorems}
	\crefname{Not}{Notation}{Notation}
	\crefname{Prf}{Proof}{Proofs}
			
		\newcommand{\fun}[1]{\bm{\mathrm{#1}}}
		
		\sbox0{$\mathbf{\xdef\mybffam{\the\fam}}$}
		\newcommand{\cat}[1]{\begingroup\fam\mybffam#1\endgroup}

		\newlist{category}{enumerate}{10}
		\setlist[category]{
			itemsep=1ex,
			topsep=0.5ex,
			parsep=0pt,
			leftmargin=3em,
			itemindent=5em,
			labelwidth=5em,
			labelsep=1em,
			align=categorylabel,
			label*=.\arabic*
		}
		\setlist[category,1]{label=\arabic*}
		\SetLabelAlign{categorylabel}{\bfseries{\scshape{#1}}:\hfil}	
	
		\newenvironment{DefI}[1][]{\ifstrempty{#1}{\Def}{\Def[#1]\index{#1}}}{\endDef}
	
		\newcommand{\DefMap}[4]{
			\begin{align*}
				\begin{array}{rcl}
					#1 & \to & #2 \\
					#3 & \mapsto & #4 
				\end{array} 
			\end{align*}
		}

		\newcommand{\emphi}[1]{\emph{#1}\index{#1}}

		\newcommand{\jeq}[2]{\stackrel{\text{#1}}{#2}}
	
	\newlist{caselist}{enumerate}{10}
	\setlist[caselist]{itemsep=1ex,topsep=1ex,parsep=0pt,leftmargin=0pt,%
	  labelwidth=*,labelsep=1ex,align=caselabel,label*=.\arabic*}
	\setlist[caselist,1]{label=\arabic*}

	\newcommand\thecaselabeltext{}
	\SetLabelAlign{caselabel}{\textsc{Case~#1}\thecaselabeltext:\hfil}

	\crefname{caselisti}{case}{cases}
	\loop\ifnum\count255<10
	   \advance\count255 1
	   \crefalias{caselist\romannumeral\count255}{caselisti}
	\repeat

	\makeatletter
	\def\ifempty#1{\def\@x{#1}\ifx\@x\@empty}
	\makeatother

	\newlist{steplist}{enumerate}{10}
	\setlist[steplist]{nolistsep,itemsep=0ex,topsep=0ex,partopsep=0px,parsep=0pt,leftmargin=0pt,%
	  labelwidth=*,labelsep=1ex,align=steplabel,label*=.\arabic*}
	\setlist[steplist,1]{label=\arabic*}

	\newcommand\thesteplabeltext{}
	\SetLabelAlign{steplabel}{\textsc{Step~#1}\thesteplabeltext:\hfil}

	\crefname{steplisti}{Step}{steplist}
	\loop\ifnum\count255<10
	   \advance\count255 1
	   \crefalias{steplist\romannumeral\count255}{steplisti}
	\repeat

	\makeatletter
	\def\ifempty#1{\def\@x{#1}\ifx\@x\@empty}
	\makeatother

	\newcommand\step[1][]{\renewcommand{\thesteplabeltext}
	  {\ifempty{#1}\else~(#1)\fi}\item}
	\DeclareMathOperator{\ash}{arsinh}				
	\DeclareMathOperator{\Aut}{Aut}					
	\DeclareMathOperator{\B}{B}						
	\DeclareMathOperator{\Cl}{C\ell}				
	\DeclareMathOperator{\Conf}{Conf}				
	\DeclareMathOperator{\C}{\mathbb{C}}			
	\DeclareMathOperator{\Diff}{Diff} 				
	\DeclareMathOperator{\Dirac}{\slashed{D}}		
	\DeclareMathOperator{\End}{End} 				
	\newcommand{\GLp}{\operatorname{GL}^{+}}		
	\newcommand{\GLtp}{\widetilde{\operatorname{GL}}^{+}}	
	\DeclareMathOperator{\GL}{GL} 					
	\DeclareMathOperator{\G}{\mathcal{G}}			
	\DeclareMathOperator{\Hom}{Hom} 				
	\DeclareMathOperator{\Isom}{Isom}				
	\DeclareMathOperator{\Iso}{Iso} 				
	\DeclareMathOperator{\K}{\mathbb{K}}			
	\DeclareMathOperator{\Lin}{span} 				
	\DeclareMathOperator{\Mon}{Mon}					
	\DeclareMathOperator{\N}{\mathbb{N}}			
	\DeclareMathOperator{\Or}{\operatorname{O}}		
	\DeclareMathOperator{\Rm}{\mathcal{R}}			
	\DeclareMathOperator{\R}{\mathbb{R}}			
	\DeclareMathOperator{\SL}{SL}					
	\DeclareMathOperator{\SO}{SO} 					
	\DeclareMathOperator{\Spin}{Spin} 				
	\DeclareMathOperator{\T}{\mathbb{T}}			
	\DeclareMathOperator{\Z}{\mathbb{Z}}			
	\DeclareMathOperator{\arsinh}{arsinh}			
	\DeclareMathOperator{\codim}{codim} 			
	\DeclareMathOperator{\dist}{dist} 				
	\DeclareMathOperator{\dom}{\mathcal{D}} 		
	\DeclareMathOperator{\dvol}{dV} 				
	\DeclareMathOperator{\ev}{ev}					
	\DeclareMathOperator{\grad}{grad} 				
	\DeclareMathOperator{\id}{id} 					
	\DeclareMathOperator{\image}{im} 				
	\DeclareMathOperator{\ind}{ind} 				
	\DeclareMathOperator{\pr}{pr} 					
	\DeclareMathOperator{\sa}{sa}					
	\DeclareMathOperator{\scal}{scal}				
	\DeclareMathOperator{\sgn}{sgn} 				
	\DeclareMathOperator{\spn}{span} 				
	\DeclareMathOperator{\spc}{\mathfrak{s}}		
	\DeclareMathOperator{\specfl}{sf}				
	\DeclareMathOperator{\spec}{spec} 				
	\DeclareMathOperator{\spin}{spin} 				
	\DeclareMathOperator{\spp}{\mathfrak{sp}}		
	
	\newcommand{\subseto}{\hspace{.1em}\mathring{\subset}\hspace{.1em}} 
	\renewcommand{\H}{\mathbb{H}}			
	\renewcommand{\O}{O} 					
	\renewcommand{\S}{S}					
	\newcommand{\D}{D}						
	\newcommand{\A}{A}						

%% file: titlepages.tex
	\begin{titlepage}
		$ $
		\vspace{8em} 
		\begin{center}
			{\huge \textbf{Dirac Eigenvalues of higher Multiplicity}}
			\begingroup
			\fontsize{12pt}{14pt}\selectfont			
				\vspace{5em}
				
				\includegraphics[scale=0.6]{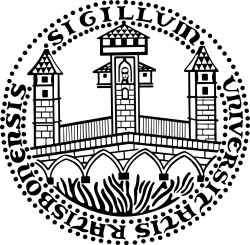} \\
				
				\vspace{6em}
				\textsc{
					Dissertation zur Erlangung des Doktorgrades \\
					der Naturwissenschaften (Dr. rer. nat.) \\
					der Fakultät für Mathematik\\
					der Universität Regensburg}
					
					\vspace{6em}
					
					\textsc{vorgelegt von}\\ 
					
					\vspace{0.5em}
					
					{\Large \textsc{ \textbf{Nikolai Nowaczyk}}}\\
					
					\vspace{0.5em}
					
					\textsc{aus Duisburg im Jahr 2014 }\\ 

			\endgroup			
		\end{center}
	\end{titlepage}
			
	\newpage
	
	\thispagestyle{empty}
	$ $
	\vspace{40em}
	{\large $ $ \\
	Promotionsgesuch eingereicht am:  25.08.2014. \\ 
	$ $\\ $ $\\ 
	Die Arbeit wurde angeleitet von: Prof. Dr. Bernd Ammann. \\ $ $ \\

	Prüfungsausschuss: \\ \\
	\begin{tabular}[H]{lll}
		& Vorsitzender: &  Prof. Dr. Helmut Abels\\
		& Erst-Gutachter: &  Prof. Dr. Bernd Ammann\\
		& Zweit-Gutachter: & Prof. Dr. Mattias Dahl, KTH Stockholm \\
		& weiterer Prüfer: &  Prof. Dr. Stefan Friedl\\
		& Ersatzprüfer: & Prof. Dr. Felix Finster\\
	\end{tabular}
	}

	\newpage
	
	\thispagestyle{empty}
	{\Large \textbf{Abstract}} \\
	Let $(M,\Theta)$ be a closed spin manifold of dimension $m \geq 3$ with fixed topological spin structure $\Theta$. For any Riemannian metric $g$, we can construct the associated Dirac operator $\Dirac^g$. The spectrum of this Dirac operator depends on $g$ of course. In 2005, Dahl conjectured that $M$ can be given a metric, for which a finite part of the spectrum consists of arbitrarily prescribed eigenvalues of arbitrary (finite) multiplicity. The only constraints one has to respect are the exception of the zero eigenvalue (due to the Atiyah-Singer index theorem) and in certain dimensions, the quaternionic structure of the eigenspaces and also the symmetry of the spectrum. Dahl also proved his conjecture in case all eigenvalues have simple multiplicities. The question, if one can prescribe arbitrary multiplicities, or if the existence of eigenvalues of higher multiplicity might somehow be topologically obstructed, has been open ever since. 
	
	In this thesis, we prove that on any closed spin manifold of dimension $m \equiv 0,6,7 \mod 8$, there exists a metric for which at least one eigenvalue is of higher multiplicity. 
	
	For the proof, we introduce a technique which ``catches'' the desired metric with a loop in the space of all Riemannian metrics.  We will construct such a loop on the sphere and transport it to a general manifold by extending some classical surgery theory results by Bär and Dahl. As a preparation, we will show that the Dirac spectrum can be described globally by a continuous family of functions on the space of Riemannian metrics and that the spinor field bundles with respect to the various metrics assemble to a continuous bundle of Hilbert spaces. These results might be useful in their own right.
	
	\vspace{3em}
	
	{\Large \textbf{Zusammenfassung}} \\
	Sei $(M,\Theta)$ eine kompakte Spin--Mannigfaltigkeit der Dimension $m \geq 3$ mit fester topologischer Spin--Struktur. Für jede Riemannsche Metrik $g$ erhalten wir einen Dirac-Operator $\Dirac^g$, dessen Spektrum von der Metrik abhängt. Dahl vermutet in einer Arbeit aus dem Jahr 2005, dass $M$ eine Metrik trägt, für die ein endlicher Teil des Dirac--Spektrums aus beliebigen vorgeschriebenen Eigenwerten beliebiger Multiplizitäten besteht. Nur der Eigenwert Null kann nicht beliebig vorgeschrieben werden (aufgrund des Atiyah-Singer Indexsatzes). Außerdem muss man in einigen Dimensionen die quaternionische Struktur der Eigenräume und die Symmetrie des Dirac--Spektrums beachten. Dahl beweist seine Vermutung für den Fall einfacher Eigenwerte. Die Frage ob auch Eigenwerte von beliebiger Multiplizität vorgeschrieben werden können oder ob die Existenz von Eigenwerten von höherer Multiplizität nicht möglicherweise topologisch obstruiert sein könnte, ist seit dem offen. 
	
	In der vorliegenden Arbeit zeigen wir, dass es auf jeder geschlossenen Spin--Mannigfaltigkeit der Dimension $m \equiv 0,6,7 \mod 8$ eine Metrik gibt, sodass der zugehörige Dirac-Operator mindestens einen Eigenwert von höherer Multiplizität besitzt. 
	
	Für den Beweis entwickeln wir eine Technik, die die gesuchte Metrik mit einer Schleife im Raum aller Riemannschen Metriken ,,einfängt''. Wir werden eine solche Schleife auf der Sphäre konstruieren und dann auf eine allgemeine Mannigfaltigkeit transportieren. Dazu erweitern wir einige klassische Resultate von Bär und Dahl aus der Chirurgietheorie. Als Vorbereitung werden wir zeigen, dass das Dirac--Spektrum vollständig durch eine stetige Familie von Funktionen auf den Riemannschen Metriken beschrieben werden kann und dass sich die Spinorfelder aller Metriken zu einem stetigen Bündel aus Hilberträumen zusammensetzen lassen. Diese Resultate könnten auch für sich genommen nützlich sein.

	\newpage

	\thispagestyle{empty}
	{\Large \textbf{Acknowledgements}} \\ $ $ \\
	I am very grateful to all the people who supported me during my time as a PhD student. I would like to thank my supervisor Bernd Ammann for introducing me to the beautiful topic of spin geometry and to Dahl's conjecture. His continuing support and our inspiring discussions gave my research the direction it needed, but also left enough freedom to pursue my own ideas. I am also grateful to my co-advisor Felix Finster for giving his support and his input. \\
	I'd like to say a big thank you to all our colleagues at the Faculty of Mathematics of the University of Regensburg for creating such a friendly and productive working atmosphere. In particular to Andreas Hermann, Jan-Hendrik Treude, Ulrich Bunke and Helmut Abels for interesting discussions. Special thanks go to Nicolas Ginoux and Olaf Müller for their strong interest in my work and for many inspiring discussions. \\
	Of course, I'd also like to thank Mattias Dahl for posing his conjecture and explaining parts of his previous work to me during his visits in Regensburg. 
	
	Many thanks also go to my parents for their lifelong support. Last but not least, I would like to thank Jesko Hüttenhain for our long lasting friendship.
	
	The research in this thesis would not have been possible without the kind support of the \emph{Studienstiftung des deutschen Volkes} and the \emph{DFG Graduiertenkolleg GRK 1692 ``Curvature, Cycles and Cohomology''}.

%% file: intro.overview.tex
\section{Motivation}
This thesis contributes to the study of the space of solutions of the Dirac equation 
\begin{align*}
	\Dirac^g \psi = \lambda \psi.
\end{align*}
Here, $\Dirac^g$ is the Dirac operator of a closed Riemannian spin manifold acting on sections $\psi$ of the spinor bundle. We are particularly interested in the \emph{spectrum of $\Dirac^g$}, i.e. the set of all $\lambda$ such that the above equation holds for $\psi \neq 0$.

Calculating explicitly the spectrum $\spec \Dirac^g$ of a Dirac operator  for a given closed Riemannian spin manifold $(M, g, \Theta^g)$ can be very difficult and this problem has been subject to extensive research. By now there are a lot of manifolds whose Dirac spectrum is well known. In many cases, estimates for some of the eigenvalues are possible. An extensive overview of available results of that kind can be found for instance in a textbook by Nicolas Ginoux, see \cite{GinouxDiracSpectrum}.

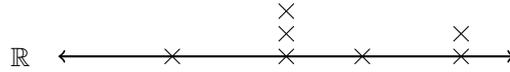
\begin{figure}[t] 
	\begin{center}
		\input{fig.diracspec}
		\caption[Typical Dirac spectrum]{A typical Dirac spectrum. Multiple crosses indicate higher multiplicites.}
		\label{FigDiracSpecTypical}
	\end{center}
\end{figure}

On the one hand, all Dirac spectra have certain properties in common. The Dirac operator is a self-adjoint elliptic differential operator of first order and as such has compact resolvent. Therefore, the Dirac spectrum $\spec \Dirac^g$ is always a closed discrete subset of the real line. The spectrum is also unbounded from both sides and all eigenvalues have finite multiplicities (see \cref{FigDiracSpecTypical} for a typical Dirac spectrum).

On the other hand, the various Dirac spectra are very different. They depend on the Riemannian metric chosen to define $\Dirac^g$. In \cite{FriedTori}, Thomas Friedrich gives an explicit formula for the Dirac spectrum of flat tori showing that it also depends on the spin structure. One might wonder, which sets can occur as Dirac spectra. To turn this into a precise mathematical question, we introduce the following problem, which is the central problem of this thesis.

\begin{Prb}
	\label{PrbDahl}
	Let $(M, \Theta)$ be a closed spin manifold of dimension $m \geq 3$ and $\Lambda_1 < \lambda_1 < \ldots < \lambda_k < \Lambda_2 \in \R$ be arbitrary real numbers and $\nu_1, \ldots, \nu_k \in \N$ be positive natural numbers. Can we find a Riemannian metric $g$ such that the complex Dirac operator $\Dirac^g$ satisfies
	\begin{align*}
		\spec \Dirac^g \cap \mathopen{]} \Lambda_1, \Lambda_2 \mathclose{[} = \{\lambda_1, \ldots, \lambda_k\}
	\end{align*}
	and for each $1 \leq j \leq k$, the eigenvalue $\lambda_j$ has multiplicity\footnote{The notion of \emph{multiplicity} in this context is more subtle than usual, see \cref{DefDahlMult} for details.} $\mu(\lambda_j)=\nu_j$? 
\end{Prb}

This problem was introduced by Mattias Dahl in \cite{DahlPresc}. In the same article, Dahl shows how to solve this problem in the case where all eigenvalues are \emph{simple} (i.e. if $\nu_j=1$ for all $1\leq j \leq k$), see \cref{ThmDahl} for the precise statement. One should remark that there are a few cases, in which \cref{PrbDahl} has no solution. These concern mainly the eigenvalue $0$ due to the Atiyah-Singer index theorem and symmetries of the spectrum in certain dimensions, see \cref{RemUnsolvableCases} for a precise classification of these cases. Apart from that, there seems to be no reason, why \cref{PrbDahl} should not be solvable for arbitrary multiplicities, which is precisely Dahl's conjecture:

\begin{Cnj}[\protect{\cite[Conjecture 3]{DahlPresc}}]
	Except for the algebraic constraints on the spinor bundle (giving quaternionic eigenspaces and symmetric spectrum) and the constraints on the zero eigenvalue coming from the Atiyah-Singer index theorem it is possible to find a Riemannian metric on any compact spin manifold with a finite part of its Dirac spectrum arbitrarily prescribed.
\end{Cnj}

Dahl also remarks that the technique he uses to solve the problem for simple eigenvalues does not work for eigenvalues with higher multiplicities. The problem of prescribing eigenvalues with higher multiplicities has been open ever since. However, the analogous problem for the Laplace operator has been solved, see \cref{SectLaplaceSchroedinger}.

\section{Main Results and Structure}
\label{SectMainResultsStructure}
The aim of this thesis is to introduce some new techniques to approach \cref{PrbDahl} in case of higher multiplicities. Ultimately, we will be able to show that the existence of higher multiplicities is at least not topologically obstructed. The precise result is as follows.

\begin{MainThm}[existence of higher multiplicities]
	\label{MainThmHigher}	
	\input{mainthm.higher}
\end{MainThm}

The restriction in the dimension stems from the fact that we require tools from real and from complex spin geometry. In dimensions $m \equiv 0,6,7 \mod 8$, complex spin geometry is the complexification of real spin geometry (see \cref{ThmCplxSpinGeomReal} for the precise result), which allows us to jump back and forth between the two.

Many constructions in spin geometry depend on the Riemannian metric $g$ on the manifold $M$: Namely, the metric spin structure $\Spin^g M$, the spinor bundle $\Sigma^g M$, the volume form and hence the space of sections $L^2(\Sigma^g M)$, the Clifford multiplication and of course the Dirac operator $\Dirac^g$ and its spectrum $\spec \Dirac^g$. Many texts in spin geometry either deal with a specific Riemannian spin manifold like a torus for instance or with an arbitrary, but fixed manifold $(M,g,\Theta^g)$ to study general geometric features. In both cases, the metric is fixed once and for all and it is even customary to drop the dependence of the metric in notation. Regarding the metric, this point of view is very static. Also, one usually considers only the complex Dirac operator, because the complex representation theory of $\Spin_m$ is easier.

To prove \cref{MainThmHigher}, we will inevitably have to deal with different metrics on the manifold. As a first step we have to study ``spin geometry in motion'', i.e. we have to precisely investigate the dependence of the constructions in spin geometry on the Riemannian metric. In \cref{ChapterReviewSpinGeometry} we will introduce these constructions in more detail. In particular, we will discuss topological spin structures and the subtleties between real and complex spin geometry. This is very important to make precise the ``right'' notion of the \emph{multiplicity} $\mu(\lambda)$ of an eigenvalue $\lambda$, see \cref{DefDahlMult}.

\cref{ChapterContDiracSpec,ChapterUSFB} deal with the dependence of the spectrum on the metric. From the point of view of \cref{MainThmHigher}, these extensions play the role of a Lemma. Nevertheless these results might be useful for other applications too.

In \cref{ChapterContDiracSpec}, we will give another formalization of the spectrum. Usually, $\spec \Dirac^g$ is viewed as a \emph{set} and to each element $\lambda \in \spec \Dirac^g$, we associate a natural number called \emph{multiplicity}. We will describe the spectrum by a non-decreasing function $\Z \to \R$ instead, or, since there is no canonical first eigenvalue, by an equivalence class of those functions, see \cref{EqDefsg,DefMonConf}. This has the advantage that the information of the set of all eigenvalues and their multiplicities is captured in a single mathematical object. This enables us to show that the spectrum in this sense depends continuously on the Riemannian metric, if one chooses the right topology, see \cref{MainThmSpec}. In particular, we obtain the following:

\begin{MainThm}[continuity of eigenvalue functions]
	\label{MainThmFun}
	\nomenclature[lambdaj]{$\lambda_j(g)$}{$j$-th eigenvalue of $\Dirac^g$}
	\input{mainthm.fun}
\end{MainThm}

This theorem generalizes the well known fact that a bounded spectral interval can be described locally by a continuous family of functions, see \cite[Prop. 7.1]{BaerMetrHarmonSpin}.
After proving this theorem, we will digress a bit and discuss to what extent these functions descend to certain \emph{moduli spaces} of the form $\Rm(M) / G$, where $G \subset \Diff(M)$, see \cref{SecModuliSpecFl}. This is closely related to the question, whether or not the indices of the eigenvalue functions are shifted, if one travels through $\Rm(M)$ along certain loops. This in turn is related to the notion of \emph{spectral flow}. We illustrate how one can use \cref{MainThmFun} to give an alternative definition of the spectral flow. Using a recent result from differential topology proven by Bernhard Hanke, Thomas Schick and Wolfgang Steimle in \cite{Hanke}, we will show that the spectral flow is non-trivial, see \cref{MainThmFlow}. The content of \cref{ChapterContDiracSpec} has been published by now in \cite{evpaper}. A more detailed overview of these results will be given in \cref{SctEvpaperIntro}.

For the proof of \cref{MainThmHigher}, controlling the eigenvalues of the Dirac operator will not be enough. We also need to control the eigenspinors. In \cref{ChapterUSFB}, we will review a well known construction introduced by Jean-Pierre Bourguignon and Paul Gauduchon in \cite{BourgGaud} to identify the spinor bundles with respect to two different metrics with one another. We will slightly generalize their results and reformulate them in a bundle theoretic language. We will show that there exists a continuous bundle of Hilbert spaces $L^2(\Sigma M) \to \Rm(M)$ such that for each Riemannian metric $g \in \Rm(M)$, the fibre over $g$ corresponds to $L^2(\Sigma^g M)$, the bundle of spinor fields for $g$, see \cref{LemUnivSpinBdleTrivs}. This will be proven by showing that the constructions in \cite{BourgGaud} of the isomorphisms to identify the spinor field bundles for two metrics are continuous in the metrics. Again, the main problem will be to define the ``right'' topologies. As an intermediate step, we will topologize $\G(\GLtp M)$, the gauge group of the topological spin structure and show that the topological spin structure $\Theta: \GLtp M \to \GLp M$ induces a nice covering map $\Theta_*:\G(\GLtp M) \to \G^{\spin}(\GLp M)$  between the corresponding gauge groups, see \cref{ThmSpinGaugeTrafosCov}. These technical preliminaries will also be useful in a later study of loops of spin diffeomorphisms. In preparation of \cref{ChapterHigher}, we also show that there are certain subsets of Riemannian metrics in $\Rm(M)$ over which the images of the spectral projection onto a fixed bounded interval $[\Lambda_1, \Lambda_2]$ assemble to a continuous vector bundle of finite rank, see \cref{ThmUSFBSpecSubs}.

\subsection{Proof Strategy for \cref{MainThmHigher}}
\label{SubSectPrfStrategyMainThmHigher}

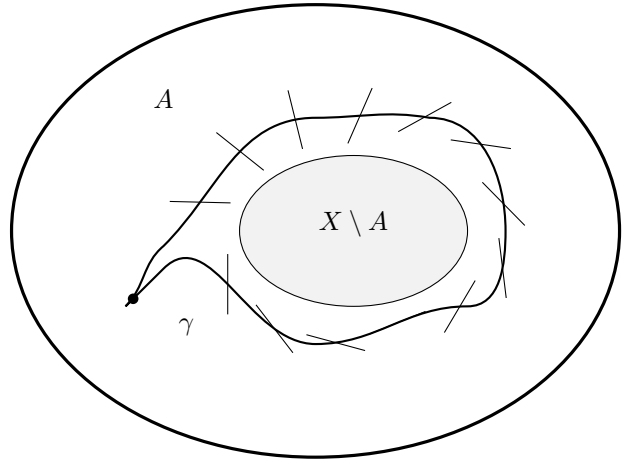
\begin{figure}[t] 
	\begin{center}
		\input{fig.lasso}
		\caption[The ``Lasso Lemma''.]{The ``Lasso Lemma''.}
		\label{FigLasso}
	\end{center}
\end{figure}

\cref{ChapterHigher} is devoted to the proof of \cref{MainThmHigher} and is the heart of this thesis. The key idea of the proof is the following simple topological reasoning to which we will refer to as the \emph{Lasso Lemma} (see also \cref{LemLasso}): Let $X$ be a simply connected topological space and let $A \subset X$ be a subspace. We want to show that $X \setminus A$ is not empty. Assume we can find a loop $\gamma:\S^1 \to A$. If $E \to A$ is a real vector bundle and $\gamma^* E \to \S^1$ is not orientable, then $\gamma$ cannot be null-homotopic. But if $\gamma$ is not null-homotopic, then $X \setminus A$ cannot be empty. Intuitively, the loop $\gamma$ ``catches'' the set $X \setminus A$. The situation is depicted in \cref{FigLasso}. We will apply this reasoning in the following way: For $X$ we take the space $(\Rm(M), \mathcal{C}^1)$ of all Riemannian metrics on $M$ endowed with $\mathcal{C}^1$-topology. The set $A$ will be $\Rm_A(M)$ (see \cref{DefPartEBundleGlobal}), which is a subspace of metrics tailor-made such that $X \setminus A \neq \emptyset$ directly implies \cref{MainThmHigher}. (The set $A$ contains the set of all metrics for which all eigenvalues are simple; we use $A$ instead of this set for technical reasons.) The bundle $E$ will be $E(M)$, which will consist of the span of the eigenspinors corresponding to a certain finite set of eigenvalues, see \cref{DefPartEBundleGlobal}. For the loop $\gamma$ we will have to construct a suitable loop $\mathbf{g}:\S^1 \to \Rm_{A}(M)$ of Riemannian metrics. 

The preliminary work of \cref{ChapterContDiracSpec} will make the definition of $\Rm_A(M)$ very easy and the preliminary work of \cref{ChapterUSFB} will enable us to quickly define the bundle $E(M)$ and to show that it is a continuous vector bundle over $\Rm_A(M)$ of finite rank, see \cref{CorPartEBGlobal}. The hard part will be to construct the loop $\mathbf{g}$ and to show that $\mathbf{g}^* E(M) \to \S^1$ is not trivial. Unfortunately, we will not be able to construct this loop directly. Therefore, we will use the following strategy: In \cref{SecLoopsSpinDiffeos}, we consider loops of spin diffeomorphisms $(f_{\alpha})_{\alpha \in \S^1}$ on $M$ and study loops of metrics induced by setting $g_{\alpha} := (f^{-1}_{\alpha})^* g$, $\alpha \in \S^1$, $g \in \Rm(M)$. We will work out a criterion when a bundle over this loop is trivial or not, see \cref{ThmS1ActionBundleTwisted}. This reduces the problem of finding a loop of metrics to finding a loop of spin diffeomorphisms. At first glance, this does not seem to help, since the construction of such a loop of diffeomorphisms is even more difficult than the construction of a loop of metrics. But in specific cases, the construction of a loop of diffeomorphisms is very easy. In \cref{SectSphere}, we will show that the family of rotations by degree $\alpha$ on the sphere $\S^m$ will suit our purpose, if we start with a metric $g_0$ that is obtained from the round metric by a small pertubation. This will give us the desired loop of metrics on the sphere $\S^m$.

\begin{figure}[t]
	\begin{center}
		\input{fig.spheresurgery}
		\caption[Connected sum with a sphere]{Connected sum with a sphere.}
		\label{FigSphereSurgery}
	\end{center}
\end{figure}
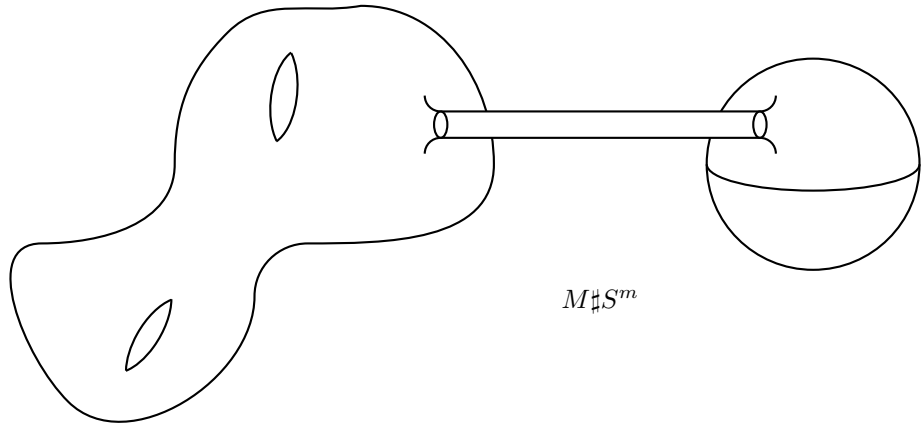

Finally, we will have to transport the loop of metrics on the sphere $\S^m$ to our original manifold $M$. Any smooth $m$-manifold $M$ is diffeomorphic to $M \sharp \S^m$, where $\sharp$ denotes a connected sum. Connected sums are a special type of \emph{surgery}. In \cref{SectSurgeryStability}, we will review the concept of surgery in the setting of Riemannian spin geometry and ultimately show that the existence of a suitable loop of metrics is stable under certain surgeries, see \cref{ThmSurgeryStabilityOdd}. Applying this to the connected sum will yield the desired result, see also \cref{FigSphereSurgery}.

\section{Comparison to Results for Laplace, Schrödinger and other Operators}
\label{SectLaplaceSchroedinger}
One should note that \cref{PrbDahl} has not only been formulated for the Dirac operator, but also for the Laplace operator and for the Schrödinger operator as well. We give a short overview of the results available for these operators. 

In 1986 Yves Colin de Verdière showed in \cite{verdiere1} that on any compact manifold $M$ of dimension $m \geq 3$ and for any positive integer $n \in \N$, there exists a Riemannian metric $g$ on $M$ such that the first eigenvalue of the Laplace operator $\Delta_g$ acting on functions has multiplicity $n$. In dimension $m=2$ this does not hold, since there are bounds on the multiplicity in terms of the genus of the surface. Shortly thereafter Colin de Verdière is able to considerably improve this result. In \cite{verdiere2} he shows that in case $m \geq 3$, for every finite sequence $0 = \lambda_1 < \lambda_2 \leq \ldots \leq \lambda_n$, there exists a metric $g$ on $M$ such that the Laplacian $\Delta_g$ acting on functions has this sequence as the first $n$ eigenvalues. Notice that a number might occur multiple times in the sequence, i.e. the eigenvalues are allowed to have higher multiplicities. The author also considers the analogous problem for Schrödinger operators $H_g = \Delta_g + V$, where $V$ is a potential such that $0$ is the smallest eigenvalue of $H_g$. He shows that on surfaces one can prescribe the first $n$ eigenvalues of $H_g$, but $n$ is bounded by the maximal number of vertices of a complete graph imbedded in $M$. Many parts of the article are formulated for more general classes of self-adjoint positive operators (notice that the Dirac operator is not positive). The author investigated the bound for the number $n$ of the Schrödinger operator in more detail in 1993, see \cite{verdiere3}.

In 2008 Pierre Jammes studied the analogous problem for the Hodge Laplacian on $p$-forms within a conformal class, see \cite{jammes1}. He is able to show that if the dimension of the manifold satisfies $m \geq 5$, then for any integer $p \in [2,m-2]$, $p \neq \tfrac{m}{2}$, any positive number $V$ and any increasing sequence $0 < \lambda_{p,1} < \ldots < \lambda_{p,n}$ there exists a Riemannian metric $g$ within a given conformal class such that the first $n$ eigenvalues of the Hodge Laplacian $\Delta_g = d \delta + \delta d$ acting on $p$-forms are simple and exactly equal to the sequence of $\lambda_j's$. In addition the volume of $(M,g)$ equals $V$. This strengthens a result obtained in 2004 by Pierre Guerini in \cite{guerini}, who studies this problem without the restriction to a conformal class on domains of $\R^m$ (with boundary conditions). In 2009 Jammes strengthens his result by allowing also double eigenvalues to occur in the sequence, see \cite{jammes3}. Finally, in 2011 Jammes solved the problem also for higher multiplicities. In \cite{jammes2} he shows that if $m \geq 6$, $p \in [2, \tfrac{m-3}{2}]$ is an integer, $V$ is any positive number and $0 < \lambda_{p,1} \leq \ldots \leq \lambda_{p,n}$, there exists a metric $g$ such that the first $n$ eigenvalues of the Hodge Laplacian $\Delta_g$ acting on co-exact $p$-forms are exactly given by the $\lambda_{p,j}$ and the volume of $(M,g)$ equals $V$. In particular, the multiplicities of the various eigenvalues can be arbitrarily large. In 2012 Jammes generalizes this result to the Witten Laplacian, see \cite{jammes4}. To any function $f \in \mathcal{C}^{\infty}(M)$ and any metric $g$, the \emph{Witten Laplacian} is defined by $\Delta_{(g,f)} := d_{f} \delta_f + \delta_f d_f$, where $d_{f} := e^{-f} d e^{f}$, $\delta_f := e^{f} \delta e^{-f}$.

Notice how the research on the problem of prescibing the eigenvalues of the Laplace operator has progressed: Jammes started with simple eigenvalues, advanced to double eigenvalues and finally considered eigenvalues of arbitrary multiplicity. Therefore, we think that a similar approach for the Dirac operator is reasonable. Also notice that it is possible to not only prescribe the eigenvalues, but also other geometric quantities like the volume. This is interesing because the geometry of the manifold and the spectrum of the Laplacian are intimitely related to each other.

A similar problem is given by the Laplace Operator $\Delta_{\Omega} = - \sum_i \partial_i^2$ on a domain $\Omega \subset \R^m$ with Dirichlet boundary conditions. The spectrum $\{\lambda_j(\Omega)\}_{j \in \N}$ of $\Delta_{\Omega}$ depends on $\Omega$, but it cannot be prescribed arbitrarily by varying $\Omega$ among all domains of $\R^m$ with a fixed volume. By the theorem of Faber-Krahn, the Ball $B$ of volume $c$ satisfies
\begin{align*}
	\lambda_1(B) = \min\{ \lambda_1(\Omega) \mid \Omega \subseto \R^m, |\Omega| = c \}.
\end{align*}
Analogously, by the theorem of Kran-Szegö, the minimum of $\lambda_2(\Omega)$ among all bounded open subsets of $\R^m$ with given volume is achieved by the union of two identical balls. A proof of these results (and many more results in this direction) can be found in \cite{henrot}.

While it is possible to prescribe eigenvalues of higher multiplicity for the Laplace operator, there are other physically motivated operators $L$ for which $Lu = \lambda u$ always implies that $\lambda$ is simple. For instance, consider the \emph{Sturm-Liouville operator} 
\begin{align*}
	L u := - \Big( \frac{d}{dx} \left( p \cdot \frac{d}{dx} \right) + q \Big)u = \lambda u 
\end{align*}
on $L^2([a,b])$ subject to the boundary conditions
\begin{align}
	\label{EqSturmLiouBoundary}
	\begin{split}
		c_a u (a) + d_a u'(a) &= 0, \\
		c_b u (b) + d_b u'(b) &= 0.
	\end{split}
\end{align}
for some fixed constants $c_a, d_a, c_b, d_b \in \R$. Here, $p$ is differentiable and positive and $q$ is continuous. As a domain for $L$ we can choose the closure of the $\mathcal{C}^2$ functions satisfying the boundary conditions \cref{EqSturmLiouBoundary} under the $L^2$-scalar product. Then $L$ is an elliptic self-adjoint operator of second order depending on the functions $p$ and $q$. However, any eigenvalue $\lambda$ of $L$ is always simple regardless of the choice of $p$ and $q$, see for instance \cite[Thm 4.1]{hartmann}.

\section{Most important Notation Conventions}
\label{SubSectNotationConventions}

The following notation conventions will be used throughout the thesis: 

\begin{tabular}[H]{ll}
	$M^m$ & a smooth closed spin manifold of dimension $m \geq 3$ \\
	$\Theta:\GLtp M \to \GLp M$ & a topological spin structure \\
	$\Rm(M)$ & space of Riemannian metrics on $M$ with $\mathcal{C}^1$-topology \\
	$\Theta^g:\Spin^g M \to \SO^g M$ & a metric spin structure \\
	$\Sigma^g_{\K} M$ & $\K$-spinor bundle w.r.t.  $g$ \\
	$\Dirac^g_{\K}: L^2(\Sigma^g_{\K} M) \to L^2(\Sigma^g_{\K} M)$ & Dirac operator \\
	$I:=[0,1]$ & the unit interval \\
\end{tabular}

\nomenclature[M]{$M$}{a smooth closed spin manifold of dimension $m \geq 3$}
\nomenclature[RmM]{$\Rm(M)$}{space of Riemannian metrics on $M$ with $\mathcal{C}^1$-topology}
\nomenclature[I]{$I$}{$I=[0,1] \subset \R$}

A comprehensive list of all the notation can be found in the \nameref{SectNomenclature}.

\begin{Rem}[real vs. complex spin geometry]
	\label{RemRealComplexNotation}
	The spinor bundle and the Dirac operator depend not only on the Riemannian metric, but also on the field $\K \in \{\R, \C\}$. In situations, where we want to stress this dependence and stress that a claim holds for both, we will index the spinor bundle and the Dirac operator with a $\K$, i.e. we will write $\Dirac^g_{\K}$ and $\Sigma^g_{\K} M$. In situations, where we want to compare real with complex spin geometry, we will put an $\R$ or a $\C$ in the index. In agreement with established literature, a spinor bundle or a Dirac operator without an index refers to the complex Dirac operator on the complex spinor bundle. 
\end{Rem}

%% file: fig.diracspec.tex
\begin{tikzpicture}

	\draw[<->,thick] (-3,0) to (3,0);

	\foreach \xp/\yp in {-1.5/0,0/0,0/0.3,0/0.6,1/0,2.3/0,2.3/0.3}{
		\draw (\xp-0.1,\yp-0.1) to (\xp+0.1,\yp+0.1);
		\draw (\xp-0.1,\yp+0.1) to (\xp+0.1,\yp-0.1);
	}
	
	\coordinate[label=$\R$] (A) at (-3.5,-0.25);

\end{tikzpicture}

%% file: mainthm.higher.tex
Let $(M, \Theta)$ be a closed spin manifold of dimension $m \equiv 0,6,7 \mod 8$. There exists a Riemannian metric $g$ on $M$ such that the complex Dirac operator $\Dirac^g_{\C}$ has at least one eigenvalue of multiplicity at least two. In addition, $g$ can be chosen such that it agrees with an arbitrary metric $\tilde g$ outside an arbitrarily small neighborhood on the manifold.

%% file: mainthm.fun.tex
Let $(M,\Theta)$ be a closed spin manifold and $\Rm(M)$ be the space of Riemannian metrics on $M$ endowed with $\mathcal{C}^1$-topology. There exists a family of continuous functions $\{ \lambda_j: \Rm(M) \to \R \}_{j \in \Z}$ such that for all $g \in \Rm(M)$, the sequence $(\lambda_j(g))_{j \in \Z}$ is non-decreasing and represents all the eigenvalues of the Dirac operator $\Dirac^g$ (counted with multiplicities). In addition, the sequence $\{ \arsinh(\lambda_j) \}_{j \in \Z}$ is equicontinuous.%

%% file: fig.lasso.tex
\begin{tikzpicture}

  \tikzset{
    decoration={
      markings,
      mark=between positions 0.125 and 0.90 step 0.0625 with {
        \xdef\maxseq{\pgfkeysvalueof{/pgf/decoration/mark info/sequence number}}
        \coordinate (pt-\maxseq)
        at (0,0);
      },
    }
  }

	
	\draw[black,very thick]
		(0,0) ellipse (4 and 3);

	\draw[gray!10,fill]
		(0.5,0) ellipse (1.5 and 1);

	\draw[black]
		(0.5,0) ellipse (1.5 and 1);

	\draw[thick,postaction={decorate}] (-2.5,-1) to (-2,-0.5)
		to [out=45,in=180] (0,-1.5)
		to [out=0,in=180] (2,-1)
		to [out=0,in=-90] (2.5,0)
		to [out=90,in=0] (1.5,1.5)
		to [out=170,in=0] (0,1.5)
		to [out=180,in=45] (-2,-0.2)
		to [out=220,in=45] (-2.5,-1);
		
	\foreach \c in {1,...,\maxseq}{
	\draw[black,rotate=1.35*360*(\c-1) / \maxseq] (pt-\c) +(0,-0.4) -- +(0,0.4);
	}
	
	\fill[black] (-2.4,-0.9) circle (0.07);
	

	\coordinate[label=above:$A$] (A) at (-2,1.5);
	\coordinate[label=above:$X \setminus A$] (A) at (0.5,-0.2);
	\coordinate[label=above:$\gamma$] (gamma) at (-1.7,-1.5);
	
\end{tikzpicture}

%% file: fig.spheresurgery.tex
\begin{tikzpicture}[scale=0.7]

\tikzset{
    partial ellipse/.style args={#1:#2:#3}{
        insert path={+ (#1:#3) arc (#1:#2:#3)}
    }
}

	\draw[thick] (6,0) circle (2);
	\draw[thick] (6,0) [partial ellipse=180:360:2 and 0.5];
	
	\draw[thick]
	(0,0)
	to [out=90, in=0] (-2.5,3)
	to [out=190, in=45] (-5,2.5)
	to [out=-135, in=90] (-6,0)
	to [out=-90, in=0] (-8.5,-1.5)
	to [out=180, in=135] (-8,-4.5)
	to [out=-45, in=-90] (-4.5,-2.5)
	to [out=90, in=180] (-3.5,-1.5)
	to [out=0, in=-90] (0,0);
	
	\draw[thick,rotate=-5] (-4.3,1) [partial ellipse=50:-68:0.5 and 1];
	\draw[thick,rotate=-6] (-3.8,0.8) [partial ellipse=115:232:0.5 and 1];
	
	\draw[thick,rotate=-30] (-4.3,-6) [partial ellipse=49:-58:0.5 and 1];
	\draw[thick,rotate=-30] (-3.7,-6.1) [partial ellipse=122:228:0.5 and 1];
	
	\fill[white] (-0.5,1) rectangle (0.5,0.5);	
	\fill[white] (4,1) rectangle (4.5,0.5);
		
	\draw[thick] (-1,1) to (5,1); 
	\draw[thick] (-1,0.5) to (5,0.5);	
	\draw[thick] [out=180,in=-90] (-1,1) to (-1.3,1.3); 
	\draw[thick] [out=180, in=90] (-1,0.5) to (-1.3,0.2); 
	\draw[thick] (5,1) [out=0, in=-90] to (5.3,1.3); 
	\draw[thick] (5,0.5) [out=0, in=90] to (5.3,0.2); 
	\draw[thick] (-1,0.75) ellipse (0.125 and 0.25); 
	\draw[thick] (5,0.75) ellipse (0.125 and 0.25);	
	
	\coordinate[label=above:$M \sharp S^m$] (A) at (2,-3);

\end{tikzpicture}

%% file: review.abstract.tex
\begin{subabstract}
In this chapter, we will review the foundations of classical spin geometry. Our aim is primarily to establish some notation and to clarify some subtleties between real and complex spin geometry. We also give a short reformulation of classical spin geometry in terms of category theory. The presentation will not be entirely self-contained and we are not aiming to give a full introduction to the subject. The reader interested in this is referred to the introductions already available, in particular \cite{FriedSpinGeom,hijazi,LM,Roe}. For readers already familiar with spin geometry, it will be sufficient to take note of the most important notation conventions, see \cref{SubSectNotationConventions}. A comprehensive list of all the notation can also be found in the \nameref{SectNomenclature}.
\end{subabstract}

%% file: review.clifford.tex
\section{Clifford Algebras}

In this section, we collect some purely algebraic facts about Clifford algebras. More details can be found in \cite[Chapter I]{LM} and \cite[Chapter 10]{greub2}. Let $K$ be a field of characteristic zero, $V$ be an $n$-dimensional vector space over $K$ and $q:V \to K$ be a quadratic form (or equivalently, a symmetric bilinear form).

\begin{Def}[Clifford algebra] \index{Clifford!algebra}
	A \emph{Clifford algebra} for $(V,q)$ is a unital $K$-algebra $A$ together with a map $\varphi:V \to A$ such that 
	\begin{align} 
		\label{EqDefCliffordProp}
		\forall v \in V: \varphi(v)^2 = -q(v) 1
	\end{align}
	and such that $(A,\varphi)$ satisfies the following universal property: For any other $K$-algebra $A'$ and any other map $\varphi':V \to A'$ satisfying \eqref{EqDefCliffordProp}, there exists a unique morphism of unital algebras $\psi:A \to A'$ such that $\psi \circ \varphi = \varphi'$, i.e. there exists a commutative diagram
	\begin{align} 
	\label{EqDefCliffordDiagram}
		\begin{split}
			\xymatrix{
				V
					\ar[r]^-{\varphi}
					\ar[d]_-{\varphi'}
				& A  
					\ar@{..>}[dl]^-{\exists! \psi}   \\
				A'.
			}
		\end{split}
	\end{align}
	A map $f:V \to A$ such that \eqref{EqDefCliffordProp} holds, is called \emph{Clifford}\index{Clifford!map}.
\end{Def}

\nomenclature[ClVq]{$\Cl(V,g)$}{Clifford algebra of $(V,g)$}
\begin{Thm}[existence and uniqueness of Clifford algebras]
	For any $(V,q)$, there exists a Clifford algebra $(A,\varphi)$, which is unique up to canonical isomorphisms (of unital algebras). We denote this object by
	\begin{align*} 
		\Cl(V) := \Cl(V,g) := A
	\end{align*}
    and for vectors $x, y \in V$, we set $x \cdot y := \varphi(x) \varphi(y) \in A$.
\end{Thm}

\begin{Prf} 
	The uniqueness part follows directly from the universal property in \cref{EqDefCliffordDiagram}. For the existence, we consider the tensor algebra 
	\begin{align*}
		T(V) := \bigoplus_{r=0}^{\infty}{V^{\otimes r}},
	\end{align*}
	where $V^{\otimes r}$, $r \in \N$, is the $r$-th tensor power of $V$ and $V^0 := K$. Let $I_q$ be the ideal generated by
	\begin{align*}
		\{ v \otimes v + q(v)1 \mid v \in V \},
	\end{align*}
	$A := T(V) / I_q$ and denote by $\pi$ the canonical projection. The map $\varphi$ is then given by
	\begin{align*}
		\xymatrix{
			V 
				\ar@{^(->}[r]
			& T(V)
				\ar@{->>}[r]^-{\pi}
			& A. }
	\end{align*}
\end{Prf}

\begin{Thm}[functoriality]
	\label{ThmCliffFunctor} 
	Let $f:(V,q) \to (V',q')$ be isometric, i.e.
	\begin{align*}
		\forall v \in V: q'(f(v)) = q(v).
	\end{align*}
	Then there exists a unique map of algebras $\tilde f=\Cl(f):\Cl(V) \to \Cl(V')$ such that
	\begin{align} 
		\label{EqCliffFunctorDef}
		\begin{split}
			\xymatrix{
				(V,q)
					\ar[d]_-{f}
					\ar[r]^-{\varphi}
				&\Cl(V,q)
					\ar@{..>}[d]^-{\tilde f} 
				\\
				(V',q')
					\ar[r]^-{\varphi'} 
				&\Cl(V',q')
				}
		\end{split}
	\end{align}
	commutes. In this sense, $\Cl$ is functorial, i.e. $\widetilde{\id_V} = \id_{\Cl(V)}$ and for any other isometric map $g:(V',q') \to (V'',q'')$, we have $\widetilde{g \circ f} = \tilde g \circ \tilde f$. 
\end{Thm}

\begin{Def}[grading]
\label{ThmCliffGrading}
	The Clifford Algebra has a natural grading constructed as follows: The map $\alpha:(V,q) \to (V,q)$, $v \mapsto -v$, is isometric. By \cref{ThmCliffFunctor}, it extends to a map $\tilde \alpha:\Cl(V,q) \to \Cl(V,q)$. This map satisfies $\tilde \alpha^2 = \id$. There exists a direct decomposition of $K$-vector spaces
	\begin{align*}
		\Cl(V,q) = \Cl^0(V,q) \oplus \Cl^1(V,q),
	\end{align*}
	where $\Cl^i(V,q)$, $i=0,1$, are the eigenspaces corresponding to the eigenvalues $(-1)^i$ of $\tilde \alpha$. We say $\Cl^0(V,q)$ is the \emph{even part} \index{even!part of the Clifford algebra} and $\Cl^1(V,q)$ is the \emph{odd part}\index{odd!part of the Clifford algebra}. 
	This decomposition satisfies
	\begin{align*}
		\forall i,j \in \{0,1\}: \Cl^i(V,q) \cdot \Cl^j(V,q) \subseteq \Cl^{i+j}(V,q), 
	\end{align*}
	where the indices are taken modulo 2. In particular, the even part is a subalgebra of $\Cl(V,q)$ (while the odd part is not).
\end{Def}

\nomenclature[Cl0V]{$\Cl^0(V,q)$}{even part of the Clifford algebra}
\nomenclature[Cl1V]{$\Cl^1(V,q)$}{odd part of the Clifford algebra}

\begin{Lem}[Clifford algebras of finite dimensional spaces]
	Let $(b_1, \ldots, b_n)$ be a basis of $V$. Then the $2^n$ vectors $1 \in \Cl(V,q)$ and
	\begin{align*}
		\{ x_{i_1} \cdot \ldots \cdot x_{i_k} \mid 1 \leq i_1 < \ldots < i_k \leq n, 1 \leq k \leq n \}
	\end{align*}
	form a basis of $\Cl(V,q)$.
\end{Lem}

\begin{Def}
    For $\K \in \{\R, \C\}$, $m \in \N$, we denote by $\Cl_{m,\K}$ the Clifford algebra of $\K^{m}$ endowed with the standard quadratic form. 
\end{Def}
\nomenclature[Clmk]{$\Cl_{m,\K}$}{Clifford algebra of $\K^{m}$}

\begin{Lem}
    \label{SpinCliffordUniversalCovering}
    For any $m \in \N$, define the group
    \begin{align*}
        \Spin_m := \{ x_{i_1} \cdot \ldots \cdot x_{i_{k}} \mid 0 \leq k \leq m, \forall 1 \leq \nu \leq k: |x_{i_\nu}| = 1 \}.
    \end{align*}
    For each $0 \neq x \in V$, let $\rho_x \in \SO_m$, be the reflection along the hyperplane $x^{\perp}$, i.e.
    \begin{align*}
        \forall y \in \R^m: \rho_x(y) = y - 2 \frac{\langle x, y \rangle }{\langle x, x \rangle }x,
    \end{align*}
    where $\langle \_, \_ \rangle$ is the Euclidean metric. Let $\vartheta_m:\Spin_m \to \SO_m$ be the unique group homomorphism satisfying 
    \begin{align*}
        \forall x \in \Spin_m: \vartheta_m(x) = \rho_x.
    \end{align*}
    Then for each $m \geq 2$, $\vartheta_m$ is a non-trivial double cover and for $m \geq 3$ it is universal.
\end{Lem}

%% file: review.spinstructs.tex
\section{Spin Structures}

\begin{Def}[frame bundles]\index{frame bundle}
	\nomenclature[GLE]{$\GL E$}{bundle of frames of $E$}
	\nomenclature[GLpE]{$\GLp E$}{bundle of positive frames on $E$}
	\nomenclature[SOE]{$\SO^h E$}{bundle of positive $h$-orthonormal frames on $E$}
	\nomenclature[bA]{$b.A$}{action of $A$ on $b$}
	Let $E \to M$ be a real vector bundle of rank $n$. We denote by $\pi_{\GL E}: \GL E \to M$ the principal $\GL_n$-bundle of frames of $E$. The action of a matrix $A \in \GL_n$ on a frame $b \in \GL E$ is denoted by $b.A \in \GL E$. In case $E$ is oriented, we denote by $\pi_{\GLp E}: \GLp E \to M$ the principal $\GLp_n$-bundle of positively oriented frames of $E$. In case $E$ is oriented and endowed with a Riemannian fibre metric $h$, we denote by $\pi_{\SO^h E}: \SO^h E \to M$ the principal $\SO_n$-bundle of positively oriented orthonormal frames on $E$. 
\end{Def}

\begin{Def}[spin structure]
	\index{spin structure!topological}
    \label{DefTopSpinStructure}
	Let $n \geq 2$ and $\vartheta_n:\GLtp_n \to \GLp_n$ be a connected non-trivial double cover of the Lie group $\GLp_n$ of invertible $n \times n$ matrices with positive determinant. Let $E \to M$ be a smooth oriented real vector bundle of rank $n$. A \emph{topological spin structure for $E$} is a $2:1$-covering map $\Theta: \GLtp E \to \GLp E$, where $\pi_{\GLtp E}: \GLtp E \to M$ is a principal $\GLtp_n$-fibre bundle such that
	\begin{align}
		\label{EqDefTopSpinStruct}
		\begin{split}
			\xymatrixcolsep{3.5em}
			\xymatrix{
				\GLtp E \times \GLtp_n
					\ar[r]
					\ar[d]^-{\Theta \times \vartheta_n}
				& \GLtp E
					\ar[d]^{\Theta}
					\ar[dr]^{\pi_{\GLtp E}}
				\\
				\GLp E \times \GLp_n
					\ar[r]
				& \GLp E
					\ar[r]_-{\pi_{\GLp E}}
				& M
			}
		\end{split}
	\end{align}
	commutes. Here, the horizontal arrows are given by the group actions. \end{Def}

\nomenclature[GLtpE]{$\GLtp E$}{topological spin structure on $E$}
\nomenclature[Theta]{$\Theta$}{a topological spin structure}
\nomenclature[thetan]{$\vartheta_n$}{the connected non-trivial double cover of $\GLp_n$}
\nomenclature[GLpn]{$\GLp_n$}{invertible $n \times n$ matrices with positive determinant}
\nomenclature[GLtpn]{$\GLtp_n$}{connected non-trivial double cover of $\GLp_n$}

\begin{Rem} $ $
\begin{enumerate}
	\item
		Since $\GLp_n$ is a Lie group, $\GLtp_n$ is a Lie group as well and $\vartheta_n$ is a homomorphism of Lie groups. 
	\item
		One can give topological criteria for existence and uniqueness of spin structures. An oriented vector bundle $E \to M$ admits a spin structure if and only if its second Stiefel-Whitney class vanishes. In that case, the equivalence classes of spin structures, see \cref{DefEquivalenceSpinStructs}, are in one-to-one correspondence with the elements of $H^1(M;\Z_2)$, see \cite[Thm. II.1.7]{LM}.
	\item
		Using the language of principal $G$-bundles, one can also express \cref{EqDefTopSpinStruct} by saying that $\Theta$ is a $\vartheta_n$-reduction of $\GLp E$ that it also a $2:1$-covering. See \cite[Kap. 2.5]{Baum} for more on reduction of principal fibre bundles. 
\end{enumerate}
\end{Rem}

\begin{Def}[spin manifold]
	\index{spin manifold}
	A \emph{spin manifold} is a tuple $(M, \Theta)$, where $M$ is a smooth oriented manifold and $\Theta$ is a topological spin structure for $TM$. We set $\GLp M := \GLp TM$ and $\GLtp M := \GLtp TM$. 
\end{Def}

\begin{Rem} $ $
	\begin{enumerate}
		\item 
			In our definition the spin structure is part of the data. Other authors define a spin manifold as a manifold, for which there exists a spin structure $\Theta$ on $TM$, but they exclude the spin structure from the data. Such a manifold could also be called a \emph{spinnable manifold}. The literature is not coherent in this matter.
		\item
			This definition of a spin manifold agrees with the definition in \cite{BaerGaudMor}. Traditionally, spin structures and spin manifolds are defined using the group $\Spin_m$ instead of $\GLtp_m$. Then one has to use metric spin structures, see \cref{DefMetricSpinStructure}, and this certainly makes the terminology more coherent. The obvious disadvantage of a metric spin structure is that it depends on the metric. The topological spin structure only needs the orientation, which makes it technically easier to handle various different metrics on a spin manifold.  
	\end{enumerate}
\end{Rem}

\begin{Def}[metric spin structure]
	\label{DefMetricSpinStructure}
	\index{spin structure!metric}
		Let $n > 2$ and let $\vartheta_n:\Spin_n \to \SO_n$ be the universal covering from \cref{SpinCliffordUniversalCovering}. Let $E \to M$ be a smooth oriented real vector bundle of rank $n$ and let $h$ be a Riemannian fibre metric for $E$. A \emph{metric spin structure for $E$} is a $2:1$-covering map $\Theta^h: \Spin^h E \to \SO^h E$, where $\pi_{\Spin^h E}: \Spin E^h \to M$ is a principal $\Spin_n$-fibre bundle such that
	\begin{align}
		\label{EqDefMetrSpinStruct}
		\begin{split}
			\xymatrixcolsep{3.5em}
			\xymatrix{
				\Spin^h E \times \Spin_n
					\ar[r]
					\ar[d]^-{\Theta^h \times \theta_n}
				& \Spin^h E
					\ar[d]^{\Theta^h}
					\ar[dr]^{\pi_{\Spin E}}
				\\
				\SO^h E \times \SO_n
					\ar[r]
				& \SO^h E
					\ar[r]_-{\pi_{\SO^h E}}
				& M
			}
		\end{split}
	\end{align}
	commutes. Here, the horizontal arrows are given by the group actions. In case $n=2$, a metric spin structure is defined analogously with $\Spin_2$ replaced by $\SO_2$ and $\vartheta_2$ replaced by the connected two-fold covering $\SO_2 \to \SO_2$. If $n=1$, a spin structure is defined to be a $2:1$-covering of $M$. 
\end{Def}

\begin{Def}[Riemannian spin manifold]
	\index{spin manifold!Riemannian}
    \label{DefRiemSpinMfd}
	A \emph{Riemannian spin manifold} is a triple $(M, g, \Theta^g)$, where $M$ is a smooth oriented manifold and $\Theta^g$ is a metric spin structure for $(TM,g)$. We set $\Spin^g M := \Spin^g TM$ and $\SO M := \SO^g TM$. 
\end{Def}

The precise relation between metric and topological spin structures is as follows.

\begin{Thm} 
	\label{ThmSpinStructTopToMetric}
	Let $E \to M$ be an oriented vector bundle of rank $n$ and $h$ be a Riemannian fibre metric for $E$. 
	\begin{enumerate}
		\item 
			Let $\Theta:\GLtp E \to \GLp E$ be a topological spin structure for $E$. Setting
			\begin{align*}
				\Spin^h E := \Theta^{-1}(\SO^h E), && \Theta^h:=\Theta|_{\Spin^h E}:\Spin^h E \to \SO^h E
			\end{align*}
			yields a metric spin structure for $E$.
		\item
			Conversely, let $\Theta^h:\Spin^h E \to \SO^h E$ be a metric spin structure for $E$. The canonical inclusion admits a lift
			\begin{align*}
				\xymatrix{
					\Spin_n
						\ar@{..>}[r]^-{\tilde \iota}
						\ar[d]^-{\theta_n}
					&\GLtp_n
						\ar[d]^-{\vartheta_n}
					\\
					\SO_n
						\ar@{^(->}[r]^-{\iota}
					& \GLp_n.
				}
			\end{align*}
			and the $\tilde \iota$-extension of $\Spin^h E$ yields a topological spin structure for $E$. 
	\end{enumerate}
\end{Thm}

\begin{Prf}
	This follows from general results about the extension and reduction of principal fibre bundles, see \cite[Kap. 2.5]{Baum}.
\end{Prf}

%% file: review.spinorbundles.tex
\section{Spinor Bundles}
\label{SectReviewSpinorBundles}
In this section, we review how to construct the spinor bundle $\Sigma ^g M$ out of the spin structure of a Riemannian spin manifold $(M, g, \Theta^g)$, see \cref{DefRiemSpinMfd}. We emphasize that one can construct a real as well as a complex spinor bundle. 
\nomenclature[K]{$\K$}{$\R$ \text{or} $\C$}

\begin{Def}[$K$-representation]
    \label{DefKReprAlg}
    \index{representation}
    Let $k$ be a field, $V$ be a $k$-vector space, $q$ be a quadratic form on $V$ and $\Cl(V,q)$ be the associated Clifford algebra. Let $K \supseteq k$ be a field. A \emph{$K$-representation} of the Clifford algebra $\Cl(V,q)$ is a $k$-algebra homomorphism
    \begin{align*}
        \rho:\Cl(V,q) \to \End_{K}(W),
    \end{align*}
    where $W$ is a vector space of finite dimension over $K$.
\end{Def}

\begin{Def}[irreducible]
    \label{DefIrreducibleRepr}
    \index{irreducible}
    Let $\rho:\Cl(V,q) \to \End_{K}(W)$ be a $K$-representation as in \cref{DefKReprAlg}. Then $\rho$ is \emph{reducible}, if there exists a decomposition $W = W_1 \oplus W_2$ over $K$ such that 
    \begin{align*}
        \forall \varphi \in \Cl(V,q): \rho(\varphi)(W_j) \subset W_j, \qquad j=1,2
    \end{align*}
    and $\dim W_j \neq 0$, $j=1,2$. A representation is \emph{irreducible}, if it is not reducible.
\end{Def}

\begin{Rem}
    These definitions agree with \cite[Def. 5.1, 5.3]{LM} and are tailor made to formulate the representation theory of the Clifford algebra. In \cref{DefKReprAlg}, we allow the field to be larger to account for the fact that certain representations of the real Clifford algebra are automatically complex or even quaternionic, see \cite[Thm. 5.8]{LM}. 
\end{Rem}

\begin{Def}[equivalence]
    \index{equivalence!of representations}
    Two $K$-representations $\rho_{j}:\Cl(V,q) \to \End_{K}(W_j)$, $j=1,2$, are \emph{equivalent}, if there exists an isomorphism $F \in \Iso_{K}(W_1,W_2)$ such that for any $\varphi \in \Cl(V,q)$, the following diagram commutes
    \begin{align*}
        \xymatrix{
            W_1
                \ar[r]^{\rho_1(\varphi)}
                \ar[d]^-{F}
            &W_1
                \ar[d]^-{F}
            \\
            W_2
                \ar[r]^-{\rho_2(\varphi)}
            &W_2.
        }
    \end{align*}
\end{Def}

\begin{Def}[representation]
    \label{DefReprGroup}
    Let $\K \in \{\R, \C\}$. A \emph{$\K$-representation}\index{K-representation} of a group $G$ on a $\K$-vector space $W$ is a group homomorphism $\rho: G \to \GL_{\K}(W)$. 
\end{Def}

\begin{Def}[spinor representation]
	\label{DefSpinorRepresentation}
    Let $\rho_{m,\K}:\Cl_{m,\K} \to \End_{\K}(\Sigma_{\K,m})$ be an irreducible representation. Then
    \begin{align*}
        \Delta_{m,\K} := \rho_{m,\K}|_{\Spin_m}:\Spin_m \to \GL_{\K}(\Sigma_{\K,m})
    \end{align*}
    is a \emph{$\K$-spinor representation}. 
\end{Def}
\nomenclature[DmK]{$\Delta_{m,\K}$}{a $\K$-spinor representation}

\begin{Rem}
    Any $\K$-spinor representation is a $\K$-representation of the group $\Spin_m$ in the sense of \cref{DefReprGroup}. 
\end{Rem}

\begin{Def}[spinor bundle]
    \index{spinor bundle}
	\label{DefSpinorBundle}
	Let $\Delta_{m,\K}:\Spin_m \to \GL_{\K}(\Sigma_{\K,m})$ be a $\K$-spinor representation and $(M^m, g, \Theta^g)$ be a Riemannian spin manifold with metric spin structure $\Theta^g:\Spin^g M \to \SO^g M$. The associated vector bundle
	\begin{align*}
		\Sigma_{\K}^g M := \Spin^g M \times_{\Delta_{m,\K}} \Sigma_{\K,m} \to M
	\end{align*}
	is a \emph{$\K$-spinor bundle of $M$ with respect to $g$}. Its elements are called \emph{spinors} and its sections are called \emph{spinor fields}. 
\end{Def}
\nomenclature[SigmaKgM]{$\Sigma_{\K}^g M$}{$\K$-spinor bundle on $M$ with respect to $g$}

\begin{Rem}
	\label{RemNoUSBByRepresentations}
	Unfortunately, there exists no analogue of \cref{ThmSpinStructTopToMetric} for spinor bundles. By \cref{ThmSpinStructTopToMetric}, a topological spin structure $\Theta:\GLtp M \to \GLp M$ contains the metric spin structures $\Theta^g$ for any metric $g$ by setting 
	\begin{align*}
		\Theta^g := \Theta|_{\Theta^{-1}(\SO^g M)}:\Spin^g M \to \SO^g M.
	\end{align*}
	Of course, one can take any representation $\tilde \rho_m:\GLtp_m \to \Aut_{\K}(W)$ and define the bundle $\GLtp M \times_{\tilde \rho_m} W$. However, no finite dimensional representation $\tilde \rho_m$ restricts to a $\K$-spinor representation. The reason for this is as follows: The $\K$-spinor representations are all faithful. By \cite[Lem. II.5.23]{LM}, for any finite-dimensional representation $\tilde \rho_m$, $m>2$, of $\GLtp_m$,  there exists a representation $\rho_m$ such that
	\begin{align*}
		\xymatrix{
			\GLtp_m
				\ar[r]^-{\tilde \rho_m}
				\ar[d]^-{\vartheta}_-{2:1}
			&\Aut_{\K}(W)
			\\
			\GLp_m 
				\ar@{..>}[ur]_-{\rho_m}
		}
	\end{align*}
	commutes. This implies that $\tilde \rho_m|_{\Spin_m}$ cannot be faithful and hence cannot agree with a $\K$-spinor representation. 
	
	The problem of forming a spinor bundle without reference to any metric is the problem of finding a \emph{universal spinor bundle}. We will come back to this issue in \cref{ChapterUSFB}. In \cref{DefUniversalSpinorBundle}, we will see that such a finite dimensional bundle does exist, if one replaces the finite dimensional manifold $M$ by the infinite-dimensional manifold $\Rm(M) \times M$.
\end{Rem}

\begin{Thm}[additional structures on the spinor bundle]
	\label{ThmSpinorBundleAddStructure}
    Let $\Sigma^g_{\K} M \to M$ be the $\K$-spinor bundle of a Riemannian spin manifold $(M,g,\Theta^g)$.
    \begin{enumerate}
        \item 
            There exists a morphism of vector bundles $\mathfrak{m}^g_{\K}:TM \otimes_{\R} \Sigma^g_{\K} M \to \Sigma^g_{\K} M$, $X \otimes \psi \mapsto X \cdot \psi$, called \emph{Clifford multiplication}\index{Clifford!multiplication}, such that
            \begin{align*}
                \forall X \in TM: \forall \psi \in \Sigma^g_{\K} M: X \cdot X \cdot \psi = - g(X,X) \psi.
            \end{align*}
        \item
            There exists a Riemannian fibre metric $\langle \_, \_ \rangle$ on $\Sigma^g_{\R} M$ respectively a Hermitian fibre metric $\langle \_, \_ \rangle$ on $\Sigma^g_{\C} M$ such that the Clifford multiplication is skew-adjoint, i.e.
            \begin{align*}
                \forall X \in TM: \forall \psi_1, \psi_2 \in \Sigma^g_{\K} M: \langle X \cdot \psi_1, \psi_2 \rangle = - \langle \psi_1, X \cdot \psi_2 \rangle.
            \end{align*}
        \item
            There exists a connection $\nabla^{g,\K}$ on $\Sigma^g_{\K} M$, called \emphi{spinorial Levi-Civita connection}, such that
            \begin{align*}
                \forall X,Y \in \mathcal{T}(M): \forall \psi \in \Gamma(\Sigma^g_{\K} M): \nabla^{g,\K}_X(Y \cdot \psi) = \nabla^g_X Y \cdot \psi + Y \cdot \nabla^{g,\K}_X \psi,
            \end{align*}
            where $\nabla^g$ denotes the Levi-civita connection on $M$.
    \end{enumerate}
\end{Thm}
\nomenclature[mK]{$\mathfrak{m}^g_{\K}$}{Clifford multiplication}
\nomenclature[nablaK]{$\nabla^{g,\K}$}{spinorial Levi-Civita connection}

\begin{Prf} $ $
	\begin{enumerate}
		\item 
			The construction is roughly as follows (see also \cite[Def. 4.2.ii)]{hijazi}): Let $x \in M$,  $X = [\Theta(\tilde b), v] \in \SO^g M_x \times_{\tau} \R^n \cong T_xM$, where $\tau:\SO_m \to \GL(\R^m)$ is the canonical representation. For any $\psi = [\tilde b, \sigma] \in \Spin^g M \times_{\Delta_{m,\K}} \Sigma_{\K,m} = \Sigma^g_{\K} M|_x$, we define
			\begin{align*}
				X \cdot \psi = [\Theta(\tilde b), v] \cdot [\tilde b, \sigma] := [\tilde b, v \cdot \sigma],
			\end{align*}
			where $v \cdot \sigma := \Delta_{m,\K}(v)(\sigma)$. Here, $\rho_{m,\K}:\Cl_m \to \GL_{\K}(\Sigma_{\K,m})$ is an irreducible representation of the Clifford algebra, see \cref{DefSpinorRepresentation}.
		\item
			Again, let $\rho_{m,\K}:\Cl_m \to \GL_{\K}(\Sigma_{\K,m})$ be as in \cref{DefSpinorRepresentation}. By averaging over the multiplicative subgroup of the Clifford algebra generated by an orthonormal basis of $\R^m$, one can show that there exists an invariant inner product of $\Sigma_{\K,m}$, for which the action of $\rho_{m,\K}$ is skew-adjoint, see \cite[Cor. 2.27]{hijazi}. This induces the desired inner product on $\Sigma^g_{\K} M$, see \cite[Def. 4.3.ii)]{hijazi}. 
		\item
			By the Fundamental Lemma of Riemannian geometry, see \cite[Thm. 5.4]{LRM}, the Levi-Civita connection $\nabla^g$ exists on the vector bundle $TM$. This gives a connection form $\omega \in \Omega^1(\SO^g M, \mathfrak{{so}_m})$ on the principal $\SO_m$-bundle $\SO^g M$. One can show that this connection form lifts to a connection form $\tilde \psi \in \Omega^1(\Spin^g M, \mathfrak{spin}_m)$, see \cite[Sect. 4.2]{hijazi}. By a general procedure, this connection form induces a connection on all vector bundles associated to $\Spin^g M$, see \cite[Satz 3.12]{Baum}. In particular, this gives a connection on the spinor bundle $\Sigma^g_{\K} M = \Spin^g M \times_{\Delta_{m,\K}} \Sigma_{\K,m}$. 
	\end{enumerate}
\end{Prf}

\begin{Rem}[local expressions]
	Let $b=(b_1, \ldots, b_m)$ be a local orthonormal frame on $M$, $\tilde b$ be its lift to $\Spin M$ and $\sigma_1, \ldots, \sigma_N$ be a basis of $\Sigma_{\K,m}$. Then the spinor fields locally given by $\psi_{\alpha} = [\tilde b, \sigma_{\alpha}]$ satisfy 
	\begin{align*}
		\forall X \in \mathcal{T}(M): \nabla^{g,\K}_X \psi_{\alpha} = \frac{1}{4} \sum_{i,j=1}^{m}{g(\nabla^g_X b_i, b_j) b_i \cdot b_j \cdot \psi_\alpha}.
	\end{align*}
	Using this formula, one can calculate $\nabla^{g, \K}_X \psi$ locally for any spinor field $\psi$. The curvature $R^{g,\K}$ of $\nabla^{g,\K}$ is related to the curvature $R^g$ of the Levi-Civita connection $\nabla^g$ by
	\begin{align*}
		\forall X, Y \in \mathcal{T}(M): \forall \psi \in  \Gamma(\Sigma^g_{\K} M): R^{g,\K}_{X,Y} \psi  = \frac{1}{4} \sum_{i,j=1}^{m}{g(R^g_{X,Y}(b_i), b_j) b_i \cdot b_j \cdot \psi}.
	\end{align*}
	More details can be found in \cite[Prop. 4.3]{hijazi}.
\end{Rem}

\begin{Def}[Dirac operator]
	\label{DefDiracOperatorSpinMfd}
	\index{Dirac operator}
	Let $(M,g,\Theta^g)$ be a Riemannian spin manifold with spinor bundle $\Sigma^g_{\K} M$. Let $\mathfrak{m}^g_{\K}$ be the Clifford multiplication, $\nabla^{g,\K}$ be the spinorial Levi-Civita connection and $\sharp^g: T^*M \to TM$ be the musical isomorphism induced by $g$. 
	The composition
	\begin{align*}
		\xymatrixcolsep{3.5em}
		\xymatrix{
			\Gamma(\Sigma_{\K}^g M)
				\ar[r]^-{\nabla^{\K,g}}
			& \Gamma(TM^* \otimes \Sigma_{\K}^g M)
				\ar[r]^-{\sharp^g \otimes \id}
			& \Gamma(TM \otimes \Sigma_{\K}^g M)
				\ar[r]^-{\mathfrak{m}_{\K}^g}
			& \Gamma(\Sigma_{\K}^g M)
			}
	\end{align*}
	is the associated \emph{Dirac operator} $\Dirac^g_{\K}$. We also think of $\Dirac^g_{\K}$ as an unbounded operator 
	\begin{align*}
		\Dirac_{\K}^g : L^2(\Sigma_{\K}^g M) \to L^2(\Sigma_{\K}^g M),
	\end{align*}
	densely defined on the first order Sobolev space $H^1(\Sigma_{\K}^g M)$. 
\end{Def}
\nomenclature[DK]{$\Dirac_{\K}^g$}{Dirac operator}

\begin{Thm}[properties of Dirac operators]
	\label{ThmSpecPropsDirac}
	The Dirac operator $\Dirac^g_{\K}$ of a closed Riemannian spin manifold $(M^m,g,\Theta^g)$ is a first order elliptic differential operator. For any local orthonormal frame $e_1, \ldots, e_m$ of $TM$
	\begin{align*}
		\Dirac^{g}_{\K} = \sum_{j=1}^{m}{e_{j} \cdot \nabla^{g,\K}_{e_j}}.
	\end{align*}
	As an operator $L^2(\Sigma^g_{\K} M) \to L^2(\Sigma^g_{\K} M)$, the Dirac operator $\Dirac^{g}_{\K}$ is self-adjoint, has compact resolvent and its spectrum $\spec \Dirac^{g}_{\K}$ is a closed discrete subset of the real line that is unbounded from both sides. For any $\lambda \in \spec \Dirac^{g,\K}$ the eigenspace $\ker (\Dirac^{g}_{\K} - \lambda)$ is finite-dimensional over $\K$.
\end{Thm}

\begin{Prf}
	The local coordinate formula follows directly from the definition. From this formula it follows that $\Dirac^g_{\K}$ is an elliptic first order differential operator. Self-adjointness follows from a calculation that is carried out for instance in \cite[Lem. 4.7 ]{hijazi} or alternatively \cite[Chapter II, Prop. 5.3]{LM} and the fact that the manifold is compact. The rest of the claim follows from the general theory for linear elliptic differential operators and abstract functional analysis, see \cite[Sect. III.§5]{LM}, or alternatively \cite[Sect. 4.2]{FriedSpinGeom} for a proof specifically for the Dirac operator.
\end{Prf}

%% file: review.realcomplex.tex
\section{Real vs. Complex Spin Geometry}    
In \cref{SectReviewSpinorBundles}, we introduced the real as well as the complex spinor bundle. In this section, we discuss some issues concerning the relationship between the two. 

\begin{Def}[complexification]
    If $W$ is an $\R$-vector space, the space $W^{\C} := W \otimes_{\R} \C$ is the \emph{complexification of $W$}.\index{complexification!of a vector space} The \emph{complexification}\index{complexification!of a representation} of a real representation of an $\R$-algebra $\rho:A \to \End_{\R}(W)$ is the map
        \DefMap{\rho^{\C}:A}{\End_{\C}(W^{\C})}{a}{(w \otimes \zeta \mapsto \rho_a(w) \otimes \zeta).}
    The analogous notion applies to the complexification of a real representation of a group.
\end{Def}
\nomenclature[WC]{$W^{\C}$}{complexification of $W$}
\nomenclature[rhoC]{$\rho^{\C}$}{complexification of $\rho$}

The following will be absolutely crucial.

\begin{Lem}
	\label{LemReprSpinGrp}
    In dimensions $m \equiv 0,6,7 \mod 8$, the complexification of the real spinor representation is equivalent to the complex spinor representation, i.e. $\Delta_{m,\R}^{\C} = \Delta_{m,\C}$.
\end{Lem}

\begin{Prf}
	It follows from the explicit classification of Clifford algebras (see for instance \cite[Thm 5.8]{LM}) that in dimensions $m \equiv 0,6,7 \mod 8$, the complexification of the real Clifford algebra is isomorphic to the complex Clifford algebra and that the irreducible complex representations are precisely the complexifications of irreducible real representations. This implies the result.
\end{Prf}

This representation theoretic fact has analogous consequences for all the spin geometric constructions: In dimensions $m \equiv 0,6,7 \mod 8$, the complexification of real spin geometry is isomorphic to the complex spin geometry. More precisely, the following hold.

\begin{Thm}
	\label{ThmCplxSpinGeomReal}
	Let $m = \dim M \equiv 0,6,7 \mod 8$. 
	\begin{enumerate}
		\item
			The complexification $(\Sigma_{\R}^g M)^{\C}$ of a real spinor bundle $\Sigma_{\R}^g M$ is a complex spinor bundle $\Sigma^g_{\C} M$.
		\item
			The complexification $(\mathfrak{m}_{\R}^g)^{\C}$ of the real Clifford multiplication $\mathfrak{m}^g_{\R}$ is the complex Clifford multiplication $\mathfrak{m}^g_{\C}$.
		\item
			The complexification $(\nabla^{g,\R})^{\C}$ of the real spinorial Levi-Civita connection on $\Sigma_{\R}^g M$ is the spinorial Levi-Civita connection $\nabla^{g,\C}$ on $\Sigma_{\C}^g M$.
		\item
			The complexification $(\Dirac_{\R}^g)^{\C}$ of the real Dirac operator $\Dirac_{\R}^g$ is the complex Dirac operator $\Dirac_{\C}^g$.
		\item
			Any $\lambda \in \R$ is an eigenvalue of $\Dirac_{\R}^g$ if and only if $\lambda$ is an eigenvalue of $\Dirac_{\C}^g$. 
	\end{enumerate}
\end{Thm}

\begin{Prf}
	\cref{LemReprSpinGrp} implies these results.
\end{Prf}

It is clear that the eigenspaces of $\Dirac^g_{\K}$ have the structure of a $\K$-vector space. In dimensions $m \equiv 2,3,4 \mod 8$, the eigenspaces of the complex Dirac operator $\Dirac^g_{\C}$ even have a quaternionic structure and hence can be thought of as vector spaces over the quaternions $\H$, see \cite[sec. 1.7]{FriedSpinGeom}. This yields to various notions of \emph{multiplicity} of an eigenvalue, which we now make precise. To that end, we first introduce a notation, which will also be useful later.

\nomenclature[SigmaA]{$L^2_A(\Sigma^g_{\K} M)$}{eigenspace of $A$}
\begin{Def}
	For any subset $A \subset \R$, we define
	\begin{align*}
		L^2_A(\Sigma^g_{\K} M)
		:= \spn \{\psi \in H^1(\Sigma^g_{\K} M) \mid \exists \lambda \in A: \Dirac^g_{\K} \psi = \lambda \psi\}.
	\end{align*}
\end{Def}	

\begin{DefI}[multiplicity]
	\label{DefDahlMult}
	For any eigenvalue $\lambda \in \spec \Dirac_{\C}^g$, the number
	\begin{align*}
		\mu(\lambda) := 
		\begin{cases}
			\dim_{\H}{L^2_{\lambda}(\Sigma_{\C}^g M)}, & m \equiv 2, 3, 4 \mod 8, \\
			\dim_{\C}{L^2_{\lambda}(\Sigma_{\C}^g M)}, & m \equiv 0, 1, 5, 6, 7 \mod 8,
		\end{cases}
	\end{align*}
	is called \emph{spin multiplicity of $\lambda$}. For an eigenvalue $\lambda \in \spec \Dirac_{\K}^g$, the number 
	\begin{align*}
		\mu_{\K}(\lambda) := \dim_{\K}{L^2_{\lambda}(\Sigma_{\K}^g M)}
	\end{align*}
	is called \emph{$\K$-multiplicity} of $\lambda$. An eigenvalue is called \emph{spin-simple}, if $\mu(\lambda) =1$ and \emph{$\K$-simple}, if $\mu_{\K}(\lambda) = 1$. \index{simple eigenvalue}
\end{DefI}
\nomenclature[mu]{$\mu(\lambda)$}{spin multiplicity of an eigenvalue $\lambda$}
\nomenclature[muK]{$\mu_{\K}(\lambda)$}{$\K$-multiplicity of an eigenvalue $\lambda$}

\begin{Rem}
	\label{RemMultK}
	In particular, any $\lambda \in \R$ satisfies 
	\begin{align*}
		\mu(\lambda) = 
		\begin{cases}
			\mu_{\C}(\lambda) = \mu_{\R}(\lambda), & m \equiv 0,6,7 \mod 8, \\
			\mu_{\C}(\lambda) = \tfrac{1}{2} \mu_{\R}(\lambda), & m \equiv 1,5 \mod 8, \\
			\tfrac{1}{2} \mu_{\C}(\lambda) = \tfrac{1}{4} \mu_{\R}(\lambda), & m \equiv 2, 3, 4 \mod 8.
		\end{cases}
	\end{align*}
	$ $
\end{Rem}

\begin{Rem}
	In closing we remark that the construction of a Dirac operator as described in \cref{DefDiracOperatorSpinMfd} admits various generalizations. First of all, one can replace the Riemannian metric by a pseudo-Riemannian metric of signature $(r,s)$. Then in principle, the construction is the same, but the representation theory of the Clifford algebra $\Cl(\R^{r,s})$ of $\R^{s+r}$ with a symmetric bilinear form of signature $(r,s)$ is much more complicated, see \cite[Chapter I. \textsection 4]{LM}. Also, the induced metric on the spinor bundle from \cref{ThmSpinorBundleAddStructure} will not be Riemannian in general. See \cite{BaumPseudo} for more details and also \cite{HarveySpinors} for the algebraic side.
	
	One can also twist the Dirac operator of the manifold with a vector bundle. Even more general, one can forget that the structures described in \cref{ThmSpinorBundleAddStructure} arise from a spin manifold and study so called \emph{Clifford bundles}, see also \cite[Def. 3.4]{Roe}. We will elaborate on this a bit further in \cref{DefCliffordBundles}. 
\end{Rem}

%% file: review.dahl.tex
\section{Dahl's Result}
Having introduced the precise notion of multiplicity, see \cref{DefDahlMult}, we are now able to state Dahl's Theorem.

\begin{Thm}[\protect{{\cite[Thm. 1]{DahlPresc}}}]
	\label{ThmDahl}
	Let $M$ be a closed spin manifold of dimension $m \geq 3$ and $L > 0$ be a real number. 
	\begin{enumerate}
		\item
			\label{ItDahlNotSym}
			Suppose that $m \equiv 3 \mod 4$ and let $l_1, \ldots, l_n$ be non-zero real numbers such that
			\begin{align*}
				-L < l_1 < \ldots < l_n < L.
			\end{align*}
			Then there exists a metric $g \in \mathcal{R}(M)$ such that
			\begin{align*}
				\spec \Dirac^g_{\C} \cap (\mathopen{]}-L, L \mathclose{[}\setminus \{ 0 \}) = \{ l_1, \ldots, l_n \}
			\end{align*}
			and all the $l_i$ satisfy $\mu(l_i)=1$. 
		\item 
			\label{ItDahlSym}
			Suppose $m \not\equiv 3 \mod 4$ and $l_1, \ldots, l_n, L \in \R$ such that
			\begin{align*}
				0 < l_1 < \ldots < l_n < L .
			\end{align*}
			Then there exists a metric $g \in \mathcal{R}(M)$ such that
			\begin{align*}
				\spec \Dirac^g_{\C} \cap (\mathopen{]}-L,L\mathclose{[} \setminus \{ 0 \}) = \{ \pm l_1, \ldots, \pm l_n \}
			\end{align*}
			and all the $\pm l_i$ satisfy $\mu(\pm l_i)=1$. 
	\end{enumerate}
\end{Thm}

\begin{Rem}[unsolvable cases]
	\label{RemUnsolvableCases}
	We also give an overview about the situations, where \cref{PrbDahl} is most certainly unsolvable and explain why.
	\begin{description}[style=multiline,leftmargin=3cm]
		\item[Quaternionic Structures]
			In dimensions $m \equiv 2,3,4 \mod 8$, the eigenspaces of the Dirac operator have a quaternionic structure, hence their complex dimension will always be even. Consequently, one cannot prescribe eigenvalues with odd complex multiplicity in these dimensions. This is precisely the reason, why we defined the multiplicity $\mu(\lambda)$ of an eigenvalue $\lambda$ in \cref{DefDahlMult} as the quaternionic dimension in this case. Using $\mu$ and not $\mu_{\C}$ in \cref{PrbDahl} already excludes this case from the problem statement.
		\item[Symmetries]
			In dimensions $m \not \equiv 3 \mod 4$, the Dirac spectrum is always symmetric, i.e. $\spec \Dirac^g = - \spec \Dirac^g$, see for instance \cite[Thm. 1.3.7.iv)]{GinouxDiracSpectrum}. Consequently, one can never prescribe a set of eigenvalues, which is unsymmetric. This is the reason why the statement of \cref{ThmDahl} is split up into two cases. If $m \not \equiv 3 \mod 4$, one can prescribe eigenvalues on the positive real line, which automatically prescribes the eigenvalues on the corresponding part of the negative real line.
		\item[Atiyah-Singer index Theorem]
			The eigenvalue $\lambda = 0$ plays an exceptional role. Its eigenspace is the kernel of the Dirac operator and its dimension is related to the topology of the manifold by the Atiyah-Singer index theorem, see for instance \cite[III.\textsection 13]{LM}. One can conclude from this theorem that for any Riemannian metric $g$ and any closed connected spin manifold $M$,
			\begin{align}
				\label{EqAtiyahBounds}
				\begin{split}
					\dim_{\C} \ker \Dirac^g_{\C} \geq 
					\begin{cases}
						|\hat A(M)|, & m \equiv 0,4 \mod 8, \\
						1, & m \equiv 1 \mod 8 \text{ and } \alpha(M) \neq 0, \\
						2, & m \equiv 2 \mod 8 \text{ and } \alpha(M) \neq 0, \\
						0, & \text{otherwise}.
					\end{cases}
				\end{split}
			\end{align}
			Therefore, one can certainly not prescribe eigenvalues with a multiplicity that violates this constraint. 
			However, one can  ask, if there exists a metric such that \cref{EqAtiyahBounds} is an equality (metrics with this property are called \emph{$D$-minimal}). In \cite{AmmannSurgery}, it is shown that the answer is affirmative on any closed connected spin manifold. Metrics satisfying this equality are even generic, i.e. they are dense in the $\mathcal{C}^{\infty}$-topology and open in $\mathcal{C}^1$-topology. \\
			On the other hand, the \emph{large kernel conjecture} states that for any $k \in \N$, one can find a metric $g_k$ such that $\dim \ker \Dirac^{g_k}_{\C} \geq k$ (if $\dim M \geq 3$), see also \cite{BaerMetrHarmonSpin}. In general, this is also an open problem. However, there are some cases for which it is solved:
			\begin{enumerate}
				\item
					for $M=\S^3$, see \cite{hitchin},
				\item
					for $m \equiv 3,7 \mod 8$ and $k=1$, see \cite{BaerMetrHarmonSpin},
				\item
					for $M=S^{2m}$, $m \geq 2$ and $k=1$, see \cite{seeger},
				\item
					for $m \geq 5$ and $k=1$, see \cite{dahlmetrinv}.					
			\end{enumerate}
			 
			All these considerations are the the reason, why the eigenvalue $0$ is excluded in the statement of \cref{ThmDahl}. 
	\end{description}
\end{Rem}

\begin{Rem}[prescribing the entire spectrum]
	Of course one could also ask, if one can prescribe the entire spectrum. We do not investigate this problem in the present thesis, since it requires very different techniques. It is clear that in addition to the constraints mentioned in \cref{RemUnsolvableCases}, there are even more restrictions. For instance, any Dirac spectrum has to satisfy the Weyl asymptotics. 
\end{Rem}

%% file: review.spinmorphisms.tex
\section{Spin Morphisms}
In the presence of a spin structure one can discuss whether or not certain morphisms lift to the spin structure. We will introduce the notion of a \emph{spin morphism} on various levels and show that a spin structure also induces a covering on the spin morphisms, see \cref{ThmSpinGaugeTrafosCov}. This part of the introduction will be carried out in detail, since this particular point of view is not so standard and we are unable to give a reference for \cref{ThmSpinGaugeTrafosCov}. We will require some basic notions from the theory of principal fibre bundles, see \cref{SectPrinGBdles,SubSectAppPrinGBdles} for notation.

\subsection{Spin Gauge Transformations}
Let $P \to M$ be a principal $G$-bundle. Recall that any gauge transformation $F \in \G(P)$ can also be described as an equivariant function $P \to G$, see \cref{LemGaugeToFunctions}. We denote by $\mathcal{C}^{\infty}_e(P,G)$ the space of all equivariant functions, see \cref{DefEquivariantGaugeFunctions}. In particular, this holds in the following situation: Let $M$ be a connected manifold and $E \to M$ be a spin vector bundle of rank $n$ with topological spin structure $\Theta:\GLtp E \to \GLp E$. Then we can apply this to the two gauge groups $\G(\GLtp E)$ and $\G(\GLp E)$.  Our aim is to show that the spin structure induces a covering between the gauge groups as well, see \cref{ThmSpinGaugeTrafosCov}. We will use the \emph{smooth compact-open-topology}, sometimes also called the \emph{weak topology} $\mathcal{C}^{\infty}_w$ (see \cref{SubSectWeakTop}) to topologize the gauge groups.

\begin{DefI}[equivariant spin functions]
	\label{DefSpinEquivFuncts}
	\nomenclature[CspineGLpeGLpn]{$\mathcal{C}^{\spin}_e(\GLp E, \GLp_n)$}{equivariant spin functions}
	Let $E \to M$ be a spin vector bundle of rank $n$ over a connected manifold $M$. A function $\sigma \in \mathcal{C}^{\infty}_e(\GLp E, \GLp_n)$ is \emph{spin}, if there is a function $\tilde \sigma \in \mathcal{C}^{\infty}_e(\GLtp E, \GLtp_n)$ such that
	\begin{align*}
		\xymatrix{
			\GLtp E
				\ar@{..>}[r]^-{\tilde \sigma}
				\ar[d]^-{\Theta}
			&\GLtp_n
				\ar[d]^-{\vartheta}
			\\
			\GLp E
				\ar[r]^-{\sigma}
			&\GLp_n
		}
	\end{align*}
	commutes. We denote by $\mathcal{C}^{\spin}_e(\GLp E, \GLp_n)$ the set of all these functions $\sigma$.
\end{DefI}

\begin{Lem} 
	\label{LemSpinGaugeFunctionsCov}
	In the situation of \cref{DefSpinEquivFuncts}, the map
		\DefMap{ \Pi: \mathcal{C}^{\infty}_e(\GLtp E, \GLtp_n)}{\mathcal{C}^{\spin}_e(\GLp E, \GLp_n)}{\tilde \sigma}{\vartheta \circ \tilde \sigma \circ \Theta^{-1}}
	is a well-defined group homomorphism with kernel $\Z_2$. If both groups are equipped with the $\mathcal{C}^{\infty}_w$-topology, this map is a $2:1$ covering.
\end{Lem}

\begin{Prf} $ $
	\begin{steplist}
		\step[well-defined]
			The definition of the map is to be understood as follows: Let $b \in \GLp E$. Then there exists $\tilde b \in \GLtp E$ such that $\Theta^{-1}(b) = \{ \pm \tilde b \}$. Here $-\tilde b = \tilde b.(-1)$, where $(-1) \in \GLtp_n$ is the non-trivial element in the fibre over the identity matrix. Since $\tilde \sigma$ is equivariant,			
			\begin{align*}
				\tilde \sigma(\pm \tilde b)
				= (\pm 1)^{-1} \tilde \sigma(\tilde b) (\pm 1)
				= \tilde \sigma( \tilde b),
			\end{align*}
			thus $\Pi$ is well-defined. It is clear that $\Pi$ is a homomorphism. 
		\step[image]
			We verify that $\sigma := \Pi(\tilde \sigma)$ is equivariant: Choose any $A \in \GLp_n$, $\tilde A \in \GLtp_n$, $\tilde b \in \GLtp E$, $b \in \GLp E$ such that 
			\begin{align*}
				\vartheta(\tilde A) = A, &&
				\Theta(\tilde b)=b, &&
				\Longrightarrow \Theta(\tilde b.\tilde A) = b.A.
			\end{align*}
			 We calculate
			\begin{align*}
				\sigma(b.A)
				= \vartheta(\tilde \sigma(\tilde b.\tilde A))
				= \vartheta(\tilde A^{-1} \tilde \sigma(\tilde b) \tilde A)
				= \vartheta(\tilde A^{-1}) \vartheta(\tilde \sigma(\tilde b)) \vartheta(\tilde A)
				= A^{-1} \sigma(b) A.
			\end{align*}
		\step[kernel]
			Furthermore, if
			\begin{align*}
				\Pi(\tilde \sigma_1) = \Pi(\tilde \sigma_2)
				\Longrightarrow  \vartheta \circ \tilde \sigma_1 \circ \Theta^{-1} = \vartheta \circ \tilde \sigma_2 \circ \Theta^{-1} 
				\Longrightarrow \vartheta \circ \tilde \sigma_1 = \vartheta \circ \tilde \sigma_2,
			\end{align*}
			then $\tilde \sigma_1$ and $\tilde \sigma_2$ are both lifts of the same map. Consequently, $\tilde \sigma_1 = \pm \tilde \sigma_2$.	
		\step[covering] 
			Let $\sigma \in \mathcal{C}^{\spin}_e(\GLp E, \GLp_n)$ be arbitrary. We have to construct a neighborhood of $\sigma$ evenly covered by $\Pi$. This will follow from the fact that $\vartheta$ is a covering.
			\begin{steplist}
				\step[construction of the neighborhoods] 
					Let $V \subset \GLp_n$ be an open connected neighborhood that is evenly covered by $\vartheta$, i.e. there exist open and disjoint subsets $\tilde V_{\pm} \subset \GLtp_n$ such that $\vartheta_{\pm} := \vartheta|_{\tilde V_{\pm}}:\tilde V_{\pm} \to V$ is a diffeomorphism. Since $\sigma$ is continuous, there exists an open connected coordinate domain $U \subset \sigma^{-1}(V) \subset \GLp E$ that is evenly covered by $\Theta$, i.e. there exist open and disjoint subsets $\tilde U_{\pm} \subset \GLtp E$ such that $\Theta_{\pm} := \Theta|_{\tilde U_{\pm}}$ is a diffeomorphism. Let $\varphi$ denote a manifold chart on $U$ and choose any compact $K \subset U$. We obtain that for any $k \in \N$,
					\begin{align*}
						B := \{ \sigma' \in \mathcal{C}^{\spin}_e(\GLp E, \GLp_n) \mid \sigma'(K) \subset V, \| \sigma' \circ \varphi^{-1}  - \sigma \circ \varphi^{-1} \|_{\mathcal{C}^k} < 1 \}
					\end{align*}
					is an open neighborhood of $\sigma$. We want to show that $B$ is evenly covered. Since $V$ and $U$ are open and connected, $\tilde U_{\pm}$ and $\tilde V_{\pm}$ are open and connected as well. Let $\tilde K_{\pm} := \Theta_{\pm}^{-1}(K) \subset \tilde U_{\pm}$. Since $\Theta_{\pm}$ is a diffeomorphism, $\tilde K_{\pm}$ is compact. Any lift $\tilde \sigma$ of $\sigma$ satisfies
					\begin{align*}
						(\vartheta \circ \tilde \sigma)(\tilde K_{\pm}) = (\sigma \circ \Theta)(\tilde K_{\pm}) = \sigma(K) \subset V.
					\end{align*}
					Since $\tilde U_{\pm}$ and $\tilde V_{\pm}$ are connected, it follows that the two lifts $\tilde \sigma_{\pm}$ of $\sigma$ are characterized by $\tilde \sigma_{\pm}(\tilde K_{\pm}) \subset \tilde V_{\mp}$.					
					We define the open neighborhood
					\begin{align*}
						\tilde B^{+} := \{ &\tilde \sigma' \in \mathcal{C}^{\spin}_e(\GLtp E, \GLtp_n) \mid \tilde \sigma'(\tilde K_{\pm}) \subset \tilde V^{\pm},  \\
						&\| \vartheta_{\pm} \circ \tilde \sigma' \circ \Theta_{\pm}^{-1} \circ \varphi^{-1} - \vartheta_{\pm} \circ \tilde \sigma_+ \circ \Theta_{\pm}^{-1} \circ \varphi^{-1} \|_{\mathcal{C}^k} < 1\}
					\end{align*}
					of $\tilde \sigma^+$ and the open neighborhood
					\begin{align*}
						\tilde B^{-} := \{ &\tilde \sigma' \in \mathcal{C}^{\spin}_e(\GLtp E, \GLtp_n) \mid \tilde \sigma'(\tilde K_{\pm}) \subset \tilde V^{\mp},  \\
						&\| \vartheta_{\mp} \circ \tilde \sigma' \circ \Theta_{\pm}^{-1} \circ \varphi^{-1} - \vartheta_{\mp} \circ \tilde \sigma_- \circ \Theta_{\pm}^{-1} \circ \varphi^{-1} \|_{\mathcal{C}^k} < 1\}
					\end{align*}			
					$\tilde \sigma^{-}$. We claim that $B$ is evenly covered by $\tilde B^{\pm}$, in particular that 
					\begin{align*}
						\Pi^{-1}(B) = \tilde B^+ \dot \cup \tilde B^{-}.
					\end{align*}
					Since $\tilde V^+$ and $\tilde V^-$ are disjoint, $\tilde B^+$ and $\tilde B^-$ are disjoint by construction.
				\step[$\Pi^{-1}(B) \subset \tilde B^+ \dot \cup \tilde B^{-}$]
					If $\tilde \sigma' \in \Pi^{-1}(B)$, then there exists $\sigma' \in B$ such that $\sigma' = \vartheta \circ \tilde \sigma' \circ \Theta^{-1}$. Again, $(\vartheta \circ \tilde \sigma')(\tilde K_{\pm}) \subset V$, thus either $\tilde \sigma'(\tilde K_{\pm}) \subset \tilde V^{\pm}$ or $\tilde \sigma'(\tilde K_{\pm}) \subset \tilde V^{\mp}$. In the former case, we calculate
					\begin{align*}
						\vartheta_{\pm} \circ \tilde \sigma' \circ \Theta_{\pm}^{-1} \circ \varphi^{-1}
						=\sigma' \circ \Theta \circ \Theta_{\pm}^{-1} \circ \varphi^{-1} 
						=\sigma' \circ \varphi^{-1}
					\end{align*}
					on $\tilde K_{\pm}$. The analogous calculation holds for $\sigma$ as well, thus
					\begin{align}
						\label{EqGaugeFunctNorm}
						\| \vartheta_{\pm} \circ \tilde \sigma' \circ \Theta_{\pm}^{-1} \circ \varphi^{-1} - \vartheta_{\pm} \circ \tilde \sigma_+ \circ \Theta_{\pm}^{-1} \circ \varphi^{-1} \|_{\mathcal{C}^k}
						=\| \sigma' \circ \varphi^{-1}  - \sigma \circ \varphi^{-1} \|_{\mathcal{C}^k}
					\end{align}
					and $\tilde \sigma' \in \tilde B^+$. In case $\tilde \sigma'(\tilde K_{\pm)} \subset \tilde V_{\mp}$, we obtain $\tilde \sigma' \in \tilde B^-$ in the same manner. 
				\step[$\tilde B^+ \dot \cup \tilde B^{-} \subset \Pi^{-1}(B)$]
					Conversely, if $\tilde \sigma' \in \tilde B^+$, then we define $\sigma' := \Pi(\tilde \sigma')$ and obtain
					\begin{align*}
						\sigma'(K)
						= (\vartheta \circ \tilde \sigma' \circ \Theta^{-1})(K)
						= (\vartheta \circ \tilde \sigma')(\tilde K_+ \cup \tilde K_-)
						\subset \vartheta(\tilde V_{\pm})
						=V
					\end{align*}
					and again by \cref{EqGaugeFunctNorm}, we obtain $\tilde \sigma' \in B$. If $\tilde \sigma' \in \tilde B^-$, we obtain $\Pi(\tilde \sigma') \in B$ in the same manner.
				\step[homeomorphism]
					By definition $\Pi|_{\tilde B^{\pm}}$ is given as pre- respectively post composition with $\Theta_{\pm}^{-1}$ respectively $\vartheta_{\pm}$. Since these maps are diffeomorphisms on their domains and composition is continuous by \cref{ThmCOkTop}, the result follows.
			\end{steplist}
	\end{steplist}
\end{Prf}

\begin{Def}[spin gauge trafo]
	\nomenclature[GspinGLpE]{$\G^{\spin}(\GLp E)$}{spin gauge transformations}
	\label{DefSpinGaugeTrafo}
	A gauge transformation $F \in \G(\GLp E)$ is \emph{spin}, if there exists a gauge transformation $\tilde F \in \G(\GLtp E)$ such that
	\begin{align*}
		\xymatrix{
			\GLtp E
				\ar[r]^-{\tilde F}
				\ar[d]^{\Theta}
			& \GLtp E
				\ar[d]^-{\Theta}
			\\
			\GLp E
				\ar[r]^-{F}
			& \GLp E
		}
	\end{align*}
	commutes. We denote by $\G^{\spin}(\GLp E)$ the set of all those gauge transformations.
\end{Def}

\begin{Thm}
	\label{ThmSpinGaugeTrafosCov}
	\nomenclature[ThetaPush]{$\Theta_*$}{covering $\G(\GLtp E) \to \G^{\spin}(\GLp E)$}
	The map
		\DefMap{\Theta_*: \G(\GLtp E)}{\G^{\spin}(\GLp E)}{\tilde F}{\Theta \circ \tilde F \circ \Theta^{-1}}
	is a well-defined group homomorphism with kernel $\Z_2$. If both groups are equipped with the $\mathcal{C}^{\infty}_w$-topology, this map is a $2:1$ covering.
\end{Thm}

\begin{Prf}
	The proof that $\Theta_*$ is a well-defined group homomorphism with kernel $\Z_2$ is analogous the proof of these claims in \cref{LemSpinGaugeFunctionsCov}.
	By construction, the diagram
	\begin{align*}
		\xymatrix{
			\G(\GLtp E)
				\ar[r]
				\ar[d]^-{\Theta_*}
			&\mathcal{C}^{\infty}_e(\GLtp E, \GLtp_n)
				\ar[d]^-{\Pi}
			\\
			\G^{\spin}(\GLp E)
				\ar[r]
			&\mathcal{C}^{\spin}_e(\GLp E, \GLp_n)
		}
	\end{align*}
	commutes. The horizontal maps are homeomorphisms by \cref{CorWeakTopGBundle}\ref{ItWeakTopSigmaHoem} and the right map is a $2:1$ covering by \cref{LemSpinGaugeFunctionsCov}.	
\end{Prf}

\subsection{Spin Morphisms of Vector Bundles}
\label{SubSectSpinMorphVB}

\begin{Def}[spin morphism]
	\index{spin morphism}
	\label{DefSpinMorphism}
	Let $\pi_j:E_j \to M_j$, $j=1,2$, be spin vector bundles of rank $n$ with spin structure $\Theta_j:\GLtp E_j \to \GLp E_j$. An orientation-preserving isomorphism of vector bundles $(f,F)$ is \emph{spin}, if there exists a morphism of principal $\GLtp_n$-bundles $\tilde F$ such that
	\begin{align}
		\begin{split}
			\xymatrixcolsep{3.5em}
			\xymatrix{
				\GLtp E_1
					\ar@{..>}[r]^{\tilde F}
					\ar[d]^{\Theta_1}
				&\GLtp E_2
					\ar[d]^-{\Theta_2}
				\\
				\GLp E_1
					\ar[r]^{F}
					\ar[d]
				& \GLp E_2
					\ar[d]
				\\
				M_1
					\ar[r]^{f}
				& M_2
			}
		\end{split}
	\end{align}	
	commutes. Here, $f$ is a diffeomorphism of smooth manifolds and $F$ is identified with the corresponding $\GLp_n$-gauge transformation
		\DefMap{\GLp E_1}{\GLp E_2}{(b_1, \ldots, b_n)}{(Fb_1, \ldots, Fb_n).}
\end{Def}

\begin{Rem}
	\label{RemSpinGaugeSpec}
	This definition is particularly important, if $E_1 = E_2$ and $\Theta_1 = \Theta_2$. In that case $(f,F)$ is spin if and only if $F \in \G^{\spin}(\GLp E)$ in the sense of \cref{DefSpinGaugeTrafo}.
\end{Rem}

\begin{Rem}
	If an isomorphism $(f,F)$ is spin and $M$ is connected, there are precisely two two spin lifts $\tilde F^\pm$ of $F$, which are related by $F^+.(-1)=F^-$, where $-1 \in \GLtp_n$.
\end{Rem}

\begin{Def}[equivalence of spin structures]
	\index{equivalence!of spin structures}
	\label{DefEquivalenceSpinStructs}
	Let $E \to M$ be an oriented vector bundle. Two spin structures $\Theta_j:P_j \to \GLp E$ are \emph{equivalent}, if there exists an isomorphism $F:P_1 \to P_2$ of principal bundles such that $(\id_M, F)$ is a spin morphism.
\end{Def}

\subsection{Spin Diffeomorphisms}
Specializing \cref{DefSpinMorphism} to the tangent bundle of a spin manifold gives rise to a more restrictive notion. 

\begin{DefI}[spin diffeomorphism]
	\label{DefSpinDiffeo}
	A diffeomorphism $f \in \Diff^+(M)$ is a \emph{spin diffeomorphism}, if $(f, f_*)$ is a spin morphism, i.e. if there exists a morphism $F$ of princial $\GLtp_m$-bundles such that
	\begin{align}
		\label{EqDefSpinDiffeo}
		\begin{split}
			\xymatrix{
				\GLtp M \ar[d]^-{\Theta}
					\ar@{..>}[r]^-{F}
				& \GLtp M
					\ar[d]^-{\Theta}
				\\
				\GLp M
					\ar[d]
					\ar[r]^-{f_*}
				& \GLp M
					\ar[d] 
				\\
				M 
					\ar[r]^-{f}
				& M
			}
		\end{split}
	\end{align}
	commutes. Here, $f_*$ is identified with
		\DefMap{f_*:\GLp M}{\GLp M}{(b_1, \ldots, b_n)}{(f_* b_1, \ldots, f_*b_n).}
	We say $F$ is a \emph{spin lift} of $f$, and define
	\begin{align*}
		\Diff^{\spin}(M)& := \{ f \in \Diff^+(M) \mid \text{$f$ is a spin diffeomorphism} \}.
	\end{align*}
\end{DefI}

\begin{Rem}[spin isometries]
	\label{RemSpinIsometries}
	If $f \in \Diff^{\spin}(M)$ and $h \in \Rm(M)$, we can set $g := f^* h$. In this case \eqref{EqDefSpinDiffeo} restricts to the analogous diagram
	\begin{align}	
		\label{EqDefSpinIsometry}
		\begin{split}
			\xymatrix{
				\Spin^g M \ar[d]^-{\Theta^g}
					\ar[r]^-{F}
				& \Spin^h M
					\ar[d]^-{\Theta^h}
				\\
				\SO^g M
					\ar[d]
					\ar[r]^-{f_*}
				& \SO^h M
					\ar[d] 
				\\
				(M,g) 
					\ar[r]^-{f}
				& (M,h)
			}
		\end{split}
	\end{align}
	of metric spin structures. We say that $f$ is a \emphi{spin isometry} in this case. Notice that this implies that $F$ induces an isometry of vector bundles
		\DefMap{\bar F:\Sigma^{g}_{\K} M }{\Sigma^h_{\K} M}{{[\tilde b, v]}}{{[F(\tilde b), v]}.}
	This isometry commutes with the Clifford multiplication and the spin connection. Therefore, $(M,g)$ and $(M,h)$ are \emphi{Dirac-isospectral}, i.e. their sets of Dirac eigenvalues, as well as their multiplicities, are equal. 
\end{Rem}

\subsection{Isotopies of Spin Diffeomorphisms}
We will need some facts about paths of diffeomorphisms and their lifts to the spin structure. Again, we use the $\mathcal{C}^{\infty}_w$-topology (see \cref{SubSectWeakTop}) on $\Diff(M,N)$.

\begin{DefI}[isotopy]
	Let $M$ and $N$ be smooth manifolds. An isotopy is a continuous path $h:I \to \Diff(M,N)$. Two diffeomorphisms $f,g \in \Diff(M,N)$ are \emph{isotopic}, if there exists an isotopy $h$ such that $h_0 = f$ and $h_1 = g$. 
\end{DefI}

\begin{Rem}
	By \cref{ThmCOkTop}, an isotopy can equivalently be defined as a continuous map $h:I \times M \to N$ such that for all $t \in I$, the induced map $h_t := h(t, \_): M \to N$ is a diffeomorphism. We do not distinguish between $h$ as a map $I \to \Diff(M,N)$ and as a map $I \times M \to N$.
\end{Rem}

\begin{Rem}
	\label{RemIsotopyOrientation}
	Let $h:I \times M \to N$ be an isotopy. Since $\det$ is continuous, either all diffeomorphisms $h_t$ are orientation-preserving or no $h_t$ is.
\end{Rem}

This phenomenon generalizes to the spin structure.

\begin{Thm}[spin isotopy invariance]
	\label{ThmSpinIsotopyInvariance}
	Let $(M,\Theta)$, $(N, \Xi)$ be spin manifolds and $h:I \times M \to N$ be an isotopy. If $h_0$ is spin in the sense of \cref{DefSpinDiffeo}, there exists a lift $\hat H$ such that
		\begin{align} 
			\label{EqSpinLiftfIsotopy}
			\begin{split}
				\xymatrix{
					I \times \GLtp M
						\ar[d]^-{\id \times \Theta}
						\ar@{..>}[r]^-{\exists \hat H} 
					& \GLtp N \ar[d]^-{\Xi}
					\\
					I \times \GLp M
						\ar[d]
						\ar[r]^-{H}
					& \GLp N
						\ar[d]
					\\
					I \times M
						\ar[r]^-{h} 
					& N.
				}
			\end{split}
		\end{align}	
	commutes. In particular, $h_t$ is spin for all $t \in I$.
\end{Thm}

\begin{Prf}
	By hypothesis, $h_0$ is spin. In particular $h_0 \in \Diff^+(M,N)$, thus $h_t \in \Diff^+(M,N)$ for all $t \in I$. Consequently, the map $H$ in \cref{EqSpinLiftfIsotopy} can be defined by 
	\begin{align*}
		\forall b \in \GLp M: \forall t \in I: H(t,b) := {h_t}_*(b).
	\end{align*}
	To show the existence of $\hat H$, we note that, since $H$ is an isotopy, it is in particular a homotopy. Consequently, $H \circ (\id \times \Theta)$ is also a homotopy between $H_0 \circ \Theta$ and $H_1 \circ \Theta $. By hypothesis, there exists a spin lift $\hat H_0$ of $h_0$. By definition, this map satisfies $\Theta \circ \hat H_0 = \Theta \circ H_0$. Since covering spaces have the homotopy lifting property (see \cref{ThmHomotopyLifting}), there exists a unique $\hat H$ such that $\Xi \circ \hat H = H \circ (\id \times \Theta)$.
\end{Prf}

\begin{Rem}
	\label{RemIsotopySpin}
	Let $h:I \times M \to N$ be an isotopy. Analogous to \cref{RemIsotopyOrientation} we can rephrase \cref{ThmSpinIsotopyInvariance} by saying that either all diffeomorphisms $h_t$ are spin or no $h_t$ is. Therefore, we will say that an isotopy $h$ \emph{is spin} if one, hence all, diffeomorphisms $h_t$ are spin.
\end{Rem}

\begin{Rem}
	\label{RemSpinLiftIsotopies}
	Recall from \cref{ThmSpinGaugeTrafosCov} that the map
		\DefMap{\Theta_*: \G(\GLtp M)}{\G^{\spin}(\GLp M)}{\hat F}{\Theta \circ \hat F \circ \Theta^{-1}}
	is a $2:1$ covering map. Notice that a diffeomorphism $f \in \Diff^+(M)$ is spin if and only if $f_* \in \G^{\spin}(\GLp M)$, see \cref{RemSpinGaugeSpec}. Let $h:I \times M \to M$ be a spin isotopy. We obtain a diagram
	\begin{align}
		\label{EqSpinIsotopyLift}
		\begin{split}
			\xymatrix{
				& \G(\GLtp M)
					\ar[d]^-{\Theta_*}_{2:1}
				\\
				I
					\ar[r]_-{h}
					\ar@{..>}[ur]^-{\hat H}
				& \G^{\spin}(\GLp M).
			}
		\end{split}
	\end{align}
\end{Rem}

%% file: review.category.tex
\section{Category Theoretic Reformulation}
In this section, we will reformulate the classical constructions in spin geometry in a modern category theoretic language. We do not strictly need this formulation later, but it gives a clear and concise overview of the structure of spin geometry and might serve as an alternative introduction to readers from other branches of mathematics like algebraic topology or algebraic geometry. We assume the reader to be familiar with the basics of category theory.

\begin{Def}[spin vector bundles]
	\nomenclature[SpinVectBdm]{$\cat{SpinVB}$}{category of spin vector bundles}
	The \emph{category of spin vector bundles}, $\cat{Spin\?VB}$, consists of
	\begin{category}
		\item[Objects]
			A \emph{spin vector bundle} is a tuple $(\pi:E \to M,\Theta)$, where $\pi:E \to M$ is a smooth oriented real vector bundle of some rank $n \in \N$ over a smooth closed manifold $M$ and $\Theta$ is a topological spin structure for $E$, see \cref{DefTopSpinStructure}.
			
		\item[Morphisms]
			A \emph{morphism $(\pi_1:E_1 \to M_1,\Theta_1) \to (\pi_2:E_2 \to M_2,\Theta_2)$ of spin vector bundles} is a triple $(f,F,\tilde F)$ such that
			\begin{align}
				\label{EqDefSpinMorph}
				\begin{split}
					\xymatrix{
						\GLtp E_1
							\ar[r]^{\tilde F}
							\ar[d]^{\Theta_1}
						&\GLtp E_2
							\ar[d]^-{\Theta_2}
						\\
						\GLp E_1
							\ar[r]^{F}
							\ar[d]
						& \GLp E_2
							\ar[d]
						\\
						M_1
							\ar[r]^{f}
						& M_2
					}
				\end{split}
		\end{align}
		commutes. Here, $f$ is a diffeomorphism between smooth manifolds, $F$ is induced by an orientation-preserving isomorphism of smooth vector bundles and $\tilde F$ is an isomorhpism of principal bundles.
	\end{category}
\end{Def}

\begin{Rem}
	Not every diffeomorphism $f:M_1 \to M_2$ admits a lift as in \cref{EqDefSpinMorph}. In case $E_j = TM_j$, $j=1,2$, and a lift exists as in \cref{EqDefSpinDiffeo}, i.e. $F=f_*$, we call such a diffeomorphism a \emph{spin diffeomorphism}. If in addition $M$ is connected, there are precisely two such lifts, the other one given by $-\tilde F := \tilde F.(-1)$, $-1 \in \GLtp_n$.
\end{Rem}

\begin{Def}[Riemannian spin vector bundles]
	\nomenclature[RiemSpinVectBdm]{$\cat{RiemSpinVB}$}{category of Riemannian spin vector bundles}
	The \emph{category of Riemannian spin vector bundles}, $\cat{Riem\?Spin\?VB}$, consists of
	\begin{category}
		\item[Objects]
			A \emph{Riemannian spin vector bundle} is a tuple $(\pi:E \to M,h, g, \nabla, \Theta^h)$, where $\pi:E \to M$ is a smooth oriented vector bundle of some rank $n \in \N$ over a smooth closed manifold $M$, $g$ is a Riemannian metric on $M$, $h$ is a Riemannian metric on $E$, $\nabla$ is a metric connection on $E$ and $\Theta^h$ is a metric spin structure for $E$, see \cref{DefMetricSpinStructure}.
			
		\item[Morphisms]
			A \emph{morphism}
			\begin{align*}
				(\pi_1: E_1 \to M_1, h_1, g_1, \nabla^1, \Theta^{h_1}) \to (\pi_2:E_2 \to M_2, h_2, g_2, \nabla^2, \Theta^{h_2})
			\end{align*}
			\emph{of Riemannian spin vector bundles} is a triple $(f,F,\tilde F)$ such that
			\begin{align*}
			\xymatrix{
				\Spin^{h_1} E_1
					\ar[r]^{\tilde F}
					\ar[d]^{\Theta_1}
				&\Spin^{h_2} E_2
					\ar[d]^-{\Theta_2}
				\\
				\SO^{h_1} E_1
					\ar[r]^{F}
					\ar[d]
				& \SO^{h_2} E_2
					\ar[d]
				\\
				M_1
					\ar[r]^{f}
				& M_2
			}
		\end{align*}
		commutes. Here, $f$ is an isometry of Riemannian manifolds, $F$ is induced by an isometry of vector bundles, which is compatible with the connections, and $\tilde F$ is an isomorphism of principal bundles.
	\end{category}
\end{Def}

\begin{Def}[spin manifolds]
	\nomenclature[SpinMfd]{$\cat{SpinMfd}$}{category of spin manifolds}
	The \emph{category of spin manifolds}, $\cat{Spin\?Mfd}$, consists of
	\begin{category}
		\item[Objects]
			A \emph{spin manifold} is a tuple $(M,\Theta)$, where $M$ is a smooth closed oriented manifold and $\Theta$ is a topological spin structure for $TM$.
		\item[Morphisms]
			A \emph{morphism of spin manifolds} $(M_1,\Theta_1) \to (M_2,\Theta_2)$ is a tuple $(f,\tilde F)$ such that $(f,f_*,F)$ is a morphism  $(TM_1 \to M_1,\Theta_1) \to (TM_2 \to M_2, \Theta_2)$ of spin vector bundles.
	\end{category}
\end{Def}

\begin{Def}[Riemannian spin manifolds]
	\nomenclature[RiemSpinMfdm]{$\cat{RiemSpinMfd}$}{category of Riemannian spin manifolds}
	The \emph{category of Riemannian spin manifolds}, $\cat{Riem\?Spin\?Mfd}$, consists of
	\begin{category}
		\item[Objects]
			A \emph{Riemannian spin manifold} is a tuple $(M,g,\Theta^g)$ such that $(M,g)$ is a smooth oriented closed Riemannian manifold and $\Theta^g$ is a metric spin structure on $(TM, g)$.
		\item[Morphisms] 
			A \emph{morphism of Riemannian spin manifolds} $(M_1,g_1,\Theta^{g_1}) \to (M_2, g_2, \Theta^{g_2})$ is a tuple $(f,F)$ such that $f:(M_1,g_1) \to (M_2, g_2)$ is an orientation-preserving isometry of Riemannian manifolds and $(f,f_*,F)$ is a morphism of spin manifolds.
	\end{category}
\end{Def}

\begin{Rem}
	These four categories are related by the following commutative diagram of categories
	\begin{align*}
		\xymatrix{
			\cat{RiemSpinMfd}
				\ar[r]^-{\fun{T}}
				\ar[d]^-{\fun{MetrTop}}
			&\cat{RiemSpinVB}
				\ar[d]^-{\fun{MetrTop}}
			\\
			\cat{SpinMfd}
				\ar[r]^-{\fun{T}}
			&\cat{SpinVB}
		}
	\end{align*}
	Here, $\fun{T}$ is the \emphi{tangent functor}, which maps a Riemannian manifold $(M,g)$ to its tangent bundle $TM$ with its Levi-Civita connection and a map $f$ to its derivative $f_*$. $\fun{MetrTop}$ is the functor that forgets the Riemannian metric and replaces the metric spin structure by a topological spin structure, see \cref{ThmSpinStructTopToMetric}. We see that spin manifolds are only a special case of spin vector bundles. 
\end{Rem}

\begin{Def}[Clifford bundles]
	\label{DefCliffordBundles}
	\nomenclature[CliffBd]{$\cat{CliffB}_{\K}$}{category of Clifford Bundles}
	The \emph{category of $\K$-Clifford Bundles}, $\cat{Cliff\?B_{\K}}$, consists of
	\begin{category}
	 \item[Objects] 
		A \emph{Clifford Bundle} is a tuple $(\pi:S \to M, h, g, \mathfrak{m}, \nabla^S)$, such that
		\begin{enumerate}
			\item 
				$(M,g)$ is a closed Riemannian manifold.
			\item
				$\pi:S \to M$ is a $\K$-vector bundle of some rank $n$ with a Riemannian respectively Hermitian metric fibre metric $h$.
			\item
				$\mathfrak{m}$ is a morphism of $\R$-vector bundles, 
					\DefMap{\mathfrak{m}: TM \otimes S}{S}{X \otimes s}{\mathfrak{m}(X \otimes s) =: X \cdot s,}
				called \emph{Clifford multiplication}, which satisfies
				\begin{align*}
					\forall X \in TM: \forall \psi \in \Gamma(S): X \cdot X \cdot \psi = - g(X,X) \psi .
				\end{align*}
			\item
				The Clifford multiplication is skew-adjoint with respect to $h$, i.e.
				\begin{align*}
					\forall X \in \mathcal{T}(M): \forall \psi_1, \psi_2 \in \Gamma(S): h(X \cdot \psi_1, \psi_2) = - h(\psi_1, X \cdot \psi_2) .
				\end{align*}
			\item
				$\nabla^S:\Gamma(S) \to \Gamma(T^*M \otimes S)$ is a connection on $S$ that is \emph{compatible} with the Levi-Civita connection $\nabla^g$ on $TM$ and the Clifford multiplication $\mathfrak{m}$, i.e.
				\begin{align*}
					\forall X, Y \in \mathcal{T}(M): \forall \psi \in \Gamma(S): \nabla^S_X(Y \cdot \psi) = \nabla^g_X Y \cdot \psi + Y \cdot \nabla^S_X \psi.
				\end{align*}
		\end{enumerate}
		
		\item[Morphisms]
			A \emph{morphism}
			\begin{align*}
				(\pi_1:S_1 \to M_1, h_1, g_1, \mathfrak{m}_1, \nabla^{S_1}) \to (\pi_2:S_2 \to M_2, h_2, g_2, \mathfrak{m}_2, \nabla^{S_2}), 
			\end{align*}
			\emph{of $\K$-Clifford bundles} is a tuple $(f,F)$, such that
			\begin{align*}
				\xymatrixcolsep{3.5em}
				\xymatrix{
					S_1
						\ar[r]^{F}
						\ar[d]^{\pi_1}
					& S_2
						\ar[d]^-{\pi_2}
					\\
					M_1
						\ar[r]^-{f}
					& M_2
				}
			\end{align*}
			commutes. Here, $f:(M_1,g_1) \to (M_2,g_2)$ is an isometry of Riemannian manifolds and $F:(S_1,h_1) \to (S_2,h_2)$ is an isometric isomorphism of $\K$-vector bundles, which satisfies $F^* \nabla^{S_1} = \nabla^{S_2}$ and $F \circ \mathfrak{m}_1 = \mathfrak{m}_2 \circ f_* \otimes F$, i.e. the diagrams
			\begin{align*}
				\xymatrixcolsep{3.5em}
				\xymatrix{
					TM_1 \otimes S_1
						\ar[r]^-{\mathfrak{m}_1}
						\ar[d]^{f_* \otimes F}
					& S_1
						\ar[d]^{F}
					& \hspace{-2em}
					& \Gamma(S_1)
						\ar[r]^-{\nabla^{S_1}}
						\ar[d]^-{\Gamma(F)}
					& \Gamma(T^* M_1 \otimes S_1)
						\ar[d]^-{\Gamma((f^{-1})^* \otimes F)}
					\\
					TM_2 \otimes S_2
						\ar[r]^-{\mathfrak{m}_2}
					& S_2
					& \hspace{-2em}
					&\Gamma(S_2)
						\ar[r]^-{\nabla^{S_2}}
					& \Gamma(T^* M_2 \otimes S_2)					
				}
			\end{align*}
			commute.	 
	\end{category}
\end{Def}

\begin{Def}[differential vector bundles]
	\nomenclature[VectDiffn]{$\cat{DiffVB}$}{category of differential vector bundles}
	The \emph{category of differential $\K$-vector bundles}, $\cat{Diff\?VB_{\K}}$, consists of
	\begin{category}
	 \item[Objects] 
		A \emph{differential vector bundle} is a tuple $(\pi: E \to M, h,g, D)$, where $\pi:(E,h) \to (M,g)$ is a $\K$-vector bundle with Riemannian respectively Hermitian fibre metric $h$ of some rank $n \in \N$ over a closed Riemannian manifold $(M,g)$ and
		\begin{align*}
			D:\Gamma(E) \to \Gamma(E)
		\end{align*}
		is a self-adjoint first order elliptic differential operator. We also regard $D$ as an unbounded operator $L^2(E) \to L^2(E)$ with dense domain $H^1(E)$, where $H^1$ is the first order Sobolev space.
	\item[Morhpisms]
		A \emph{morphism}
		\begin{align*}
			(\pi_1: E_1 \to M_1, h_1,g_1, D_1) \to (\pi_2: E_2 \to M_2, h_2,g_2, D_2)
		\end{align*}
		\emph{of differential $\K$-vector bundles} is a tuple $(f,F)$ such that $f:(M_1,g_1) \to (M_2,g_2)$ is an isometry of Riemannian manifolds and $F:(E_1,h_1) \to (E_2,h_2)$ is an isometric isomorphism of $\K$-vector bundles such that
		\begin{align*}
			\xymatrix{
				\Gamma(E_1)
					\ar[r]^-{D_1}
					\ar[d]^{\Gamma(F)}
				& \Gamma(E_1)
					\ar[d]^{\Gamma(F)}
				\\
				\Gamma(E_2)
					\ar[r]^{D_2}
				& \Gamma(E_2)
			}
		\end{align*}
		commutes.		
	\end{category}
\end{Def}

\begin{Def}[Dirac functor]
	\nomenclature[DiracFun]{$\fun{\Dirac}$}{Dirac functor}
	The \emph{Dirac functor} $\fun{\Dirac_{\K}}: \cat{CliffB_{\K}} \to \cat{DiffVB_{\K}}$ is defined
	\begin{category}
		\item[on objects by] 
		\begin{align*}
			\fun{\Dirac_{\K}}((\pi:S \to M, h, g, \mathfrak{m}, \nabla^S)) := (\pi:S \to M, h, g, \Dirac^g_{\K}),
		\end{align*}
		where the \emph{Dirac operator} $\Dirac^g_{\K}$ is defined to be the composition
		\begin{align*}
			\xymatrixcolsep{3.5em}
			\xymatrix{
				\Gamma(S)
					\ar[r]^-{\nabla^S}
				& \Gamma(T^*M \otimes S)
					\ar[r]^-{\sharp \otimes \id}
				& \Gamma(TM \otimes S)
					\ar[r]^-{\mathfrak{m}}
				& \Gamma(S).
			}
		\end{align*}
		\item[on morphisms by]
			\begin{align*}
				\fun{\Dirac}((f,F)) := (f,F). 
			\end{align*}			
	\end{category}
\end{Def}

\begin{Def}[spin geometry functor]
	\label{DefSpinGeometryFunctor}
	\nomenclature[SigmaFun]{$\fun{\Sigma}_{\K}$}{spin geometry functor}
	For any $n \in \N$, let $\rho_n:\Spin_n \to \GL(\Sigma_n)$ be a $\K$-spinor representation.
	The \emph{spin geometry functor} 
	\begin{align*}
		\fun{\Sigma_{\K}}: \cat{RiemSpinVB} \to \cat{CliffB_{\K}}
	\end{align*}
	is defined 
	\begin{category}
	 \item[on objects by] 
	 \begin{align*}
		\fun{\Sigma_{\K}}(\pi:E \to M,h, g, \Theta^h) := (\pi^{\Sigma}_g: \Sigma^g_{\K} E \to M, h, g, \mathfrak{m}, \nabla^S).
	 \end{align*}
	 Here, $\Sigma^g_{\K} E := \Spin E \times_{\rho_n} \Sigma_n$, where $n$ is the rank of $E$, is the $\K$-spinor bundle, see also \cref{DefSpinorBundle}. The Clifford multiplication $\mathfrak{m}$, the fibre metric $h$ and the spinorial Levi-Civita connection $\nabla^S$ are constructed as in \cref{ThmSpinorBundleAddStructure} (the theorem is formulated for spin manifolds, but the construction for spin vector bundles is similar).
	\item[on morphisms by]
		\begin{align*}
			\fun{\Sigma_{\K}}(f,F) := (f,\bar F),
		\end{align*}
		where
			\DefMap{\bar F: \Sigma^{g_1}_{\K} E_1}{\Sigma^{g_2}_{\K} E_2}{{[s,\sigma]}}{{[F(s),\sigma]}.}
	\end{category}
\end{Def}

\begin{Rem}
	\label{RemSpinReprChoice}
	The definition of the functor $\fun{\Sigma}_{\K}$ depends on the choice of the $\K$-spinor representations $\{\rho_n\}_{n \in \N}$. But once such a choice is fixed, this is a well-defined functor. 
\end{Rem}

\begin{Rem}
	The results of this section can be expressed by saying that there is a diagram
	\begin{align*}
		\xymatrixcolsep{3.5em}
		\xymatrix{
			\cat{RiemSpinVB}
				\ar[r]^-{\fun{\Sigma_{\K}}}
				\ar[d]^-{\fun{MetrTop}}
			&\cat{CliffB_{\K}}
				\ar[r]^-{\fun{\Dirac_{\K}}}
			&\cat{DiffVB_{\K}}
			\\
			\cat{SpinVB}
		}
	\end{align*}
	of categories and functors. Unfortunately, one cannot complete the lower row in a similar fashion, see \cref{RemNoUSBByRepresentations}.
\end{Rem}

%% file: evpaper.abstract.tex
\begin{subabstract}
	It is well known that on a bounded spectral interval, the Dirac spectrum can be described locally by a sequence of continuous functions of the Riemannian metric, which evaluated at any metric is non-decreasing, see \cref{ThmBaerSpecC1Conti}. In this chapter, we extend this result to a global version. We view the spectrum of a Dirac operator as a function $\Z \to \R$ and endow the space of all spectra with an $\arsinh$-uniform metric. We prove that the spectrum of the Dirac operator depends continuously on the Riemannian metric. As a corollary, we obtain the existence of a non-decreasing family of functions on the space of all Riemannian metrics, which represents the entire Dirac spectrum at any metric. We also show that, due to spectral flow, these functions do not descend to the space of Riemannian metrics modulo spin diffeomorphisms in general. The content of this chapter has been published by now in a similar version, see \cite{evpaper}.
\end{subabstract}

%% file: evpaper.intro.tex
\section{Introduction and Statement of the Results} 
\label{SctEvpaperIntro}

In this chapter, we prove that the Dirac spectrum depends continuously on the Riemannian metric. More precisely, we will prove:

\begin{recall}[\cref{MainThmFun}.]
	\input{mainthm.fun}
\end{recall}

Before we start, we give a motivation for this particular formalization. By \cref{ThmSpecPropsDirac}, $\spec \Dirac_{\C}^g$ is a subset of the real line, which is closed, discrete and unbounded from both sides. The elements of $\spec \Dirac_{\C}^g$ consist entirely of eigenvalues of finite multiplicity. However, treating $\spec \Dirac_{\C}$ as a function from the Riemannian metrics to the subsets of the real line is inconvenient, because it is unclear what the continuity assertion should mean and because $\spec \Dirac_{\C}^g$ as a set does not contain any information about the multiplicities of the eigenvalues. Both problems can be solved by describing the spectrum as a sequence, i.e. a function $\Z \to \R$. Intuitively, we would like to enumerate the eigenvalues from $- \infty$ to $+ \infty$ as a non-decreasing sequence indexed by $\Z$. The problem is that this is not well-defined, because it is unclear which eigenvalue should be the ``first'' one. Formally, we can avoid this problem via the following definition.

\begin{Def}
	\label{EqDefsg}
	\nomenclature[sg]{$\spc^g$}{ordered spectral function of $\Dirac_{\C}^g$}
	For any $g \in \Rm(M)$, let $\spc^g:\Z \to \R$ be the unique non-decreasing function such that $\spc^{g}(\Z)=\spec \Dirac_{\C}^{g}$, 
	\begin{align*}
		\forall \lambda \in \R: \dim_{\C} \ker(\Dirac_{\C}^g - \lambda) = \sharp (\spc^g)^{-1}(\lambda),
	\end{align*}
	 and $\spc^g(0)$ is the first eigenvalue $\geq 0$ of $\Dirac_{\C}^g$. 
\end{Def}

\begin{figure}[t] 
	\begin{center}
		\input{fig.specenum}
		\caption[Ordered spectral function]{Setting $\lambda_j:=\spc^{g}(j)$, we obtain this enumeration of $\spec \Dirac^g$. }
		\label{FigOrdSpecFuncDisc}
	\end{center}
\end{figure}

The situation is depicted in \cref{FigOrdSpecFuncDisc}. (We will give an analogous definition for slightly more general operators later, see \cref{DefOrderedSpecFunT}). 
Now, $\spc^g$ is well-defined, but as it will turn out, the requirement that $\spc^g(0)$ be the first eigenvalue $\geq 0$ has some drawbacks. Namely, the map $g \mapsto \spc^g(j)$, $j \in \Z$, will in general not be continuous, see \cref{RemIntuitionDisc}. To obtain a more natural notion, we introduce the following definition.

\begin{Def}[$\Mon$ and $\Conf$]
	\label{DefMonConf}
	\nomenclature[Mon]{$\Mon$}{functions $\Z \to \R$ that are monotonous and proper}
	\nomenclature[Conf]{$\Conf$}{$\Mon / \tau$}
	Define
	\begin{align*} 
		\Mon := \{ u:\Z \to \R \mid u \text{ is non-decreasing and proper} \} \subset  \R^{\Z}.
	\end{align*} 
	The group $(\Z,+)$ acts canonically on $\Mon$ via shifts, i.e.
	\begin{align} \label{EqDeftau}
	\begin{split}
		\begin{array}{rcl}
		\Mon \times \Z & \to & \Mon \\
		(u,z) & \mapsto & (j \mapsto (u.z)(j):=u(j+z))
		\end{array} 
	\end{split}
	\end{align}
	and the quotient
	\begin{align*}
		\Conf := \Mon / \Z
	\end{align*} 
	is called the \emph{configuration space}. Let $\pi:\Mon \to \Conf$, $u \mapsto \bar u$, be the quotient map.
\end{Def}

By construction $\spc^g \in \Mon$ and $\overline{\spc}^g := \pi (\spc^g) \in \Conf$. This defines maps
\begin{align*}
	\begin{array}{rclcrcl}
		\spc:\Rm(M) & \to & \Mon, & \hspace{4em} & \overline{\spc}:\Rm(M) & \to & \Conf, \\
		g & \mapsto & \spc^g, & & g & \mapsto & \overline{\spc}^{g}.
	\end{array} 
\end{align*}
We would like to say that $\overline{\spc}$ is continuous. To make formal sense of this, we introduce a topology on $\Mon$ and $\Conf$.
\nomenclature[s]{$\spc$}{discontinuous Dirac spectrum $\Rm(M) \to \Mon$}
\nomenclature[sbar]{$\overline{\spc}$}{Dirac spectrum $\Rm(M) \to \Conf$}

\begin{Def} [$\arsinh$-topology]
	\label{DefArsinhTopology}
	\nomenclature[da]{$d_a$}{$\arsinh$-metric}
	\nomenclature[dabar]{$\bar d_a$}{$\arsinh$-quotient-metric}
	The topology induced by the metric $d_a$ defined by
	\begin{align*}
		\forall u,v \in \R^{\Z}: d_a(u,v) := \sup_{j \in \Z}{|\arsinh(u(j)) - \arsinh(v(j))}| \in [0,\infty]
	\end{align*}
	on $\R^{\Z}$ is called the \emphi{$\arsinh$-topology}. The group $\Z$ acts by isometries with respect to $d_a$ and the quotient topology on $\Conf$ is induced by the metric $\bar d_a$ described by
	\begin{align} \label{EqDefbarda}
		\forall u \in \bar u, v \in \bar v \in \Conf: \bar d_a(\bar u, \bar v) = \inf_{j \in \Z}{d_a(u,v.j)}.
	\end{align}
	The use of this metric on the quotient is common in metric geometry, cf. \cite[Lemma 3.3.6]{Burago}.
\end{Def}

This allows us to formulate the following theorem, which is the major technical part of this chapter.
\begin{Thm}
	\label{MainThmSpec}
	\nomenclature[shat]{$\widehat{s}$}{continuous lift of $\overline{s}$}
	The map $\overline{\spc} = \pi \circ \spc$ admits a lift $\widehat{\spc}$ against $\pi$ such that
	\begin{align}
        \label{EqMainThmSpc}
        \begin{split}
            \xymatrix{
                & (\Mon,d_a)
                    \ar[d]^-{\pi} 
                \\
                (\Rm(M),\mathcal{C}^1)
                    \ar[r]^-{\overline{\spc}}
                    \ar[ur]^-{\widehat{\spc}}
                & (\Conf, \bar d_a)}
        \end{split}
	\end{align}
	is a commutative diagram of topological spaces.
\end{Thm}
First of all, we quickly convince ourselves that this theorem implies \cref{MainThmFun}. 

\nomenclature[ev]{$\ev_j$}{$\R^{\Z} \to \R$, $u \mapsto u(j)$}
\begin{Prf}[\cref{MainThmFun}]
	Clearly, the evaluation $\ev_j:(\Mon,d_a) \to \R$, $u \mapsto u(j)$, is a continuous map for any $j \in \Z$. Consequently, by \cref{MainThmSpec}, the functions $\{ \lambda_j := \ev_j \circ \widehat{\spc} \}_{j \in \Z}$ have the desired properties. The last assertion follows from \cref{RemMultK}.
\end{Prf}

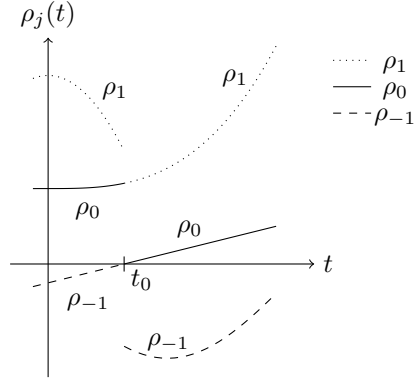
\begin{figure}[t] 
	\begin{center}
		\input{fig.jumps}
		\caption[Zeroes of eigenvalue functions.]{A zero of $\rho_{-1}$ at $t_0$ can cause discontinuities at $t_0$ in all $\rho_j$.}
		\label{FigJumps}
	\end{center}
\end{figure}

\begin{Rem}[intuitive explanation]
	\label{RemIntuitionDisc}
	It is a more subtle problem than one might think to choose functions $\{ \lambda_j \}_{j \in \Z}$, which depend continuously on $g \in \Rm(M)$ and represent the entire Dirac spectrum. The functions induced by $\spc$ (we call these $\rho_j := \ev_j \circ \spc$ for the moment) are not continuous in general. To see this, imagine a continuous path of metrics $(g_t)_{t \in \R}$ and consider $\rho_j:\R \to \R$ as functions of $t$, see \Cref{FigJumps}. Since $\rho_0(t)$ is the first eigenvalue $\geq 0$ of $\Dirac_{\C}^{g_t}$, this function will have a jump at a point $t_0$, where $\rho_0(t_0)>0$ and $\rho_{-1}(t_0)=0$. This can cause discontinuities in all the other functions $\rho_j$ as well. \\
	However, for any $k \in \Z$, the sequence $\rho_j' := \rho_{j+k}$, $j \in \Z$, gives another enumeration of the spectrum. Intuitively, \cref{MainThmSpec} states that if one uses this freedom in the enumeration of the eigenvalues at each metric in the ``right'' way, one obtains a globally well-defined family of continuous functions representing all the Dirac eigenvalues.
\end{Rem}

The rest of this \textname is organized as follows: After a slight generalization of our notation in \Cref{SctEvpaperNotSpect}, the main part will be \Cref{SctEvpaperProofMainThm}, which is devoted to build up the technical results needed in the proof of \cref{MainThmSpec}. Finally, in \Cref{SecModuliSpecFl} we will investigate to what extent the functions $\widehat{\spc}$ and $\overline{\spc}$ descend to certain quotients of $\Rm(M)$ called \emph{moduli spaces}. Our central result will be that there exists an obstruction, the \emph{spectral flow}, for $\widehat{\spc}$ to descend onto $\Rm(M) / \Diff^{\spin}(M)$. This will be made precise in \cref{DefSpinDiffeo} and \cref{LemSpecFlowDef} and the central result will be stated in \cref{MainThmFlow}.

Using these results, the actual proof of \cref{MainThmSpec} becomes very short.

\begin{Prf}[\cref{MainThmSpec}]
	By \cref{ThmSpecContArsinh}, the map $\bar \spc:(\Rm(M),\mathcal{C}^1) \to (\Conf, \bar d_a)$ is continuous. By \cref{ThmPiCoveringArsinh}, the map $\pi:(\Mon,d_a) \to (\Conf,\bar d_a)$ is a covering map. Since $\Rm(M)$ is path-connected, locally path-connected and simply connected, the result follows from the Lifting Theorem of Algebraic Topology, see \cref{ThmLifting}.
\end{Prf}

\begin{Rem}[uniqueness]
	From this proof, we see that the lift $\widehat{\spc}$ is not unique. In fact there are $\Z$ possibilities of how to lift $\overline{\spc}$ against $\pi$. One can use this freedom to arrange that $\widehat{\spc}^{g_0} = \spc^{g_0}$ for one fixed $g_0 \in \Rm(M)$.
\end{Rem}

\begin{Rem}[supremum norm]
	Can one replace the metric $d_a$ from \cref{DefArsinhTopology} by the ordinary supremum norm? Of course one can also topologize $\Mon$ using the supremum norm instead, but then no equicontinuity result like \cref{MainThmFun} can hold. This can be quickly seen as follows: Take any metric $g \in \Rm(M)$. For any $t > 0$, the tensor field $t g$ is also a Riemannian metric and $tg \to g$ as $t \to 1$ in $\Rm(M)$. By conformal covariance, $\spec \Dirac_{\C}^{tg} = t^{-1/2} \spec \Dirac_{\C}^g$, see for instance \cite[Prop. 5.13]{hijazi}. So if the $\{\lambda_j\}_{j \in \Z}$ themselves were equicontinuous, for $\varepsilon = 1$ there exists $\delta > 0$ such that
	\begin{align*}
		\forall t \in I_{\delta}(1): \forall j \in \Z : |t^{-1/2}-1||\lambda_j(g)| = |\lambda_j(tg) - \lambda_j(g)|  < 1.
	\end{align*}
	Since $\lim_{j \to \infty}{\lambda_j(g)} = \infty$, this simply cannot be true.
\end{Rem}

\begin{Rem}[$\arsinh$]
	However, we will see later in the proof \cref{MainThmSpec} that it is not so important to choose precisely the function $\arsinh$ in \cref{DefArsinhTopology} to obtain the result. The most important feature of $\arsinh$ that will enter in the proof is that $\xi \mapsto \arsinh'(\xi) |\xi|$ is a bounded function on $\R^n$. Therefore, one could replace $\arsinh$ by a function with similar growth properties.
\end{Rem}

Finally, we remark that the following well known theorem implies that a bounded spectral interval of the Dirac operator can be described locally by continuous functions. Consequently, \cref{MainThmSpec} can be thought of as a global analogue of this local result.

\begin{Thm}[\protect{\cite[Prop. 7.1]{BaerMetrHarmonSpin}}]
	\label{ThmBaerSpecC1Conti}
	Let $(M,g)$ be a closed Riemannian spin manifold with Dirac operator $\Dirac_{\C}^g$ having spectrum $\spec \Dirac_{\C}^g$. Let $\Lambda > 0$ such that $-\Lambda,\Lambda \notin \spec \Dirac_{\C}^g$ and enumerate
	\begin{align*}
		\spec \Dirac_{\C}^g \cap \mathopen{]} -\Lambda,\Lambda \mathclose{[} = \{ \lambda_1 \leq \lambda_2 \leq \cdots \leq \lambda_n \}.
	\end{align*}
	For any $\varepsilon > 0$, there exists a $\mathcal{C}^1$-neighborhood $U$ of $g$ such that for any $g' \in U$ 
	\begin{enumerate}
		\item $\spec \Dirac_{\C}^{g'} \cap \mathopen] -\Lambda, \Lambda \mathclose[ = \{ \lambda'_1 \leq \cdots \leq \lambda'_n \}$,
		\item $\forall 1 \leq i \leq n: |\lambda_i - \lambda_i'| < \varepsilon$.
	\end{enumerate}
\end{Thm}

\begin{Rem}
	\label{RemEvpaperRealComplexSpinMult}
	In this chapter, all results will be formulated for the complex Dirac operator and complex multiplicities. However, in dimensions $m \equiv 2,3,4 \mod 8$, where the eigenspaces have quaternionic dimensions, we can take $\{\lambda_j\}_{j \in 2\Z}$ and obtain a family of continuous functions that represents the entire Dirac spectrum at every metric counted with spin multiplicity (see \cref{DefDahlMult}). For the real Dirac operator, the proof of \cref{MainThmFun} does not directly go through, because it relies on Kato's pertubation theory, see \cref{ThmAnalyticEVKato}, which is formulated for holomorphic families of operators. However, the $\R$-multiplicities of the eigenvalues of the real Dirac operator $\Dirac^g_{\R}$ are the same as the $\C$-multiplicities of its complexification $(\Dirac^g_{\R})^{\C}$. Therefore, we can apply \cref{ThmDiscFamA} to the complexification of the real Dirac operator and obtain a version of \cref{MainThmFun} for the real Dirac operator and $\R$-multiplicities. Of course, in dimensions $m \equiv 0,6,7 \mod 8$, it makes no difference anyway, which multiplicity is chosen, see \cref{RemMultK}.
\end{Rem}

%% file: fig.specenum.tex
\begin{tikzpicture}

	\draw[<->,thick] (-4,0) to (4,0);
	\draw[thick] (0,-0.2) to (0,0.2);

	\coordinate[label=$0$] (A) at (0,-0.7);

	\foreach \xp/\yp in {1/0,2/0,2/0.3,2.8/0}{
		\draw (\xp-0.1,\yp-0.1) to (\xp+0.1,\yp+0.1);
		\draw (\xp-0.1,\yp+0.1) to (\xp+0.1,\yp-0.1);
	}

	\foreach \xp/\yp in {-1.3/0,-2.4/0,-2.4/0.3}{
		\draw (\xp-0.1,\yp-0.1) to (\xp+0.1,\yp+0.1);
		\draw (\xp-0.1,\yp+0.1) to (\xp+0.1,\yp-0.1);
	}
	
	\coordinate[label=$\lambda_0$] (A) at (1.1,-0.7);
	\coordinate[label=$\lambda_1$] (A) at (2.1,-0.7);
	\coordinate[label=$\lambda_2$] (A) at (2.1,-1.1);
	\coordinate[label=$\lambda_3$] (A) at (2.9,-0.7);
	\coordinate[label=$\ldots$] (A) at (3.7,-0.5);

	\coordinate[label=$\lambda_{-1}$] (A) at (-1.4,-0.7);
	\coordinate[label=$\lambda_{-2}$] (A) at (-2.4,-0.7);
	\coordinate[label=$\lambda_{-3}$] (A) at (-2.4,-1.1);
	\coordinate[label=$\ldots$] (A) at (-3.3,-0.5);

\end{tikzpicture}

%% file: fig.jumps.tex
\begin{tikzpicture}
	\draw[->] (-0.5,0) -- (3.5,0) node[right] {$t$};
	\draw[->] (0,-1.5) -- (0,3) node[above] {$\rho_j(t)$};
	\draw[color=black,  domain=-0.2:1] plot[id=fa] function{0.07 * x**3 +1 };
	
	\draw[color=black, dotted, domain=1:3] plot[id=fb] function{0.07 * x**3 +1 };
	
	\draw[color=black, dashed,  domain=-0.2:1] plot[id=fc] function{0.25*x-0.25};

	\draw[color=black, domain=1:3] plot[id=fd] function{0.25*x-0.25};
	
	\draw[color=black, dotted, domain=-0.2:1] plot[id=fe] function{-x*x + 2.5};
	
	\draw[color=black, dashed, domain=1:3] plot[id=ff] function{-sin(x)-0.25};
	
	\node [label={[xshift=15, yshift=-25]$\rho_{-1}$}] {};
	\node [label={[xshift=45, yshift=-40]$\rho_{-1}$}] {};
	\node [label={[xshift=15, yshift=10]$\rho_{0}$}] {};
	\node [label={[xshift=53, yshift=3]$\rho_{0}$}] {};
	\node [label={[xshift=25, yshift=55]$\rho_{1}$}] {};
	\node [label={[xshift=70, yshift=60]$\rho_{1}$}] {};
	
	\draw[color=black] (1,0.1) -- (1,-0.1);
	\node [label={[xshift=35, yshift=-17]$t_0$}]  {};
	
	\draw[color=black, dotted] (3.75,2.7) -- (4.25,2.7);
	\draw[color=black] (3.75,2.35) -- (4.25,2.35);
	\draw[color=black, dashed] (3.75,2) -- (4.25,2);
	
	\node [label={[xshift=130, yshift=65]$\rho_{1}$}]  {};
	\node [label={[xshift=130, yshift=55]$\rho_{0}$}] {};
	\node [label={[xshift=130, yshift=45]$\rho_{-1}$}] {};

\end{tikzpicture}

%% file: evpaper.famdiscrete.tex
\section{Families of Discrete Operators}
\label{SctEvpaperNotSpect}
For the proof of \cref{MainThmSpec} we need a version of the function $\spc^g$ defined for operators slightly more general than Dirac operators. In this section, we introduce the necessary definitions and notational conventions. Let $X,Y$ be complex Banach spaces. We denote by $C(X,Y)$ the space of closed unbounded operators $T:X \supset \dom(T) \to Y$. Let $B(X,Y)$ denote the bounded operators $X \to Y$. We set $C(X):=C(X,X)$ and $B(X):=B(X,X)$. The spectrum of $T$ is denoted by $\spec T \subset \C$.

\nomenclature[CXY]{$C(X, Y)$}{closed unbounded operators $X \to Y$}
\nomenclature[CX]{$C(X)$}{closed unbounded operators on $X$}
\nomenclature[BXY]{$B(X,Y)$}{bounded operators $X \to Y$}
\nomenclature[BX]{$B(X)$}{bounded operators on $X$}
\nomenclature[spec]{$\spec T$}{spectrum of $T$}

\begin{Def}[discrete operator]
	\label{DefDiscreteOperator} \index{discrete!operator}
	An operator $T \in C(X)$ is \emph{discrete}, if $T$ has compact resolvent and $\spec T \subset \R$ is unbounded from both sides. (It follows  that $\spec T$ is closed and consists solely of eigenvalues of finite multiplicity.)
\end{Def}

\nomenclature[spcT]{$\spc_T$}{spectral function of $T$}
\begin{DefI}[ordered spectral function] 
	\label{DefOrderedSpecFunT}
	Let $T \in C(X)$ be discrete. The sequence $\spc_T \in \R^{\Z}$, uniquely defined by the properties
	\begin{enumerate}
		\item $0 = \min \{ j \in \Z \mid \spc_T(j) \geq 0 \}$.
		\item $\forall i,j \in \Z: i \leq j \Longrightarrow \spc_T(i) \leq \spc_T(j)$.
		\item $\forall \lambda \in \R: \sharp (\spc_T)^{-1}(\lambda)= \dim_{\C} \ker (T - \lambda) $.
	\end{enumerate}
	is the \emph{ordered spectral function of $T$}. 
\end{DefI}

\nomenclature[spp]{$\spp_T(I)$}{spectrum of $T$ in $I$}
\begin{DefI}[spectral parts] 
	Let $T \in C(X)$ be discrete. To denote parts of the ordered spectrum, we introduce the following notation: If $I \subset \R$ is an interval, then $(\spc_T)^{-1}(I) = \{k, k+1, \ldots, l\}$ for some $k, l \in \Z$, $k \leq l$. The sequence
	\begin{align*} 
		\spp_T(I) := (\spc_T(i))_{k \leq i \leq l},
	\end{align*}
	is the \emph{spectral part of $T$ in $I$}. 
\end{DefI}

\begin{Def}[discrete family]
	\index{discrete!family}
	Let $E$ be any set. A map $T:E \to C(X)$ is a \emph{discrete family}, if for any $e \in E$, the operator $T_e$ is discrete in the sense of \cref{DefDiscreteOperator}. We obtain a function
	\DefMap{\spc_T:E}{\Mon}{e}{\spc_T^e:=\spc_{T(e)}.}
	Analogously, we set $\spp_T^e := \spp_{T(e)}$.
\end{Def}

\begin{Rem}[family of Dirac operators]
	In view of \cref{ThmSpinorIdentification}, we can apply the above in particular to Dirac operators. Namely, we fix $g \in \Rm(M)$ and set $X:=L^2(\Sigma^g_{\C} M)$ and $E:=\Rm(M)$. Then $h \mapsto \Dirac^h_g$ is a discrete family. We will suppress its name in notation and just write $\spc^h_g$ for the ordered spectral function of $\Dirac^h_g$. Since $\Dirac^h_{\C}$ and $\Dirac^h_g$ are isospectral, we can ignore the reference metric entirely and simply write $\spc^h$.
\end{Rem}

%% file: evpaper.prfmain.tex
\section{Proof of \cref{MainThmSpec}}
\label{SctEvpaperProofMainThm}

In this section, we carry out the details of the proof of \cref{MainThmSpec}. The idea to construct the $\arsinh$-topology was inspired by a paper of John Lott, cf. \cite[Theorem 2]{Lott}. The arguments require some basic notions from analytic pertubation theory. A modified version of some results by Kato is needed, cf. \cite{kato}. Applying analytic pertubation theory to families of Dirac operators is a technique that is also used in other contexts, cf. \cite{BaerGaudMor}, \cite{BourgGaud}, \cite{AndreasDiss}. \\
Let $X$,$Y$ be complex Banach spaces, and let $X'$ be the topological dual space of $X$. For any operator $T$, we denote its adjoint by $T^*$. Let $\Omega \subset \C$ be an open and connected subset. Recall that a function $f:\Omega \to X$ is \emphi{holomorphic}, if for all $\zeta_0 \in \Omega$, the limit
\begin{align*}
	f'(\zeta_0) := \lim_{\zeta \to \zeta_0}{\tfrac{f(\zeta)-f(\zeta_0)}{\zeta-\zeta_0}}
\end{align*}
exists in $(X,\|\_\|_X)$. A family of operators $T:\Omega \to B(X,Y)$ is called \emph{bounded holomorphic}\index{holomorphic!bounded}, if $T$ is a holomorphic map in the above sense. To treat the unbounded case, the following notions are crucial.

\begin{Def}[holomorphic family of type (A)]
	\index{holomorphic!family of type (A)}
	A family of operators $T:\Omega \to C(X,Y)$, $\zeta \mapsto T_\zeta$, is \emph{holomorphic of type (A)}, if the domain $\dom(T_\zeta) =: \dom(T)$ is independent of $\zeta$, and for any $x \in \dom(T)$ the map $\Omega \to Y$, $\zeta \mapsto T_\zeta x$, is holomorphic.
\end{Def}

\begin{Def}[self-adjoint holomorphic family of type(A)]
	\index{holomorphic!self-adjoint of type (A)}
	A family $T:\Omega \to C(H)$ is \emph{self-adjoint holomorphic of type (A)}, if it is holomorphic of type (A), $H$ is a Hilbert space, $\Omega$ is symmetric with respect to complex conjugation, and
	\begin{align*}
		\forall \zeta \in \Omega: T^*_{\zeta} = T_{\bar \zeta}.
	\end{align*}
\end{Def}

These families are particularly important for our purposes due to the following useful theorem.

\begin{Thm}[\protect{\cite[VII.\textsection 3.5, Thm. 3.9]{kato}}]
	\label{ThmAnalyticEVKato}
	Let $T:\Omega \to C(H)$ be a self-adjoint holomorphic family of type (A), and let $I \subset \Omega \cap \R$ be an interval. Assume that $T_\zeta$ has compact resolvent for all $\zeta \in \Omega$. Then, there exists a family of functions $\{\lambda_n \in \mathcal{C}^\omega(I,\R)\}_{n \in \N}$ and a family of functions $\{u_n \in \mathcal{C}^\omega(I,H)\}_{n \in \N}$ such that for all $t \in I$, the $(\lambda_n(t))_{n \in \N}$ represent all the eigenvalues of $T_t$ (counted with multiplicity), $T_tu_n(t)=\lambda_n(t)u_n(t)$, and the $(u_n(t))_{n \in \N}$ form a complete orthonormal system of $H$.
\end{Thm}

Derivatives of holomorphic families can be estimated using the following theorem.

\begin{Thm}[\protect{\cite[VII.\textsection 2.1, p.375f]{kato}}]
	\label{ThmHolAFamilyDerGrowth}
	Let $T:\Omega \to C(X,Y)$ be a holomorphic family of type (A). For any $\zeta \in \Omega$ define the operator
	\begin{align*}
		T'_{\zeta}: \dom(T) \to Y, && u \mapsto T'_{\zeta}u := \tfrac{d}{d\zeta}(T_{\zeta}u).
	\end{align*}
	Then $T'$ is a map from $\Omega$ to the unbounded operators $X \to Y$ (but $T'_{\zeta}$ is in general not closed). For any compact $K \subset \Omega$, there exists $C_K>0$ such that 
	\begin{align*}
		\forall \zeta \in K: \; \forall u \in \dom(T): \; \|T'_{\zeta}u\|_Y \leq C_K(\|u\|_X + \|T_{\zeta}u\|_Y).
	\end{align*}
	If $\zeta_0 \in K$ is arbitrary, $Z:=\dom(T)$ and $\|u\|_Z := \|u\|_X + \|T_{\zeta_0}u\|_Y$, then $C_K := \alpha_K^{-1} \beta_K$ does the job, where
	\begin{align} \label{EqHolAFamilyDerGrowthConstantSpec}
		\alpha_K := \inf_{\zeta \in K}{ \inf_{\|u\|_Z=1}{\|u\|_X + \|T_{\zeta}u\|_Y}}, &&
		\beta_K := \sup_{\zeta \in K}{\|T'_{\zeta}\|_{B(Z,Y)}}.
	\end{align}
\end{Thm}

This can be used to prove the following result about the growth of eigenvalues. 

\begin{Thm}[\protect{\cite[VII.\textsection 3.4, Thm. 3.6]{kato}}]
	\label{ThmHolAFamilyEVGrowth}
	Let $T:\Omega \to C(H)$ be a self-adjoint holomorphic family of type (A).  Let $I \subset \Omega \cap \R$ be a compact interval, and let $J \subset I$ be open. Assume that $\lambda \in \mathcal{C}^\omega(J,\R)$ is an \emph{eigenvalue function}, i.e. for all $ t \in J$ the value $\lambda(t)$ is an eigenvalue of $T_t$. Then
	\begin{align} \label{EqHolAFamilyConstantExp}
		\forall t,t_0 \in J: \; |\lambda(t) - \lambda(t_0)| \leq (1 + |\lambda(t_0)|)(\exp(C_I |t-t_0|) - 1),
	\end{align}
	where $C_I$ is the constant from \cref{ThmHolAFamilyDerGrowth}. 
\end{Thm}

The preceding \cref{ThmHolAFamilyEVGrowth} provides two key insights into the growth of eigenvalue functions. First of all, it is remarkable that the constant $C_I$ in \eqref{EqHolAFamilyConstantExp} does not depend on the eigenvalue function $\lambda$. In particular, if we consider a family of eigenvalue functions $\{\lambda_n\}_{n \in \N}$ as in \cref{ThmHolAFamilyDerGrowth}, the constant $C_I$ is uniform in $n$. This will be crucial later in the proof of \cref{CorAnalyticGlobalGrowth}. Secondly, we see that, due to the factor $1 + |\lambda(t_0)|$ in \eqref{EqHolAFamilyConstantExp}, an eigenvalue function grows faster the larger it is. This is the reason why we cannot expect the continuity result of \cref{MainThmSpec} to hold for the ordinary supremum norm. However, as we will show in the next corollary, we can get rid of this factor by reformulating \eqref{EqHolAFamilyConstantExp} in terms of the $\arsinh$-topology.

\begin{Cor}[growth of eigenvalues]
	\label{CorHolAFamilyEVGrowth}
	In the situation of \cref{ThmHolAFamilyEVGrowth}, the following holds in addition: For any $t_0 \in I$ and $\varepsilon > 0$, there exists $\delta > 0$ such that for all $t \in I_{\delta}(t_0) \cap J$ and all eigenvalue functions $\lambda \in \mathcal{C}^\omega(J,\R)$ we have
	\begin{align} \label{EqHolAFamilyDerGrowthArsinh}
		|\arsinh(\lambda(t)) - \arsinh(\lambda(t_0))| < \varepsilon.
	\end{align}
	There exist universal constants (i.e. independent of the family $T$) $C_1,C_2 > 0$ such that 
	\begin{align} \label{EqHolAFamilyDerGrowthDelta}
		\delta := C_I^{-1} \ln(\min(C_1, \varepsilon C_2)  + 1 )
	\end{align} 
	does the job.
\end{Cor}

\begin{Prf}
	The key observation needed is that $\arsinh(t)$ grows slower the larger $|t|$ gets. This follows simply from the formula $\arsinh'(t) = (1+t^2)^{-1/2}$. We will show that this neutralizes the $(1+|\lambda(t_0)|)$-factor in \eqref{EqHolAFamilyConstantExp} when the growth of $\lambda$ is measured in the $\arsinh$-metric. The $\exp$-term in \eqref{EqHolAFamilyConstantExp} can be estimated by a standard continuity argument.

	\begin{steplist}
	\step[$\exp$-term]
		The function $\alpha:\R \to \R$, $t \mapsto \exp(C_I |t-t_0|) - 1$, is continuous and satisfies $\alpha(t_0)=0$. Notice that for $b > 0$ 
		\begin{align}
			\label{EqHolAFamilyEVGrowthDeltaCalc}
			|\alpha(t)| < b \Longleftrightarrow |t-t_0| < C_I^{-1} \ln(b+1).
		\end{align}
		In particular there exists $\delta_1 > 0$ such that 
		\begin{align} \label{EqHolAFamilyConstantExpDeltaQuarter}
			\forall t \in I_{\delta_1}(t_0): |\alpha(t)| < \tfrac{1}{4}.
		\end{align}
		So let $t \in I_{\delta_1}(t_0)$. 
		
	\step[preliminary estimate] Setting $\lambda_0:=\lambda(t_0)$ we can reformulate \eqref{EqHolAFamilyConstantExp} as
	\begin{align} \label{HolAFamilyEVGrowthReform}
	\lambda_0 - (1 + |\lambda_0|)\alpha(t) < \lambda(t) < \lambda_0 +(1+|\lambda_0|) \alpha(t).
	\end{align}
	Since
	\begin{align*}
		\lim_{|R| \to \infty}{\tfrac{|R|}{1+|R|}} = 1 ,
	\end{align*}
	and the convergence is monotonously increasing, there exists $R > 0$ such that
	\begin{align} \label{EqEVfunArsinhExQuotient}
		\forall |\eta| \geq R: \tfrac{1}{2} < \tfrac{|\eta|}{1 + |\eta|}.
	\end{align}
	Now assume $|\lambda_0| \geq R$. In case $\lambda_0 \geq R > 0$, we calculate
	\begin{align}
		\label{EqHolAFamilyEVGrowthBoundBelowPositive}
		\tfrac{\lambda_0}{1 + \lambda_0} > \tfrac{1}{2} \jeq{\eqref{EqHolAFamilyConstantExpDeltaQuarter}}{\geq} 2 \alpha(t)		
		\Longrightarrow \tfrac{1}{2} \lambda_0 \geq \alpha(t)(1+\lambda_0)
		\Longrightarrow \lambda_0 - \alpha(t)(1 +|\lambda_0|) \geq  \tfrac{1}{2} \lambda_0 .
	\end{align}
	Analogously, if  $\lambda_0 \leq -R < 0$, we calculate
	\begin{align}
		\label{EqHolAFamilyEVGrowthBoundBelowNegative}
		\lambda_0 + \alpha(t) (1 + |\lambda_0|) < \tfrac{1}{2} \lambda_0.
	\end{align}

	\step[$\arsinh$-metric]
		Define the constants
		\begin{align*}
			C_0 := \sup_{t \in \R}{\tfrac{1 + |t|}{\sqrt{1 + t^2}}}, &&
			C_1 := \tfrac{1}{4}, &&
			C_2 := \min \left( \tfrac{1}{R+1}, \tfrac{1}{2C_0} \right),
		\end{align*}
		and set 
		\begin{align*}
		\delta_2 := C_I^{-1} \ln(\min(C_1, \varepsilon C_2)  + 1 ).
		\end{align*}
		By \eqref{EqHolAFamilyEVGrowthDeltaCalc} this implies 
		\begin{align} \label{EqHolAFamilyEVGrowthDeltaInp}
			\forall t \in I_{\delta_2}(t_0): \alpha(t) < \min(C_1,\varepsilon C_2) \leq \varepsilon C_2.
		\end{align}
		So let $t \in I_{\delta_2}(t_0)$ be arbitrary and set $c_\pm := \lambda_0 \pm (1 + |\lambda_0|)\alpha(t) $. It follows from the Taylor series expansion of $\arsinh$ that there exists $\xi \in [\lambda_0,c_+]$ such that
		\begin{align} \label{EqHolAFamilyEVGrowthFinal} 
			\arsinh(c_+) - \arsinh(\lambda_0) 
			= \arsinh'(\xi) (1 + |\lambda_0|)\alpha(t) 
			=\frac{(1 + |\lambda_0|)}{\sqrt{1 + \xi^2}} \alpha(t).
		\end{align}
		Now in case $|\lambda_0| \leq R$, we continue this estimate by
		\begin{align*}
			\eqref{EqHolAFamilyEVGrowthFinal}
			\leq (1 + |\lambda_0|) \alpha(t)
			\leq (1 + R) \alpha(t) 
			\jeq{\eqref{EqHolAFamilyEVGrowthDeltaInp}}{<} \varepsilon.
		\end{align*}
		In case $\lambda_0 \geq R$, we continue this estimate by
		\begin{align*}
			\eqref{EqHolAFamilyEVGrowthFinal}
			\leq \frac{(1 + |\lambda_0|)}{\sqrt{1 + \lambda_0^2}} \alpha(t)
			\leq C_0 \alpha(t)  
			\jeq{\eqref{EqHolAFamilyEVGrowthDeltaInp}}{<} \varepsilon.
		\end{align*}
		In case $\lambda_0 \leq -R$, we continue this estimate by
		\begin{align*}
			\eqref{EqHolAFamilyEVGrowthFinal}
			\leq \frac{(1 + |\lambda_0|)}{\sqrt{1 + c_+^2}} \alpha(t)
			\jeq{\eqref{EqHolAFamilyEVGrowthBoundBelowNegative}}{\leq} \frac{(1 + |\lambda_0|)}{\sqrt{1 + \tfrac{1}{4} \lambda_0^2}} \alpha(t)
			\leq 2 C_0 \alpha(t)  
			\jeq{\eqref{EqHolAFamilyEVGrowthDeltaInp}}{<} \varepsilon.
		\end{align*}
		Consequently, since $\arsinh$ is strictly increasing, in all cases we obtain 
		\begin{align*}
			\arsinh(\lambda(t))
			\jeq{\eqref{HolAFamilyEVGrowthReform}}{<}\arsinh(\lambda_0 + (1 + |\lambda_0|)\alpha(t))
			< \arsinh(\lambda_0) + \varepsilon.
		\end{align*}
		By an analogous argument, we obtain 
		\begin{align*}
			\arsinh(\lambda(t))
			> \arsinh(\lambda_0 - (1+|\lambda_0|)\alpha(t))
			\geq \arsinh(\lambda_0) - \varepsilon .
		\end{align*}
	\end{steplist}
	This proves the claim.	
\end{Prf}

In the next step, we will apply the preceding result to discrete families.

\begin{Not} 
	\nomenclature[daxy]{$d_a(x, y)$}{$|\ash(x)-\ash(y)|$}
	\nomenclature[Iex]{$I_{\varepsilon}(x)$}{$\varepsilon$-neighborhood of $x$}
	\nomenclature[IeS]{$I_{\varepsilon}(S)$}{$\varepsilon$-hull of $S$}
	\label{NotTechnical}
	Since the following proof is somewhat technical, we abbreviate $d_a(x,y):=|\ash(x)-\ash(y)|$ for $x,y \in \R$. For any $\varepsilon > 0$, the $\varepsilon$-neighborhoods of $x \in \R$ and $\varepsilon$-hulls of a set $S \subset \R$ will be denoted, respectively, by 
	\begin{align*}
		I_\varepsilon(x) := \{ t \in \R \mid |t-x| < \varepsilon \}, &&
		I_\varepsilon(S) := \bigcup_{x \in S}{I_\varepsilon(x)}.
	\end{align*}
\end{Not}

\begin{Cor}[spectral growth]
	\label{CorAnalyticGlobalGrowth}
	Let $\Omega \subset \C$ be open, let $I \subset \Omega \cap \R$ be an interval, and let $T:\Omega \to C(H)$ be a discrete and self-adjoint holomorphic family of type (A). For any $t_0 \in I$ and $\varepsilon > 0$, there exists $\delta > 0$ such that 
	\begin{align} \label{EqAnalyticGlobalGrowth}
		\forall t \in I_\delta(t_0) \cap I: \exists k \in \Z: \forall j \in \Z: d_a(\spc^{t_0}_T(j),\spc^{t}_T(j+k)) < \varepsilon.
	\end{align}
\end{Cor}

\begin{Prf} $ $ 
\begin{steplist}
	\step[apply \cref{ThmAnalyticEVKato}]
		Certainly the family of eigenfunctions from \cref{ThmAnalyticEVKato} can be $\Z$-reindexed to a family $\{ \lambda_j \in \mathcal{C}^\omega(I,\R) \}_{j \in \Z}$ satisfying $\lambda_j(t_0)=\spc_T^{t_0}(j)$, $j \in \Z$. By \cref{CorHolAFamilyEVGrowth} 
		\begin{align*}
			\exists \delta > 0: \forall t \in I_\delta(t_0) \cap I: \forall j \in \Z: |\ash(\lambda_j(t_0)) - \ash(\lambda_j(t))| < \varepsilon .
		\end{align*}
		Fix any $t \in I_\delta(t_0) \cap I$ and let $\sigma:\Z \to \Z$ be the bijection satisfying $\spc_T^t(\sigma(j)) = \lambda_j(t)$. This implies 
		\begin{align} \label{EqAnalyticGlobalGrowthBij}
			\forall j \in \Z: d_a(\spc_T^{t_0}(j),\spc_T^t(\sigma(j)) < \varepsilon,
		\end{align}
		which is almost \eqref{EqAnalyticGlobalGrowth}, except that $\sigma$ might not be given by a translation.
		
	\step[general idea]
		We will show that we may replace the bijection $\sigma$ by an increasing bijection $\tau$, which still satisfies \eqref{EqAnalyticGlobalGrowthBij}. Since every increasing bijection $\Z \to \Z$ is given by a translation $\tau^{(k)}:\Z \to \Z$, $z \mapsto z + k$, for some $k \in \Z$, this implies the claim. To that end we will first show how to modify $\sigma$ on finite subsets and then use the pigeonhole principle to conclude the argument.

	\step[on finite subsets]
		For any $n \in \N$ set $I_n := \{-n, \ldots, n\}$ and consider the function $\sigma_n:=\sigma|_{I_n}:I_n \to \Z$. This function is injective and satisfies \eqref{EqAnalyticGlobalGrowthBij} for all $-n \leq j \leq n$. Furthermore, setting
		\begin{align*}
			\spp^T_{t_0}([\lambda_{-n}(t_0), \lambda_n(t_0)]) &=: (\lambda_{-n}, \ldots, \lambda_{n})
		\end{align*}
		we obtain numbers $n',m' \in \Z$ such that the eigenvalues $\mu_j := \spc_T^{t}(j)$ satisfy
		\begin{align} \label{AnalyticLocalGrowthEVNumber}
			I_\varepsilon(\ash(\spp^T_{t}([\lambda_{-n},\lambda_{n}]))) & = (\ash(\mu_{n'}), \ldots, \ash(\mu_{m'}) ), 
		\end{align}
		and we have the estimate
		\begin{align} \label{EqAnalyticGlobalGrowthLocalized} 
			\forall -n \leq j \leq n: & |\ash(\lambda_j) - \ash(\mu_{\sigma_n(j)}) | < \varepsilon.
		\end{align}
		We will show that $\sigma_n$ can be modified to an increasing injection $\tilde \sigma_n$, which satisfies $\image(\sigma_n) = \image(\tilde{\sigma}_n)$ and \eqref{EqAnalyticGlobalGrowthLocalized}. To that end, choose any $-n \leq i < j \leq n$ and assume that $\sigma_n(j)<\sigma_n(i)$. Notice that by construction
		\begin{align}
			\label{EqAnalyticLocalGrowthMonotoneLambdaMu}
			i < j \Longrightarrow  \lambda_i \leq \lambda_j , && \sigma_n(j) < \sigma_n(i) \Longrightarrow \mu_{\sigma_n(j)} \leq \mu_{\sigma_n(i)}.
		\end{align}
		Define the function $\tilde \sigma_n$ by setting 
		\begin{align*}
			\tilde \sigma_n|_{\{-n, \ldots, n\} \setminus \{i,j\}} := \sigma_n, &&
		\tilde \sigma _n(i)=\sigma_n(j), &&
		\tilde \sigma_n(j)=\sigma_n(i).
		\end{align*}
		It is clear that $\tilde \sigma_n$ is still injective and $\image(\tilde \sigma_n) = \image(\sigma_n)$. To show that it still satisfies \eqref{EqAnalyticGlobalGrowthLocalized}, we distinguish two cases (see \Cref{FigCones}). First consider the case that $\lambda_i = \lambda_j$. Then, it follows automatically that
		\begin{align*}
			|\ash(\lambda_i) - \ash(\mu_{\tilde \sigma_n(i)})| = |\ash(\lambda_j) - \ash(\mu_{\sigma_n(j)})| < \varepsilon,
		\end{align*}
		and the same for $j$. In case $\lambda_i \neq \lambda_j$, it follows that $\lambda_i < \lambda_j$. This implies
		\begin{align*}
			\ash(\lambda_i) - \varepsilon 
			& < \ash(\lambda_j) - \varepsilon 
			< \ash(\mu_{\sigma_n(j)}) \\
			& \jeq{\eqref{EqAnalyticLocalGrowthMonotoneLambdaMu}}{\leq} \ash(\mu_{\sigma_n(i)}) 
			< \ash(\lambda_i) + \varepsilon 
			< \ash(\lambda_j) + \varepsilon,
		\end{align*}
		hence
		\begin{align*}
			\mu_{\sigma_n(i)},\mu_{\sigma_n(j)} \in I_\varepsilon(\ash(\lambda_i)) \cap I_\varepsilon(\ash(\lambda_j)).
		\end{align*}
		In particular, this intersection is not empty. Consequently, $\tilde \sigma_n$ satisfies \eqref{EqAnalyticGlobalGrowthLocalized}. By repeating this procedure for all index pairs $(i,j)$, $-n \leq i \leq n$, $i<j \leq n$, it follows that $\sigma_n$ can be modified finitely many times in this manner to obtain an increasing injection having the same image, which still satisfies \eqref{EqAnalyticGlobalGrowthLocalized}. For simplicity, denote this function also by $\tilde \sigma_n$ and, define
		\DefMap{\tilde \tau_n:\Z}{\Z}{j}{\begin{cases}
		\tilde \sigma_n(j), & -n \leq j \leq n, \\
		\sigma(j), & \text{otherwise}.
		\end{cases}}
		This function is still bijective, still satisfies \eqref{EqAnalyticGlobalGrowthBij} and is increasing on $I_n$. Define $J_n:=\tilde{\tau}_n(I_n)$.
		
	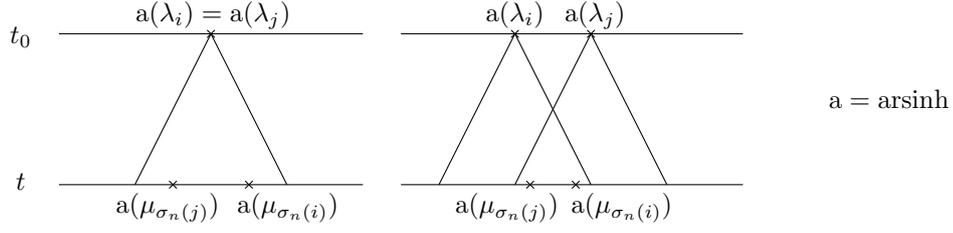
\begin{figure}[t] 
		\begin{center}
			\input{fig.cones}
			\caption[Proof of \cref{CorAnalyticGlobalGrowth}.]{The two possibilities for $\lambda_j$.}
			\label{FigCones}
		\end{center}
	\end{figure}

	\step[pigeonhole principle] 
		Unfortunately, it might happen that $\tilde{\tau}_{n+1}|_{I_n} \neq \tilde{\tau}_n$. Due to \eqref{EqAnalyticGlobalGrowthLocalized} however, there exists $n_1$ such that all $n \geq n_1$ satisfy $\tilde{\tau}_n(I_1) \subset J_{n_1}$. Since there are only finitely many functions $I_{1} \to J_{n_1}$, there must be at least one such function occurring infinitely often in the sequence $\{\tilde{\tau}_n|_{I_1}\}_{n \in \N}$. Thus, there exists an infinite subset $\N_1 \subset \N_0:=\N$ such that $\tilde{\tau}_n|_{I_1}$ is the same for all $n \in \N_1$. \\
		Now, the same holds for $I_2$: There exists $n_2 \geq n_1$ such that all $n \geq n_2$ satisfy $\tau_n(I_2) \subset J_{n_2}$. Again, since there are only finitely many functions $I_2 \to J_{n_2}$, one of them must occur infinitely often in the sequence $\{\tau_n|_{I_2}\}_{n \in \N_1}$. Consequently, there exists an infinite subset $\N_2 \subset \N_1$ such that $\tau_n|_{I_2}$ is the same for all $n \in \N_2$. This process can be continued indefinitely for all the intervals $I_\nu$, $\nu \in \N$. Finally, the function
		\DefMap{\tau:\Z}{\Z}{j}{\tilde{\tau}_n(j), \; j \in I_\nu, n \in \N_{\nu}}
		does the job: It is well-defined, satisfies \eqref{EqAnalyticGlobalGrowthBij}, remains injective and is surjective: Since the sets $\{I_n\}_{n \in \N}$ exhaust all of $\Z$, and since the $\tilde \tau_n$ are bijective and increasing on $I_n$, it follows that the $J_n$ are also sets of subsequent numbers in $\Z$. Thus, by injectivity of the $\tilde{\tau}_n$, the $\{J_n\}_{n \in \N}$ exhaust all of $\Z$.
	\end{steplist}
\end{Prf}

The preceding \cref{CorAnalyticGlobalGrowth} is almost the result we need to conclude the proof of \cref{MainThmSpec}, except that it is formulated only for paths of operators. As a last step we provide a framework, which allows us to pass from paths of operators to families of operators parametrized by more general spaces.

\begin{Def}[discrete family of type (A)]
	\label{DefDiscFamA} \index{discrete!family of type(A)}
	Let $H$ be a complex Hilbert space. A discrete family $T:E \to C(H)$ is \emph{self-adjoint of type (A)}, if 
	\begin{enumerate}
		\item There exists a dense subspace $Z \subset H$, such that all $e \in E$ satisfy $\dom T_e = Z$. We set $\dom T:=Z$. 
		\item For all $e \in E$, the operator $T_e$ is self-adjoint.
		\item There exists a norm $| \_ |$ on $Z$ such that, for all $e \in E$, the operator $T_e:(Z,| \_ |) \to (H, \| \_ \|_H)$ is bounded and the graph norm of $T_e$ is equivalent to $| \_ |$.
		\item $E$ is a topological space. 
		\item The map $E \to B(Z,H)$, $e \mapsto T_e$, is continuous.
	\end{enumerate}
\end{Def}

\begin{Thm}
	\label{ThmDiscFamA}
	Let $T:E \to C(H)$ be a discrete self-adjoint family of type (A). For any $e_0 \in E$ and any $\varepsilon > 0$ there exists an open neighborhood $U \subset E$ of $e_0$ such that 
	\begin{align*}
		\forall e \in U: \exists k \in \Z: \forall j \in \Z: d_a(\spc^{e_0}_T(j),\spc^e_T(j+k)) < \varepsilon.
	\end{align*}
\end{Thm}

\begin{Prf}
	Let $\varepsilon > 0$ and $e_0 \in E$ be arbitrary. As in \cref{DefDiscFamA}, let $\| \_ \|_H$ be the norm in $H$, $Z:=\dom T$, and let $\| \_ \|_Z$ be the graph norm of $T_{e_0}$ on $Z$. Finally, let $\| \_ \|$ be the associated operator norm in $B(Z,H)$ (which is then also equivalent to the operator norm induced by $| \_ |$).
	
	\begin{steplist}
	\step[setup and strategy]
		By construction, for any $e_1 \in E$ 
		\begin{align*}
			D_{e_1}(\zeta) := \zeta T_{e_1} + (1-\zeta) T_{e_0} = T_{e_0} + \zeta (T_{e_1} - T_{e_0}), && \zeta \in \C,
		\end{align*}
		defines a discrete self-adjoint holomorphic family $D_{e_1}:\C \to C(H)$ of type (A) with domain $Z$. The idea is to prove the claim using \cref{CorAnalyticGlobalGrowth}. The only problem is that for any two $e_1, e_2 \in E$, the families $D_{e_1}$ and $D_{e_2}$ are different. Hence their constants $C_{I,e_1}$, $C_{I,e_2}$ from \cref{ThmHolAFamilyDerGrowth} for the interval $I$ could differ. Consequently, their associated deltas $\delta_{e_1}$, $\delta_{e_2}$ from \cref{CorAnalyticGlobalGrowth} could also differ. We will show that there exists an open neighborhood $U$ around $e_0$ sufficiently small such that for all $e_1 \in U$, the $\delta = \delta_{e_1}$ is $\geq 1$ , if $\zeta_0$ is always set to $\zeta_0:=t_0:=0$. This will prove the claim.
	\step[preliminary estimate] 
		Recall that by \eqref{EqHolAFamilyDerGrowthDelta} there are $C_1,C_2 > 0$ such that
		\begin{align*} 
			\delta_{e_1} = C_{I,e_1}^{-1} \ln(\min(C_1, \varepsilon C_2)  + 1 ).
		\end{align*} 
		Since $\lim_{t \to 0}{e^t} = 1$, there exists $\varepsilon_1 > 0$ such that 
		\begin{align} \label{PreqDiracSpecContArsinhExpCont}
			\forall |t| \leq 2 \varepsilon_1: \exp(t) - 1 \leq \min(C_1, \varepsilon C_2).
		\end{align}
	\step[construction of $U$]
		Since $T$ is discrete and self-adjoint of type (A), the map $E \to B(Z,H)$, $e \mapsto T_{e}$, is continuous. Consequently, there exists an open neighborhood $U$ of $e_0$ such that
		\begin{align} \label{EqPreqDiracSpecContArsinhDCont}
			\forall e_1 \in U: \|T_{e_1} - T_{e_0}\| < \min\left( \tfrac{1}{2}, \varepsilon_1 \right).
		\end{align}
		Now for any $e_1 \in U$, $t \in [0,1]$ and  $\varphi \in Z$ we have
		\begin{align*}
			\|D_{e_1}(t) \varphi \|_H
			\geq \|T_{e_0} \varphi \|_H - \| T_{e_1} - T_{e_0} \| \|\varphi \|_Z.
		\end{align*}
		Therefore, applying  \eqref{EqHolAFamilyDerGrowthConstantSpec} to $D_{e_1}$,  we obtain 
		\begin{align*}
			\alpha_{I,e_1}
			&=\inf_{t \in I}{\inf_{\|\varphi\|_Z=1}}{\|\varphi\|_H + \|D_{e_1}(t) \varphi\|_H} 
			\geq 1 - \| T_{e_0} - T_{e_0} \| 
			\jeq{\eqref{EqPreqDiracSpecContArsinhDCont}}{>} \tfrac{1}{2}, \\
			\beta_{I,e_1} &= \sup_{t \in I}{\|D'_{e_1}(t)\|} = \|T_{e_1} - T_{e_0}\| < \varepsilon_1.
		\end{align*}
		Altogether, we achieved for any $e_1 \in U$
		\begin{align*}
			C_{I,e_1} = \alpha_{I,e_1}^{-1} \beta_{I,e_1} < 2 \varepsilon_1.
		\end{align*}
		By \eqref{PreqDiracSpecContArsinhExpCont}, this implies
		\begin{align*}
			\exp(C_{I,e_1}) - 1 \leq \min(C_1, \varepsilon C_2)
			\Longrightarrow \delta_{e_1} = C_{I,e_1}^{-1}\ln(\min(C_1, \varepsilon C_2)+1) \geq 1,
		\end{align*}
	\end{steplist}
		which proves the claim.	
\end{Prf}

Finally, we apply all our results to Dirac operators.

\begin{Thm}
	\label{ThmSpecContArsinh}
	The map
	\begin{align*}
		\bar \spc:(\Rm(M),\mathcal{C}^1) \to (\Conf,\bar d_a), \; \; g \mapsto \overline{\spc}^g,
	\end{align*}	
	is continuous. 
\end{Thm}

\begin{Prf}
	Let $g_0 \in \Rm(M)$ and $\varepsilon > 0$ be arbitrary. By definition of $\bar d_a$, cf. \eqref{EqDefbarda}, it suffices to find an open neighborhood $U \subset \Rm(M)$ such that
	\begin{align} \label{EqSpecContArsinh}
		\forall g' \in U: \exists k \in \Z: \forall j \in \Z: d_a(\spc^{g}(j), \spc^{g'}(j+k)) < \varepsilon.
	\end{align}
	By \cref{ThmSpinorIdentification}, the map $\Rm(M) \to B(H^1(\Sigma_{\C}^{g_0}M),L^2(\Sigma_{\C}^{g_0}M))$, $h \mapsto \Dirac^h_{g_0}$, is a discrete self-adjoint family of type (A). Consequently, by \cref{ThmDiscFamA}, there exists $U$ such that \eqref{EqSpecContArsinh} holds. 
\end{Prf}

To apply the Lifting Theorem (\cref{ThmLifting}), we quickly verify that $\pi$ is a covering map.

\begin{Thm}
	\label{ThmPiCoveringArsinh}
	The map $\pi:(\Mon,d_a) \to (\Conf, \bar d_a)$ is a covering map with fibre~$\Z$. 
\end{Thm}

\begin{figure}[t] 
	\begin{center}
		\input{fig.cont}
		\caption[An evenly covered neighborhood for $u$.]{An evenly covered neighborhood for $u$.}
		\label{FigCont}
	\end{center}
\end{figure}
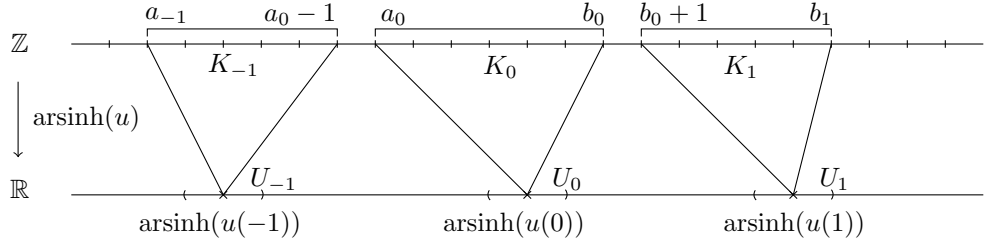

\begin{Prf}
	In this proof we also use \cref{NotTechnical}. By the definition of the group action, see \eqref{EqDeftau}, $\Z$ acts on $\Mon$ by isometries. In particular, the group action is continuous. We will show that, for each $u \in \Mon$, there exists an open neighborhood $V$ such that
	\begin{align} \label{EqTauPropDisc}
		\pi^{-1}(\pi(V)) = \dot{\bigcup}_{j \in \Z}{V.j}.
	\end{align}
	To see this, note that the function $\ash \circ u$ is non-decreasing and proper. The set $K_0:=(\ash \circ u)^{-1}(\ash(u(0)))$ is of the form $K_0 = \{a_0, \ldots, b_0\}$ for some $a_0 \leq b_0$, $a_0,b_0 \in \Z$. For the same reason, there exist $b_1$, $a_{-1} \in \Z$ such that (see also \Cref{FigCont})
	\begin{align*}
		(\ash \circ u)^{-1}(\ash(u(0))) &= \{a_0, \ldots, b_0\} = K_0, \\
		(\ash \circ u)^{-1}(\ash(u(b_0 + 1))) &= \{b_0+1, \ldots, b_1 \} =: K_1, \\
		(\ash \circ u)^{-1}(\ash(u(a_0 - 1))) &= \{ a_{-1}, \ldots, a_0 - 1\} =:K_{-1}.
	\end{align*}
	Since $\ash(u(\Z))$ is discrete, there exists $\varepsilon > 0$ such that 
	\begin{align*}
		I_{\varepsilon}(\ash(u(0))) \cap I_{\varepsilon}(\ash(u(b_0+1))) = \emptyset, &&
		I_{\varepsilon}(\ash(u(0))) \cap I_{\varepsilon}(\ash(u(a_0-1))) = \emptyset, 
	\end{align*}
	Thus, we obtain open sets
	\begin{align*}
		U_0 := I_{\varepsilon}(\ash(u(0))), && 
		U_1 := I_{\varepsilon}(\ash(u(b_0 + 1))), && 
		U_{-1} := I_{\varepsilon}(\ash(u(a_0-1))), && 
	\end{align*}
	which are mutually disjoint. To see that $V := B_\varepsilon(u)$ satisfies \eqref{EqTauPropDisc}, suppose to the contrary that there exists  $v \in V$ and $j \in \Z$ such that $v.j \in V$. Assume $j > 0$ (the proof for $j<0$ is entirely analogous). By hypothesis, this implies that $\ash(v(b_0)) \in U_0$ and $\ash(v(b_0 + j)) = \ash((v.j(b_0))) \in U_0 $. But $\ash \circ v$ is non-decreasing, so $\ash(v(b_0 + j)) \geq \ash(v(b_0 + 1)) \in U_1$. This implies that $\ash(v(b_0+j)) \notin U_0$, which is a contradiction. \\
	Finally, to see that $\pi$ is a covering map, let $u \in [u] \in (\Conf,\bar d_a)$ be arbitrary. Let $V$ be an open neighborhood of $u$ satisfying \eqref{EqTauPropDisc}. Then, $\bar V := \pi(V)$ is evenly covered. The map $\pi$ is open due to \cref{EqTauPropDisc} and the fact that $\Z$ acts by isometries. Thus, $\pi$ is a covering map.
\end{Prf}

%% file: fig.cones.tex
\begin{tikzpicture}
	\draw (3,0) -- (-1.5,0) node[left] {};
	\draw (-2,0) -- (-6,0) node[left] {};
	\draw (3,-2) -- (-1.5,-2) node[left] {};
	\draw (-2,-2) -- (-6,-2) node[left] {};

	\draw (-5,-2) -- (-4,0) -- (-3,-2) {};
	
	\draw (-1,-2) -- (0,0) -- (1,-2) {};
	
	\draw (0,-2) -- (1,0) -- (2,-2) {};	
	
	\node [label={[xshift=-185, yshift=-12]$t_0$}]  {};
	\node [label={[xshift=-185, yshift=-66]$t$}]  {};
	\node [label={[xshift=-113, yshift=-5]$\operatorname{a}(\lambda_{i}) = \operatorname{a}(\lambda_j)$}]  {};
	\node [label={[xshift=0, yshift=-5]$\operatorname{a}(\lambda_{i})$}]  {};
	\node [label={[xshift=30, yshift=-5]$\operatorname{a}(\lambda_{j})$}]  {};
	\node [label={[xshift=140, yshift=-35]$\operatorname{a} = \operatorname{arsinh}$}]  {};
	
	\node [label={[xshift=-130, yshift=-78]$\operatorname{a}(\mu_{\sigma_n(j)})$}] {};
	\node [label={[xshift=-85, yshift=-78]$\operatorname{a}(\mu_{\sigma_n(i)})$}] {};
	\node [label={[xshift=-3, yshift=-78]$\operatorname{a}(\mu_{\sigma_n(j)})$}] {};
	\node [label={[xshift=40, yshift=-78]$\operatorname{a}(\mu_{\sigma_n(i)})$}] {};
	
	\foreach \m in {-4,0,1}{
		\draw (\m-0.05,-0.05) -- (\m+0.05,+0.05) {};
		\draw (\m-0.05,0.05) -- (\m+0.05,-0.05) {};
	}
	
	\foreach \m in {-4.5,-3.5,0.2,0.8}{
		\draw (\m-0.05,-2-0.05) -- (\m+0.05,-2+0.05) {};
		\draw (\m-0.05,-2+0.05) -- (\m+0.05,-2-0.05) {};
	}

\end{tikzpicture}

%% file: fig.cont.tex
\begin{tikzpicture}
	\draw (6,0) -- (-6,0) node[left] {};
	\draw (6,-2) -- (-6,-2) node[left] {};
	
	\draw[->] (-6.7,-0.5) -- (-6.7,-1.5) {};
		
	\foreach \z in {-5.5, -5, ..., 5.5}{
		\draw (\z,-0.05) -- (\z,+0.05) {};
	}
	
	\foreach \m in {-4, 0, 3.5}{
		\draw (\m-0.05,-2-0.05) -- (\m+0.05,-2+0.05) {};
		\draw (\m-0.05,-2+0.05) -- (\m+0.05,-2-0.05) {};
	}
	
	\draw (-5,0) -- (-4,-2) -- (-2.5,0) {};
	\draw (-2,0) -- (0,-2) -- (1,0) {};
	\draw (1.5,0) -- (3.5,-2) -- (4,0) {};	
	
	\draw (-5,0.1) -- (-5,0.2) -- (-2.5,0.2) -- (-2.5,0.1) {};
	\draw (-2,0.1) -- (-2,0.2) -- (1,0.2) -- (1,0.1) {};
	\draw (1.5,0.1) -- (1.5,0.2) -- (4,0.2) -- (4,0.1) {};
		
	\node [label={[xshift=-190, yshift=-10]$\mathbb{Z}$}]  {};
	\node [label={[xshift=-190, yshift=-65]$\mathbb{R}$}]  {};
	\node [label={[xshift=-165, yshift=-40]$\ash(u)$}]  {};
	
	\node [label={[xshift=-135, yshift=0]$a_{-1}$}]  {};
	\node [label={[xshift=-85, yshift=0]$a_{0}-1$}]  {};
	\node [label={[xshift=-50, yshift=0]$a_{0}$}]  {};
	\node [label={[xshift=25, yshift=0]$b_{0}$}]  {};
	\node [label={[xshift=55, yshift=0]$b_{0}+1$}]  {};
	\node [label={[xshift=110, yshift=0]$b_{1}$}]  {};
	
	\node [label={[xshift=-110, yshift=-20]$K_{-1}$}]  {};
	\node [label={[xshift=-10, yshift=-20]$K_{0}$}]  {};
	\node [label={[xshift=80, yshift=-20]$K_{1}$}]  {};
	
	\node [label={[xshift=-115, yshift=-80]$\ash(u(-1))$}]  {};
	\node [label={[xshift=-5, yshift=-80]$\ash(u(0))$}]  {};
	\node [label={[xshift=100, yshift=-80]$\ash(u(1))$}]  {};
	
	\node [label={[xshift=-95, yshift=-63]$U_{-1}$}]  {};
	\node [label={[xshift=15, yshift=-62]$U_{0}$}]  {};
	\node [label={[xshift=115, yshift=-62]$U_{1}$}]  {};
	
	\draw (-4.5,-1.94) arc (145:225:0.1) (-4.5,-2.05) ;
	\draw (-3.5,-2.06) arc (-40:40:0.1) (-3.5,-2.10) ;
	\draw (-0.5,-1.94) arc (145:225:0.1) (-0.5,-2.05) ;
	\draw (0.5,-2.06) arc (-40:40:0.1) (0.5,-2.10) ;
	\draw (3,-1.94) arc (145:225:0.1) (3,-2.05) ;
	\draw (4,-2.06) arc (-40:40:0.1) (4,-2.10) ;
	
\end{tikzpicture}

%% file: evpaper.specfl.tex
\section{Moduli Spaces and Spectral Flow}
\label{SecModuliSpecFl}

Let $\Diff(M)$ be the diffeomorphism group of $M$. This group acts canonically on the space of Riemannian metrics via 
	\DefMap{\Rm(M) \times \Diff(M)}{\Rm(M)}{(g,f)}{f^*g.} 
For any subgroup $G \subset \Diff(M)$, the quotient space $\Rm(M) / G$ is called a \emphi{moduli space}. We want to investigate when the map $\bar \spc$ (respectively $\widehat{\spc}$) from \cref{MainThmSpec} (and hence the family of functions $(\lambda_j)_{j \in \Z}$ from \cref{MainThmFun}) descends to a continuous map on a moduli space. More precisely, we ask whether or not, there exists a diagram
\begin{align}
	\label{EqPassModuliSpace}
	\xymatrixcolsep{3.5em}
	\begin{split}
		\xymatrix{
			\Rm(M)
				\ar[r]^-{\widehat{\spc}}
				\ar@{->>}[d]
				\ar[dr]^/-1.5em/{\overline{\spc}}|{\phantom{*}}
			&\Mon
				\ar[d]^-{\pi}
			\\
			\Rm(M) / G
				\ar@{..>}[r]_-{\exists ?}
				\ar@{..>}[ur]^/-1.5em/{\exists ?}
			&\Conf
		}
	\end{split}
\end{align}
and if the dashed map goes to $\Conf$ or even to $\Mon$. We give an answer to this question, when $G$ is one of the groups 
\begin{align*}
	\Diff^+(M), &&
	\Diff^0(M), &&
	\Diff^{\spin}(M), &&
	\ker \specfl.
\end{align*}
Here, $\Diff^+(M)$ denotes the subgroup of orientation-preserving diffeomorphisms and $\Diff^0(M)$ are the diffeomorphisms that are isotopic to the identity. The group $\Diff^{\spin}(M)$ was defined in \cref{DefSpinDiffeo}. The group $\ker \specfl$ is the kernel of the \emph{spectral flow}, see \cref{RemKerSpecfl} for details. The results are collected in \Cref{FigPassing}. We inserted a ``$\checkmark$'' into the column $\Mon$, if the arrow in \cref{EqPassModuliSpace} always exists and goes to $\Mon$, and a ``$\times$'' if this is not the case in general. In the column $\Conf$, we use the same notation. Since
\begin{align*}
	\Diff^0(M) \subset \ker \specfl \subset \Diff^{\spin}(M) \subset \Diff^+(M),
\end{align*}
the overall picture is not surprising. However, proving ``$\times$'' for $\Diff^{\spin}(M)$ turns out to be quite difficult, see \cref{MainThmFlow}.

\begin{figure}[t]
	\begin{center}
		\begin{tabular}{|c|c|c|c|}
		\hline
		$G$ & $\Mon$ & $\Conf$ & \text{Details}\\
		\hline
		$\Diff^0(M)$ & $\checkmark$ &  $\checkmark$ &  \cref{ThmPassDiff0} \\
		$\ker \specfl$ & $\checkmark$ & $\checkmark$ &  \cref{RemKerSpecfl} \\
		$\Diff^{\spin}(M)$ & $\times$ & $\checkmark$ & \cref{MainThmFlow}, \cref{ThmPassDiffSpin} \\
		$\Diff^+(M)$ & $\times$ & $\times$ & \cref{RemPassSPecflTorus}\\
		\hline
		\end{tabular}
	\end{center}
	\caption[Spectral functions and moduli spaces.]{Does the ordered spectral function pass to a moduli space?}
	\label{FigPassing}
\end{figure}

The preliminary considerations in \cref{RemSpinIsometries} immediately imply the following.

\begin{Thm} 
	\label{ThmPassDiffSpin}
	There exists a commutative diagram
	\begin{align}
		\label{EqDiffSpinPass}
		\begin{split}
			\xymatrix{
				\Rm(M) 
					\ar[d]
					\ar[r]^-{\overline{\spc}} 
				& (\Conf, \bar d_a) 
				\\
				\Rm(M) / \Diff^{\spin}(M)
					\ar@{..>}[ur]^-{\exists \spc^{\spin}}
			}
		\end{split}
	\end{align}	
\end{Thm}

\nomenclature[spcspin]{$\spc^{\spin}$}{passing of $\overline{s}$ to $\Rm(M)/\Diff^{\spin}(M)$}
\nomenclature[spc0]{$\spc^{0}$}{passing of $\widehat{s}$ to $\Rm(M)/\Diff^{0}(M)$}

With only a little more work, we get an even stronger statement for $\Diff^0(M)$.

\begin{Thm} 
	\label{ThmPassDiff0}
	There exists a commutative diagram 
	\begin{align}
		\label{EqPassDiff0Mon}
		\begin{split}
			\xymatrix{
				\Rm(M)
					\ar[d]
					\ar[r]^-{\widehat{\spc}} 
				& (\Mon, d_a)
				\\
				\Rm(M) / \Diff^0(M) 
					\ar@{..>}[ur]^-{\exists \spc^0}
			}
		\end{split}
	\end{align}
\end{Thm}

\begin{Prf}  
	The claim will follow from the universal property of the topological quotient, if we can show that 
	\begin{align*}
		\forall g_0 \in \Rm(M): \forall f \in \Diff^0(M): \widehat{\spc}^{g_0} = \widehat{\spc}^{f^*g_0}.
	\end{align*}
	Let $h:I \times M \to M$ be an isotopy from $h_0 = \id$ to $h_1 = f$. Clearly, $\id \in \Diff^{\spin}(M)$, thus we obtain a lift 
	\begin{align} 
		\label{EqSpinLiftfF}
		\begin{split}
			\xymatrix{
				I \times \GLtp M 
					\ar[d]^-{\id \times \Theta}
					\ar@{..>}[r]^-{\exists \tilde H} 
				& \GLtp M \ar[d]^-{\Theta}
				\\
				I \times \GLp M 
					\ar[d]
					\ar[r]^-{H}
				& \GLp M
					\ar[d] \\
				I \times M
					\ar[r]^-{H} 
				& M.
			}
		\end{split}
	\end{align}
	by \cref{ThmSpinIsotopyInvariance}. Specializing to any $t \in I$
	\begin{align*}
		\xymatrix{
			\Spin^{g_t} M
				\ar[d]^-{\Theta^{g_t}}
				\ar[r]^-{\tilde H_t}
			& \Spin^{g_0} M 
				\ar[d]^-{\Theta^{g_0}}
			\\
			\SO^{g_t} M
				\ar[d]
				\ar[r]^-{(H_t)_*}
			& \SO^{g_0} M 
				\ar[d]
			\\
			(M,g_t) 
				\ar[r]^-{H_t}
			& (M,g_0)
		}
	\end{align*}
	we obtain that for all $t \in I$, the map $h_t$ is a spin isometry in the sense of \eqref{EqDefSpinIsometry}. Therefore, $(M,g_t)$ and $(M,g_0)$ are Dirac isospectral for all $t \in I$ by \cref{RemSpinIsometries}. This implies $\widehat{\spc}^{g_0} = \widehat{\spc}^{g_1}$.
\end{Prf}

\begin{Rem}[a counter-example on the torus]
	\label{RemPassSPecflTorus}
	It remains to discuss the group $\Diff^+(M)$, and one might ask if \eqref{EqDiffSpinPass} still holds, if $\Diff^{\spin}(M)$ is replaced by $\Diff^+(M)$. This is false in general. A counter-example is provided by the standard torus $\T^3 = \R^3 / \Z^3$ equipped with the induced Euclidean metric $\bar g$. It is well known that the (equivalence classes of) spin structures on $\T^3$ are in one-to-one correspondence with tuples $\delta \in \Z_2^3$, see for instance \cite{FriedTori}. We denote by $\Spin^{\bar g}_{\delta} \T^3$ the spin structure associated to $\delta$. The map	
	\begin{align*}
		f := 
		\begin{pmatrix}
			1 & 1 & 0\\
			0 & 1 & 0 \\
			0 & 0 & 1 \\
		\end{pmatrix} :\R^3 \to \R^3
	\end{align*}	
	preserves $\Z^3$ and satisfies $\det(f) = 1$. Hence it induces a diffeomorphism $\bar f \in \Diff^+(\T^3)$. One checks that there is a commutative diagram
	\begin{align*}
		\xymatrix{
			\Spin_{(1,1,0)}^{\bar g} \T^3 
				\ar[r]^-{F}
				\ar[d]^{\Theta^{\bar g}} 
			&\Spin^{\bar f^* \bar g}_{(1,0,0)} \T^3
				\ar@{..>}[r]^-{\nexists}
				\ar[d]^{\Theta^{\bar f^* \bar g}}
			& \Spin_{(1,0,0)}^{\bar g} \T^3
				\ar[d]^{\Theta^{\bar g}}
			\\
			\SO^{\bar g} \T^3 
				\ar[r]^-{\bar f^{-1}_*}
				\ar[d]
			&\SO^{\bar f^* \bar g} \T^3
				\ar[r]^-{\bar f_*}
				\ar[d]
			& \SO^{\bar g} \T^3
				\ar[d]
			\\
			(\T^3, \bar g)
				\ar[r]^-{\bar f^{-1}}
			&(\T^3,\bar f^* \bar g)
				\ar[r]^-{\bar f}
			& (\T^3, \bar g). 
		}
	\end{align*}
	The map in the right upper row cannot exist; if it did, the spin structures corresponding to $(1,1,0)$ and $(1,0,0)$ would be equivalent. The left part of the above diagram is a spin isometry analogous to \eqref{EqDefSpinIsometry}. Therefore $\Dirac^{f^* \bar g}_{(1,0,0)}$ and $\Dirac^{\bar g}_{(1,1,0)}$ are isospectral, but the spectra of $\Dirac^{\bar g}_{(1,1,0)}$ and $\Dirac^{\bar g}_{(1,0,0)}$ are already different as a set. This follows from the explicit computation of the spectra of Euclidean tori, see also \cite{FriedTori}. Consequently, $\spec \Dirac^{\bar f^*g}_{(1,0,0)} \neq \spec \Dirac^{\bar g}_{(1,0,0)}$, and no diagram analogous to \eqref{EqDiffSpinPass} can exist for $\Diff^+(\T^3)$.
\end{Rem}

\begin{Rem}
	Notice that in \eqref{EqPassDiff0Mon} the map $\spc^0$ goes from the moduli space for $\Diff^0(M)$ to $\Mon$, whereas in \eqref{EqDiffSpinPass} the corresponding map $\spc^{\spin}$ goes to $\Conf$. Therefore one might ask, if one could improve \eqref{EqDiffSpinPass} by lifting $\spc^{\spin}$ to a map $\widehat{\spc}^{\spin}$ such that
\begin{align}
	\label{EqLiftSpinMonDef}
	\begin{split}
		\xymatrix{
			& \Mon
				\ar[d]
			\\
			\Rm(M) / \Diff^{\spin}(M)
				\ar[r]^-{\spc^{\spin}}
				\ar@{..>}[ur]^-{\widehat{\spc}^{\spin}}
			& \Conf
			}
	\end{split}
\end{align}
	commutes. This question is not so easy to answer, and the rest of this section is devoted to the proof that this is not possible in general. To see where the problem lies, it will be convenient to introduce the following terminology.
\end{Rem}

\begin{Lem}[spectral flow] $ $
    \label{LemSpecFlowDef} \index{spectral flow}
	\nomenclature[sfgf]{$\specfl_g(f)$}{spectral flow of $f$ at $g$}
	\nomenclature[sf]{$\specfl(f)$}{$\Rm(M) \to \Z$ spectral flow of $f$}
	\nomenclature[sfg]{$\specfl(\mathbf{g})$}{spectral flow of $\mathbf{g}$}
    \begin{enumerate}
        \item 
            For any $f \in \Diff^{\spin}(M)$ and $g \in \Rm(M)$, there exists a unique $\specfl_g(f) \in \Z$ such that
            \begin{align}
                \label{EqDefSpectralFlow}
                \forall j \in \Z: \widehat{\spc}^{g}(j) = \widehat{\spc}^{f^*g}(j - \specfl_g(f)).
            \end{align}
            The induced map $\specfl(f):\Rm(M) \to \Z$ is called the \emph{spectral flow of $f$}. 
        \item
            Let $\mathbf{g}:[0,1] \to \Rm(M)$, $t \mapsto g_t$, be a continuous path of metrics. Let $\spc:\Rm(M) \to \Mon$ be the ordered spectral function for the associated Dirac operators. Take a lift $\widehat{\spc}:\Rm(M) \to \Mon$ of $\overline{\spc}$ such that $\widehat{\spc}^{g_0} = \spc^{g_0}$ as in \eqref{EqMainThmSpc}. There exists a unique integer $\specfl(\mathbf{g}) \in \Z$ such that
            \begin{align*}
                \forall j \in \Z: \widehat{\spc}^{g_1}(j) = \spc^{g_1}(j+\specfl(\mathbf{g})).
            \end{align*}
            The integer $\specfl(\mathbf{g})$ is called the \emph{(Dirac) spectral flow along $\mathbf{g}$}. 
        \item
            \label{ItSpecFlowHomotpyInvariant}
            For any $f \in \Diff^{\spin}(M)$ and any family $\mathbf{g}$ joining $g_0$ and $f^* g_0$, we have $\specfl_{g_0}(f)=\specfl(\mathbf{g})$.
    \end{enumerate}
\end{Lem}

\begin{Prf} $ $
    \begin{enumerate}
        \item 
            By \cref{ThmPassDiffSpin}, the map $\overline{\spc}:\Rm(M) \to \Conf$ descends to a quotient map 
			\begin{align*}
				\spc^{\spin}:\Rm(M) / \Diff^{\spin}(M) \to \Conf.
			\end{align*}
			This precisely means that $\widehat{\spc}^g$ and $\widehat{\spc}^{f^*g}$ are equal in $\Conf$. By the definition of $\Conf$, this implies the existence of $\specfl_g(f)$ as required. 
        \item 
            This follows directly from \eqref{EqMainThmSpc} and the fact that $\spc$ and $\widehat{\spc}$ are equal in $\Conf$.
        \item Set $g_1:=f^*g_0$, and let $\mathbf{g}$ be a family of metrics joining $g_0$ and $g_1$. By \cref{ThmPassDiffSpin}, we obtain $\spc^{g_0} = \spc^{g_1}$. Take a lift $\widehat{\spc}$ satisfying $\widehat{\spc}^{g_0} = \spc^{g_0}$. This implies for all $j \in \Z$
			\begin{align*}
				\spc^{g_0}(j)
				&=\widehat{\spc}^{g_0}(j)
				=\widehat{\spc}^{g_1}(j - \specfl_{g_0}(f)) 
				=\spc^{g_1}(j - \specfl_{g_0}(f) + \specfl(\mathbf{g})) \\
				&=\spc^{g_0}(j - \specfl_{g_0}(f) + \specfl(\mathbf{g})),
			\end{align*}
			which implies $\specfl(\mathbf{g}) - \specfl_{g_0}(f)  = 0$, since $\spc^{g_0}$ is non-decreasing and all eigenvalues are of finite multiplicity.
    \end{enumerate}
\end{Prf}

\begin{Rem}[spectral flow]
    \label{RemSpecFlowPhillips}
	Intuitively, the spectral flow $\specfl(\mathbf{g})$ of a path $\mathbf{g}:[0,1] \to \Rm(M)$ counts the signed number of eigenvalues of the associated path $\Dirac^{g_t}$ of Dirac operators that cross $0$ from below when $t$ runs from $0$ to $1$. The sign is positive, if the crossing is from below, whereas it is negative, if the crossing is from above.\\
    The concept of spectral flow is well known in other contexts. A good introduction can be found in a paper by Phillips, see \cite{phillips}. Phillips introduces the spectral flow for continuous paths $[0,1] \to \mathcal{F}_*^{\sa}$, where $\mathcal{F}_*^{\sa}$ is the non-trivial component of the space of self-adjoint Fredholm operators on a complex separable Hilbert space $H$. In this general setup, the definition of spectral flow is a little tricky, see \cite[Prop. 2]{phillips}. But for paths of Dirac operators, it coincides with the definition given in \cref{LemSpecFlowDef} above (by \cref{ThmSpinorIdentification} we can think of all Dirac operators $\Dirac^{g_t}$, $t \in [0,1]$, of a path $\mathbf{g}$ as defined on the same Hilbert space). Therefore, we have found a convenient alternative for describing the spectral flow in this case using the continuous function $\widehat{\spc}$. 
\end{Rem}

\begin{Rem}
	\label{RemKerSpecfl}
    By \cite[Prop. 3]{phillips}, the spectral flow of a path of operators depends only on the homotopy class of the path. Consequently, since $\Rm(M)$ is simply connected, $\specfl(\mathbf{g})$ depends only on $g_0$ and $g_1$. It follows that $\specfl:\Diff^{\spin}(M) \to \Z$ is a group homomorphism.
	The map $\widehat{\spc}$ certainly descends to
	\begin{align*}
		\xymatrix{
			\Rm(M)
				\ar[r]^-{\widehat{\spc}}
				\ar[d]
			& \Mon
			\\
			\Rm(M) / \ker \specfl
				\ar@{..>}[ur]^-{\exists}
		}
	\end{align*}
	and $\ker \specfl$ is the largest subgroup of $\Diff^{\spin}(M)$ with this property. Rephrased in these terms, we conclude that the map $\spc^{\spin}$ from \eqref{EqDiffSpinPass} lifts to a map $\widehat{\spc}^{\spin}$ as in \eqref{EqLiftSpinMonDef} if and only if $\specfl(f)=0$ for all $f \in \Diff^{\spin}(M)$. Consequently, we must show the following theorem.%
\end{Rem}

\begin{Thm}
	\label{MainThmFlow}
	There exists a spin manifold $(M,\Theta)$ and a diffeomorphism $f \in \Diff^{\spin}(M)$ such that $\specfl(f) \neq 0$. 
\end{Thm}

\begin{Prf}
	The general idea is to use a recent result from differential topology, which implies the existence of a fibre bundle $P \to S^1$ such that $P$ is spin, $\widehat{A}(P) \neq 0$, and $P$ has a $2$-connected fibre type $M$. Setting $S^1 = [0,1] / \{0 \sim 1\}$, we can view $M$ as a fibre $M=P_{[0]}$. The bundle $P$ will be isomorphic to a bundle $P_f$, obtained from the trivial bundle $[0,1] \times M \to [0,1]$ by identifying $(1,x)$ with $(0,f(x))$, $x \in M$, for a suitable diffeomorphism $f \in \Diff^{\spin}(M)$. A lift of $f$ to a spin morphism yields a spin structure on $P_f$. To show that $f$ has nontrivial spectral flow, we will cut open the bundle along $[0]$, obtain a trivial fibre bundle $[0,1] \times M \to [0,1]$, and glue in two infinite half-cylinders at both sides (see \Cref{FigCut}). This gives a bundle of the form $\R \times M \to \R$, and in each fibre we get a Dirac operator. Using various index theorems, we will show that the $\widehat{A}$-genus of $P$ equals the index of $\R \times M$, which in turn equals the spectral flow of the associated family of Dirac operators in the fibres, which finally equals the spectral flow of $f$. The technical details of this proof rely on several other theorems, which are collected in \cref{SectEvpaperAppendixFundRes} (one might want to take a look at these first). All Dirac operators will be complex, so we drop the subscript $\C$ in notation.
	
	\begin{steplist}
	\step[construct a bundle]
		Apply \cref{ThmHanke} to $(k,l) = (1,2)$ and obtain a fibre bundle $P \to S^1$ with some fibre type $M$, where $\dim P = 4n$, $n$ odd, and 
		\begin{align}
			\label{EqIndexHanke}
			\widehat{A}(P) \neq 0.
		\end{align}
		Since $M$ is $2$-connected, $M$ has a unique spin structure (up to equivalence). It follows that $m := \dim M = 4n-1 \equiv 3 \mod 4$, and also $m \equiv 3 \mod 8$, since $n$ is odd. Therefore, by \cref{ThmAmmannSurgery}, there exists a metric $g_0$ on $M$ such that the associated Dirac Operator $\Dirac^M$ is invertible. By \cref{LemClassBundlesS1}, $P$ is isomorphic to $P_f = [0,1] \times M / f$ for some $f \in \Diff^{\spin}(M)$. Define $g_1 := f^* g_0$, and connect $g_0$ with $g_1$ in $\Rm(M)$ by the linear path $g_t := t g_1 + (1-t) g_0$, $t \in [0,1]$. Endow $[0,1] \times M$ with the generalized cylinder metric $dt^2 + g_t$. Denote by $\pi:[0,1] \to S^1$ the canonical projection. We obtain a commutative diagram
		\begin{align*}
			\xymatrix{
				[0,1] \times M
					\ar[r]
					\ar[d]
				& P_f
					\ar[r]^-{\cong}
					\ar[d]
				& P
					\ar[d]
				\\
				[0,1]
					\ar[r]^-{\pi}
				& S^1
					\ar[r]^-{\id}
				& S^1.
				}
		\end{align*}
		By construction, we can push forward the metric $dt^2 + g_t$ on $[0,1] \times M$ to a metric on $P_f$, and then further to $P$ such that the above row consists of local isometries. The right map is actually an isometry along which we can pull back the spin structure on $P$ to $P_f$. This map is then a spin isometry, and therefore we will no longer distinguish between $P$ and $P_f$. The left map is an isometry, except that it identifies $\{0\} \times M$ with $\{1\} \times M$. Notice that, since $[0,1] \times M$ is simply connected, the  spin structure on $[0,1] \times M$ obtained by pulling back the spin structure on $P_f$ along $\pi$ is equivalent to the canonical product spin structure on $[0,1] \times M$.
		
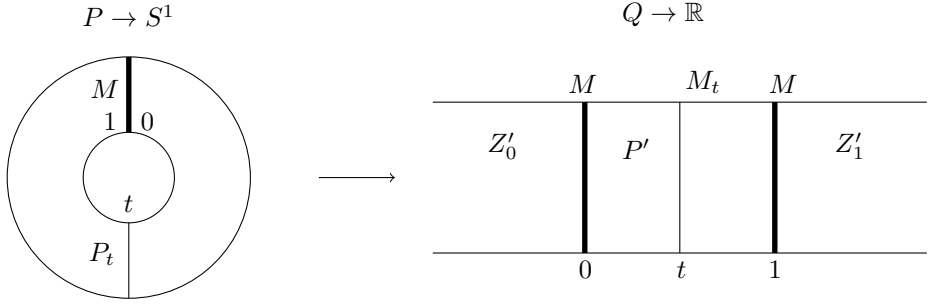
\begin{figure}[t] 
	\begin{center}
		\input{fig.cut}
		\caption[Proof of \cref{MainThmFlow}.]{$P$ is cut open at $[0]$, and we obtain a bundle $P' \to [0,1]$. Then we glue in two half-cylinders $Z_0'$ and $Z_1'$ to obtain a bundle $Q \to \R$.}
		\label{FigCut}
	\end{center}
\end{figure}		

	\step[trivialize]
		The Riemannian manifold $P':= ([0,1] \times M,dt^2 + g_t)$ has two isometric boundary components. Geometrically, $P'$ is obtained from $P$ by cutting $M=P_{[0]}$ out of $P$ and adding two boundaries $P'_0$ and $P'_1$, i.e. $P' = (P \setminus P_{[0]}) \coprod P'_0 \coprod P'_1$, where $P'_0$, $P'_1$ are two isometric copies of $P_{[0]}$. By \cref{ThmBB}, we obtain
		\begin{align}
			\label{EqIndexTrivializedBundle}
			\ind(\Dirac_+^P) = \ind(\Dirac_+^{P'}).
		\end{align}		
	
	\step[index of half-cylinders]
		Now, set
		\begin{align*}
			\begin{array}{ll}
				Z_0' := (\mathopen{]} -\infty, 0 \mathclose{]} \times M, dt^2 + g_0), 
				&Z_1' := (\mathopen{[} 1,\infty \mathclose{[} \times M, dt^2 + g_1), \\
				Z_1'' := (\mathopen{[} 0,\infty \mathclose{[} \times M, dt^2 + g_0), 
				&Z := (\R \times M, dt^2+g_0), \\
				Z' := Z_0' \textstyle \coprod Z_1', 
				&Z'' := Z_0' \textstyle \coprod Z_1''.
			\end{array}
		\end{align*}
		Since $Z$ is a Riemannian product, it follows that
		\begin{align}
			\label{EqInvertibleOnProduct}
			\forall \psi \in \Gamma_c(\Sigma Z): \| \Dirac^Z \psi \|_{L^2(\Sigma Z)}^2 \geq \lambda_{\min}^2 \|\psi\|_{L^2(\Sigma Z)}^2,
		\end{align}
		where $\lambda_{\min}$ is the eigenvalue of $\Dirac^{g_0}$ of minimal absolute magnitude. By construction, $\Dirac^{g_0}$ is invertible, thus $\lambda_{\min} > 0$. Therefore $\Dirac^{Z}$ is invertible and \emph{coercive at infinity} (see \cref{ThmBB} for the definition). By \cref{ThmBB}, this implies
		\begin{align}
			\begin{split}
				\label{EqIndexHalfZylinderZero}
				0
				&= \ind(\Dirac^{Z}_+)
				= \ind(\Dirac^{Z''}_+)
				= \ind(\Dirac^{Z_0'}_+) + \ind(\Dirac^{Z_1''}_+) \\
				&= \ind(\Dirac^{Z_0'}_+) + \ind(\Dirac^{Z_1'}_+)
				= \ind(\Dirac^{Z'}_+),
			\end{split}
		\end{align}
		where we used the fact that $Z_1'$ and $Z_1''$ are spin isometric.
		
	\step[glue in the half-cylinders]
		Now glue $Z'$ to $P'$ ($Z'_0$ at $\{0\} \times M$ and $Z'_1$ at $\{1\} \times M$) and obtain a bundle $Q =( \R \times M, dt^2 + g_t)$ where $g_t = g_0$ for $t \leq 0$ and $g_t = g_1$ for $t \geq 1$. Since $\Dirac^{g_1}$ is invertible as well, it follows that  $\Dirac^{Z'}$ satisfies an estimate analogous to \eqref{EqInvertibleOnProduct}. We therefore see that $\Dirac^{Q}$ is also coercive at infinity (take $K:=P'$ as the compact subset). By \cref{ThmBB}, we obtain
		\begin{align}
			\label{EqIndexHalfCylinderGlued}
			\ind(\Dirac_+^{Q})
			= \ind(\Dirac_+^{P'}) + \ind(\Dirac_+^{Z'})
			\jeq{\eqref{EqIndexHalfZylinderZero}}{=} \ind(\Dirac_+^{P'}).
		\end{align}		
		
	\step[apply hypersurface theory]
		Each $Q_t = \{t\} \times M$ is a hypersurface in $Q$, and $\tilde \partial_t \in \mathcal{T}(Q)$ (horizontal lift of $\partial_t$) provides a unit normal field for all $Q_t$. Therefore, we can apply some standard results about the Dirac operator on hypersurfaces, see \cite{BaerGaudMor}: Since $m$ is odd
		\begin{align*}
			\Sigma Q|_{Q_t} = \Sigma^+ Q|_{Q_t} \oplus \Sigma^- Q|_{Q_t} = \Sigma^+ M_t \oplus \Sigma^- M_t, && M_t := (M,g_t),
		\end{align*}
		where $\Sigma^+ M_t = \Sigma^- M_t = \Sigma M_t$ as Hermitian vector bundles. The Clifford multiplication ``$\cdot$'' in $\Sigma Q$ is related to the Clifford multiplication ``$\bullet_t^\pm$'' in $\Sigma^\pm M_t$ by $X \bullet_t^\pm \psi = \pm \tilde \partial _t \cdot X \cdot \psi$. Setting
		\begin{align}
			\label{EqDiracCopy}
			\tilde \Dirac^{M_t} := (\Dirac^{M_t} \oplus (-\Dirac^{M_t})),
		\end{align}
		we obtain the Dirac equation on hypersurfaces
		\begin{align}
			\label{EqDiracHyper}
				\tilde \partial_t \cdot \Dirac^{Q} = \tilde \Dirac^{M_t} +\tfrac{m}{2}H_t -  \nabla^{\Sigma Q}_{\tilde \partial_t}: \Gamma(\Sigma Q|_{Q_t}) \to  \Gamma(\Sigma Q|_{Q_t})
		\end{align}
		for all $t \in \R$. 
	
	\step[identification of the Spinor spaces]
		Let $\psi \in \Gamma(\Sigma Q)$ be a spinor field. For any $t \in \R$ this defines a section $\psi_t \in \Gamma(\Sigma Q|_{Q_t})$.  Therefore, we can also think of $\psi$ as a ``section''  of
		\begin{align*}
			\bigcup_{t \in \R}{\Gamma(\Sigma Q|_{Q_t})} \to \R.
		\end{align*}
		and \eqref{EqDiracHyper} tells us how $\tilde \partial_t \cdot \Dirac^Q$ acts on these sections under this identification. We would like to apply \cref{ThmSalamon}, and therefore have to solve the problem that for various $t$ the Hilbert spaces $L^2(\Sigma Q|_{Q_t})$ are different. As discussed in \cite{BaerGaudMor}, we will use the following identification: For any $x \in M$, consider the curve $\gamma^x:\R \to \R \times M$, $t \mapsto (t,x)$. Each spinor field $\psi \in \Gamma(\Sigma Q)$ determines a section $\psi^x$ in $\Sigma Q$ along $\gamma_x$. Using the connection $\nabla^{\Sigma Q}$, we obtain a parallel translation $\tau^{t}_{0}:\Sigma_x M_{t} \to \Sigma_x M_{0}$. Notice that the Clifford multiplication ``$\cdot$'', the vector field $\tilde \partial_t$ and the volume form $\omega$ that determines the splitting $\Sigma Q = \Sigma^+Q \oplus \Sigma^-Q$ are all parallel.
		Identifying $M$ with $M_0$, we obtain a map
			\DefMap{\tau(\psi):=\bar \psi: \R}{\Gamma(\Sigma^+ M) \oplus \Gamma(\Sigma^- M)}{t}{(x \mapsto \tau_{t}^0(\psi^+_{(t,x)})+\tau_{t}^0(\psi^-_{(t,x)})).}
		This identification defines an isometry $\tau: L^2(\Sigma Q) \to L^2(\R, L^2(\Sigma^+ M \oplus \Sigma^- M))$. The operator in \eqref{EqDiracHyper} can be pulled back via a commutative diagram
		\begin{align*}
			\xymatrix{
				L^2(\Sigma Q)
					\ar[r]
					\ar[d]^-{\tau}
				&L^2(\Sigma Q)
					\ar[d]^-{\tau}
				\\
				L^2(\R, L^2(\Sigma^+ M \oplus \Sigma^- M))
					\ar[r]
				&L^2(\R, L^2(\Sigma^+ M \oplus \Sigma^- M)),
				}
		\end{align*}
		and the upper row is given by \eqref{EqDiracHyper}. We calculate what this equation looks like in the lower row: Since Clifford multiplication and $\tilde \partial_t$ are parallel,
		\begin{align*}
			 \tau \circ (\tilde \partial_t \cdot \Dirac^{Q}) \circ \tau^{-1}
			 = \tilde \partial_t \cdot (\tau \circ \Dirac^{Q} \circ \tau^{-1})
			 =: \tilde \partial_t \cdot \Dirac^{Q}_{M}.
		\end{align*}
		Since there is a splitting $(\Sigma Q, \nabla^Q) = (\Sigma^+Q,\nabla^{+} \oplus \nabla^{-})$, it suffices to check the following for a $\psi \in \Gamma(\Sigma^+ Q)$: Let $x \in M$ be arbitrary, and let $D_x$ be the covariant derivative induced by $\nabla^{\Sigma Q}$ along $\gamma^x$. For any $t_0 \in \R$ 
		\begin{align*}
			\nabla_{\tilde \partial t}^{\Sigma Q}{\psi}|_{{(t_0,x)}}
			&=D_x(\psi^x)(t_0)
			=\lim_{t \to t_0}{\frac{\tau_t^{t_0}(\psi^x(t)) - \psi^x(t_0)}{t-t_0}},
		\end{align*}
		and consequently
		\begin{align*}
			\overline{\nabla_{\tilde \partial t}^{\Sigma Q}{\psi}}(t_0)|_x
			=&\tau_{t_0}^0(\nabla_{\tilde \partial t}^{\Sigma Q}{\psi}|_{{(t_0,x)}})
			=\tau_{t_0}^0\left( \lim_{t \to t_0}{\frac{\tau_t^{t_0}(\psi^x(t)) - \psi^x(t_0)}{t-t_0}} \right) \\
			=& \lim_{t \to t_0}{\frac{\tau_{t_0}^0(\tau_t^{t_0}(\psi^x(t))) - \tau_{t_0}^0(\psi^x(t_0))}{t-t_0}} 
			= \lim_{t \to t_0}{\frac{\tau_t^{0}(\psi_{(t,x)}) - \tau_{t_0}^0(\psi|_{(t_0,x)})}{t-t_0}} \\
			=& \lim_{t \to t_0}{\frac{\bar \psi(t)|_x - \bar \psi(t_0)|_x}{t-t_0}} 
			=\tfrac{d\bar \psi}{dt}(t_0)|_x.
		\end{align*}
		This implies for any $\bar \psi \in L^2(\R, H^1(\Sigma^+ M \oplus \Sigma^- M))$
		\begin{align*}
			(\tau \circ \nabla_{\tilde \partial t}^{\Sigma Q} \circ \tau^{-1})(\bar \psi) = \tfrac{d}{dt} \bar \psi.
		\end{align*}		
		All in all \eqref{EqDiracHyper} transforms under $\tau \circ \_ \circ \tau^{-1}$ into
	\begin{align}
			\label{EqDiracHyperTransform}
				\tilde \partial_t \cdot \Dirac^{Q}_M = \tilde \Dirac^{Q}_M + \tfrac{m}{2}H - \tfrac{d}{dt},
		\end{align}
		where $H:\R \to \R$, $t \mapsto H_t$, and $\tilde \Dirac^Q_M:\R \to \Gamma(\Sigma M \oplus \Sigma M)$ is given by $\tilde \Dirac^Q_M(t) = \bar \Dirac(t) \oplus (- \bar \Dirac(t))$, $\bar \Dirac(t) = \tau \circ \Dirac^{M_t} \circ \tau^{-1}$. 
	
	\step[apply \cref{ThmSalamon}] Set $H := L^2(\R, L^2(\Sigma M))$, $W := L^2(\R, H^1(\Sigma M))$ and $A(t) = \bar \Dirac(t)$. Using \cref{ThmSalamon}, we obtain 
	\begin{align}
		\label{EqIndexSpecFlowSalamon}
		\begin{split}
			\ind(\Dirac^Q_+)
			&=\ind((\tilde \partial_t \cdot \Dirac^Q_M - \tfrac{m}{2} H)_+) 
			\jeq{\eqref{EqDiracHyperTransform}}{=}\ind((\Dirac^{Q}_M  - \tfrac{d}{dt})_+) \\
			&=\ind((\bar \Dirac  - \tfrac{d}{dt})_+) 
			=\ind(\tfrac{d}{dt} - \bar \Dirac)
			=\specfl(\bar \Dirac).
		\end{split}
	\end{align}
	Now the spectral flow $\specfl(\bar \Dirac)$ in the sense of Salamon, cf. \cite{Salamon}, coincides with the spectral flow in the sense of \cref{LemSpecFlowDef}.
	
	\step[final argument]
		By the classical Atiyah-Singer index theorem, cf. \cite[Thm. III.13.10]{LM}, we obtain 
		\begin{align}
			\label{EqIndexAtiyahSingerClassical}
			\widehat{A}(P) = \ind(\Dirac_+^P). 	
		\end{align}
		Consequently, we can put all the steps together to obtain
		\begin{align*}
			0 \jeq{\eqref{EqIndexHanke}}{\neq} \widehat{A}(P)
			&\jeq{\eqref{EqIndexAtiyahSingerClassical}}{=} \ind(\Dirac_+^P)
			\jeq{\eqref{EqIndexTrivializedBundle}}{=} \ind(\Dirac_+^{P'}) \\
			&\jeq{\eqref{EqIndexHalfCylinderGlued}}{=} \ind(\Dirac_+^{Q}) 
			\jeq{\eqref{EqIndexSpecFlowSalamon}}{=} \specfl(\bar \Dirac)
			=\specfl_{g_0}(f).
		\end{align*}
	\end{steplist}
\end{Prf}

%% file: fig.cut.tex
\begin{tikzpicture}
	\def\Mx{0}%
	\def\My{0}%
	\def\r{0.6}
	\def\inc{0.5}
	\draw (\Mx,\My) circle (\r);
	\draw (\Mx,\My) circle (\r+2*\inc);
	\draw[line width=2] (\Mx,\My+\r) -- (\Mx,\My+\r+2*\inc);
	\draw (\Mx,\My-\r) -- (\Mx,\My-\r-2*\inc);
	
	\draw[->] (\Mx+2.5,\My) -- (\Mx+3.5,\My);
	
	\def\Z{6}
	\def\O{8.5}
	\def\vert{1}
	
	\draw (\Z-2,\vert) -- (\O+2,\vert);
	\draw (\Z-2,-\vert) -- (\O+2,-\vert);	
	
	\draw[line width=2] (\Z,-\vert) -- (\Z,\vert);
	\draw[line width=2] (\O,-\vert) -- (\O,\vert);
	\draw (0.5*\Z+0.5*\O,-\vert) -- ( 0.5*\Z+0.5*\O,\vert);
	
	\node [label={[xshift=0, yshift=50]$P \to S^1$}]  {};
	\node [label={[xshift=-9, yshift=23]$M$}]  {};	
	\node [label={[xshift=-10, yshift=-40]$P_t$}]  {};
	\node [label={[xshift=7, yshift=11]$0$}]  {};
	\node [label={[xshift=-7, yshift=11]$1$}]  {};
	\node [label={[xshift=0, yshift=-20]$t$}]  {};
	
	\node [label={[xshift=171, yshift=-45]$0$}]  {};
	\node [label={[xshift=242, yshift=-45]$1$}]  {};
	\node [label={[xshift=207, yshift=-45]$t$}]  {};
	
	\node [label={[xshift=200, yshift=50]$Q \to \mathbb{R}$}]  {};
	\node [label={[xshift=215, yshift=25]$M_t$}]  {};
	\node [label={[xshift=170, yshift=25]$M$}]  {};
	\node [label={[xshift=245, yshift=25]$M$}]  {};
	\node [label={[xshift=190, yshift=0]$P'$}]  {};
	\node [label={[xshift=140, yshift=0]$Z_0'$}]  {};
	\node [label={[xshift=270, yshift=0]$Z_1'$}]  {};
	
\end{tikzpicture}

%% file: usb.abstract.tex
\begin{subabstract}
	For any two Riemannian metrics, the associated Dirac operators cannot be compared directly, because they are defined on different spaces. In this chapter, we will review a well known construction to identify the various spinor bundles with one another. This gives the statement that the Dirac operator depends continuously on the Riemannian metric a precise meaning and a proof. We will then reformulate these results and prove that the collection of all spinor fields formed with respect to the various metrics assemble to a continuous Hilbert bundle over the space of all Riemannian metrics. We will also show that on certain subsets of the Riemannian metrics, the spans of eigenspinors belonging to a finite number of eigenvalues assemble to a continuous vector bundle of finite rank. All these are important technical preparations for the proof of \cref{MainThmHigher}.
\end{subabstract}

%% file: usb.intro.tex
\section{Introduction and Statement of the Results}

For any two Riemannian metrics $g,h \in \Rm(M)$, the two Dirac operators
\begin{align*}
	\Dirac^g:H^1(\Sigma^g M) \subset L^2(\Sigma^g M) \to L^2(\Sigma^g M), &&
	\Dirac^h:H^1(\Sigma^h M) \subset  L^2(\Sigma^h M) \to L^2(\Sigma^h M),
\end{align*}
cannot be compared directly, because not only the operator $\Dirac^g$ depends on the metric $g$, but also its domain $H^1(\Sigma^g M)$. Therefore, the expression $\Dirac^g - \Dirac^h$ makes no sense, because these operators are defined on different spaces. 

A solution to this problem is to systematically construct identification isomorphisms $\bar \beta_{g,h}:L^2(\Sigma^g M) \to L^2(\Sigma^h M)$ for any two metrics $g$ and $h$ and use these maps to pull back one Dirac operator to the domain of definition of the other. 

The precise construction of $\bar \beta_{g,h}$ is rather complicated and will be carried out in several steps. We will construct:
\begin{enumerate}
	\item
		$a_{g,h}, b_{g,h} \in \Gamma(\Iso M)$ as an auxiliary tool, see \cref{Lemagh,Lembgh},
	\item
		$c_{g,h} :\SO^g M \to \SO^h M$ to identify the frame bundles, see \cref{Lemcgh}, 
	\item
		$\gamma_{g,h}: \Spin^g M \to \Spin^h M$ to identify the spin structures, see \cref{LemGammagh},
	\item
		$\beta_{g,h}: \Sigma^g_{\K} M \to \Sigma^h_{\K} M$ to identify the spinor bundles, see \cref{LemBetaDef},
	\item
		$\bar \beta_{g,h}: L^2(\Sigma^g_{\K} M) \to L^2(\Sigma^h_{\K} M)$ to identify the spinor fields, see \cref{CorBetabargh}.
\end{enumerate}

The ideas of these constructions originate from \cite{BourgGaud}, who construct the identification isomorphisms via the parallel transport of the lift of a partial connection. Following \cite{MaierGen}, we do not introduce a connection and use the Lifting theorem instead. In \cite{BourgGaud,MaierGen} these identifications are applied to the calculus of variations, where this problem becomes a pressing concern. One should remark that for Lorentz metrics there also exists an alternative approach by \emph{generalized cylinders}, see \cite{BaerGaudMor}.

We give a self-contained review of all the technical details of the construction mentioned above. We will require this insight to prove that 
	\DefMap{L^2(\Sigma_{\K} M) = \coprod_{g \in \Rm(M)}{L^2(\Sigma^g_{\K} M)}}{\Rm(M)}{\psi \in L^2(\Sigma^g_{\K} M)}{g,}
is a continuous bundle of Hilbert spaces, for which the $\bar \beta_{\_,g}$ are global trivializations, see \cref{LemUnivSpinBdleTrivs}. This is equivalent to showing that $\bar \beta_{g,h}$ depends continuously on $h$. Since $\bar \beta_{g,h}$ is constructed step by step according to the above list, we will show that in each of these steps, the objects depend continuously on $h$. To make sense of this claim, we will find spaces, which are defined independently of the metric, in which these objects live and explain how to topologize these spaces. We will use the weak topology from \cref{SubSectWeakTop} and its application in the topologization of the gauge groups from \cref{ThmSpinGaugeTrafosCov}.

Once the identifications $\bar \beta_{g,h}$ are constructed, one can use them to pull back all the Dirac operators $\Dirac^h$ to a single spinor bundle $L^2(\Sigma^g_{\K} M)$  via
\begin{align*}
	\Dirac^h_g := \bar \beta_{h,g} \circ \Dirac^h_{\K} \circ \bar \beta_{g,h}:L^2(\Sigma^g_{\K} M) \to L^2(\Sigma^g_{\K} M).
\end{align*}
The operator $\Dirac^h_g$ has virtually the same properties as $\Dirac^h_{\K}$, but is defined on $L^2(\Sigma^g_{\K} M)$. In particular $\Dirac^h_g$ is isospectral to $\Dirac^h_{\K}$. One can compute an explicit local coordinate formula for $\Dirac^h_g$, see \cref{ThmDiracghLocal}, which makes it possible to compare this operator with $\Dirac^g_{\K}$. It follows from this comparison that the operators $\Dirac^h_g$ depend continuously on $h$. This is a fact we already used, see \cref{ThmSpinorIdentification}. 

Finally, we will show that there are subsets $U \subset \Rm(M)$ of Riemannian metrics over which $L^2_{[\Lambda_1, \Lambda_2]}(\Sigma M)$, i.e. the eigenspinors with respect to eigenvalues $\Lambda_1 < \lambda < \Lambda_2$, is a continuous vector bundle of finite rank, see \cref{ThmUSFBSpecSubs}. This will be crucial to construct the vector bundle $E$ in our application of the Lasso Lemma (see \cref{LemLasso}) in the proof of \cref{MainThmHigher}.

In this chapter, $M$ still denotes a closed spin manifold of dimension $m \geq 3$ with a fixed topological spin structure $\Theta:\GLtp M \to \GLp M$. We denote by $\Iso^+ M$ the space of orientation-preserving isomorphism fields on $M$ thought of as a subset of $\End M = T^* M \otimes TM$. The corresponding space of sections $\Gamma(\End M)$ is given the $\mathcal{C}^1$-topology. All results hold for the real as well as for the complex Dirac operator.

\nomenclature[IsopM]{$\Iso^+ M$}{isomorphism fields on $M$}
\nomenclature[GammaIsopM]{$\Gamma(\Iso^+ M)$}{space of smooth sections in $\Iso^+ M$ endowed with $\mathcal{C}^1$-topology}

%% file: usb.usb.tex
\section{Construction of the Identification Maps}

\subsection{Construction of $a$ and $b$}

\begin{Lem}
	\label{Lemagh}
	\nomenclature[agh]{$a_{g,h}$}{an auxiliary map}
	For any $g,h \in \Rm(M)$, there exists a unique $a_{g,h} \in \Gamma(\Iso^+ M)$ such that
	\begin{align}
		\label{EqAlphaghDef}
		\forall x \in M: \forall v,w \in T_xM: g(a_{g,h}(v),w) = h(v,w).
	\end{align}	
	These maps have the following properties:
	\begin{enumerate}
		\item
			$a_{g,h}$ is an isomorphism and $a_{g,h}^{-1} = a_{h,g}$.
		\item
			$a_{g,h}$ is self-adjoint with respect to $g$ and $h$.
		\item
			$a_{g,h}$ is positive definite\footnote{We say that $f \in \End(V,g)$ is \emph{positive definite}, if $f$ is self-adjoint and 
			\begin{align*}
				\forall 0 \neq v \in V: g(f(v),v)>0.	
			\end{align*}
			 } with respect to $g$ and $h$.  \item $a$ is \emph{cocyclic in $\Rm(M)$}, i.e. if $k \in \Rm(M)$ is another Riemannian metric
			 \begin{align}
				\label{EqaghCocyclic}
				a_{g,h} \circ a_{h,k} = a_{g,k}, &&
				a_{g,g} = \id.
			 \end{align}
	\end{enumerate}	
\end{Lem}

\begin{Prf}
	First, we show uniqueness: Assume that $a_{g,h}, a'_{g,h} \in \Gamma(\End M)$ both satisfy \eqref{EqAlphaghDef}. This immediately implies for any $x \in M$, $v, w \in T_xM$
	\begin{align*}
		g(a_{g,h}(v) - a'_{g,h}(v),w) = h(v, w) - h(v, w) =  0,
	\end{align*}
	thus $a_{g,h}(v) = a'_{g,h}(v)$, since $g$ is positive definite. To show existence, it suffices to construct $a_{g,h}$ locally. Let $(c_1, \ldots, c_m)$ be a local $h$-ONF on $U$. By the Gram-Schmidt process
	\begin{align*}
		b_1 := \frac{c_1}{\|c_1\|_g}, && b_i' := c_i - \sum_{l=1}^{i-1}{g(b_i,c_l)c_l}, && b_i := \frac{b_i'}{\|b_i'\|_g},
	\end{align*}
	we obtain a $g$-ONB $b_1, \ldots, b_m$. Let $\xi_{g,h}$ be the local isomorphism field mapping the $c_i$'s to the $b_i's$. This map satisfies
	\begin{align*}
		\forall x \in M: \forall v,w \in T_xM: g(\xi_{g,h}(v),\xi_{g,h}(w)) = h(v,w),
	\end{align*}	
	i.e. it is an isometry between $(TU, h)$ and $(TU, g)$. Consequently, denoting by $\xi^{*_g}_{g,h}$ the $g$-adjoint of $\xi_{g,h}$, we find that $a_{g,h} := \xi^{*_g}_{g,h} \circ \xi_{g,h}$ satisfies \eqref{EqAlphaghDef}. The properties of $a_{g,h}$ are shown as follows.
	\begin{enumerate}
		\item This follows from the fact that $\xi_{g,h}$ is an isomorphism and $\xi_{g,h}^{-1} = \xi_{h,g}$.
		\item Locally, $a_{g,h}=\xi^{*_g}_{g,h} \circ \xi_{g,h}$, thus $a_{g,h}$ is $g$-self-adjoint. We calculate for any $v,w \in T_xM$,
		\begin{align*}
			h(a_{g,h}(v),w)
			&=g(a^2_{g,h}(v),w) 
			=g(v,a^2_{g,h}(w)) 
			=g(a^2_{g,h}(w), v) \\
			&=h(a_{g,h}(w), v)
			=h(v,a_{g,h}(w)).
		\end{align*}
		\item Since $h$ is positive definite, for any $0 \neq v \in TM$
		\begin{align*}
			g(a_{g,h}(v),v) = h(v,v) > 0,
		\end{align*}
		so $a_{g,h}$ is $g$-positive definite. Since $g$ is positive definite and $a_{g,h}$ is a $g$-self-adjoint isomorphism
		\begin{align*}
			h(a_{g,h}(v),v) = g(a_{g,h}^2(v),v) = g(a_{g,h}(v),a_{g,h}(v)) > 0.
		\end{align*}
		\item Since $a_{g,h}$ is uniquely determined by \eqref{EqAlphaghDef}, we just have to verify
		\begin{align*}
			g((a_{g,h} \circ a_{h,k})(v),w) 
			=h(a_{h,k}(v),w)
			=k(v,w), &&
			g(\id(v),w)) = g(v,w).
		\end{align*}
	\end{enumerate}
\end{Prf}

\begin{Cor}
	\label{Coraghcont}
	The maps $a_{g,h}$ from \cref{Lemagh} define a continuous map
		\DefMap{a:\Rm(M)^2}{\Gamma(\Iso^+ M)}{(g,h)}{a_{g,h}.}
	Here, $\Gamma(\Iso^+ M)$ is topologized as an open subset of $\Gamma(\End M) = \Gamma(T^*M \otimes TM)$ endowed with $\mathcal{C}^1$-topology.
\end{Cor}

\begin{Prf}
	This follows from the explicit construction of $a_{g,h}$ in the proof of \cref{Lemagh} and the fact that the Gram-Schmidt process depends $\mathcal{C}^1$-continuously on the involved metric. We can start with a fixed local frame $e_1, \ldots, e_m$ and apply the Gram-Schmidt process to it with respect to $h$. We obtain a well-defined local $h$-ONF $c_1, \ldots, c_m$ that depends $\mathcal{C}^1$-continuously on $h$. From this, we can construct the map $\xi_{g,h}$ which then depends $\mathcal{C}^1$-continuously on $g$ and $h$. 
\end{Prf}

\begin{Lem}
	\label{Lembgh}
	\nomenclature[bgh]{$b_{g,h}$}{an isomorphism field to identify bases}
	For any $g,h \in \Rm(M)$, there exists $b_{g,h} \in \Gamma(\Iso^+ M)$ such that
	\begin{align*}
		\forall x \in M: b_{g,h}(x)^2 = a_{g,h}(x)^{-1}.
	\end{align*}
	These endomorphisms satisfy:
	\begin{enumerate}
		\item
			$b_{g,h}$ is an isomorphism and $b_{g,h}^{-1} = b_{h,g}$.
		\item
			$b_{g,g} = \id$.
		\item
			$b_{g,h}$ is self-adjoint with respect to $g$ and $h$.
		\item
			$b_{g,h}$ is positive definite with respect to $g$ and $h$.
		\item
			$b_{g,h} \circ b_{h,k} = b_{g,k}$ if and only if $a_{g,h}$ and $a_{h,k}$ commute.
	\end{enumerate}	
\end{Lem}

\begin{Prf}
	Since $a_{g,h}$ is positive definite, so is $a_{g,h}^{-1}$. Consequently, there exists a unique square root, which we call $b_{g,h}$. The rest of the claims follow directly from \cref{Lemagh}. For the cocycle condition, notice that $a_{g,h}$ and $a_{h,k}$ commute if and only if their square roots commute.\footnote{It is clear that the endomorphisms commute, if their square roots commute. For the converse, notice that both endomorphisms are self-adjoint. Hence, their coordinate matrices are symmetric, hence diagonalizable. Now, two diagonalizable matrices $A$ and $B$ are simultanously diagonalizable if and only if $AB=BA$, see for instance \cite[4.3.6]{Fischer}. Now, if $D = SAS^{-1}$ is a diagonal matrix, $S \in \GL_n$, then $\sqrt{A} = S \sqrt{D} S^{-1}$ and analogously for $B$. This implies $\sqrt{A}\sqrt{B}=\sqrt{B}\sqrt{A}$ and thus the claim.} 
	Since $a_{g,h}$ and $a_{h,k}$ are both invertible, these isomorphisms commute if and only if
	\begin{align*}
		(a_{g,h}^{-1/2} \circ a_{h,k}^{-1/2})^2 = a_{g,h}^{-1} \circ a_{h,k}^{-1}.
	\end{align*}
	Therefore, by \cref{EqaghCocyclic} 
	\begin{align*}
		a_{g,k} = a_{h,k} \circ a_{g,h}
		&\Longleftrightarrow a_{g,k}^{-1} = a_{g,h}^{-1} \circ a_{h,k}^{-1} \\
		&\Longleftrightarrow a_{g,k}^{-1} = (a_{g,h}^{-1/2} \circ a_{h,k}^{-1/2})^2 \\
		&\Longleftrightarrow b_{g,k}^2 = (b_{g,h} \circ b_{h,k})^2 \\
		& \Longleftrightarrow b_{g,k} = b_{g,h} \circ b_{h,k}.
	\end{align*}
\end{Prf}

\begin{Cor}
	\label{CorbghCont}
	The maps $b_{g,h}$ from \cref{Lembgh} define a continuous map
		\DefMap{b:\Rm(M)^2}{\Gamma( \Iso^+ M)}{(g,h)}{b_{g,h}.}
\end{Cor}

\begin{Prf}
	Since taking square roots of a positive endomorphism is smooth, $b_{g,h} \in \Gamma(\Iso^+ M)$ is a smooth field, which depends continuously on $g$ and $h$ by \cref{Coraghcont}.
\end{Prf}

\subsection{Identification of the Frame Bundles}

\begin{Lem}
	\label{Lemcgh}
	\nomenclature[cgh]{$c_{g,h}$}{a gauge transformation mapping $\SO^g M$ to $\SO^h M$}
	The isomorphism field $b_{g,h}$ from \cref{Lembgh} induces an isomorphism of principal $\GLp_m$-bundles
	\DefMap{c_{g,h}:\GLp M }{\GLp M}{(e_i)_{i=1}^m}{(b_{g,h}(e_i))_{i=1}^m}
	which restricts to an isomorphism $c_{g,h}:\SO^g M \to \SO^h M$ of principal $\SO_m$-bundles. In addition, the isomorphisms $c_{g,h}$ satisfy $c_{g,h}^{-1} = c_{h,g}$ and $c_{g,g} = \id$. 
\end{Lem}

\begin{Prf} Let $(e_1, \ldots, e_m) \in \GLp M$ and $A \in \GL_m^+$ be arbitrary.
	\begin{steplist}
		\step[orientation]
			If $\Lambda \in \Omega^m(M)$ is an orientation form
			\begin{align*}
				\Lambda(b_{g,h}(e_1), \ldots, b_{g,h}(e_n)) = \det(b_{g,h}) \Lambda(e_1, \ldots, e_n).
			\end{align*}			
			Since $b_{g,h} \in \Gamma(\Iso^+M)$, we obtain $\det(b_{g,h})> 0$ and therefore $c_{g,h}$ preserves the orientation of a basis.
		\step[equivariance]
			We calculate
			\begin{align*}
				b_{g,h}(e_i.A)
				=\sum_{j=1}^{m}{b_{g,h}(A^j_i e_j)}
				=\sum_{j=1}^{m}{A^j_i b_{g,h}( e_j)}
				=b_{g,h}(e_i).A.
			\end{align*}
			This automatically implies that $c_{g,h}$ is an isomorphism, see \cref{ThmGBundleMorphismIso}.
		\step[restriction]
			Let $(e_1, \ldots, e_m) \in \SO^g M$. For any $1 \leq i,j \leq m$, we calculate (using \cref{Lembgh})
			\begin{align*}
				h(b_{g,h}(e_i),b_{g,h}(e_j))
				&=h(b^2_{g,h}(e_i),e_j) 
				=h(a^{-1}_{g,h}(e_i), e_j)
				=g(e_i,e_j)
				=\delta_{ij}.
			\end{align*}
			Consequently, $(b_{g,h}(e_i), \ldots, b_{g,h}(e_j))$ is an $h$-orthogonal basis. 
	\end{steplist}
\end{Prf}

\begin{Rem}
	The reason we defined $c_{g,h}$ on the larger space $\GLp M$ is that $\GLp M$ does not depend on $g$ and $h$. 
\end{Rem}

\begin{Lem}
	\label{ThmIsoIsoMor}
	The map
		\DefMap{\Phi: (\Gamma(\Iso M), \mathcal{C}^1) }{(\G( \GL M), \mathcal{C}^1_w)}{F}{\hat F,}
	where
	\begin{align}
		\label{EqIsoIsoGaugeDef}
		\forall b=(b_1, \ldots, b_m) \in \GL M: \hat F(b_1, \ldots, b_m) := (F b_1, \ldots, F b_m),
	\end{align}
	is a homeomorphism. Here, $\mathcal{C}^1_w$ is the weak topology, see \cref{DefWeakTopology}. It restricts to a map
	\begin{align*}
		\Phi: \Gamma(\Iso^+ M) \to \G(\GLp M).
	\end{align*}
\end{Lem}

\begin{Prf} $ $
	\begin{steplist}
		\step[equivariance]
			To see that $\Phi(F)$ is a gauge transformation, we take any $A \in \GL_m$, $b \in \GL M$ and calculate
			\begin{align*}
				\Phi(F)(b.A) 
				&= \Phi(F)(A^{i_1}_1 b_1, \ldots, A^{i_m}_m b_m)
				= (F(A^{i_1}_1 b_1), \ldots, F(A^{i_m}_m b_m)) \\
				&= (A^{i_1}_1 F(b_1), \ldots, A^{i_m}_m F(b_m))
				=\Phi(F)(b).A.
			\end{align*}
		\step[inverse]
			Define $\Phi^{-1}$ as follows: Let $\alpha \in \G(\GL M)$ and choose any basis $b=(b_1, \ldots, b_m) \in \GL_x M$, $x \in M$. By construction $\alpha(b)=:c=(c_1, \ldots, c_m) \in \GL_x M$. Let $\Phi^{-1}(\alpha)|_{\Iso_x M}$ be the isomorphism $F_x$ mapping $b_i$ to $c_i$, $1 \leq i \leq m$. Then the definition of $F_x \in \Iso(T_xM)$ does not depend on the choice of the basis $b$: Let  $b'$ be any other basis, $c':=\alpha(b)$ and $F_x'$ be defined analogously. There exists a unique $A \in \GL_m$ such that $b=b'.A$. Since $\alpha$ is a gauge transformation, we obtain $c = \alpha(b) = \alpha(b').A = c'.A$.  Therefore, for any $1 \leq j \leq m$,
			\begin{align*}
				F'_x(b_j)
				=F'(A^{\nu}_j b'_\nu)
				=A^{\nu}_j F'(b'_\nu)
				=A^{\nu}_j  c'_{\nu}
				=A^{\nu}_j (A^{-1})^{\mu}_\nu c_\mu
				=(A^{-1}A)^\mu_j c_\mu
				=c_j
				=F_x(b_j),
			\end{align*}
			thus $F'_x = F_x$. Obviously, $\Phi^{-1}$ is an inverse to $\Phi$ and $\hat F$ preserves the orientation, if $F$ does.
		\step[continuity]
			To see that $\Phi$ is a homeomorphism, we compose $\Phi$ with the homeomorphism $\sigma$, see \cref{CorWeakTopGBundle}\ref{ItWeakTopSigmaHoem}, and obtain that the map
			\begin{align*}
				\xymatrix{
					\Gamma(\Iso M)
						\ar[r]^-{\Phi}
					& \G(\GL M)
						\ar[r]^-{\sigma}
					& \mathcal{C}^{\infty}_e(\GL M, \GL_m)
				}
			\end{align*}
			is given by $F \mapsto \sigma_{\Phi(F)}$, where for any $b \in \GL M$, $\Phi(F)(b) = b. \sigma_{\Phi(F)}(b)$. By definition, $\sigma_{\Phi(F)}(b)$ is simply the coordinate matrix of $F$ with respect to $b$, thus $\sigma \circ \Phi$ is a homeomorphism, which implies the claim.
	\end{steplist} 
\end{Prf}

\begin{Cor}
	\label{CorcghCont}
	The map
		\DefMap{c:\Rm(M)^2}{\G(\GLp M)}{(g,h)}{c_{g,h},}
	where $c_{g,h}$ is defined in \cref{Lemcgh}, is continuous. 
\end{Cor}

\begin{Prf}
	By \cref{CorbghCont}, the map $b$ is continuous, by \cref{ThmIsoIsoMor}, $\Phi$ is continuous and by definition $c = \Phi \circ b$. 
\end{Prf}

\subsection{Identification of the Spin Structures}

\begin{Thm}
	\label{LemGammagh}
	\nomenclature[gammagh]{$\gamma_{g,h}$}{a gauge transformation mapping $\Spin^g M$ to $\Spin^h M$}
	For any $g,h \in \Rm(M)$, there exists an isomorphism $\gamma_{g,h}$ of principal $\GLtp_m$-bundles such that
	\begin{align}
		\label{EqGammaghTopLift}
		\begin{split}
			\xymatrixcolsep{3.5em}
			\xymatrix{
				\GLtp M
					\ar[r]^-{\gamma_{g,h}}
					\ar[d]^-{\Theta}
				&\GLtp M
					\ar[d]^-{\Theta}
				\\
				\GLp M
					\ar[r]^-{c_{g,h}}
				&\GLp M
			}
		\end{split}
	\end{align}
	commutes. In particular, the restriction of this map to $\Spin^g M$ yields a commutative diagram
	\begin{align} 
		\label{EqGammaghProp}
		\begin{split}
			\xymatrixcolsep{3.5em}
			\xymatrix{
				\Spin^g M
					\ar[r]^-{\gamma_{g,h}}
					\ar[d]^-{\Theta^g}
				&\Spin^h M
					\ar[d]^-{\Theta^h}
				\\
				\SO^g M
					\ar[r]^-{c_{g,h}}
				&\SO^h M
			}
		\end{split}		
	\end{align}
	and $\gamma_{g,h}:\Spin^g M \to \Spin^h M$ is an isomorphism of principal $\Spin_m$-bundles.
\end{Thm}

\begin{Prf} Recall that $I=[0,1]$.
	\begin{steplist}
	\step[construction]
		Setting $g_I :=(g_s)_{s \in I}$, where
		\begin{align*}
			\forall s \in I: g_s := (1-s)g + sh \in \Rm(M),
		\end{align*}
		defines a path of Riemannian metrics. The map
			\DefMap{F:I \times \GLtp M}{\GL^+(M)}{(s, \tilde b)}{c_{g,g_s}(\Theta(\tilde b)),}
		is a homotopy. We would like to apply the Homotopy Lifting Theorem (see \cref{ThmHomotopyLifting}) to show that there exists $G$ such that
		\begin{align}
			\label{EqLemmaghLiftFG}
			\begin{split}
				\xymatrixcolsep{3.5em}
				\xymatrix{
					& & \GLtp M
						\ar[d]^-{\Theta} \\
					\Spin^g M
						\ar@/^2pc/[rru]_-{\iota}
						\ar@/_2pc/[rr]^-{\Theta}
						\ar@{^(->}[r]^-{\{0\} \times \iota }
					&I \times \GLtp M 
						\ar@{..>}[ru]^-{G}
						\ar[r]^-{F}
					& \GLp M
				}
			\end{split}
		\end{align}
		commutes. Here, $\iota$ is the canonical inclusion. In order to do this, we just have to verify
		\begin{align*}
			\forall \tilde b \in \Spin^g M: (\Theta \circ \iota)(\tilde b) &= \Theta(\tilde b) = c_{g,g}(\Theta(\tilde b)) = F(\tilde b,0).
		\end{align*}
		Consequently, there exists $G$ such that $\Theta \circ G = F$ and $G(\_,0)=\iota$. We define
			\DefMap{\gamma_{g,h}:\GLtp M}{\GLtp M}{\tilde b}{G(\tilde b,1)}
		and obtain
		\begin{align*}
			\forall \tilde b \in \GLp M: \Theta (\gamma_{g,h}(\tilde b)) = \Theta(G(\tilde b,1)) = F(\tilde b,1) = c_{g,h}(\Theta(\tilde b)),
		\end{align*}
		which proves \cref{EqGammaghProp}.
		
	\step[equivariance] This follows directly from \cref{LemLiftEquiv}.
	\end{steplist}
\end{Prf}

\begin{Rem}
	One can show that $\gamma_{g,h}^{-1} = \gamma_{h,g}$, $\gamma_{g,g} = \id$.
\end{Rem}

\begin{Cor}
	\label{CorGammaCont}
	The map $c$ from \cref{Lemcgh} admits a lift against the covering $\Theta_*$ from \cref{ThmSpinGaugeTrafosCov}, i.e. there exists a continuous map $\gamma$ such that
	\begin{align*}
		\xymatrix{
			& \G(\GLtp M)
				\ar[d]^{\Theta_*}
			\\
			\Rm(M)^2
				\ar[r]^-{c}
				\ar@{..>}[ur]^{\gamma}
			&\G^{\spin}(\GLp M)
		}
	\end{align*}
	commutes. One can chose $\gamma$ such that for all $g, h \in \Rm(M)$, $\gamma(g, h) = \gamma_{g,h}$, where $\gamma_{g,h}$ is as in \cref{LemGammagh}.
\end{Cor}

\begin{Prf}
	By \cref{LemGammagh}, the map $c$ is actually a map to $\G^{\spin}(\GLp M)$, see \cref{DefSpinGaugeTrafo}. By \cref{ThmSpinGaugeTrafosCov}, the map $\Theta_*$ is a covering. Since $\Rm(M)^2$ is simply connected, the existence of the lift follows from the Lifting Theorem of Algebraic Topology, see \cref{ThmLifting}. By definition of $\Theta_*$ and \cref{EqGammaghTopLift}, we can ensure that $\gamma(g,h)=\gamma_{g,h}$.
\end{Prf}

\subsection{Identification of the Spinor Bundles}

\begin{Lem}
	\label{LemBetaDef}
	\nomenclature[betagh]{$\beta_{g,h}$}{isometry between $\Sigma^g M$ and $\Sigma^h M$}
	The isomorphism $\gamma_{g,h}$ from \cref{LemGammagh} induces an isomorphism of vector bundles
		\DefMap{\beta_{g,h}:\Sigma^g_{\K} M}{\Sigma^h_{\K} M}{{[s,\sigma]}}{[\gamma_{g,h}(s),\sigma],}
	which is a fibrewise isometry with respect to the metrics on $\Sigma^g_{\K} M$, $\Sigma^h_{\K} M$ satisfying $\beta_{g,h}^{-1} = \beta_{h,g}$. Furthermore, Clifford multiplication satisfies
	\begin{align}
		\label{EqBetaDefClifford}
		\forall X \in TM: \forall \psi \in \Sigma^g_{\K} M: \beta_{g,h}(X \cdot \psi) = b_{g,h}(X) \cdot \beta_{g,h}(\psi).
	\end{align}
\end{Lem}

\begin{Prf} $ $
	\begin{steplist}
	\step
		By \cref{LemGBundleIsoAssocVB}, the map $\beta_{g,h}$ is well-defined and an isometry. By \cref{LemGammagh} $\beta_{g,h}^{-1} = \beta_{h,g}$. 
	\step[Clifford multiplication]
		Take any
		\begin{align*}
			X = [b,x] \in TM \cong \SO^g M \times_\rho \R^m, &&
			\psi = [\tilde b, \sigma] \in \Sigma^g_{\K} M, &&
			\Theta(\tilde b)=b.
		\end{align*}
		By definition of the Clifford multiplication,
		\begin{align}
			\label{EqDefCliffRecall}
			[b,x] \cdot [\tilde b, \sigma] = [\tilde b, x \cdot \sigma].
		\end{align}
		By \eqref{EqGammaghProp}, this implies
		\begin{align} \label{EqBetaDefMulti}
			[c_{g,h}(b),x] \cdot [\gamma_{g,h}(\tilde b), \sigma] = [\gamma_{g,h}(\tilde b), x \cdot \sigma].
		\end{align}
		Consequently,
		\begin{align*}
			\beta_{g,h}(X \cdot \psi)
			&= \beta_{g,h}([b,x] \cdot [\tilde b, \sigma])
			\jeq{\eqref{EqDefCliffRecall}}{=} \beta_{g,h}([\tilde b, x \cdot \sigma])
			= [\gamma_{g,h}(\tilde b), x \cdot \sigma] \\
			&\jeq{\eqref{EqBetaDefMulti}}{=}[c_{g,h}(b),x] \cdot [\gamma_{g,h}(\tilde b), \sigma] 
			=b_{g,h}([b,x]) \cdot \beta_{g,h}([\tilde b, \sigma]) 
			=b_{g,h}(X) \cdot \beta_{g,h}(\psi).
		\end{align*}
	\end{steplist}
\end{Prf}

\begin{Cor}
	\nomenclature[betaggprime]{$\beta^{g,g'}$}{trivialization for $\Gamma(\Sigma M)$}
	For any fixed $g, g' \in \Rm(M)$, the map
		\DefMap{\beta^{g, g'}:\Rm(M)}{\Gamma(\Isom(\Sigma^g_{\K} M,\Sigma^{g'}_{\K} M))}{h}{\beta_{h,g'} \circ \beta_{g,h}.}
	is continuous. Here, $\beta_{g,h}$ is defined as in \cref{LemBetaDef} and the isometry fields are topologized as a subset of the homomorphism fields $\Gamma(\Hom(\Sigma^g_{\K} M, \Sigma^{g'}_{\K} M) = \Gamma((\Sigma^g_{\K} M)^*  \otimes \Sigma^{g'}_{\K} M)$ endowed with $\mathcal{C}^1$-topology.
\end{Cor}

\begin{Prf}
	By \cref{ThmCOkTop}\ref{ItCOkTopComp} and \cref{CorGammaCont}, the map 
		\DefMap{\gamma_{g',g} \circ \gamma^{g, g'}:\Rm(M)}{\G(\Spin^g M)}{h}{\gamma_{g',g} \circ \gamma_{h, g'} \circ \gamma_{g, h}}
	is continuous. Since $\Sigma^g_{\K} M = \Spin^g M \times_{\rho} W$ is an associated vector bundle, the induced map 
		\DefMap{\Rm(M)}{\Gamma(\Isom(\Sigma^g_{\K} M))}{h}{{([s, w] \mapsto [(\gamma_{g',g} \circ \gamma^{g, g'})(s), w])}}
	is continuous by \cref{ThmContGaugeAssocIsos}. This map equals $\beta^{g,g'}$ up to a post-composition with $\beta_{g',g}$. Since composition is continuous by \cref{ThmCOkTop}\ref{ItCOkTopComp}, this implies the claim. 
\end{Prf}

\subsection{Identification of the Spinor Fields}

Although the isomorphism $\beta_{g,h}:\Sigma^g_{\K} M \to \Sigma^h_{\K} M$ from \cref{LemBetaDef} is an isometry of vector bundles, this map does not induce an isometry $L^2(\Sigma^g_{\K} M) \to L^2(\Sigma^h_{\K} M)$ between the Hilbert spaces, because the Riemannian volume forms $\dvol_{g}$ and $\dvol_{h}$ on $M$ might differ. In order to compensate this, we need to introduce the following. 

\begin{Lem}
	\label{LemfghDef}
	\nomenclature[fgh]{$f_{g,h}$}{function relating the volume forms $\dvol_{h}$ and $\dvol_{g}$}
	For any $g,h \in \Rm(M)$, there exists a unique function $f_{g,h} \in \mathcal{C}^\infty(M,\R_+)$ such that
	\begin{align}
		\label{EqDeffg}
		\dvol_{h} = f_{g,h}^2 \dvol_{g}.
	\end{align}
	These functions satisfy:
	\begin{enumerate}
		 \item
			Locally, $f_{g,h}$ can be expressed as follows: Let $\varphi$ be a chart on $U \subset M$, let $\det(g) := \det(g_{ij})>0$ and $g_{ij} = g(\partial \varphi_i, \partial \varphi_j)$. Then we have
			\begin{align*}
				f_{g,h} = \sqrt[4]{\frac{\det(h)}{\det(g)}}.
			\end{align*}
		 \item
			$f_{g,g} = 1$.
		 \item
			$f_{g,h}^{-1} = f_{h,g}$.
		 \item
			For any $k \in \Rm(M)$, we obtain $f_{h,k} \cdot f_{g,h} = f_{g,k}$.
	\end{enumerate}
\end{Lem}

\begin{Prf}
	Since $\dvol_{g}, \dvol_{h} \in \Omega^m(M)$ are both smooth positive Riemannian volume forms and $\Omega^m M$ is a smooth line bundle, we immediately obtain the existence of $f_{g,h} \in \mathcal{C}^\infty(M,\R_+)$. Define $\omega := \partial \varphi \wedge \ldots \wedge \partial \varphi^m$. By
	\cite[13.23]{LSM}
	\begin{align*}
		\dvol_{g} = \sqrt{\det(g)} \omega, && \dvol_{h} = \sqrt{\det(h)} \omega.
	\end{align*}
	Consequently,
	\begin{align*}
		 \sqrt[4]{\frac{\det(h)}{\det(g)}}^2 \dvol_{g}
		 =\frac{\sqrt{\det(h)}}{\sqrt{\det(g)}} \sqrt{\det(g)} \omega
		 =\dvol_{h},
	\end{align*}
	thus $f_{g,h}$ has the local representation as claimed. The rest of the claim follows from this local representation.
\end{Prf}

\begin{Cor}
	\label{CorBetabargh}
	\nomenclature[barbetagh]{$\bar \beta_{g,h}$}{isomorphism $L^2(\Sigma^g M) \to L^2(\Sigma^h M)$}
	\nomenclature[barbetaggprim]{$\bar \beta^{g,g'}$}{trivialization for $L^2(\Sigma M)$}
	Let $g, g', h \in \Rm(M)$, $f_{g,h}$ be as in \cref{EqDeffg} and $\beta_{g,h}$ be as in \cref{LemBetaDef}.
	\begin{enumerate}
		\item
			For any $g,h \in \Rm(M)$ the map
			\begin{align*}
				\bar \beta_{g,h} := \frac{1}{f_{g,h}} \beta_{g,h}: L^2(\Sigma^g_{\K} M) \to L^2(\Sigma^h_{\K} M)
			\end{align*}
			is an isometry of Hilbert spaces. 
		\item
			The restriction $\bar \beta|_{H^1(\Sigma^g_{\K} M)}$ is a continuous isomorphism $H^1(\Sigma^g_{\K} M) \to H^1(\Sigma^h_{\K} M)$.
		\item
			For any fixed $g, g' \in \Rm(M)$, this defines a continuous map
				\DefMap{\bar \beta^{g,g'}:\Rm(M)}{\Isom(L^2(\Sigma^g_{\K} M), L^2(\Sigma^{g'}_{\K} M))}{h}{\bar \beta_{h,g'} \circ \bar \beta_{g,h}.}
		\item
			\label{ItH1L2Adapted}
			This map induces a continuous map
				\DefMap{\bar \beta^{g, g'}:\Rm(M)}{\Iso(H^1(\Sigma^g_{\K} M), H^1(\Sigma^{g'}_{\K}M))}{h}{\bar \beta^{g, g'}(h)|_{H^1(\Sigma^g_{\K} M)}.}
	\end{enumerate}
\end{Cor}

\begin{Prf} $ $ 
	\begin{enumerate}
		\item 
			By \cref{LemBetaDef}, $\beta_{g,h}:\Sigma^g_{\K} M \to \Sigma^h_{\K} M$ defines a field of isometries. Consequently, for any $\psi \in L^2(\Sigma^g_{\K} M)$,
			\begin{align}
				\label{EqbarBetaIsometry}
				\begin{split}
					\|\bar \beta_{g,h}(\psi)\|_{L^2(\Sigma^h_{\K} M)}^2
					&=\int_{M}{|f_{g,h}^{-1} \beta_{g,h}(\psi)|^2 \dvol_{h}} 
					= \int_{M}{|\psi|^2 f_{g,h}^{-2} \dvol_{h}} \\
					&\jeq{\cref{EqDeffg}}{=} \int_{M}{|\psi|^2 \dvol_{g}} 
					=\|\psi\|_{L^2(\Sigma^g_{\K} M)}^2.
				\end{split}
			\end{align}
		\item  
			For any $\psi \in H^1(\Sigma^g_{\K} M)$, 
			\begin{align*}
				\|\nabla^h(\bar \beta_{g,h}(\psi))\|_{L^2(\Sigma^h_{\K} M)}
				&=\|\nabla^h(\bar \beta_{g,h})(\psi) + \bar \beta_{g,h}(\nabla^g \psi)\|_{L^2(\Sigma^h_{\K} M)} \\
				&\leq \|\nabla \bar \beta_{g,h} \|_{\mathcal{C}^0} \|\psi\|_{L^2(\Sigma^g_{\K} M)} + \| \bar \beta_{g,h}(\nabla^g \psi)\|_{L^2(\Sigma^h_{\K} M)} \\
				&\leq \|\nabla \bar \beta_{g,h} \|_{\mathcal{C}^0} \|\psi\|_{H^1(\Sigma^g_{\K} M)} + \| \psi\|_{H^1(\Sigma^g_{\K} M)} \\
				&= (\|\nabla \bar \beta_{g,h} \|_{\mathcal{C}^0} + 1) \|\psi\|_{H^1(\Sigma^g_{\K} M)}.
			\end{align*}
			Combining this with \cref{EqbarBetaIsometry} yields $\bar \beta_{g,h}(\psi) \in H^1(\Sigma^h_{\K} M)$. Consequently, $\bar \beta:H^1(\Sigma^g_{\K} M) \to H^1(\Sigma^h_{\K} M)$. By the same argument, $\bar \beta_{g,h}^{-1} = \bar \beta_{h,g}:H^1(\Sigma^h_{\K} M) \to H^1(\Sigma^g_{\K} M)$. Thus, $\bar \beta_{g,h}$ is an isomorphism $H^1(\Sigma^g_{\K} M) \to H^1(\Sigma^h_{\K} M)$ (however, it is not an isometry with respect to the Sobolev norms).
		\item
			By \cref{LemBetaDef}, we already have a continuous map
			\begin{align*}
				\beta^{g,g'}: \Rm(M) \to \Gamma(\Isom(\Sigma^g_{\K} M, \Sigma^{g'}_{\K} M)).
			\end{align*}
			The metrics $g$ and $g'$ induce a fibre metric on $\Hom(\Sigma^g_{\K} M, \Sigma^{g'}_{\K} M)$. Therefore, the $\mathcal{C}^1$-topology can be described by the induced metric on the sections. 
			With this in mind, we directly verify the continuity claim: Let $\varepsilon > 0$ and $h \in \Rm(M)$. Certainly, there exists a constant
			\begin{align*}
				C_{g',g } := \sup_{x \in M}{f_{g',g}(x)} = \|f\|_{\mathcal{C}^0} \; \in \; ]0,\infty[.
			\end{align*}
			Since $\beta^{g,g'}$ is continuous, there exists an open neighborhood $U \subseto \Rm(M)$ of $h$ such that for all $h' \in U$, $\beta^{g,g'}(h') \in B_{\varepsilon/C_{g',g}}^{\mathcal{C}^1}(\beta^{g,g'}(h))$. This implies that for all $h' \in U$ and all $\psi \in L^2(\Sigma^{g}_{\K} M)$
			\begin{align}
				\label{EqBarBetaggpCont}
				\begin{split}
					\|(\bar \beta^{g,g'}(h) - \bar \beta^{g,g'}(h'))\psi\|_{L^2(\Sigma^{g'}_{\K} M)}^2
					& =\int_{M}{ |f_{g,h}^{-1} f_{h,g'}^{-1} \beta^{g,g'}(h)(\psi) - f_{g,h}^{-1} f_{h,g'}^{-1} \beta^{g,g'}(h')\psi|^2 d_{g'}V} \\
					& =\int_{M}{ f_{g',g}^2 |(\beta^{g,g'}(h) - \beta^{g,g'}(h'))\psi|^2 d_{g'}V} \\
					& \leq C_{g',g}^2 \|\beta^{g,g'}(h) - \beta^{g,g'}(h') \|^2_{\mathcal{C}^1} \| \psi\|_{L^2(\Sigma^{g'}_{\K} M)}^2 \\
					&\leq \varepsilon^2 \|\psi\|_{L^2(\Sigma^{g'}_{\K} M)}^2.
				\end{split}
			\end{align}
			Thus $\bar \beta(h') \in B_\varepsilon(\bar \beta^{g, g'}(h))$.
		\item
			Analogously, define $C'_{g',g} := \|f_{g',g}\|_{\mathcal{C}^1}$ and $C := \max(C'_{g',g}, C_{g',g})$. Since $\beta^{g,g'}$ is continuous, for any $h \in \Rm(M)$, we obtain an open neighborhood $U \subset \Rm(M)$ such that for all $h' \in U$, $\beta^{g,g'}(h') \in B_{\varepsilon / 3 C}^{\mathcal{C}^1}(h)$. Since
			\begin{align*}
				\|\nabla (\bar \beta^{g,g'}(h))(\psi)\|_{L^2(\Sigma^{g'}_{\K} M)}
				&=\|\nabla (f_{g',g} \beta^{g,g'}(h))(\psi)\|_{L^2(\Sigma^{g'}_{\K} M)} \\
				&\leq C \|\beta^{g,g'}(h)(\psi)\|_{L^2(\Sigma^{g'}_{\K} M)} 
				+C \|(\nabla \beta^{g,g'}(h))(\psi)\|_{L^2(\Sigma^{g'}_{\K} M)}, \\
				&\leq  2 C \| \beta^{g,g'}(h)\|_{\mathcal{C}^1} \|\psi\|_{L^2(\Sigma^{g'}_{\K} M)},
			\end{align*}
			we obtain for any $h' \in U$,
			\begin{align*}
				\|\nabla^{g'}((\bar \beta^{g, g'}(h) - \bar \beta^{g, g'}(h'))(\psi))\|_{L^2(\Sigma^{g'}_{\K} M)}
				&\leq  \|\nabla(\bar \beta^{g, g'}(h) - \bar \beta^{g, g'}(h')) (\psi) \|_{L^2(\Sigma^{g}_{\K} M)} \\
				& +\|\beta^{g, g'}(h) - \beta^{g, g'}(h') \|_{\mathcal{C}^0} \|\nabla^g \psi\|_{L^2(\Sigma^g_{\K} M)} \\
				&\leq 3 C \| \beta^{g, g'}(h) - \beta^{g, g'}(h') \|_{\mathcal{C}^1} \|\psi\|_{H^1(\Sigma^g_{\K} M)} \\
				&\leq \varepsilon \|\psi\|_{H^1(\Sigma^g_{\K} M)}.
			\end{align*}
			Combining this with \cref{EqBarBetaggpCont}, yields the result.
	\end{enumerate}
\end{Prf}

\section{Construction of the Bundle}

\begin{Def}[universal spinor bundle]
	\nomenclature[SigmaM]{$\Sigma_{\K} M$}{universal spinor bundle}
	\label{DefUniversalSpinorBundle}
	For any $g \in \Rm(M)$ and $x \in M$ let $\Sigma_{\K} M _{(g,x)} := \Sigma^g_{\K} M|_x$ be the spinor space at $x$ with respect to $g$. The bundle
		\DefMap{\pi: \Sigma_{\K} M := \coprod_{(g,x) \in \Rm(M) \times M}{\Sigma^g_{\K} M|_x}}{\Rm(M) \times M}{\psi \in \Sigma^g_{\K} M|_x}{(g, x)}
	is called \emphi{universal spinor bundle}. 
\end{Def}

\begin{Lem} 
	The universal spinor bundle $\pi: \Sigma_{\K} M \to \Rm(M) \times M$ carries a unique vector bundle topology such that for any $g \in \Rm(M)$ and any local trivialization $\alpha:\Sigma^g_{\K} M|_U \to U \times \R^n$ of $\Sigma^g_{\K} M$,
		\DefMap{\Phi_g: \Sigma_{\K} M|_U}{\Rm(M) \times U \times \R^n}{\psi \in \Sigma^h_{\K} M|_x}{(h,x,\alpha(\beta_{h,g}(\psi))),}
	is a local trivialization. Here, $\beta_{g,h}$ is the isomorphism defined in \cref{LemBetaDef}. 
\end{Lem}

\begin{Prf}
	For any other reference metric $g' \in \Rm(M)$ and any other local trivialization $\alpha'$,
		\DefMap{\Phi_{g'} \circ \Phi_g^{-1}: \Rm(M) \times U \times \R^n}{\Rm(M) \times U \times \R^n}{(h,x,v)}{(h,x,\alpha'(\beta_{h,g'}(\beta_{g,h}(\alpha^{-1}(v))))),}
	which is continuous by \cref{LemBetaDef}. 
\end{Prf}

\begin{Def}[universal spinor field bundle]
	\nomenclature[GammaSigmaM]{$\Gamma(\Sigma_{\K} M)$}{smooth universal spinor field bundle}
	\nomenclature[L2SigmaM]{$L^2(\Sigma_{\K} M)$}{universal spinor field bundle}
	\label{DefUniSpinorFieldBdle} 
	The smooth section bundle
		\DefMap{\Gamma(\Sigma_{\K} M) := \coprod_{g \in \Rm(M)}{\Gamma(\Sigma^g_{\K} M)}}{\Rm(M)}{\psi \in \Gamma(\Sigma^g_{\K} M)}{g}
	is called the \emph{(smooth) universal spinor field bundle}. Analogously we define
		\DefMap{L^2(\Sigma_{\K} M) := \coprod_{g \in \Rm(M)}{L^2(\Sigma^g_{\K} M)}}{\Rm(M)}{\psi \in L^2(\Sigma^g_{\K} M)}{g}
	to be the \emphi{universal spinor field bundle}. 
\end{Def}

\begin{Thm}[topologization of universal spinor field bundles] $ $
	\label{LemUnivSpinBdleTrivs}
	\begin{enumerate}
		\item 
			The smooth universal spinor field bundle $\Gamma(\Sigma_{\K} M) \to \Rm(M)$ carries a unique topology as an infinite dimensional vector bundle (see \cref{DefVectBanHilBundle}) such that for any fixed $g \in \Rm(M)$,
				\DefMap{\beta_{\_,g}:\Gamma(\Sigma_{\K} M)}{\Rm(M) \times \Gamma(\Sigma^g_{\K} M)}{\psi \in \Gamma( \Sigma^h_{\K} M)}{(h,\beta_{h,g}(\psi)),}
			is a global trivialization. 
		\item
			Analogously, the universal spinor field bundle $L^2(\Sigma_{\K} M) \to \Rm(M)$ carries a unique Hilbert bundle topology such that for any fixed $g \in \Rm(M)$, the map
				\DefMap{\bar \beta_{\_,g}: L^2( \Sigma_{\K} M)}{\Rm(M) \times L^2(\Sigma^g_{\K} M)}{\psi \in L^2(\Sigma^h_{\K} M)}{(h,\bar \beta_{h,g}(\psi)),}
			is a global trivialization. 
		\item
			The subspaces 
			\begin{align*}
				H^1(\Sigma_{\K} M) := \coprod_{g \in \Rm(M)}{H^1(\Sigma^g_{\K} M)} \subset L^2(\Sigma_{\K} M) \to \Rm(M)
			\end{align*}
			are a domain subbundle (see \cref{DefDomainSubbundle}) of $L^2(\Sigma_{\K} M)$ and $\bar \beta_{\_, g}$ is a trivialization adapted to $H^1(\Sigma^g_{\K} M)$.
	\end{enumerate}
\end{Thm}

\vspace{1em}

\begin{Prf} $ $ 
	\begin{enumerate}
		\item 
			We just define the topology by such a global trivialization. We have to show that for any two reference metrics $g,g'$ the transition functions
				\DefMap{ \beta_{\_,g'} \circ \beta_{\_,g}^{-1} : \Rm(M) \times \Gamma(\Sigma^g_{\K} M)}{\Rm(M) \times \Gamma(\Sigma^{g'}_{\K} M)}{(h,\psi)}{(h,\beta_{h,g'}( \beta_{g,h}(\psi))),}
			are continuous. But this was already done in \cref{LemBetaDef}.
		\item
			Analogously,  
				\DefMap{ \bar \beta_{\_,g'} \circ \bar \beta_{\_,g}^{-1} : \Rm(M) \times L^2(\Sigma^g_{\K} M)}{\Rm(M) \times L^2(\Sigma^{g'}_{\K} M)}{(h,\psi)}{(h,\bar \beta_{h,g'}(\bar \beta_{g,h}(\psi))),}
			is continuous by \cref{CorBetabargh}. 
		\item
			$\bar \beta_{\_, g}$ is an adapted trivialization (as required by \cref{DefDomainSubbundle}) by \cref{CorBetabargh}\ref{ItH1L2Adapted}. 
	\end{enumerate}
\end{Prf}

The results of this section can also be described in a category theoretic language.

\begin{Def}
	The \emph{Universal Spinor Field Bundle} is a functor
	\begin{align*}
		\fun{USBF}_{\K}: \cat{SpinMfd} \to \cat{HilB}_{\K}
	\end{align*}
	defined	
	\begin{category}
		\item[on objects by] 
		\begin{align*}
			\fun{USBF}_{\K}(M,\Theta) := (L^2(\Sigma_{\K} M) \to \Rm(M)),
		\end{align*}
		where $L^2(\Sigma_{\K} M) \to \Rm(M)$ is the universal spinor field bundle as in \cref{DefUniSpinorFieldBdle}.
		\item[on morphisms by] For any morphism of spin manifolds
			\begin{align*}
				\xymatrixcolsep{3.5em}
				\xymatrix{
					\GLtp M_1
						\ar[r]^-{\tilde F}
						\ar[d]^-{\Theta_1}
					&\GLtp M_2
						\ar[d]^-{\Theta_2}
					\\
					\GLp M_1
						\ar[r]^-{f_*}
						\ar[d]
					&\GLp M_2
						\ar[d]
					\\
					M_1
						\ar[r]^-{f}
					&M_2,
				}
			\end{align*}
			we set
			\begin{align*}
				\fun{USBF}_{\K}((f,\tilde F)) := ((f^{-1})^*,\bar F),
			\end{align*}
			where 
				\DefMap{\bar F:L^2(\Sigma_{\K} M_1)}{L^2(\Sigma_{\K} M_2)}{{\psi=[\tilde b, v])}}{{[\tilde F(\tilde b),v]}}
			and obtain a morphism
			\begin{align*}
				\xymatrixcolsep{3.5em}
				\xymatrix{
					L^2(\Sigma_{\K} M_1)
						\ar[r]^-{\bar F}
						\ar[d]
					&L^2(\Sigma_{\K} M_2)
						\ar[d]
					\\
					\Rm(M_1)
						\ar[r]^-{(f^{-1})^*}
					&\Rm(M_2)
				}
			\end{align*}
			of Hilbert bundles.
	\end{category}	
\end{Def}

%% file: usb.univdirac.tex
\section{The Universal Dirac Operator}
\label{SecUnivDiracOperator}

Using the identification maps $\bar \beta_{g,h}$ constructed in \cref{CorBetabargh} one can pull back the Dirac operator $\Dirac^h_{\K}$ to $L^2(\Sigma^g_{\K} M)$. 

\begin{Thm}[\protect{\cite[Thm. 20]{BourgGaud}}]
	\label{ThmDiracghLocal}
	The operator
	\begin{align*}
		\Dirac^{h}_{g} := \bar \beta_{h,g} \circ \Dirac^h_{\K} \circ \bar \beta_{g,h}: \Gamma(\Sigma^g_{\K} M) \to \Gamma(\Sigma^g_{\K} M)
	\end{align*}
	has the following local coordinate representation: If $(e_1, \ldots, e_n) \in \SO^g U $ and $\psi \in \Gamma(\Sigma^g_{\K} M)$, then on $U$
	\begin{align} 
		\label{EqDiracghLocal}
		\begin{split}
			\Dirac^{h}_{g} \psi & = \sum_{i=1}^{n}{e_i \cdot \nabla^g_{b_{g,h}e_i}\psi}
			+\frac{1}{4} \sum_{i,j=1}^{n}{e_i \cdot e_j \cdot \left( b_{h,g} (\nabla^h_{b_{g,h}(e_i)}(b_{g,h}e_j)) - \nabla^g_{b_{g,h}(e_i)}e_j \right) \cdot \psi} \\
			&- \frac{1}{f_{g,h}} b_{g,h}(\grad^g(f_{g,h})) \cdot \psi
		\end{split}
	\end{align}
\end{Thm}

This classical result can also be reformulated in a bundle theoretic language. 

\begin{Thm} 
	\nomenclature[Dirac]{$\Dirac$}{universal Dirac operator}
	The Dirac operator is an unbounded morphism
	\DefMap{\Dirac: L^2(\Sigma_{\K} M)}{L^2(\Sigma_{\K} M)}{\psi \in H^1(\Sigma^g_{\K} M)}{\Dirac^g_{\K} \psi}
	with domain subbundle $H^1(\Sigma_{\K} M)$. (See \cref{SecHilbertBundles} for the definition of an unbounded morphism and a domain subbundle.)
\end{Thm}

\begin{Prf}
	As required by \cref{DefHilbertBundleUnboundedMorphism}, we have to show that $\Dirac \in \Gamma(\B(H^1(\Sigma_{\K} M), L^2(\Sigma_{\K} M)))$. By definition, it suffices to check that the map
		\DefMap{\Rm(M)}{\B(H^1(\Sigma^g_{\K} M),L^2(\Sigma^g_{\K} M))}{h}{\Dirac^h_g}
	is continuous. But this follows from \eqref{EqDiracghLocal} by using the fact that all terms on the right hand side depend $\mathcal{C}^1$-continuously on $h$. For $b_{g,h}$ this was shown in \cref{LemBetaDef} and for $\nabla^h$ this follows from the fact that the Levi-Civita connection of a metric depends only on the coordinate functions of the metric and their first order derivatives.
\end{Prf}

%% file: usb.specdec.tex
\section{Spectral Decomposition}
\label{SecUsbSpecProps}

\begin{Rem} 
	We recall that the Dirac operator $\Dirac^g_{\K}$ is a self-adjoint elliptic differential operator. Therefore, we obtain a decomposition
	\begin{align}
		\label{EqL2DecGeneral}
		L^2(\Sigma^g_{\K} M) = \underset{\lambda \in \spec \Dirac^g_{\K}}{\overline{\bigoplus}}{L^2_{\lambda}(\Sigma^g_{\K} M),}
	\end{align}
	where $L^2_{\lambda}(\Sigma^g_{\K} M)$ is the eigenspace of $\lambda$ and $\overline{\oplus}$ is the direct sum of Hilbert spaces (see for instance \cite[Thm. 5.8]{LM}). This also gives a pointwise decomposition of the universal spinor field bundle, but it does not give a well-defined decomposition of the Hilbert bundle $L^2(\Sigma_{\K} M)$ into continuous subbundles over $\Rm(M)$ of eigenspaces. The problem is that if $g$ varies, the various eigenvalues might split up or unite, hence changing the dimension of the corresponding eigenspaces. However, the following holds.
\end{Rem}
	
\begin{Thm}
	\label{ThmUSFBSpecSubs}
	Let $\Lambda_1, \Lambda_2 \in \R$, $\Lambda_1 < \Lambda_2$, and let $X \subset \Rm(M)$ be a connected subspace of Riemannian metrics satisfying
	\begin{align}
		\label{PartEigenBundleRMHyp}
		\forall g \in X: \Lambda_1, \Lambda_2 \notin \spec \Dirac^{g}_{\K}.
	\end{align}
	Then there exists an open neighborhood $U \subset \Rm(M)$ of $X$ such that the \emph{partial eigenbundle}
	\begin{align*}
		\pi^{\Sigma}_{[\Lambda_1, \Lambda_2]}: L^2_{[\Lambda_1, \Lambda_2]}(\Sigma_{\K} M) \to U
	\end{align*}
	is a continuous vector bundle of finite rank. 
\end{Thm}

\begin{Prf}
	We will use our results about continuity of Dirac spectra, see \cref{MainThmFun}, to show that the bundle has constant rank. To see that it is continuous, we represent it as the image of a continuous family of eigenprojections, see \cref{ThmSpectralSeparation}.
\begin{steplist}
	\step[finding the rank]
		By \cref{MainThmFun}, there exist continuous functions $\lambda_j:\Rm(M) \to \R$, $j \in \Z$, such that for each $g \in \Rm(M)$, the sequence $(\lambda_j(g))_{j \in \Z}$ is non-decreasing and represents the entire Dirac spectrum (counted with $\K$-multiplicities). By continuity, $\lambda_j(X) \subset \R$ is connected and therefore an interval for any $j \in \Z$. Consequently, \cref{PartEigenBundleRMHyp} implies that either $\lambda_j(X) \subset \mathopen{]} \Lambda_1, \Lambda_2 \mathclose{[}$ or $\lambda_j(X) \cap \mathopen{]} \Lambda_1, \Lambda_2 \mathclose{[}  = \emptyset$. Since all Dirac spectra are discrete and unbounded, the former is satisfied by only finitely many functions. After shifting the enumeration if necessary, we can assume that these functions are given by $\lambda_1, \ldots, \lambda_k$ for some suitable $k \in \N$. These functions represent the Dirac spectrum of all the metrics $g \in X$ in $[\Lambda_1, \Lambda_2]$. We claim that this $k$ is the rank of $\pi^{\Sigma}_{[\Lambda_1, \Lambda_2]}$.
	\step[finding $U$]
		Now, fix any $g \in X$. Since the Dirac spectrum is discrete and $\Lambda_1, \Lambda_2 \notin \spec \Dirac^g_{\K}$, there exists $\varepsilon > 0$ such that 
		\begin{align}
			\label{EqDiracDiscreteExploit}
			I_{\varepsilon}(\spec \Dirac^g_{\K} \cap [\Lambda_1, \Lambda_2]) \subset \mathopen{]} \Lambda_1, \Lambda_2 \mathclose{[} , &&
			I_{\varepsilon}(\spec \Dirac^g_{\K} \setminus [\Lambda_1, \Lambda_2]) \subset \R \setminus [\Lambda_1, \Lambda_2].
		\end{align}
		By continuity, there exists an open neighborhood $U_g$ of $g$ in $\Rm(M)$ such that 
		\begin{align*}
			\forall h \in U_g: & \forall j \in \{1, \ldots, k\}: \lambda_j(h) \in [\Lambda_1, \Lambda_2] \\
			& \forall j \in \Z \setminus \{1, \ldots, k\}: \lambda_j(h) \notin [\Lambda_1, \Lambda_2]
		\end{align*}
		Consequently, at each metric $h \in U_g$, the vector space
		\begin{align*}
			L^2_{[\Lambda_1, \Lambda_2]}(\Sigma^h_{\K} M) = \sum_{j=1}^k{L^2_{\lambda_j(h)}(\Sigma^h_{\K} M)}
		\end{align*}
		have the same rank $k$. Define $U := \bigcup_{g \in X}{U_g}$.
	\step[continuity]
		To see that these vector spaces vary continuously with $h$, take the global trivialization
		\begin{align*}
			\bar \beta_{\_, g}: L^2(\Sigma_{\K} M) \to \Rm(M) \times L^2(\Sigma^g_{\K} M)
		\end{align*}
		from \cref{LemUnivSpinBdleTrivs}. Since $\spec \Dirac^{h}_g = \spec \Dirac^h_{\K}$ by \cref{CorMaierIdSpectralProp} below, the image $\bar \beta_{\_, g}(L^2_{[\Lambda_1, \Lambda_2]}(\Sigma_{\K} U)) =: E_{\K}$ satisfies
		\begin{align*}
			\forall h \in U: E_{\K}|_h
			=\{h\} \times \Lin \{\psi \in H^1(\Sigma^g_{\K} M) \mid \exists \lambda \in [\Lambda_1, \Lambda_2]: \Dirac^h_g \psi = \lambda \psi \} 
		\end{align*}
		and it suffices to check that $E_{\K} \to U$ is a continuous vector bundle of rank $k$. 
		
		Now, for $\K=\C$, we argue as follows: By discreteness of $\spec \Dirac^g_{\C}$, there exists a simple closed curve $\gamma: [0,1] \to \C$ such that $\spec \Dirac^g_{\C} \cap [\Lambda_1, \Lambda_2]$ is in the interior of $\gamma$ and $\spec \Dirac^g_{\C} \setminus [\Lambda_1, \Lambda_2]$ is in the exterior of $\gamma$.  By \cref{ThmSpectralSeparation} and \cref{CorSpectralEigenProjection}, the map
			\DefMap{P:U}{L^2(\Sigma^g_{\C} M)}{h}{P_h := \oint_{\gamma}{(z-\Dirac^h_g)^{-1}dz}}
		is a projection onto $E_{\K}|_h$ of rank $k$ for each $h \in U$. Since $P$ is a continuous map, the result now follows from \cref{ThmImageBundle}.
	\step[real case]
		In case $\K=\R$, we apply \cref{ThmSpectralSeparation,CorSpectralEigenProjection,ThmImageBundle} to the complexification of the real Dirac operator $(\Dirac^g_{\R})^{\C}$. It follows that $E_{\R}^{\C}$ is a continuous bundle and therefore, $E_{\R}$ is a continuous bundle as well.
\end{steplist}
\end{Prf}

\begin{Lem}[spectral properties] 
	\label{CorMaierIdSpectralProp}
	The operators $\Dirac^h_g$ from \cref{ThmDiracghLocal} satisfy
	\begin{align*}
		\spec \Dirac^g_h = \spec \Dirac^h_{\K}, &&
		\bar \beta_{h,g} ( \ker (\Dirac^h_{\K} - \lambda) ) = \ker (\Dirac^h_g - \lambda).
	\end{align*}
\end{Lem}

\begin{Prf}
	For any $\psi \in \Gamma( \Sigma^h_{\K} M)$,
	\begin{align*}
		&\Dirac^h_{\K} \psi = \lambda \psi \\
		\Longrightarrow &(\Dirac^h_{\K} \circ \bar \beta_{g,h} \circ \bar \beta_{h,g})(\psi) = \lambda (\bar \beta_{g,h} \circ \bar \beta_{h,g})(\psi) \\
		\Longrightarrow &(\bar \beta_{h,g} \circ \Dirac^h_{\K} \circ \bar \beta_{g,h} \circ \bar \beta_{h,g})(\psi) = \lambda (\bar \beta_{h,g} \circ \bar \beta_{g,h} \circ \bar \beta_{h,g})(\psi) \\		
		\Longrightarrow & \Dirac^g_h( \bar \beta_{h,g}(\psi)) = \lambda \bar \beta_{h,g}(\psi),
	\end{align*}
	and vice versa.
\end{Prf}

%% file: higher.abstract.tex
\begin{subabstract}
	In this chapter, we will prove the existence of metrics, for which the Dirac operator has at least one eigenvalue of multiplicity at least two (in dimensions $m \equiv 0,6,7 \mod 8$), see \cref{MainThmHigher}. We will prove this by ``catching'' the desired metric in the space of all Riemannian metrics using the ``Lasso Lemma'', see \cref{LemLasso}. We will show how such a metric can be ``caught'' on the sphere with a loop of metrics induced by a family of rotations. Finally, we will transport this loop to an arbitrary manifold of suitable dimension by extending some known results concerning the surgery of spin manifolds. 
\end{subabstract}

%% file: higher.intro.tex
\section{Introduction and Statement of the Results}
\label{SectHigherIntro}

In this chapter, we will carry out the details of the proof of \cref{MainThmHigher} and use the strategy explained in \cref{SubSectPrfStrategyMainThmHigher}. We start by stating the \emph{Lasso Lemma}, see also 
\cref{FigLasso}. 

\begin{Lem}[Lasso Lemma]
	\label{LemLasso}
	Let $X$ be a simply connected topological space and $A \subset X$ be any subspace. Let $\gamma:\S^1 \to A$ be a loop and $E \to A$ be a real vector bundle such that $\gamma^* E \to \S^1$ is not orientable. Then $X \setminus A$ is not empty. 
\end{Lem}

\begin{Prf} 	
	Let $\iota:A \hookrightarrow X$ be the inclusion. Since $X$ is simply connected, $\iota \circ \gamma$ is null-homotopic. Consequently, there exists an extension $\hat \gamma:\bar \D^2 \to X$ such that $\hat \gamma|_{\S^1} = \iota \circ \gamma$. Now, if $\hat \gamma(\D^2)$ were contained in $A$, then $\hat \gamma$ would be null-homotopic in $A$, since $\bar \D^2$ is contractible. But then $\gamma$ would also be null-homotopic, thus $\gamma^* E \to \S^1$ would be trivial, hence orientable. Contradiction! Consequently, there exists $x \in \bar \D^2$ such that $\hat \gamma(x) \in X \setminus A$. In particular, $X \setminus A \neq \emptyset$. 
\end{Prf}

Recall that we want to apply this reasoning in the case where $X=(\Rm(M),\mathcal{C}^1)$ is the space of Riemannian metrics and $A$ is a subspace of metrics $\Rm_A(M)$, which will be specified in \cref{DefPartEBundleGlobal}. The space $\Rm_A(M)$ is a complicated subspace built in a way such that $\Rm(M) \setminus \Rm_A(M) \neq \emptyset$ will directly imply \cref{MainThmHigher}. For the vector bundle, we will use a bundle $E=E(M)$ consisting of a certain span of eigenspinors corresponding to a finite set of eigenvalues, see \cref{CorPartEBGlobal}.

As we said before, the construction of the loop $\gamma$, which is then a loop $\mathbf{g}:\S^1 \to \Rm_A(M)$ of Riemannian metrics will be the hard part. To that end, we will introduce a notion of an \emph{odd loop of metrics $\mathbf{g}:\S^1 \to \Rm_A(M)$}, see \cref{DefOddLoopMetrics}. These loops will have the property that $\mathbf{g}^*E(M) \to \S^1$ is not trivial. Hence, \cref{LemLasso} will imply that $\Rm(M) \setminus \Rm_A(M) \neq \emptyset$ and $M$ will \emph{admit metrics of higher multiplicities}, see \cref{DefMfdHigherMults}. In \cref{ThmOddLoopsSphere} we will show that the sphere admits an odd loop of metrics and in \cref{ThmSurgeryStabilityOdd} we will show that this property is stable under certain surgeries. 

%% file: higher.orientors.tex
\section{Triviality of Vector Bundles over $\S^1$}
\label{SectHigherOrientors}

Ultimately, we want to apply \cref{LemLasso} and therefore, we will have to verify that a certain vector bundle over $\S^1$ is not trivial. As a preparation, we will give a neat criterion to check this, see \cref{LemOrientorFrame}. The question whether or not a given vector bundle is trivial has of course been studied in much more generality in the context of characteristic classes and $K$-Theory. However, the set of isomorphism classes of vector bundles of rank $n$ over $\S^1$ has only two elements, see for instance \cite[p.25]{hatcher}. One class represents the trivial, hence orientable bundles, the other class consists of vector bundles that are non-orientable, hence non-trivial. Therefore, it turns out that for classifying vector bundles over $\S^1$, it is simpler if we go without the general abstract machinery. As a byproduct, we give an alternative proof of the classification of vector bundles over $\S^1$. We will also have to show that the class of such a vector bundle remains constant for bundles ``nearby'', see \cref{ThmSignStability}.

It will be convenient to view $\S^1$ as $\S^1 = [0,1] / \sim$, where $\sim$ just identifies $0$ with $1$. We denote by $[0]=[1]$ the equivalence class of $\{0,1\}$ in $\S^1$. The projection will always be denoted by%
\nomenclature[S1]{$\S^1$}{$[0,1] / \sim $}%
\nomenclature[piS1]{$\pi_{\S^1}$}{$\pi_{\S^1}: [0,1] \to \S^1$}%
\begin{align*}%
		\pi_{\S^1}:[0,1] \to \S^1.
\end{align*} %

\begin{Lem}[frame curves]
	\nomenclature[sgnPsi]{$\sgn(\Psi)$}{sign of a frame curve}
	\nomenclature[sgnE]{$\sgn(E)$}{sign of a vector bundle over $\S^1$}
	\label{LemOrientorFrame}%
	Let $\pi_{E}: E \to \S^1$ be a real vector bundle of rank $n$. There exists a global section $\Psi \in \Gamma(\pi_{\S^1}^* \GL E)$ and for any such $\Psi$ there exists a unique matrix $A \in \GL_n$ such that $\Psi(1)=\Psi(0).A$. This matrix satisfies	
	\begin{align*}
		\det(A) > 0 && 
		\Longleftrightarrow &&  \text{$E$ is orientable}
		&&\Longleftrightarrow && \text{$E$ is trivial}.
	\end{align*}
	Such a section $\Psi$ will be called a \emphi{frame curve} for $E$, $A$ is called its \emph{sign matrix} and we set $\sgn(E) := \sgn(\Psi) := \sgn(\det(A))$. \index{sign!of a vector bundle over $\S^1$}\index{sign!of a frame curve}
\end{Lem}

\begin{Prf} 
	Consider the pull-back diagram
	\begin{align}
		\label{EqFrameCurveExists}
		\begin{split}
			\xymatrix{
				\pi_{\S^1}^*( \GL E)
					\ar[r]^-{F}
					\ar[d]
				& \GL E
					\ar[d]^-{\pi_{\GL E}}
				\\
				[0,1]
					\ar[r]^-{\pi_{\S^1}}
					\ar@/^/[u]^{\Psi}
				& \S^1.
			}
		\end{split}
	\end{align}
	Since every fibre bundle over the contractible space $[0,1]$ is trivial (see for instance \cite[Chpt. 4, Cor. 10.3]{husemoller}), $\pi_{\S^1}^*(\GL E)$ admits a global section $\Psi$. Existence and uniqueness of $A$ follow from the fact that the action of $\GL_n$ is simply-transitive on the fibres of the principal $\GL_n$-bundle $\pi_{\S^1}^*\GL E$ and that by construction, $F(\Psi(0))$ and $F(\Psi(1))$ are in the same fibre. 
	
	If $E$ is trivial, then $E$ is orientable. If $E$ is orientable, then $\pi_{\S^1}^* \GL^+ E$ admits a global section $s$ and the curve $\Psi := F \circ s$ has a sign matrix $A \in \GLp_n$. Finally, let $\Psi \in \Gamma(\pi_{\S^1}^* \GL E)$ be a section with sign matrix $A_0 \in \GLp_n$. Since $\GLp_n$ is connected, there exists a continuous curve $A:[0,1] \to \GLp_n$ such that $A(0)=I_n$ and $A(1)=A^{-1}_0$. The curve $\Psi' := \Psi.A$ satisfies
	\begin{align*}
		\Psi'(0)=\Psi(0).A(0) = \Psi(0), &&
		\Psi'(1)=\Psi(1).A(1) = \Psi(0).
	\end{align*}
	Consequently, $F \circ \Psi'$ descends to a global section $\bar \Psi:\S^1 \to \GLp E$ and $E$ is trivial.
\end{Prf}

\begin{Rem}
	In case $E$ carries a metric, it will sometimes be convenient to study curves of orthonormal frames. In that case $\Psi \in \Gamma(\pi_{S^1}^* \Or E)$, where $\Or E$ is the principal $\Or_n$-bundle of orthonormal frames on $E$. The sign matrix $A$ of $\Psi$ will then be in $\Or_n$ and the bundle is trivial if and only if $A \in \SO_n$.
\end{Rem}

\begin{Rem}
	For any principal $G$-bundle $P \to M$ and any smooth map $f:N \to M$, we can consider the pullback $Q:=f^*P$ and the corresponding diagram
	\begin{align*}
		\xymatrix{
			Q
				\ar[r]^-{F}
				\ar[d]^-{\pi_{Q}}
			&P
				\ar[d]^-{\pi_P}
			\\
			N
				\ar[r]^-{f}
			&M.	
		}
	\end{align*}
	By construction of the pullback, for any section $s \in \Gamma(Q)$, we obtain a map $s':=F \circ s:N \to P$, which satisfies
	\begin{align}
		\label{EqPullBackCharSects}
		\pi_P \circ s' = f.
	\end{align}
	Conversely, any map $s':N \to P$, which satisfies \cref{EqPullBackCharSects}, gives a section $s \in \Gamma(Q)$ by setting $s := (\id_N, s')$. This correspondence is one-to-one and we will identify sections of $\pi_{\S^1}^*(\GL E)$ with curves $\Psi:[0,1] \to \GL E$ satisfying 
	\begin{align}
		\label{DefOrientor}
		\pi_{\GL E} \circ \Psi = \pi_{\S^1}.
	\end{align}
	$ $
\end{Rem}

To make precise the claim that the sign of a vector bundle as defined in \cref{LemOrientorFrame} does not change for bundles ``nearby'', we study vector bundles that are subbundles of a Hilbert bundle (see \cref{DefVectBanHilBundle} for a precise definition of a Hilbert bundle).

\begin{Thm}[sign stability]
	\label{ThmSignStability}
	Let $E,\tilde E \to \S^1$ be two $k$-dimensional subbundles of a Hilbert bundle $\mathcal{H} \to \S^1$ with induced metric. Denote by $S \tilde E \to \S^1$ the bundle of unit spheres in $\tilde E$. If
	\begin{align}
		\label{EqSignStabDistHyp}
		\forall \alpha \in \S^1: d(E_{\alpha}, S \tilde E_{\alpha}) < 1,
	\end{align}
	then $\sgn(E) = \sgn(\tilde E)$.
\end{Thm}

\begin{Prf} $ $
	\begin{steplist}
		\step[construct frame curves]
			Let $\tilde \psi \in \Gamma(\pi_{S^1}^* \Or \tilde E)$ be a curve of orthonormal frames for $\tilde E$. We define
			\begin{align*}
				\forall 1 \leq j \leq k: \forall \alpha \in S^1: \psi^{\alpha}_j := P_{\alpha}(\tilde \psi^{\alpha}_j) \in E_{\alpha},
			\end{align*}
			where $P_{\alpha}:\mathcal{H}_{\alpha} \to \mathcal{H}_{\alpha}$ denotes the orthogonal projection onto $E_{\alpha}$. This defines continuous sections $\psi_j \in \Gamma(\mathcal{H})$, $j=1, \ldots, k$. We claim that for each $\alpha \in S^1$, the $\psi^{\alpha}_1, \ldots, \psi^{\alpha}_k$ are a basis for $E_{\alpha}$. So let $c_1, \ldots, c_k \in \R$ such that
			\begin{align*}
				0 = \sum_{j=1}^{k}{c_j \psi^{\alpha}_j}
				=\sum_{j=1}^{k}{c_j P_{\alpha}(\tilde \psi^{\alpha}_j)}
				=P_{\alpha} \Big{(} \underbrace{\sum_{j=1}^{k}{c_j \tilde \psi^{\alpha}_j}}_{=:\tilde v} \Big{)}
			\end{align*}
			We claim that $\tilde v = 0$. Otherwise (by dividing the equation above by the norm of $\tilde v$) we can assume that $\tilde v$ has unit length. Since $P_{\alpha}$ is the orthogonal projection onto $E_{\alpha}$, the equation $P_{\alpha}(\tilde v) = 0$ simply means that $\tilde v$ is perpendicular to $E_{\alpha}$. This implies
			\begin{align*}
				d(E_{\alpha}, \tilde v) 
				=\|P_{\alpha}(\tilde v) - \tilde v\|
				=\|\tilde v\|
				=1,
			\end{align*}
			which contradicts our assumption \cref{EqSignStabDistHyp}. Consequently, we obtain $\tilde v = 0$. Since the $\tilde \psi_j^{\alpha}$ are a basis, this implies $c_1=\ldots=c_k=0$. This implies that $\Psi \in \Gamma{\pi_{S^1}^* \GL E}$ is a frame curve for $E$.
		\step[compare signs]
			Let $\tilde A, A \in \GL_k$ be the sign matrices of $\tilde \psi$ respectively $\psi$, i.e.
			\begin{align*}
				\tilde \psi(1) = \tilde \psi(0).\tilde A, && \psi(1) = \psi(0).A.
			\end{align*}
			By construction,
			\begin{align*}
				\psi(1)
				= P_{1}(\tilde \psi (1))
				= P_{1}(\tilde \psi (0). \tilde A)
				= P_{0}(\tilde \psi (0)).\tilde A
				= \psi(0).\tilde A,
			\end{align*}
			thus $A = \tilde A$ and $\sgn(E)=\sgn(\tilde E)$ as claimed.
	\end{steplist}
\end{Prf}

\begin{Rem}
	As indicated in the proof, condition \cref{EqSignStabDistHyp} is equivalent to the fact that no non-zero element in $\tilde E_{\alpha}$ is  perpendicular to $E_{\alpha}$.
\end{Rem}

%% file: higher.spindiffeos.tex
\section{Loops of Spin Diffeomorphisms}
\label{SecLoopsSpinDiffeos}

In this section, we introduce a technique to produce certain loops of metrics via loops of spin diffeomorphisms. On all mapping spaces like $\Diff M$ or $\G(\GLp M)$, we will always use the $\mathcal{C}^{\infty}_w$-topology, see \cref{DefWeakTopology}.

\subsection{Loops}

\begin{DefI}[spin loop]
	A \emph{loop of spin diffeomorphisms} is a continuous map $f:\S^1 \to \Diff^{\spin}(M)$. 
\end{DefI}

\begin{Rem}
	By \cref{ThmCOkTop}, we can identify such a spin loop $f$ with a continuous map $f:\S^1 \times M \to M$, such that for each $\alpha \in \S^1$, $f_{\alpha} := f(\alpha, \_) \in \Diff^{\spin}(M)$, and vice versa. We do not distinguish between $f$ as a loop $\S^1 \to \Diff^{\spin}(M)$ and as a map $\S^1 \times M \to M$.
\end{Rem}

\begin{Rem}[loops as isotopies]
	\index{associated isotopy}
	Any loop $f: \S^1 \to \Diff^{\spin}(M)$ defines an \emph{associated isotopy} $h:= f \circ \pi_{\S^1}: I \to \Diff^{\spin}(M)$. Recall from \cref{RemSpinLiftIsotopies} that spin isotopies can be lifted to the spin structure via
	\begin{align}
		\label{EqSpinIsotopyLiftVia}
		\begin{split}
			\xymatrix{
				& \G(\GLtp M)
					\ar[d]^-{\Theta_*}_{2:1}
				\\
				I
					\ar[r]_-{H}
					\ar@{..>}[ur]^-{\hat H}
				& \G^{\spin}(\GLp M),
			}
		\end{split}
	\end{align}
	where $H_t = {h_t}_*$. Any lift $\hat H:I \to \G(\GLtp M)$ satisfies $\Theta_*(\hat H_1) = \Theta_*(\hat H_0)$, thus $\hat H_1 = \pm \hat H_0$. In particular, $\hat H$ does \underline{not} have to be a loop again. 
\end{Rem}

\begin{Def}[sign of a loop]
	\label{DefSgnLoopDiff}
	Let $f:\S^1 \to \Diff^{\spin}(M)$ be a loop, $h := f \circ \pi_{\S^1}: I \to \Diff^{\spin}(M)$ be its associated isotopy and $\hat H$ be a lift as in \cref{EqSpinIsotopyLiftVia}. The unique number $\sgn(f) \in \{\pm 1\}$ such that 
	\begin{align}
		\label{EqRotoidOddOrEven}
		 \hat H_1 = {\sgn(f)} \hat H_0
	\end{align}
	is called the \emph{sign}\index{sign!of a spin loop} of $f$. We say $f$ is \emph{even}\index{even!loop of diffeomorphisms}, if $\sgn(f)=+1$ and \emph{odd}\index{odd!loop of diffeomorphisms}, if $\sgn(f)=-1$. 
\end{Def}
\nomenclature[sgnf]{$\sgn(f)$}{sign of a spin loop}

Notice that $\sgn(f)$ does not depend on the lift $\hat H$ chosen to define it. The sign has the following abstract characterization.

\begin{Lem}
	\label{LemSgnAbstract}
	The map $\Theta_*:\G(\GLtp M) \to \G^{\spin}(\GLp M)$ is a principal $\Z_2$-bundle. The connecting homomorphism $\delta$ from its long exact homotopy sequence
	\begin{align*}
		\xymatrix{
			\pi_1(\G(\GLtp M), \id_{\GLtp M})
							\ar[r]^-{{\Theta_*}_{\sharp}}
			&\pi_1(\G^{\spin}(\GLp M), \id_{\GLp M})
				\ar[r]^-{\delta}
			&\pi_0(\{ \pm {\id}_{\GLtp M} \}) 
		}
	\end{align*}
	satisfies $\delta(f_*) = \sgn(f) \id_{\GLtp M}$ for any loop $f:\S^1 \to \Diff^{\spin}(M)$.
\end{Lem}

\begin{Prf}
	Since $\Theta_*$ is a $2:1$-covering, it is a normal covering, hence a principal $\Z_2$-fibre bundle. Recall the definition of $\delta$ in this case\footnote{For the general case, see for instance \cite[Def. 17.3]{steenrod}}: Let $\gamma \in [\gamma] \in \pi_1(\G^{\spin}(\GLp M), \id_{\GLp M})$ and denote by $\iota:\S^0 \hookrightarrow \D^1=[0,1]$, $+1, -1 \mapsto 0,1$, the canonical inclusion. Let $\hat \gamma:I \to \G(\GLtp M)$ be the lift of $\gamma \circ \pi_{\S^1}:I \to \G^{\spin}(\GLp M)$ as in \cref{EqSpinIsotopyLiftVia}:
	\begin{align*}
		\xymatrix{
			(\S^0, +1)
				\ar@{..>}[rr]^-{\delta(\gamma)}
				\ar[d]^-{\iota}
			&& \G(\GLtp M, \id_{\GLtp M})
				\ar[d]^-{\Theta_*}_-{2:1}
			\\
			(I, \partial I)
				\ar[r]_-{\pi_{\S^1}}
				\ar@{..>}[urr]^-{\hat \gamma}
			&(\S^1, [0])
				\ar[r]_-{\gamma}
			&\G^{\spin}(\GLp M, \id_{\GLp M})
		}
	\end{align*}
	 The map $\pi_{\S^1} \circ \iota$ is constant. Therefore, $\Theta_* \circ \hat \gamma \circ \iota = \gamma \circ \pi_{\S^1} \circ \iota$ is also constant. Consequently, the image of $\hat \gamma \circ \iota$ is in the fibre and there exists a map of pointed spaces $\delta(\gamma):(\S^0,+1) \to (\{ \pm \id_{\GLtp M} \}, \id_{\GLtp M})$ such that $\widehat \gamma \circ \iota = \delta(\gamma)$. Now, there are only two such maps $(\S^0,+1) \to (\{ \pm \id_{\GLtp M} \}, \id_{\GLtp M})$ characterized by mapping $-1$ to $\pm \id_{\GLtp M}$. By construction, $\delta(\gamma)(-1) = \sgn(\gamma) \id_{\GLtp M}$.
\end{Prf}
 
\begin{Lem}
	The sign induces a group homomorphism
	\begin{align*}
		\sgn:\pi_1(\Diff^{\spin}(M), \id_M) \to \Z_2
	\end{align*}
	and $\ker \sgn$ are precisely the homotopy classes of even loops.
\end{Lem}

\begin{Prf}
	It follows from \cref{LemSgnAbstract} that $\sgn$ is well-defined on homotopy classes. To see that $\sgn$ is a group homomorphism, let $f^{(1)}, f^{(2)} \in \pi_1(\Diff^{\spin}(M), \id_M)$ and consider
		\DefMap{f := f^{(2)} * f^{(1)}:\S^1}{\Diff^{\spin}(M)}{t}{ 
			\begin{cases}
				f^{(1)}(2t), & 0 \leq t \leq \tfrac{1}{2}, \\
				f^{(2)}(2t-1), & \tfrac{1}{2} \leq t \leq 1.
			\end{cases}
		}
	Let $\hat f^{(1)}:I \to \G(\GLtp M)$ be a lift of $f^{(1)}_*:I \to \G^{\spin}(\GLp M)$ such that $\hat f^{(1)}(0)=\id_{\GLtp M}$ and $\hat f^{(1)}(1)=\sgn(f^{(1)}) \id_{\GLtp M}$. Let $\hat f^{(2)}$ be a lift of $f^{(2)}_*$ such that $\hat f^{(2)}(0)=\hat f^{(1)}(1)$. Then $\hat f = \hat f^{(2)} * \hat f^{(1)}$ is a lift of $f_*$ such that $\hat f(0)=\id_{\GLtp M}$. Consequently, 
	\begin{align*}
		\hat f(1)
		&=\hat f^{(2)}(1)
		=\sgn(f^{(2)}) \hat f^{(2)}(0) 
		=\sgn(f^{(2)}) \hat f^{(1)}(1) \\
		&=\sgn(f^{(2)}) \sgn(f^{(1)}) f^{(1)}(0)
		=\sgn(f^{(2)}) \sgn(f^{(1)}) \id_{\GLtp M}
	\end{align*}
	and $\sgn(f) = \sgn(f^{(2)}) \sgn(f^{(1)})$. Clearly, the constant map $\S^1 \to \Diff^{\spin}(M)$, $\alpha \mapsto \id_M$, lifts to $\id_{\GLtp M}$, so $\sgn$ is a group homomorphism as claimed.
\end{Prf}

\begin{Rem}
	Although one cannot just replace the isotopy $h$ associated to the loop $f$ by the loop $f$ itself in \cref{EqSpinIsotopyLiftVia}, there always exists a lift $\hat f$ such that 
	\begin{align*}
		\begin{split}
			\xymatrix{
				\S^1
					\ar[r]^-{\hat f}
					\ar[d]^-{\cdot 2}
				& \G(\GLtp M)
					\ar[d]^-{\Theta_*}_{2:1}
				\\
				\S^1
					\ar[r]_-{f_*}
					\ar@[gray]@{-->}[ur]^-{\color{gray}\hat F}
				& \G^{\spin}(\GLp M)
			}
		\end{split}
	\end{align*}
	commutes. Here $\cdot 2$ denotes the double cover of $\S^1$. This follows from the fact that $f \circ \cdot 2 = f * f$, thus $\sgn(f \circ \cdot 2) = \sgn(f)^2 = +1$ and thus, $f \circ \cdot 2$ is even. Nevertheless, the lift $\hat F$ exists if and only if $f$ itself is even.
\end{Rem}

\begin{Rem}
	By construction, one can always complete \cref{EqSpinIsotopyLiftVia} to
	\begin{align}
		\label{EqSpinLoopLift}
		\begin{split}
			\xymatrix{
				I
					\ar[r]^-{\hat H}
					\ar[d]^-{\pi_{\S^1}}
				& \G(\GLtp M)
					\ar[d]^-{\Theta_*}_{2:1}
				\\
				\S^1
					\ar[r]_-{f_*}
					\ar@[gray]@{-->}[ur]^-{\color{gray}\hat F}
				& \G^{\spin}(\GLp M).
			}
		\end{split}
	\end{align}
	Again, the map $\hat F$ exists if and only if $f$ is even, but $\hat H$ exists always.
\end{Rem}

\begin{Rem}
	In the special case where $f:\S^1 \to \Diff^{\spin}(M)$ is a group action, the notions of odd and even in the sense of \cref{DefSgnLoopDiff} coincide with the notions of an odd respectively even group action in the sense of \cite[IV.\textsection 3, p.295]{LM}.
\end{Rem}

\subsection{Associated Loops of Metrics}

Next, we study how loops of diffeomorphisms induce loops of Riemannian metrics.

\begin{Def}[associated loops of metrics]
	\label{DefAssocSpinLoop}
	Let $f:\S^1 \to \Diff^{\spin}(M)$ be a loop of spin diffeomorphisms and $g \in \Rm(M)$ be any Riemannian metric. The family of metrics 
		\DefMap{\mathbf{g}:\S^1}{\Rm(M)}{\alpha}{g_\alpha := (f_\alpha^{-1})^*g}
	is called an \emph{associated loop of metrics}. \index{loop of metrics}
\end{Def}

\begin{Rem}
	\label{RemAssocMetricsSpinIsom}
	In the situation of \cref{DefAssocSpinLoop}, for any $\alpha \in \S^1$, the map 
	\begin{align*}
		f_\alpha:(M,g) \to (M,g_\alpha)
	\end{align*}
	is a Riemannian spin isometry by construction (recall \cref{RemSpinIsometries}). As a consequence, all $(M, g_{\alpha})$, $\alpha \in \S^1$, are Dirac isospectral and 
	\begin{align*}
		\spec \Dirac^{\mathbf{g}}_{\K} := \spec \Dirac^{g_\alpha}_{\K}, \alpha \in \S^1,
	\end{align*}
	is well-defined. Furthermore, we get an isometry
		\DefMap{\bar H_t: \Sigma^{g_0}_{\K} M}{\Sigma^{g_t}_{\K} M}{{\psi=[s,v]}}{{[\hat H_t(s),v]}}
	between all the spinor bundles for any $t \in [0,1]$. Here, $\hat H$ is as in \cref{EqSpinLoopLift}. The induced map on sections, denoted by $\bar H_{t}$, satisfies
	\begin{align}
		\label{EqSpinDiffeoDiracCommute}
		\Dirac^{g_{t}}_{\K} \circ \bar H_{t} = \bar H_{t} \circ \Dirac^{g_0}_{\K}
	\end{align}
	and therefore maps eigenspinors to eigenspinors.
	\nomenclature[specmbfg]{$\spec \Dirac^{\mathbf{g}}$}{spectrum of a loop}
\end{Rem}

The following will be crucial to verify the hypothesis of \cref{LemLasso}.

\begin{Thm}
	\label{ThmS1ActionBundleTwisted}
	Let $A \subset \Rm(M)$ be any subset, $f:\S^1 \to \Diff^{\spin}(M)$ be a loop of spin diffeomorphisms and $g \in \Rm(M)$ be such that the associated loop of metrics $\alpha \mapsto (f_{\alpha}^{-1})^* g$ is a map $\mathbf{g}:\S^1 \to A$. Furthermore, let $E \subset L^2(\Sigma_{\R} M) \to A$ be a vector bundle of rank $n$. Let $\bar H_t$ be the map induced by $f$ as in \cref{RemAssocMetricsSpinIsom} and assume $H_t(E) \subset E$ for any $t \in [0,1]$.
	Then
	\begin{align*}
		\sgn(\mathbf{g}^*E) =
		\begin{cases}
			-1, & \text{$f$ is odd and $n$ is odd}, \\
			+1, & \text{otherwise}.
		\end{cases}
	\end{align*}
\end{Thm}

\begin{Prf} 
	For any basis $(0,(\psi_1, \ldots, \psi_n)) \in \mathbf{g}^*E|_0$, the curve 
		\DefMap{\Psi:[0,1]}{\GL(\mathbf{g}^*E)}{t}{(t,(\bar H_{t}(\psi_1), \ldots, \bar H_{t}(\psi_n)))}
	is a frame curve for $\mathbf{g}^*E \to \S^1$ (see \cref{LemOrientorFrame}). For any $p \in M$ and any $\psi_j|_p =: [\tilde b,v] \in \Sigma_{\R}^{g_0} M|_p$, we calculate
	\begin{align*}
		\bar H_1(\psi_j|_p) 
		=\bar H_1([\tilde b,v])
		=[\hat H_1(\tilde b),v]
		\jeq{\eqref{EqRotoidOddOrEven}}{=} [\sgn(f) \tilde b,v]
		=\sgn(f) \psi_j|_p.
	\end{align*}
	Consequently, the sign matrix $A \in \GL_n$ of $\Psi$ (see \cref{LemOrientorFrame}) is given by $A = \sgn(f) I_n$, which has determinant $\sgn(f)^{n}$. By definition,
	\begin{align*}
		\sgn(\mathbf{g}^*E)
		=\sgn(\Psi) 
		=\sgn(\det(A))
		= \sgn(\sgn(f)^n),
	\end{align*}
	which implies the result.
\end{Prf}

\subsection{Metrics of higher Multiplicities}

We are now in a position to establish a relationship between loops of spin diffeomorphisms and the existence of metrics $g$ for which $\Dirac^g_{\K}$ has an eigenvalue of higher multiplicity. 

\begin{Def}[manifolds admitting higher multiplicities]
	\label{DefMfdHigherMults}
	A spin manifold $(M,\Theta)$ \emph{admits a metric of higher $\K$-multiplicities}, if there exists a metric $g$ on $M$ such that $\Dirac^g_{\K}$ has at least one eigenvalue $\lambda$ of multiplicity $\mu_{\K}(\lambda) \geq 2$. \index{metric of higher multiplicities}
\end{Def}

\begin{Def}
	\label{DefPartEBundleGlobal}
	Let $\{ \lambda_j: \Rm(M) \to \R \}_{j \in \Z}$ be a non-decreasing sequence of eigenvalue functions as in \cref{MainThmFun}. For a fixed $n \geq 1$, we define 
	\begin{align*}
		\Rm_A(M) &:= \{g \in \Rm(M) \mid \exists 1 \leq j \leq n: \mu_{\R}(\lambda_j(g)) \in 2 \N +1, \lambda_{0}(g) < \lambda_1(g), \lambda_{n}(g) < \lambda_{n+1}(g) \}, \\
		E^g(M) &:=  \Lin \{ \psi \in H^1(\Sigma^g_{\R} M) \mid \exists 1 \leq j \leq n: \Dirac^g_{\R} \psi = \lambda_j(g) \psi \} , \qquad
		E(M) := \bigcup_{g \in \Rm_A(M)}{E^g(M)}.
	\end{align*}
\end{Def}

\nomenclature[RAM]{$\Rm_{A}(M)$}{a subset of Riemannian metrics}
\nomenclature[EM]{$E(M)$}{a finite dimensional bundle over $\Rm_{A}(M)$}

\begin{Cor}
	\label{CorPartEBGlobal}
	In the situation of \cref{DefPartEBundleGlobal}, the map $\pi:E(M) \to \mathcal{R}_{A}(M)$, $\psi \in E^g \mapsto g$, is a continuous vector bundle of rank $n$.
\end{Cor}

\begin{Prf}
	Let $g \in \Rm_{A}(M)$. By construction, there exist $\Lambda_1, \Lambda_2 \in \R$ such that
	\begin{align*}
		\lambda_{0}(g) < \Lambda_1 < \lambda_1(g) \leq \ldots \leq \lambda_n(g) < \Lambda_2 < \lambda_{n+1}(g).
	\end{align*}
	By continuity, there exists an open and connected neighborhood $U$ of $g$ in $\Rm(M)$ such that 
	\begin{align*}
		\forall h \in U: \spec \Dirac^h_{\R} \cap [\Lambda_1, \Lambda_2] = \{\lambda_1(h) \leq \ldots \leq \lambda_n(h)\}.
	\end{align*}
	By \cref{ThmUSFBSpecSubs}, $L^2_{[\Lambda_1, \Lambda_2]}(\Sigma_{\R} M) \to U$ is a continuous vector bundle of rank $n$ and by construction $L^2_{[\Lambda_1, \Lambda_2]}(\Sigma_{\R} M)|_{U \cap \Rm_{A}(M)} = E(M)|_{U \cap \Rm_{A}(M)}$. 
\end{Prf}

\begin{Def}[odd loop]
	\label{DefOddLoopMetrics}
	A closed spin manifold $(M, \Theta)$ \emph{admits an odd loop of metrics}\index{odd!loop of metrics}, if there exists an odd number $n \in \N$, $\Lambda_1 < \Lambda_2$ and a continuous map
	\begin{align*}
		\mathbf{g}:(\S^1, \tau_{\S^1}) \to (\Rm_{A}(M), \mathcal{C}^2)
	\end{align*}
	such that
	\begin{enumerate}
		\item $\forall \alpha \in \S^1: \lambda_0(g_{\alpha}) < \Lambda_1 < \lambda_1(g_{\alpha}) \leq \lambda_n(g_{\alpha}) < \Lambda_2 < \lambda_{n+1}(g_{\alpha})$ and
		\item the bundle $E(M) \to \Rm_{\A}(M)$ from \cref{CorPartEBGlobal} satisfies $\sgn(\mathbf{g}^*E(M)) = -1$.
	\end{enumerate}
	Here, $\{ \lambda_j\}_{j \in \Z}$ is the same enumeration of eigenvalue functions as in \cref{DefPartEBundleGlobal}.
\end{Def}

\begin{Thm}[existence of Dirac eigenvalues of higher multiplicity] 
	\label{ThmExistsHigherMult} 
	Let $(M,\Theta)$ be a closed spin manifold. If $(M, \Theta)$ admits an odd loop of metrics, then $(M, \Theta)$ admits a metric of higher $\R$-multiplicities. In case $m \equiv 0,6,7 \mod 8$, $M$ also admits a metric of higher $\C$-multiplicities.
\end{Thm}

\begin{Prf}
	By hypothesis, $\sgn(\mathbf{g}^*E(M)=-1)$, so $\mathbf{g}^*E(M) \to \S^1$ is not trivial by \cref{LemOrientorFrame}. By \cref{LemLasso}, 
	\begin{align*}
		\emptyset \neq \Rm_A(M)^c 
		& = \{g \in \Rm(M) \mid \forall 1 \leq j \leq n: \mu_{\R}(\lambda_j(g)) \in 2 \N \text{ or } \lambda_{0}(g) = \lambda_1(g) \text{ or } \lambda_{n}(g) = \lambda_{n+1}(g) \} \\
		& \subset \{g \in \Rm(M) \mid \exists 1 \leq j \leq n: \mu_{\R}(\lambda_j(g)) \geq 2 \}.
	\end{align*}
	This implies the claim for $\R$-multiplicities. This in turn implies the claim for $\C$-multiplicities by \cref{RemMultK}.
\end{Prf}
 
\begin{Cor}
	\label{CorOddSpinLoopMetrics}
	Let $f:\S^1 \to \Diff^{\spin}(M)$ be an odd loop of spin diffeomorphisms. If $M$ admits a metric $g$ such that one eigenvalue $\lambda$ of $\Dirac^g_{\R}$ is of odd $\R$-multiplicity, then $M$ admits an odd loop of metrics. In particular, $M$ admits higher $\R$-multiplicities. In case $m \equiv 0,6,7 \mod 8$, $M$ also admits higher $\C$-multiplicities.
\end{Cor}

\begin{Prf}
	Let $\mathbf{g}:\S^1 \to \Rm(M)$ be the loop of metrics induced by $f$ and $g$. Let $(\lambda_j)_{j \in \Z}$ be a sequence of eigenvalue functions as in \cref{MainThmFun}. By shifting this sequence if necessary, we can assume that $\lambda_0(g_0) < \lambda_1(g_0) = \ldots = \lambda_n(g_0) < \lambda_{n+1}(g_0)$, for some odd $n \in \N$. Since all $(M,g_{\alpha})$, $\alpha \in \S^1$, are isospectral, the same relation holds for all $g_{\alpha}$. Therefore, $\mathbf{g}$ is actually a map $\S^1 \to \Rm_A(M)$, where $\Rm_{A}(M)$ is as in \cref{DefPartEBundleGlobal}. Let $E(M) \to \Rm_{A}(M)$ be the vector bundle from \cref{DefPartEBundleGlobal}. Since $n$ is odd and $f$ is odd, we obtain $\sgn(\mathbf{g}^*E(M))=-1$ by \cref{ThmS1ActionBundleTwisted}. Consequently, $\mathbf{g}$ is an odd loop of metrics. This implies the result by \cref{ThmExistsHigherMult}.
\end{Prf}

%% file: higher.sphere.tex
\section{The Sphere}
\label{SectSphere}
In this section, we consider one of the most obvious loops of spin diffeomorphisms, namely the rotations $R_{2 \pi \alpha}$ of angle $2 \pi \alpha$, $\alpha \in \S^1$, on the sphere $\S^m$, $m \geq 3$. In particular, we illustrate that the notions developed in \cref{SecLoopsSpinDiffeos} are not empty. In \cref{LemRotationsOddLoop}, we will show that $f:\S^1 \to \Diff^{\spin}(\S^m)$, $\alpha \mapsto R_{2 \pi \alpha}$, is an odd loop of spin diffeomorphisms. Since $R_{2 \pi \alpha}$ is an isometry with respect to the round metric $g\degree$, the induced loop of metrics is constant, hence not very interesting. But in \cref{ThmNbhdRoundSphereOdd}, we will prove that small pertubations of the round metric yield metrics $g$ such that the induced loop of metrics $(R_{2 \pi \alpha}^{-1})^* g$ is indeed an odd loop of metrics in the sense of \cref{DefOddLoopMetrics}. We first fix this as the main result of this section and then attend to the proof of the auxiliary claims, see \cref{SubSectRotSpinLoop,SubSectNbhdRoundSphere}.

\begin{Thm}[odd loops on the sphere]
	\label{ThmOddLoopsSphere}
	The sphere $(\S^m, \Theta)$, $m \equiv 0,6,7 \mod 8$, admits an odd loop of metrics.
\end{Thm}

\begin{Prf} 
	By \cref{ThmNbhdRoundSphereOdd} and \cref{RemMultK}, there exists a metric $g$ on $\S^m$ for which at least one eigenvalue of $\Dirac^g_{\R}$ has odd $\R$-multiplicity. By \cref{LemRotationsOddLoop}, the rotations are an odd loop of spin diffeomorphisms on $\S^m$. \cref{CorOddSpinLoopMetrics} implies the result.
\end{Prf}

\begin{Rem}
	\cref{ThmOddLoopsSphere} immediately implies that $\S^m$ admits metrics of higher multiplicities by \cref{ThmExistsHigherMult}. But the importance of \cref{ThmOddLoopsSphere} is not at all that it implies this result. It is well known, see \cref{ThmDiracSpecRoundSphere} below, that lots of Dirac eigenvalues of $(\S^m, g\degree)$ are of higher multiplicity. The importance is that \cref{ThmExistsHigherMult} provides a sufficient condition for the existence of metrics of higher multiplicities and \cref{ThmOddLoopsSphere} shows that this condition is satisfied on the sphere. In \cref{SectSurgeryStability}, we will show that odd loops of metrics are stable under certain surgeries, in particular under the connected sum operation. Therefore, if some manifold $N$ admits an odd loop of metrics, then $M \sharp N$ also admits such a loop. Since $M \sharp \S^m \cong M$, we will prove that any manifold $M$ admits higher $\C$-multiplicities for $m \equiv 0,6,7 \mod 8$.
\end{Rem}

\begin{Thm}[Dirac spectrum of the round sphere, see for instance {\cite[Thm. 1]{BaerSphere}}]
	\label{ThmDiracSpecRoundSphere}
	The Dirac operator on the round sphere $(\S^m,g\degree)$, $m \geq 2$, satisfies 
	\begin{align*}
		\spec \Dirac^{g\degree}_{\C} = \{ \lambda_k^\pm := \pm \left( \frac{m}{2} + k \right) \mid k \in \N \}
	\end{align*}
	and each $\lambda_k^\pm$ has $\C$-multiplicity $\mu_k := 2^{\lfloor \frac{m}{2} \rfloor} \binom{m+k-1}{k}$.
\end{Thm} 

\subsection{The Loop of Rotations on the Sphere}
\label{SubSectRotSpinLoop}

\begin{DefI}[rotation]
	\nomenclature[Ralpha]{$R_{\alpha}$}{rotation by an angle $\alpha$}
	\label{DefRotations}
	For any $\alpha \in \R$, $m \in \N$, we define the rotation
	\begin{align*}
		R_{\alpha}:=
		\begin{pmatrix}
			I_{m-1} & 0 & 0 \\
			0 & \cos(\alpha) & -\sin(\alpha) \\
			0 & \sin(\alpha) & \cos(\alpha) 
		\end{pmatrix}:\R^{m+1} \to \R^{m+1}.
	\end{align*}
\end{DefI}

\begin{Lem}
	\label{LemRotationsOddLoop}
	The map $f:\S^1 \to \Diff^{\spin}(\S^m)$, $\alpha \mapsto R_{2 \pi \alpha}|_{\S^m}$,  is an odd loop of spin diffeomorphisms.
\end{Lem}

\begin{Prf}
	Chose the round metric $g\degree$ on $\S^m$. In this case, the spin structure on $\S^m$ is simply given by
	\begin{align*}
		\Theta := \vartheta_{m+1}: \Spin_{m+1} \to \SO_{m+1},
	\end{align*}
	see \cref{ThmSpinStructureRoundSphere} for details. By a tedious calculation carried out in \cref{LemRoundSphereLift}, we see that the map
		\DefMap{\hat f: \R}{\G(\Spin(\S^m))}{\alpha}{v \mapsto (\cos(\tfrac{\alpha}{2}) + \sin(\tfrac{\alpha}{2})e_{m-1}e_m)v}
	is a lift of $\R \to \Diff^{\spin}(\S^m)$, $\alpha \mapsto R_{\alpha}|_{\S^m}$. It follows from this explicit formula that
	\begin{align*}
		\hat f(2 \pi)
		&= \cos(\pi) + \sin(\pi)e_{m-1}e_m \\
		&= -1 
		=- \Big(  \cos(0) + \sin(0)e_{m-1}e_m \Big) 
		=- \hat f(0).
	\end{align*}
	This implies that $f$ is odd.
\end{Prf}

\subsection{The Odd Neighborhood Theorem}
\label{SubSectNbhdRoundSphere}

\begin{Thm}[odd neighborhood theorem]
\label{ThmNbhdRoundSphereOdd}
	Let $(M,\Theta)$ be a closed spin manifold of dimension $m \equiv 0,6,7 \mod 8$ and $g_0$ be any Riemannian metric on $M$. In every $\mathcal{C}^1$-neighborhood of $g_0 \in \Rm(M)$, there exists $g \in \Rm(M)$ such that $\Dirac^g_{\C}$ has an eigenvalue $\lambda$ of odd multiplicity $\mu_{\C}(\lambda)$. 
\end{Thm}
		
	\begin{figure}[t] 
		\begin{center}
			\input{fig.oddnbhd}
			\caption{Finding an odd metric near $g_0$.}
			\label{FigOddnbhd}
		\end{center}
	\end{figure}
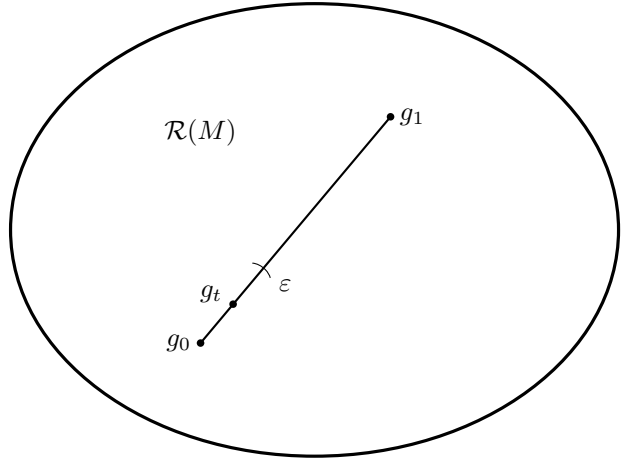

\begin{figure}[t] 
	\begin{center}
		\input{fig.nbhdsphere}
		\caption[Reduction to a finite index set.]{On $I_{\delta}(t_0)$, the $\lambda_i$, $i \in N_n(t_0)$, intersect only with one another. }
		\label{FigRoundSphereOdd}
	\end{center}
\end{figure}
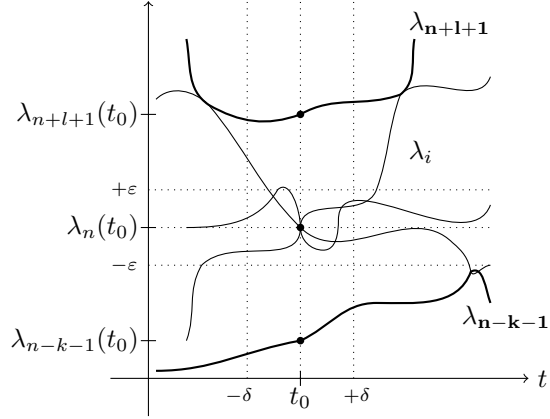

\begin{Prf}
	The idea of this proof is as follows: By \cref{ThmDahl} and \cref{RemMultK}, there exists a metric $g_1 \in \Rm(M)$ such that $\Dirac^{g_1}_{\C}$ has an eigenvalue of $\C$-multiplicity $1$, which is odd. Connect the metric $g_1$ with $g_0$, i.e. define the path $g_t := tg_1 + (1-t)g_0$, $t \in [0,1]$. This path is real-analytic. By \cite[Lem. A.0.16]{AndreasDiss}, the Dirac operators $\Dirac^{g_0, g_t}_{\C}$ are the restriction of a self-adjoint holomorphic family of type (A) onto $I$. Therefore, the eigenvalues of $\Dirac^{g_t}_{\C}$ can be described by a real-analytic family of functions $\{ \lambda_n:[0,1] \to \R \}_{n \in \N}$, see \cref{ThmAnalyticEVKato}. This means that for any $t \in [0,1]$, the sequence $(\lambda_n(t))_{n \in \N}$ represents all the eigenvalues of $\Dirac^{g_t}_{\C}$ counted with $\C$-multiplicity (but possibly not ordered by magnitude).
	To prove the claim, the argument will be by contradiction and roughly go as follows: If this claim is wrong, there exists an open neighborhood around $0$ in which all metrics $g_t$ have even multiplicities. Since the eigenvalue functions $\lambda_n$ are real analytic, this behavior extends to all of $[0,1]$ and therefore, all the eigenvalues of $g_1$ have even multiplicity as well. This contradicts the choice of $g_1$. To make this argument precise, we need a technical Lemma, which is discussed in \cref{LemRealAnalyticPartners} below. 
	\begin{steplist}
	\step[precise claim]
		We want to show: For any $\varepsilon > 0$, there exists $0 < t < \varepsilon$ such that $\Dirac^{g_t}_{\C}$ has an eigenvalue of odd multiplicity, see \cref{FigOddnbhd}. By contradiction assume that there exists $\varepsilon > 0$ such that for all $0 < t < \varepsilon $ all eigenvalues of $\Dirac^{g_t}_{\C}$ are of even multiplicity, i.e. if $\mu_n(t)$  is the $\C$-multiplicity of $\lambda_n(t)$, we obtain
		\begin{align}
			\label{EqRounsSphereOddMultHyp}
			\forall n \in \N: \forall t \in \mathopen{]} 0,\varepsilon \mathclose{[}: \mu_n(t) \in 2 \N. 
		\end{align}

	\step[verify hypothesis of Lemma] 
		We want to apply \cref{LemRealAnalyticPartners} and therefore have to verify its hypothesis. It remains only to check that we can replace the infinite index set $\N$ in \cref{EqRounsSphereOddMultHyp} by a finite index set. To that end, we define for any $n \in \N$, the index sets
			\DefMap{N_n:I}{\mathcal{P}(\N)}{t}{\{ i \in \N \mid \lambda_i(t) = \lambda_n(t) \}.}
		Fix any $n \in \N$ and $t_0 \in \mathopen{]} 0, \varepsilon \mathclose{[}$. We want to show that the functions $\{ \lambda_i \mid i \in N_n(t_0) \}$ satisfy the hypothesis of \cref{LemRealAnalyticPartners} on a small neighborhood around $t_0$. By construction, this is an even number of functions. Renumber the eigenvalues such that $(\lambda_j(t_0))_{j \in \N}$ in non-decreasing. We have to show that their multiplicity within the index set 
		\begin{align*}
			N_n(t_0) =: \{n-k, \ldots, n \ldots, n+l\}
		\end{align*}
		is also even, i.e.
		\begin{align*}
			\exists \delta > 0: \forall t \in I_\delta(t_0): \tilde \mu_n(t) := | \{ i \in N_n(t_0) \mid \lambda_i(t) = \lambda_n(t) \} | \in 2 \N.
		\end{align*}
		 Since Dirac eigenvalues are discrete, there exists $\varepsilon > 0$ such that the sets
		\begin{align*}
			I_{\varepsilon}(\lambda_{n-k-1}(t_0)), && 
			I_{\varepsilon}(\lambda_n(t_0)), &&
			I_{\varepsilon}(\lambda_{n+l+1}(t_0))
		\end{align*}
		are mutually disjoint. By continuity of the eigenvalue functions, see \cref{MainThmFun}, there exists $\delta > 0$ such that
		\begin{align*}
			\forall i \in N_n(t_0): \forall t \in I_{\delta}(t_0): \tilde \mu_i(t) = \mu_i(t) \in 2 \N,
		\end{align*}
		i.e. on $I_{\delta}(t_0)$ the eigenvalue functions $\lambda_{n-k}(t), \ldots, \lambda_{n+l}(t)$ intersect only with each other, see \cref{FigRoundSphereOdd}. Therefore, \cref{LemRealAnalyticPartners} is applicable.

	\step[final argument]
		By the conclusion of \cref{LemRealAnalyticPartners}, there exists a smallest index $m \neq n$ such that $\lambda_m = \lambda_{n}$ on an open neighborhood around $t_0$. But since all these functions are real analytic, this implies $\lambda_m \equiv \lambda_{n}$. This in turn implies that
		\begin{align*}
		\forall n \in \N: \forall t \in [0,1]: \mu_n(t) \in 2 \N.
		\end{align*}
		This is a contradiction to the fact that $\Dirac^{g_1}_{\R}$ has an eigenvalue of multiplicity $1 \notin  2 \N$.
	\end{steplist}
\end{Prf}

\begin{Lem}[real analytic partners]
	\label{LemRealAnalyticPartners}
	Let $\lambda_1, \ldots , \lambda_N \in \mathcal{C}^{\omega}(I,\R)$, $N$ even, where $I$ is an open interval containing $0$. Assume that 
	\begin{align*}
		\forall 1 \leq m \leq N: & \forall t \in I: \mu_m(t) := | \{ i \in \N \mid 1 \leq i \leq N, \lambda_i(t) = \lambda_m(t) \}| \in 2 \N.
	\end{align*}
	Then there exists $2 \leq m \leq N$ such that $\lambda_1 \equiv \lambda_m$. This claim is also true, if $I$ is replaced by any interval $t_0 + I$, $t_0 \in \R$. 
\end{Lem}

\begin{Prf}
	The idea is to use a Taylor series expansion of the functions. Since these are real analytic, it suffices to check that
	\begin{align*}
		\exists 2 \leq m \leq N: \forall k \in \N: \lambda^{(k)}_m(0) := \frac{d^k \lambda_m}{dt^k}(0) = \lambda^{(k)}_1(0).
	\end{align*}
	By contradiction assume this is false, i.e.
	\begin{align} \label{EqRealAnalyticContr}
		\forall 2 \leq m \leq N: \exists k \in \N: \lambda^{(k)}_m(0) \neq \lambda^{(k)}_1(0).
	\end{align}
	Let $2 \leq m \leq N$ be arbitrary and choose the minimal $k$ such that \eqref{EqRealAnalyticContr} holds. Expand all the $\lambda_m$ into a Taylor series 
	\begin{align*}
		\lambda_m(t) = \sum_{\nu=0}^{k-1}{a_{m,\nu}t^\nu} + a_{m,k}t^k + R_m(t), && a_{m,\nu} \in \N.
	\end{align*}
	By assumption 
	\begin{align*}
		\forall 0 \leq \nu \leq k-1: a_{m,\nu} = a_{1,\nu}, && |a_{m,k} - a_{1,k}| =: \delta > 0.
	\end{align*}
	The residue $R_m$ satisfies
	\begin{align*}
		\exists \varepsilon_m > 0: \forall t \in I \cap I_{\varepsilon_m}(0): \left| \frac{R_m(t)}{t^k} \right| < \frac{\delta}{4}.
	\end{align*}
	We can assume that the same estimate holds for $R_1$ as well. In case $a_{1,k} - a_{m,k} > 0$, we calculate for any $t \in I$, $0 < t < \varepsilon_m$,
	\begin{align*}
		\lambda_1(t) - \lambda_m(t)
		&=(a_{1,k} - a_{m,k})t^k + R_1(t) - R_m(t) \\
		& \geq \delta t^k - |R_1(t)| - |R_m(t)| \\
		&=(\delta - |R_1(t)t^{-k}| - |R_m(t)t^{-k}|)t^k \\
		&\geq \frac{\delta}{2} t^k > 0.
	\end{align*}
	In case $a_{m,k} - a_{1,k} > 0$, we do the same calculation with $\lambda_m(t) - \lambda_1(t)$. In any case, we obtain that $\lambda_1(t) \neq \lambda_m(t)$. \\
	All in all, defining
	\begin{align*}
		\varepsilon := \min_{2 \leq m \leq N}{\varepsilon_m}, &&
		U_\varepsilon := I \cap \R_{>0} \cap I_{\varepsilon}(0) \neq \emptyset,
	\end{align*}
	we conclude
	\begin{align*}
		\forall 2 \leq m \leq N: \forall t \in U_\varepsilon:  \lambda_1(t) \neq \lambda_m(t).
	\end{align*}
	Consequently, any $t \in U_\varepsilon$ satisfies $\mu_1(t) = 1 \notin 2 \N$. Contradiction!
\end{Prf}

%% file: fig.oddnbhd.tex
\begin{tikzpicture}

\tikzset{
	partial ellipse/.style args={#1:#2:#3}{
		insert path={+ (#1:#3) arc (#1:#2:#3)}
	}
}

\tikzset{
	decoration={
		markings,
			mark=between positions 0.3 and 0.90 step 1 with {
			\draw[thin] (0,0) [partial ellipse=50:-50:0.125 and 0.2];
			\coordinate[label=right:$\varepsilon$] (eps) at (-0,-0.2);
			\fill[black] (-0.5,0) circle (0.05);
			\coordinate[label=left:$g_t$] (t) at (-0.4,0.1);
		},
	}
}

	\draw[black,very thick]
		(0,0) ellipse (4 and 3);

	\draw[thick,postaction={decorate}] (-1.5,-1.5) to (1,1.5);

	\fill[black] (-1.5,-1.5) circle (0.05);
	\fill[black] (1,1.5) circle (0.05);

	\coordinate[label=above:$\mathcal{R}(M)$] (RM) at (-1.5,1);
	\coordinate[label=left:$g_0$] (A) at (-1.5,-1.5);
	\coordinate[label=right:$g_1$] (gamma) at (1,1.5);

\end{tikzpicture}

%% file: fig.nbhdsphere.tex
\begin{tikzpicture}
	\draw[->] (-0.5,0) -- (5,0) node[right] {$t$};
	
	\draw[->] (0,-0.5) -- (0,5) node[above] {};
	
	\draw[color=black] (2,0.1) -- (2,-0.1);
	\coordinate[label=below:$t_0$] (t0) at (2,0);
	\coordinate[label=below:${}_{+\delta}$] (t0) at (2.75,0);
	\coordinate[label=below:${}_{-\delta}$] (t0) at (1.2,0);

	\draw[color=black] (0.1,2) -- (-0.1,2);
	\coordinate[label=left:$\lambda_n(t_0)$] (t0) at (0,2);
	\coordinate[label=left:${}_{+\varepsilon}$] (t0) at (0,2.5);
	\coordinate[label=left:${}_{-\varepsilon}$] (t0) at (0,1.5);
	\fill (2,2) circle (0.05);
	\coordinate[label=right:$\lambda_i$] (t0) at (3.3,3);
	
	\draw[color=black] (0.1,0.5) -- (-0.1,0.5);
	\coordinate[label=left:$\lambda_{n-k-1}(t_0)$] (t0) at (0,0.5);
	\fill (2,0.5) circle (0.05);	
	\coordinate[label=right:$\mathbf{\lambda_{n-k-1}}$] (t0) at (4,0.8);
	
	\draw[color=black] (0.1,3.5) -- (-0.1,3.5);
	\coordinate[label=left:$\lambda_{n+l+1}(t_0)$] (t0) at (0,3.5);
	\fill (2,3.5) circle (0.05);
	\coordinate[label=right:$\mathbf{\lambda_{n+l+1}}$] (t0) at (3.3,4.7);

	\draw[dotted] (2,0) -- (2,5);
	\draw[dotted] (1.3,0) -- (1.3,5);
	\draw[dotted] (2.7,0) -- (2.7,5);
	
	\draw[dotted] (0,2) -- (4.5,2);
	\draw[dotted] (0,2.5) -- (4.5,2.5);
	\draw[dotted] (0,1.5) -- (4.5,1.5);
	
	\coordinate (a1) at (0.5,4.5);
	\coordinate (a2) at (1,3);
	\coordinate (b1) at (0.1,3.7);
	\coordinate (b2) at (2,3.5);
	\coordinate (c1) at (intersection of a1--a2 and b1--b2);
	
	\coordinate (a3) at (3,2.5);
	\coordinate (a4) at (3.5, 4.5);
	\coordinate (b3) at (2,3.5);
	\coordinate (b4) at (4.5, 4);
	\coordinate (c2) at (intersection of a3--a4 and b3--b4);
	
	\coordinate (a5) at (3,1);
	\coordinate (a6) at (4.5,1.5);
	\coordinate (b5) at (4,1.8);
	\coordinate (b6) at (4.5,1);
	\coordinate (c3) at (intersection of a5--a6 and b5--b6);
	
	\draw[thick] (a1)
	to [out=-80, in=145] (c1)
	to [out=-30, in=200] (2,3.5)
	to [out=30, in=200] (c2)
	to [out=40, in=-100] (a4);
	
	\draw
	(b1)
	to [in=135] (c1)
	to [out=-50, in=135] (2,2) 
	to [out=-45, in=135] (4,1.8)
	to [out=-40, in=100] (c3)
	to [out=-70, in=180] (4.5,1.5);
	
	\draw
	(0.5, 2)
	to [out=0, in=-120] (1.7,2.5)
	to [out=45, in=90] (2,2)
	to [out=-90, in=180] (2.3,1.7)
	to [out=0, in=-90] (2.5, 2.1)
	to [out=90, in=-110] (4.5, 2.3);	
	
	\draw (0.5, 0.5)
	to [out=70, in=-120]  (0.7,1.5)
	to [out=45, in=-90] (2,2)
	to [out=90, in=-110] (3,2.5)
	to [out=70, in=-110] (c2)
	to [out=45, in=-110] (b4);

	\draw[thick]
	(0.1,0.1) 
	to [out=0, in=190] (2,0.5)
	to [out=25, in=180] (a5)
	to [out=0, in=-110] (c3)
	to [out=30, in=110] (b6);

\end{tikzpicture}

%% file: higher.surgery.tex
\section{Surgery Stability}
\label{SectSurgeryStability}
The main objective of this section is to prove that an odd loop of metrics in the sense of \cref{DefOddLoopMetrics} is stable under certain \emph{surgeries}, see \cref{ThmSurgeryStabilityOdd}. To that end, we will introduce some basic notions concerning the surgery theory of spin manifolds and recall some well known results by Bär and Dahl published in \cite{BaerDahlSurgery}. 

\begin{Def}
	For any $0 < r < r'$, we define
	\begin{align*}
		\D^m(r) &:= \{x \in \R^m \mid |x| < r\} \\
		\A^m(r, r') &:= \{x \in \R^m \mid r' \leq |x| \leq r\} \\
		\S^m(r) &:= \{x \in \R^{m+1} \mid |x| = r\},
	\end{align*}
	where $| \_ |$ denotes the Euclidean norm. We also set $\S^m := \S^m(1)$ and $\D^m := \D(1)$. If $(M,\dist)$ is a metric space and $S \subset M$, we also set
	\begin{align}
		\label{EqDefUsrUsrrp}
		\begin{split}
			U_S(r) &:= \{x \in M \mid \dist(x, S) < r \} \\
			A_S(r, r') & = \{x \in M \mid r \leq \dist(x,S) \leq r' \}.
		\end{split}
	\end{align}
\end{Def}	
	\nomenclature[Dr]{$\D^m(r)$}{Euclidean ball of radius $r$}
	\nomenclature[D]{$\D^m$}{Euclidean unit ball}
	\nomenclature[Sr]{$\S^m(r)$}{sphere of radius $r$}
	\nomenclature[S]{$\S^m$}{Euclidean unit sphere}
	\nomenclature[USr]{$U_S(r)$}{$r$-neighborhood of $S$}
	\nomenclature[ASrr]{$U_S(r, r')$}{annulus around $S$ of radii $r$ and $r'$}

\begin{DefI}[surgery]
	\label{DefNormalSurgery}
	Let $M$ be a smooth $m$-manifold, let 
	\begin{align*}
		f:\S^k \times \overline{\D^{m-k}} \to M
	\end{align*}
	be a smooth embedding and set $S:=f(\S^k \times \{0\})$, $U := f(\S^k \times \D^{m-k})$. The manifold 
	\begin{align*}
		\tilde M := \Big( (M \setminus U) \amalg (\overline{\D^{k+1}} \times \S^{m-k-1})  \Big)/ \sim,
	\end{align*}
	where $\sim$ is the equivalence relation generated by 
	\begin{align*}
		\forall x \in \S^k \times S^{m-k-1}: x \sim f(x) \in \partial U,
	\end{align*}
	\emph{is obtained by surgery in dimension $k$ along $S$} from $M$. The number $m-k$ is the \emph{codimension} of the surgery. The map $f$ is the \emph{surgery map}\index{surgery!map} and $S$ is the \emph{surgery sphere}\index{surgery!sphere}.
\end{DefI}

\begin{figure}[H] 
	\begin{center}
		\input{fig.surgerydec}
		\caption[Manifold after surgery.]{The manifold after surgery. Notice that $\partial \tilde U$ is identified with $\partial U \subset (M \setminus U)$.}
		\label{FigSurgeryDec}
	\end{center}
\end{figure}
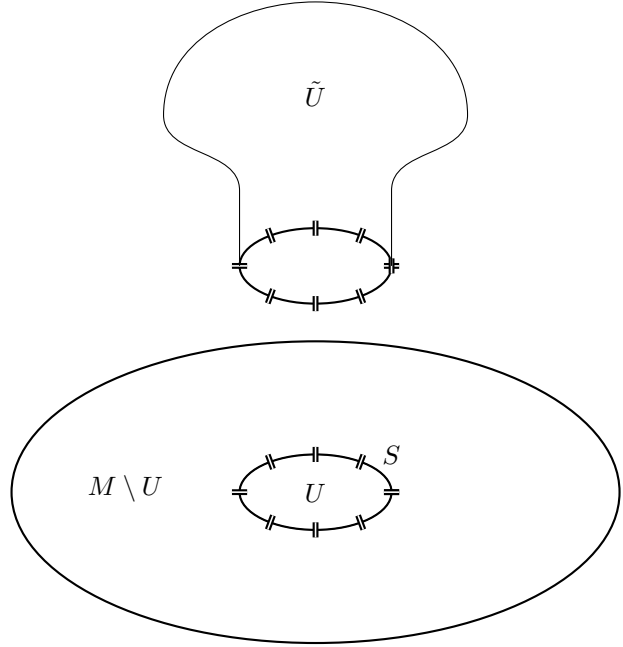

\begin{Rem}
	The space $\tilde M$ is again a smooth manifold (see for instance \cite[IV.1]{kosinski} for a very detailled discussion of the connected sum). The manifold $\tilde M$ is always of the form
	\begin{align}
		\label{EqSurgeryDecompDahl}
		\tilde M = \left( M \setminus U \right) \dot \cup \tilde U,
	\end{align}
	where $\tilde U \subset \tilde M$ is open. Here, by slight abuse of notation, $(M \setminus U) \subset M$ also denotes the image of $M \setminus U$ in the quotient $\tilde M $, see \cref{FigSurgeryDec}.
\end{Rem}

\begin{Rem}[spin structures and surgery]
	It can be shown that if one performs surgery in codimension $m-k \geq 3$, the spin structure on $M$ always extends uniquely (up to equivalence) to a spin structure on $\tilde M$, if $k \neq 1$. In case $k=1$, the boundary $S^k \times S^{m-k-1}$ has two different spin structures, but only one of them extends to $\D^{k+1} \times \S^{m-k-1}$. Adopting the convention from \cite[p. 56]{BaerDahlSurgery}, we assume that the map $f$ is chosen such that it induces the spin structure that extends. Also, we will only perform surgeries in codimension $m-k \geq 3$.
\end{Rem}

It is natural to ask how the Dirac spectra of a spin manifold before and after surgery are related to one another. For a precise statement, the following notion is useful.

\begin{Def}[$(\Lambda_1, \Lambda_2, \varepsilon)$-spectral close]
	\index{spectral close}
	Let $T: H \to H$ and $T':H' \to H'$ be two densely defined operators on Hilbert spaces $H$ and $H'$ (over $\K$) with discrete spectrum. Let $\varepsilon > 0$ and $\Lambda_1, \Lambda_2 \in \R$, $\Lambda_1 < \Lambda_2$. Then $T$ and $T'$ are \emph{\emph{$(\Lambda_1, \Lambda_2, \varepsilon)$}}-spectral close if
	\begin{enumerate}
		\item 
			$ \Lambda_1, \Lambda_2 \notin (\spec T \cup \spec T')$.
		\item
			The operators $T$ and $T'$ have the same number $n$ of eigenvalues in $\mathopen{]} \Lambda_1, \Lambda_2 \mathclose{[}$, counted with $\K$-multiplicities.
		\item
			If $\{\lambda_1 \leq \ldots \leq \lambda_n \}$ are the eigenvalues of $T$ in $\mathopen{]} \Lambda_1, \Lambda_2 \mathclose{[}$ and $\{\lambda_1' \leq \ldots \leq \lambda_n' \}$ are the eigenvalues of $T'$ in $\mathopen{]} \Lambda_1, \Lambda_2 \mathclose{[}$, then
			\begin{align*}
				\forall 1 \leq j \leq n: |\lambda_j - \lambda_j'| < \varepsilon.
			\end{align*}
	\end{enumerate}
\end{Def}

Using this terminology, a central result is the following

\begin{Thm}[\protect{\cite[Thm. 1.2]{BaerDahlSurgery}}]
	\label{ThmBaerDahlSurgeryPimped}
	Let $(M^m, g, \Theta^g)$ be a closed Riemannian spin manifold, let $0 \leq k \leq m-3$ and $f:\S^k \times \overline{\D}^{m-k} \to M$ be any surgery map with surgery sphere $S$ as in \cref{DefNormalSurgery}. For any $\varepsilon > 0$ (sufficiently small) and any $\Lambda > 0$, $\pm \Lambda \notin \spec \Dirac^g$, there exists a Riemannian spin manifold $(\tilde M^{\varepsilon}, \tilde g^{\varepsilon})$, which is obtained from $(M, g)$ by surgery such that $\Dirac^g$ and $\Dirac^{\tilde g^{\varepsilon}}$ are $(-\Lambda, \Lambda, \varepsilon)$-spectral close. This manifold is of the form $\tilde M^{\varepsilon} = (M \setminus U_{\varepsilon}) \dot \cup \tilde U_{\varepsilon}$, where $U_{\varepsilon}$ is an (arbitrarily small) neighborhood of $S$ and the metric $\tilde g^{\varepsilon}$ can be chosen such that 
	\begin{align*}
		\tilde g|_{M \setminus U_{\varepsilon}} = g|_{M \setminus U_{\varepsilon}}.
	\end{align*}
\end{Thm}

We will require not only the statement of \cref{ThmBaerDahlSurgeryPimped}, but also some arguments from the proof, which relies on estimates of certain Rayleigh quotients and these are very useful in their own right. One of the technical obstacles here is that the spinors on $M$ and on $\tilde M^{\varepsilon}$ cannot be compared directly, since they are defined on different manifolds. A simple yet effective tool to solve this problem are cut-off functions adapted to the surgery.

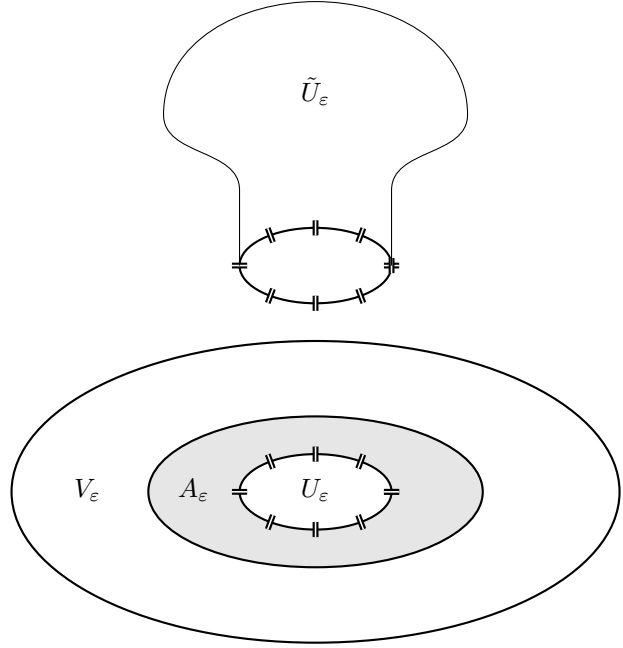
\begin{figure}[t] 
	\begin{center}
		\input{fig.surgeryprep}
		\caption[Preparing a manifold for surgery.]{Preparing a manifold $M = U_{\varepsilon} \cup A_{\varepsilon} \cup V_{\varepsilon}$ for surgery.}
		\label{FigSurgeryPrep}
	\end{center}
\end{figure}

\begin{Def}[adapted cut-off functions]
	\label{DefSurgeryCuttOffFunctions}
	In the situation of \cref{ThmBaerDahlSurgeryPimped}, assume that for each $\varepsilon > 0$ (sufficiently small), we have a decomposition of $M$ into $M = U_{\varepsilon} \dot \cup A_{\varepsilon} \dot \cup V_{\varepsilon}$, where $U_{\varepsilon} = U_S(r_{\varepsilon})$ and $A_{\varepsilon} = A_S(r_{\varepsilon}, r'_{\varepsilon})$ for some $r_{\varepsilon}, r'_{\varepsilon} > 0$ as in \cref{EqDefUsrUsrrp}. A family of cut-off functions $\chi^{\varepsilon} \in \mathcal{C}^\infty_c(M)$ is \emph{adapted} to these decompositions, if
	\begin{enumerate}
		\item 
			$0 \leq \chi^{\varepsilon} \leq 1$,
		\item
			$\chi^{\varepsilon} \equiv 0 \text{ on a neighborhood of } \bar U_{\varepsilon}$, 
		\item
			$\chi^{\varepsilon} \equiv 1 \text{ on } V_{\varepsilon}$, 
		\item
			$|\nabla \chi^{\varepsilon}| \leq \tfrac{c}{r_{\varepsilon}} \text{ on }M$ for some constant $c>0$.
	\end{enumerate}
	In case $\tilde M^{\varepsilon} = (M \setminus U_{\varepsilon}) \dot \cup \tilde U_{\varepsilon}$ is obtained from $M$ by surgery, the restriction $\chi^{\varepsilon}|_{M \setminus U_{\varepsilon}}$ can be extended smoothly by zero to a function $\chi^{\varepsilon} \in \mathcal{C}_c^\infty(\tilde M^{\varepsilon})$. The situation is depicted in \cref{FigSurgeryPrep}.
\end{Def}

\begin{Rem}[cutting off spinor fields]
	\label{RemCuttingOffSpinorFields}
	We can use the cut-off functions from \cref{DefSurgeryCuttOffFunctions} to transport spinor fields from $(M,g)$ to $(\tilde M^{\varepsilon}, \tilde g^{\varepsilon})$ and vice versa: For any $\psi \in L^2(\Sigma^g_{\K} M)$, we can think of $\chi^{\varepsilon} \psi$ as an element in $L^2(\Sigma^{\tilde g}_{\K} \tilde M^{\varepsilon})$ by extending $\chi^{\varepsilon} \psi$ to all of $\tilde M^{\varepsilon}$ by zero. Analogously, for any $\tilde \psi \in L^2(\Sigma^{\tilde g}_{\K} \tilde M^{\varepsilon})$, we can think of $\chi^{\varepsilon} \tilde \psi$ as an element in $L^2(\Sigma^g_{\K} M)$ by extending $\chi^{\varepsilon} \tilde \psi|_{M \setminus U_{\varepsilon}}$ by zero to $M$. This correspondence is not an isomorphism, but one does not loose ``too much'': It preserves smoothness and it is shown in \cite{BaerDahlSurgery} that under the assumptions of \cref{ThmBaerDahlSurgeryPimped}, for each eigenspinor $\tilde \psi^{\varepsilon} \in L^2(\Sigma^{\tilde g^{\varepsilon}} \tilde M^{\varepsilon})$ to an eigenvalue $\tilde \lambda^{\varepsilon} \in ]-\Lambda, \Lambda[$, the spinor field $\psi^{\varepsilon} := \chi^{\varepsilon} \tilde \psi^{\varepsilon} \in L^2(\Sigma^g M)$ satisfies
	\begin{align}
		\label{EqBaerDahlEstimate}
		\| \Dirac^g \psi^{\varepsilon}\|_{L^2(\Sigma^g M)} &< (\Lambda + \tfrac{\varepsilon}{2}) \|\tilde \psi^{\varepsilon}\|_{L^2(\Sigma^{\tilde g} \tilde M^{\varepsilon})}, \\
		\label{EqBaerDahlBelow}
		\|\psi^{\varepsilon}\|_{L^2(\Sigma^g M)} & \geq \frac{\Lambda + \tfrac{\varepsilon}{2}}{\Lambda + \varepsilon}  \|\tilde \psi^{\varepsilon}\|_{L^2(\Sigma^{\tilde g} \tilde M^{\varepsilon})},
	\end{align}
	The proof of \cref{EqBaerDahlEstimate,EqBaerDahlBelow} is an integral part of the proof of \cref{ThmBaerDahlSurgeryPimped}, see \cite[p. 69]{BaerDahlSurgery}. In combination, they imply the following crucial estimate for the Rayleigh quotient
	\begin{align}
		\label{EqDahlSurgeryRRBound}
		\frac{\| \Dirac^g \psi^{\varepsilon}\|_{L^2(\Sigma^g M)}^2}{\|\psi^{\varepsilon}\|_{L^2(\Sigma^g M)}^2} < (\Lambda + \varepsilon)^2.
	\end{align}
	This estimate is then used to apply the Min-Max principle, see \cref{ThmMinMaxPrinciple}.
\end{Rem}

We will require the following version of \cref{ThmBaerDahlSurgeryPimped} that is slightly more general. 

\begin{Thm}
	\label{ThmSurgeryFamilyPimped}
	Let $(M,\Theta)$ be a closed spin manifold of dimension $m \geq 3$, let $0 \leq k \leq m-3$ and $f$ be a surgery map with surgery sphere $S \subset M$ of dimension $k$ as in \cref{DefNormalSurgery}. Let $(A,\tau_A)$ be a compact topological space, 
	\begin{align*}
		\mathbf{g}:(A,\tau_A) \to (\Rm(M), \mathcal{C}^2)
	\end{align*}
	be a continuous family of Riemannian metrics and let $\Lambda_1, \Lambda_2 \in \R$, $\Lambda_1 < \Lambda_2$, such that
	\begin{align*}
		\forall \alpha \in A: \Lambda_1, \Lambda_2 \notin \spec \Dirac^{g_{\alpha}}_{\K}.
	\end{align*}
	For any $\varepsilon > 0$ (sufficiently small), there exists a spin manifold $(\tilde M^{\varepsilon}, \tilde \Theta^{\varepsilon})$, and a continuous family of Riemannian metrics
	\begin{align*}
		\mathbf{\tilde g^{\varepsilon}}:(A, \tau_A) \to (\Rm(\tilde M^{\varepsilon}), \mathcal{C}^2)
	\end{align*}
	such that for each $\alpha \in A$, the manifold $(\tilde M^{\varepsilon}, \tilde g^{\varepsilon}_{\alpha})$ is obtained from $(M, g_{\alpha})$ by surgery along $S$ and such that for all $\alpha \in A$, the operators $\Dirac^{g_{\alpha}}_{\K}$ and $\Dirac^{\tilde g^{\varepsilon}_{\alpha}}_{\K}$ are $(\Lambda_1, \Lambda_2, \varepsilon)$-spectral close. For any open neighborhood $V \subset M$ of the surgery sphere, one can choose $\mathbf{\tilde g^{\varepsilon}}$ such that
	\begin{align*}
		\forall \alpha \in A: \tilde g^{\varepsilon}_{\alpha}|_{M \setminus V} = g_{\alpha}|_{M \setminus V}.
	\end{align*}
	Moreover, if $\tilde \psi^{\varepsilon}_{\alpha} \in L^2_{[\Lambda_1, \Lambda_2]}(\Sigma^{\tilde g^{\varepsilon}_{\alpha}}_{\K} \tilde M^{\varepsilon})$, $\alpha \in A$, and $\chi^{\varepsilon}_{\alpha}$ is the cut-off function from \cref{DefSurgeryCuttOffFunctions}, the spinor field $\psi^{\varepsilon}_{\alpha} := \chi^{\varepsilon}_{\alpha} \tilde \psi^{\varepsilon}_{\alpha}$ satisfies
	\begin{align}
		\label{EqDahlSurgeryRRBoundPimped}
		\frac{\| (\Dirac^{g_{\alpha}}_{\K}-c) \psi^{\varepsilon}_{\alpha}\|_{L^2(\Sigma_{\K}^{g_{\alpha}} M)}^2}{\|\psi^{\varepsilon}_{\alpha}\|_{L^2(\Sigma^{g_{\alpha}}_{\K} M)}^2} < (l + \varepsilon)^2,
	\end{align}		
	where $c:= \tfrac{1}{2}(\Lambda_1 + \Lambda_2)$, $l := \tfrac{1}{2}|\Lambda_2 - \Lambda_1|$.
\end{Thm}

\begin{Rem}
	\cref{ThmSurgeryFamilyPimped} generalizes \cref{ThmBaerDahlSurgeryPimped} in the following ways.
	\begin{enumerate}
		\item 
			The metric $g$ is replaced by a compact $\mathcal{C}^2$-continuous family of metrics. It has already been observed by Dahl in a later paper, see \cite[Thm. 4]{DahlPresc}, that the proof of \cref{ThmBaerDahlSurgeryPimped} goes through in this case.
		\item
			The interval $[-\Lambda, \Lambda]$ is replaced by the interval $[\Lambda_1, \Lambda_2]$, which might not be symmetric around zero. This is why one has to introduce $c$ and $l$ in \cref{EqDahlSurgeryRRBoundPimped}.
		\item
			The field is $\K \in \{\R, \C\}$. This simply makes no difference in the proof.
	\end{enumerate}
	In \cref{SecSurgeryExtension}, we provide some more technical details on how to modify the proof of \cref{ThmSurgeryFamilyPimped} to obtain \cref{ThmBaerDahlSurgeryPimped}.
\end{Rem}

We are now able to prove the main result of this section.

\begin{Thm}[surgery stability]
	\label{ThmSurgeryStabilityOdd}
	Let $(M^m, \Theta)$ be a closed spin manifold, let $0 \leq k \leq m-3$ and assume that $M$ admits an odd loop of metrics as in \cref{DefOddLoopMetrics}. Let $(\tilde M, \tilde \Theta)$ be obtained from $(M, \Theta)$ by surgery in dimension $k$. Then $(\tilde M, \tilde \Theta)$ also admits an odd loop of metrics.
\end{Thm}

\begin{Prf} 
	Let $n \in \N$, $\{ \lambda_j:\Rm(M) \to \R \}_{j \in \Z}$, $\Lambda_1 < \Lambda_2$ and $\mathbf{g}:\S^1 \to \Rm_A(M)$, $\alpha \mapsto g_{\alpha}$, be the odd loop on $M$ as in \cref{DefOddLoopMetrics}. The idea is to apply \cref{ThmSurgeryFamilyPimped} and show that the resulting loop of metrics does the job. We can assume that $M$ is connected (otherwise we apply the following to a connected component of $M$).

	\begin{steplist}
		\step[apply the surgery theorem]
		By \cref{ThmSurgeryFamilyPimped}, for any $\varepsilon > 0$ (sufficiently small) and any $\alpha \in \S^1$, there exists a manifold $(\tilde M^{\varepsilon}, \tilde g^{\varepsilon}_{\alpha})$ obtained from $(M,g_{\alpha})$ by a surgery such that 
		\begin{align*}
			\mathbf{\tilde g}^{\varepsilon}:(\S^1, \tau_{\S^1}) \to (\Rm(\tilde M^{\varepsilon}), \mathcal{C}^2), \qquad 
			\alpha \mapsto \tilde g^{\varepsilon}_{\alpha},
		\end{align*}
		is continuous and such that for each $\alpha \in \S^1$, the operators $\Dirac^{g_{\alpha}}_{\R}$ and $\Dirac^{\tilde g_{\alpha}}_{\R}$ are $(\Lambda_1, \Lambda_2, \varepsilon)$-spectral close. Consequently, there exists an enumeration $\{\tilde \lambda_j^{\varepsilon}:\Rm(\tilde M^{\varepsilon}) \to \R\}_{j \in \Z}$ of the eigenvalue functions for $\tilde M^{\varepsilon}$ such that 
		\begin{align*}
			\forall \alpha \in \S^1: \spec \Dirac^{\tilde g^{\varepsilon}_{\alpha}} \cap \mathclose] \Lambda_1, \Lambda_2\mathopen[ & = \{\tilde \lambda^{\varepsilon}_1(\tilde g^{\varepsilon}_{\alpha}) \leq \ldots \leq \tilde \lambda^{\varepsilon}_n(\tilde g^{\varepsilon}_{\alpha}) \}, \\
			\forall 1 \leq j \leq n: |\lambda_j(g_{\alpha}) - \tilde \lambda_j^{\varepsilon}(\tilde g^{\varepsilon}_{\alpha})| &< \varepsilon.
		\end{align*}
		Using the same $n \in \N$ and the functions $\{\tilde \lambda^{\varepsilon}_j\}_{j \in \Z}$, we define the bundle $E(\tilde M^{\varepsilon}) \to \Rm_A(\tilde M^{\varepsilon})$ as in \cref{DefPartEBundleGlobal} for $\tilde M^{\varepsilon}$. We obtain that $\mathbf{\tilde g}^{\varepsilon}$ is actually a map $\S^1 \to \Rm_{A}(\tilde M^{\varepsilon})$. We claim that for $\varepsilon$ small enough, $\mathbf{\tilde g}^{\varepsilon}$ is an odd loop of metrics on $\tilde M^{\varepsilon}$. To prove this, it remains only to show that $\sgn((\mathbf{\tilde g}^{\varepsilon})^* E(\tilde M^{\varepsilon})) = -1$. 
		
		\step[passing from $\tilde M^{\varepsilon}$ to $M^{\varepsilon}$] 
			Let $L^2(\Sigma_{\R} M) \to \Rm(M)$ be the real universal spinor field bundle, see \cref{DefUniSpinorFieldBdle}, $E(M)$ be the bundle from \cref{DefPartEBundleGlobal} and define
			\begin{align*}
				\mathcal{H} &:= \mathbf{g}^*(L^2(\Sigma_{\R} M)) \to \S^1, \\
				E &:=  \mathbf{g}^*(E(M)) \to \S^1, \\
				\tilde E^{\varepsilon} &:= (\mathbf{\tilde g^{\varepsilon}})^*(E(\tilde M^{\varepsilon})) \to \S^1.
			\end{align*}
			Let 
			\begin{align*}
				\tilde \Psi^{\varepsilon} = (\tilde \psi^{\varepsilon}_1, \ldots, \tilde \psi^{\varepsilon}_n) \in \Gamma(\pi_{S^1}^* (\Or \tilde E^{\varepsilon}))
			\end{align*}
			be an orthonormal frame curve (see \cref{LemOrientorFrame}) for the bundle $\tilde E^{\varepsilon}$. For any $t \in I$, let $\chi^{\varepsilon}_t := \chi^{\varepsilon}_{\pi_{S^1}(t)}$, be the canonical cut-off functions from \cref{DefSurgeryCuttOffFunctions}. These functions can be chosen such that $\chi^{\varepsilon}_t$ depends continuously on $t$. We define $\psi_j^{\varepsilon}(t) := \chi^{\varepsilon}_t \tilde \psi_j^{\varepsilon}(t) \in \Gamma(\Sigma^{g^{\varepsilon}_t}_{\R} M)$, $1 \leq j \leq n$. By the unique continuation property, see \cref{LemWeakUCPLinIndep}, $\psi_1^{\varepsilon}(t), \ldots, \psi^{\varepsilon}_n(t)$ are linearly independent for all $t \in I$. Therefore, $\Psi^{\varepsilon} := (\psi^{\varepsilon}_1, \ldots, \psi^{\varepsilon}_n):I \to \mathcal{H}$ is a frame curve for the subbundle $E^{\varepsilon} \to \S^1$ spanned by $\Psi^{\varepsilon}$ in $\mathcal{H}$.  We obtain a diagram
			\begin{align*}
				\xymatrix{
					E^{\varepsilon}
						\ar[dr]
					& \tilde E^{\varepsilon}
						\ar[l]_-{\cdot \chi^{\varepsilon}}
						\ar[d]
						\ar[r]
					& E(\tilde M^{\varepsilon})
						\ar[d]
					\\
					E
						\ar[r]
					&\S^1
						\ar[r]_-{\mathbf{\tilde g}^{\varepsilon}}
					&\Rm_{A}(\tilde M^{\varepsilon}).
				}
			\end{align*}			
		\step[compare signs]
			By definition, 
			\begin{align*}
				\tilde \Psi^{\varepsilon}(1) & = \tilde \Psi^{\varepsilon}(0).\tilde A^{\varepsilon}, &&
				\Psi^{\varepsilon}(1)  = \Psi^{\varepsilon}(0). A^{\varepsilon}, 
			\end{align*}
			for some matrices $\tilde A^{\varepsilon}, A^{\varepsilon} \in \GL_n$. In coordinates, the first equation reads as
			\begin{align*}
				\forall 1 \leq j \leq n: \tilde \psi^{\varepsilon}_j(1) = \sum_{\nu=1}^{n}{\tilde \psi^{\varepsilon}_{\nu}(0) \tilde A^{\varepsilon}_{\nu j}}.
			\end{align*}
			This implies for any $1 \leq j \leq n$,
			\begin{align}
				\label{EqSignMatrixTildePsiEps}
				\psi^{\varepsilon}_j(1)
				& = \chi^{\varepsilon}_1 \tilde \psi^{\varepsilon}_j(1)
				 = \sum_{\nu=1}^{n}{\chi^{\varepsilon}_0 \tilde \psi^{\varepsilon}_{\nu}(0) \tilde A^{\varepsilon}_{\nu j}} 
				= \sum_{\nu=1}^{n}{\psi^{\varepsilon}_{\nu}(0) \tilde A^{\varepsilon}_{\nu j}},
			\end{align}
			thus $\tilde A^{\varepsilon} = A^{\varepsilon}$ and therefore $\sgn(\tilde E^{\varepsilon}) = \sgn(E^{\varepsilon})$.  By hypothesis $\sgn(E)=-1$, so it suffices to check that $\sgn(E^{\varepsilon}) = \sgn(E)$.
			
		\step[apply sign stability]
			We consider the $\lambda_j$'s as functions on $I$ by pulling them back via $\mathbf{g} \circ \pi_{S^1}$ (and analogously for $\tilde \lambda_j^{\varepsilon}$). For any $t \in I$, if $\lambda_1(t), \ldots, \lambda_n(t)$ are the eigenvalues of $\Dirac^{g_t}_{\R}$ in $[\Lambda_1, \Lambda_2]$, then $(\lambda_1(t)-c)^2, \ldots, (\lambda_n(t) - c)^2)$ are the eigenvalues of $(\Dirac^{g_t}_{\R}-c)^2$ in $[0,l^2]$, where $c := \tfrac{1}{2} \left( \Lambda_1 + \Lambda_2 \right)$ and $l := \tfrac{1}{2} |\Lambda_2 - \Lambda_1|$. The span of their collective eigenspinors is the same space $E_t$. We set $\varphi^{\varepsilon}_j(t) := \psi^{\varepsilon}_j(t) / \|\psi^{\varepsilon}_j(t)\|_{L^2(\Sigma^{g_t}_{\R}M) }$. Then $\Phi^{\varepsilon} = (\varphi^{\varepsilon}_1, \ldots, \varphi^{\varepsilon}_n)$ is still a frame curve for $E^{\varepsilon}$ and it follows from \cref{EqDahlSurgeryRRBoundPimped} that the Rayleigh quotients satisfy
			\begin{align*}
				\forall t \in I: \frac{\| (\Dirac^{g_t}_{\R} -c) \varphi_j^{\varepsilon}(t)\|_{L^2(\Sigma^{g_t}_{\R} M)}^2 }{\|\varphi_j^{\varepsilon}(t)\|_{L^2(\Sigma^{g_t}_{\R} M)}^2} < (l + \varepsilon)^2.
			\end{align*}
			Now, choose $\varepsilon$ small enough such that for all $t \in I$, we have $l + \varepsilon < \mu_{n+1}(t)$, where $\mu_{n+1}(t)$ is the $(n+1)$-th eigenvalue of $(\Dirac^{g_t}_{\R} - c)^2$. This is possible due to the continuity of $\mathbf{g}$ and since
			\begin{align*}
				\forall t \in I: \lambda_0(t) < \Lambda_1 < \lambda_1(t) \leq \ldots \leq \lambda_n(t) < \Lambda_2 < \lambda_{n+1}(t)
			\end{align*}
			by hypothesis. By \cref{ThmRayleighRitzDistance}, we obtain
			\begin{align*}
				\forall t \in I: \forall 1 \leq j \leq n: d(E_t,\varphi^{\varepsilon}_j(t))^2 \leq \frac{l + \varepsilon}{\mu_{n+1}(t)} < 1.
			\end{align*}
			By \cref{ThmSignStability}, we obtain $\sgn(E) = \sgn(E^{\varepsilon})$ and the theorem is proven.
	\end{steplist}
\end{Prf}

\begin{Lem}[weak UCP and linear independence]
	\label{LemWeakUCPLinIndep}
	Let $(M, g)$ be a closed connected Riemannian spin manifold, $\Lambda_1, \Lambda_2 \in \R$, $\Lambda_1 < \Lambda_2$, and let $ \psi_1, \ldots,  \psi_n \in L^2_{[\Lambda_1, \Lambda_2]}(\Sigma_{\K}^g M)$ be any  basis. Assume that $V \subset M$ is open and $\psi_1', \ldots, \psi_n' \in L^2(\Sigma_{\K}^g M)$ such that $\psi'_j|_V = \psi_j|_V$, $1 \leq j \leq n$. Then $\psi_1', \ldots, \psi_n'$ are linearly independent.
\end{Lem}

\begin{Prf}
	By applying a suitable isomorphism, we can assume that $ \psi_j$ is an eigenspinor to the eigenvalue $\lambda_j$, $1 \leq j \leq n$. Let $\pi_j:L^2(\Sigma_{\K}^{ g} M) \to L^2_{\lambda_j}(\Sigma_{\K}^{ g} M)$ be the orthogonal projection and let $c_1, \ldots, c_k \in \R$ such that $\psi' := \sum_{j=1}^n{c_j \psi_j'}=0$. Define $\psi := \sum_{j=1}^{n}{c_j \psi_j}$. Then $(\Dirac^{ g} - \sum_{j=1}^{n}{\lambda_j \pi_j})  \psi = 0$ and $\psi|_V = \psi'|_{V} = 0$. Since $V$ is open, we obtain $0 \equiv  \psi = \sum_{i=j}^{n}{c_j  \psi_j}$ by the weak unique continuation property for Dirac type operators, see \cite[Rem. 2.3c)]{BoossWeakUCP}. By hypothesis, this implies $c_j = 0$ for all $1 \leq j \leq n$.
\end{Prf}

%% file: fig.surgerydec.tex
\begin{tikzpicture}[
  decoration={
    markings,
    mark=between positions 0 and 1 step 0.125 with {\draw[double] (0,-0.1)--(0,0.1); },
  }
]

	\draw[black,thick]
		(0,0) ellipse (4 and 2);

	\draw
		(0,0) ellipse (1 and 0.5);
	\draw[black,thick,postaction={decorate}]
		(0,0) ellipse (1 and 0.5);

	\draw[black,thick,postaction={decorate}]
		(0,3) ellipse (1 and 0.5);		
		
	\draw (-1,3) to (-1,4)
	to [out=90, in=-90] (-2,5)
	to [out=90,in=180] (0,6.5)
	to [out=0,in=90] (2,5)
	to [out=-90, in=90] (1,4)
	to (1,3);

	\coordinate[label=above:$\tilde U$] (Ut) at (0,5);
	\coordinate[label=above:$M \setminus U$] (V) at (-2.5,-0.25);
	\coordinate[label=above:$U$] (U) at (0,-0.25);
	\coordinate[label=above:$S$] (S) at (1,0.25);

\end{tikzpicture}

%% file: fig.surgeryprep.tex
\begin{tikzpicture}[
  decoration={
    markings,
    mark=between positions 0 and 1 step 0.125 with {\draw[double] (0,-0.1)--(0,0.1); },
  }
]

	\draw[black,thick]
		(0,0) ellipse (4 and 2);
		
	\draw[fill,gray!20]
		(0,0) ellipse (2.2 and 1);
	\draw[thick,black]
		(0,0) ellipse (2.2 and 1);

	\draw[fill,white]
		(0,0) ellipse (1 and 0.5);
	\draw[black,thick,postaction={decorate}]
		(0,0) ellipse (1 and 0.5);

	\draw[black,thick,postaction={decorate}]
		(0,3) ellipse (1 and 0.5);		
		
	\draw (-1,3) to (-1,4)
	to [out=90, in=-90] (-2,5)
	to [out=90,in=180] (0,6.5)
	to [out=0,in=90] (2,5)
	to [out=-90, in=90] (1,4)
	to (1,3);

	\coordinate[label=above:$\tilde U_{\varepsilon}$] (Ut) at (0,5);
	\coordinate[label=above:$V_{\varepsilon}$] (V) at (-3,-0.25);
	\coordinate[label=above:$A_{\varepsilon}$] (A) at (-1.6,-0.25);
	\coordinate[label=above:$U_{\varepsilon}$] (U) at (0,-0.25);

\end{tikzpicture}

%% file: higher.prfmainthm.tex
\section{Proof of \cref{MainThmHigher}}

After all these preparations, the actual proof of \cref{MainThmHigher} becomes very short.

\begin{recall}[\cref{MainThmHigher}]
	\input{mainthm.higher}
\end{recall}

\begin{Prf}[\cref{MainThmHigher}] 
	By \cref{ThmOddLoopsSphere}, the sphere $(\S^m, \Theta)$ admits an odd loop of metrics. By \cref{ThmSurgeryStabilityOdd} applied to $M \amalg \S^m$ and $k=0$, the resulting manifold $\tilde M = M \sharp \S^m$ admits an odd loop of metrics. By \cref{ThmExistsHigherMult}, $\tilde M$ admits higher multiplicities. But $\tilde M$ is diffeomorphic to $M$. On can arrange this diffeomorphism to be the identity outside a neighborhood of the surgery sphere and inside this neighborhood, the spin structure is unique (up to equivalence). Therefore, this diffeomorphism is a spin diffeomorphism. Declaring this to be an isometry yields the desired result. 
\end{Prf}

%% file: review.appendix.tex
\section{Principal $G$-Bundles}
\label{SectPrinGBdles}

\begin{Def}
	Let $M$ and $P$ be smooth manifolds, $G$ be a Lie group acting from the right on $P$. A smooth surjection $\pi:P \to M$ is a \emph{principal $G$-bundle}, if
	\begin{enumerate}
		\item 
			Near each point in $M$ there exists an open neighborhood $U$ and an equivariant diffeomorphism $\phi:\pi^{-1}(U) \to U \times G$ such that $\pr_1 \circ \phi = \pi$. (Here, $G$ acts on $U \times G$ via right multiplication on $G$).
		\item
			The group action is simply-transitive on the fibres $P_x := \pi^{-1}(x)$, $x \in M$, of $\pi$.
	\end{enumerate}
\end{Def}

\begin{Def}
	\nomenclature[GP]{$\G(P)$}{gauge transformations of $P$}
	Let $\pi_P: P \to M$ and $\pi_Q:Q \to N$ be two principal $G$-bundles. A \emph{morphism of $G$-bundles} is a map $f:M \to N$ together with a $G$-equivariant map $F:P \to Q$ such that the diagram
	\begin{align*}
		\xymatrix{
			P
				\ar[d]^-{\pi_P}
				\ar[r]^-{F}
			&Q 
				\ar[d]^-{\pi_Q} \\
			M 	\ar[r]^-{f}
			& N
			}
	\end{align*}
	commutes. In case $P=Q$ and $f=\id$, we say $F$ is a \emph{gauge transformation}. Denote by $\G(P)$ the set of all gauge transformations. (This is a group under composition ``$\circ$''.)
\end{Def}

Morphisms of principal bundles have the following remarkable property.

\begin{Thm}[see also \protect{\cite[Thm. 4.3.2]{husemoller}}]
	\label{ThmGBundleMorphismIso}
	Every morphism
	\begin{align*}
		\xymatrix{
			P
				\ar[d]_-{\pi_P}
				\ar[r]^-{F}
			&
			Q
				\ar[d]^-{\pi_Q}\\
			M
				\ar[r]^-{\id}
			&M
		}
	\end{align*}
	of principal $G$-bundles is an isomorphism.
\end{Thm}

\begin{Prf} $ $
	\begin{steplist}
		\step[surjectivity]
			Let $q \in Q$ be arbitrary. Then $x:=\pi_Q(q) \in M$. Take any point $p' \in P_x$. We obtain 
			\begin{align*}
				\pi_Q(F(p')) = \pi_P(p') = x = \pi_Q(q).
			\end{align*}
			Consequently, there exists $g \in G$ such that $F(p').g = q$. The point $p:=p'.g$ satisfies $F(p)=F(p').g=q$.
		\step[injectivity] 
			Let $p_1,p_2 \in P$ such that $F(p_1) = F(p_2)$. This implies
			\begin{align*}
				\pi_P(p_1) = \pi_Q(F(p_1)) = \pi_Q(F(p_2)) = \pi_P(p_2).
			\end{align*}
			Consequently, there exists $g \in G$ such that $p_2 = p_1.g$. But this implies
			\begin{align*}
				F(p_1) = F(p_2) = F(p_1).g,
			\end{align*}
			thus $g = e$, i.e. $p_2 = p_1$.
		\step[smoothness] 
			We have shown that $F$ is bijective, hence there exists an inverse $F^{-1}:Q \to P$. To see that this map is smooth, consider a neighborhood $U \subseto M$ over which both bundles are trivial. We obtain the commutative diagram
			\begin{align*}
				\xymatrix{
					P_U
						\ar[d]^-{\varphi}
						\ar[r]^-{F}
					& Q_U
						\ar[d]^-{\psi} \\
					U \times G
						\ar[r]^-{F_U} 
					& U \times G.
				}
			\end{align*}
			Here, $\psi$ and $\psi$ are local trivializations and $F_U := \psi \circ F \circ \varphi^{-1}$. Because $F$ preserves the fibres and $Q$ is a smooth bundle, there exists $f \in \mathcal{C}^{\infty}(U \times G, G)$ such that
			\begin{align*}
				\forall (x,g) \in U \times G: F_U(x,g) = (x,f(x,g) g).
			\end{align*}
			Define the function
				\DefMap{G_U: U \times G}{U \times G}{(x,g)}{(x,f(x,g)^{-1} g).}
			Clearly $F_U \circ G_U = G_U \circ F_U = \id_{U \times G}$. Therefore
			\begin{align*}
				F^{-1}|_U = \varphi^{-1} \circ G_U \circ \psi \in \mathcal{C}^{\infty}(Q_U,P_U).
			\end{align*}
	\end{steplist}
\end{Prf}

\begin{Cor}
	\label{CorGBundleMorphismIso} Let
	\begin{align*}
		\xymatrix{
			P
				\ar[d]^-{\pi_P}
				\ar[r]^-{F}
			& Q 
				\ar[d]^-{\pi_Q} \\
			M
				\ar[r]^-{f} 
			& N 
		}
	\end{align*}
	be a morphism of principal $G$-bundles. If $f$ is a diffeomorphism, then $F$ is an isomorphism.
\end{Cor}

\begin{Lem}
	\label{LemGBundleIsoAssocVB}  Let $F:P \to Q$ be an isomorphism of principal $G$-bundles over $M$ and $\rho:G \to \Aut_{\K}(V)$ be a representation. Then
		\DefMap{\bar F:P \times_\rho V}{Q \times_\rho V}{{[p,v]}}{[F(p),v]}
	is an isomorphism between the associated vector bundles. If $V$ carries a $G$-invariant metric $\langle \_, \_ \rangle$, this metric descends to a fibre metric on $P \times_\rho V$ and $Q \times_\rho V$ and $F$ is an isometry of vector bundles with respect to these metrics.
\end{Lem}

\begin{Prf}
	Since $F$ is equivariant, $\bar F$ is well-defined and its inverse is given by 
	\DefMap{\bar F^{-1}:Q \times_\rho V}{P \times_\rho V}{{[q,v]}}{[F^{-1}(q),v].}
	Since $\langle \_, \_ \rangle$ is $G$-invariant, the equation
	\begin{align*}
	\forall p \in P: \forall v,w \in V: ([p,v],[p,w]) := \langle v, w \rangle
	\end{align*}
	gives a well-defined fibre metric $( \_, \_, )$ on $P \times_\rho V$. The analogous equation gives a metric on $Q \times_\rho V$ and by construction, $F$ is an isometry.
\end{Prf}

%% file: usb.appendix.tex
\section{Topologies on Mapping Spaces}

\subsection{Compact-Open Topology}

We recall some facts concerning the compact-open topology, see for instance \cite[VII.2]{Bredon} for more details.

\begin{DefI}[compact-open topology]
	\nomenclature[VKU]{$V(K,U)$}{set of all functions $f$ such that $f(K) \subset U$}
	Let $X_1$ be a locally compact Hausdorff space, $X_2$ be any Hausdorff space and $\mathcal{C}(X_1,X_2)$ be the set of all continuous functions $f:X_1 \to X_2$. For any compact $K \subset X_1$ and any open $U \subset X_2$, define
	\begin{align*}
		V(K,U) := \{ f \in \mathcal{C}(X_1,X_2) \mid f(K) \subset U \}.
	\end{align*}
	Declare all the sets $V(K,U)$ to be a subbase for a topology on $\mathcal{C}(X_1,X_2)$. This topology is called the \emph{compact-open topology}.
\end{DefI}

\begin{Thm}
	\label{ThmCOCompEvCont}
	Let $X_1$, $X_2$, $X_3$ be Hausdorff spaces.
	\begin{enumerate}
		\item 
			If $X_1$ and $X_2$ are locally compact, then the composition
				\DefMap{\circ: \mathcal{C}(X_2,X_3) \times \mathcal{C}(X_1,X_2) }{\mathcal{C}(X_1,X_3)}{(g,f)}{g \circ f}
			is continuous.
		\item
			If $X_1$ is locally compact, the evaluation 
				\DefMap{\ev: \mathcal{C}(X_1,X_2) \times X_1}{X_2}{(f,x_1)}{f(x_1)}
			is continuous.
		\item
			If $X_1$ and $X_2$ are locally compact, the \emph{exponential map}
				\DefMap{\mathcal{C}(X_1, \mathcal{C}(X_2, X_3))}{\mathcal{C}(X_1 \times X_2, X_3)}{f}{(x_1,x_2) \mapsto f(x_1)(x_2)}
			is a homeomorphism.
	\end{enumerate}
\end{Thm} 

\subsection{Weak Topology}
\label{SubSectWeakTop}

We need a version of the compact-open topology that also captures smoothness and derivatives. We stick to \cite[Chpt. 2.1]{HirschDiffTop} and introduce the following notion. 

\begin{DefI}[weak topology]
	\label{DefWeakTopology}
	\nomenclature[Ckw]{$\mathcal{C}^{k}_{w}$}{weak $\mathcal{C}^k$-topology}
	Let $M$ and $N$ be two (possibly non-compact) smooth manifolds and $f \in \mathcal{C}^{k}(M,N)$, $k \in \N \cup \{\infty\}$. Let $\varphi$ be a chart on $U \subseto M$, $\psi$ be a chart on $V \subseto N$ and $K \subset U$ be compact such that $f(K) \subset V$. For any $0 \leq k < \infty$ and any $0 < \varepsilon \leq \infty$, the set
	\begin{align}
		\label{EqWeakTopSubbasicNbhd}
			\begin{split}
			B & := B(f, \varphi, U, \psi, V, K, \varepsilon, k)  \\
						 & := \{g \in \mathcal{C}^{k}(M,N) \mid g(K) \subset V, \|\psi \circ f \circ \varphi^{-1} - \psi \circ g \circ \varphi^{-1}\|_{\mathcal{C}^k(\varphi(K), \psi(V))} < \varepsilon\}
			\end{split}
	\end{align}
	is a \emph{subbasic neighborhood of $f$}. The \emph{compact-open $\mathcal{C}^k$-topology}, $k < \infty$, is the topology on $\mathcal{C}^k(M,N)$ generated by all subbasic neighborhoods. The \emph{compact-open $\mathcal{C}^{\infty}$-topology} on $\mathcal{C}^{\infty}(M,N)$ is the union over all compact-open $\mathcal{C}^k$-topologies. These topologies are also called the \emph{weak topologies} for any $k \in \N \cup \{\infty\}$ and are denoted by $\mathcal{C}^k_{w}$.
\end{DefI}

\begin{Rem}
	A neighborhood subbasis of $f$ in the $\mathcal{C}^k_w$-topology is already given by all the sets of the form $B_i := B(f, \varphi_i, U_i, \psi_i, V_i, K, \varepsilon, k)$, where $i \in I$, $I$ a countable index set, $\varepsilon > 0$, $\{ U_i \}_{i \in I}$ is an open cover of $M$, and $K_i \subset U_i$ is compact such that $f(K_i) \subset V_i$ and $\{ K_i \}_{i \in I}$ is still a cover of $M$.
\end{Rem}

This topology has properties analogous to the properties of the compact-open topology mentioned in \cref{ThmCOCompEvCont}.

\begin{Thm}
	\label{ThmCOkTop}
	Let $X$, $Y$, $Z$ be smooth manifolds and $k \in \N \cup \{\infty\}$. The following hold with respect to any $\mathcal{C}^k_w$-topology.
	\begin{enumerate}
		\item 
			\label{ItCOkTopComp}
			The composition 
				\DefMap{\circ: \mathcal{C}^{k}(Y,Z) \times \mathcal{C}^{k}(X,Y) }{\mathcal{C}^{k}(X,Z)}{(g,f)}{g \circ f}
			is continuous.
		\item
			\label{ItCokEvalCont}
			The evaluation
				\DefMap{\ev: \mathcal{C}^{k}(X,Y) \times X}{Y}{(f,x)}{f(x)}
			is continuous.
		\item
			Let $f:X \to \mathcal{C}^{k}(Y, Z)$ and $F:X \times Y \to Z$ be related by 
			\begin{align}
				\label{EqWeakTopExpLawPre}
				\forall (x, y) \in X \times Y: F(x, y) = f(x)(y).
			\end{align}
			If $f$ is continuous, then $F$ is continuous and for each $x \in X$, the map $F(x, \_):Y \to Z$ is in $\mathcal{C}^{k}$. Converseley, if $F \in \mathcal{C}^k(X \times Y, Z)$, then $f$ is continuous.
	\end{enumerate}
\end{Thm}

The weak topology has the convenient property that the usual operations in differential geometry are continuous.

\begin{Lem}
	Let $E \to M$ be a smooth vector bundle. The inclusion
	\begin{align*}
		(\Gamma(E), \mathcal{C}^k) \to (\mathcal{C}^{\infty}(M, E), \mathcal{C}^k_w)
	\end{align*}
	is continuous.
\end{Lem}

\begin{Thm}
	\label{ThmTFunctorCont}
	The tangential functor
		\DefMap{T:(\mathcal{C}^{\infty}(M), \mathcal{C}^{k+1}_w)}{(\Gamma(\Hom TM), \mathcal{C}^k)}{f}{f_*}
	is continuous. 
\end{Thm}

\begin{Cor}
	The group of diffeomorphisms
	\begin{align*}
		\Diff(M) \subset \mathcal{C}^{\infty}(M,M)
	\end{align*}
	is a topological group with respect to $\mathcal{C}^k_w$, $k \in \N \cup \{\infty\}$.
\end{Cor}

\subsection{Applications to Principal $G$-Bundles}
\label{SubSectAppPrinGBdles}

The weak topology can also be used to topologize certain spaces that arise naturally when working with principal $G$-bundles.

\begin{Def}
	\nomenclature[Cinftye]{$\mathcal{C}^{\infty}_e(P,G)$}{space of equivariant functions}
	\label{DefEquivariantGaugeFunctions}
	Let $\pi_P:P \to M$ be a principal $G$-bundle. Define the space
	\begin{align*}
		\mathcal{C}^{\infty}_e(P,G) := \{\sigma \in \mathcal{C}^{\infty}(P,G) \mid \forall p \in P :\forall g \in G: \sigma(p.g) = g^{-1} \sigma(p) g\}
	\end{align*}
	of $G$-equivariant maps from $P$ to $G$. This space is a topological group under pointwise multiplication ``$\cdot$'' in any $\mathcal{C}^k_{w}$-topology.
\end{Def}

\begin{Lem}
	\label{LemGaugeToFunctions}
	Let $\pi_P:P \to M$ be a principal $G$-bundle. For any $F \in \mathcal{G}(P)$ and any $p \in P$, there exists a unique $\sigma_F(p) \in G$ satisfying
	\begin{align}
		\label{EqDefGaugeSigma}
		F(p) = p.\sigma_F(p).
	\end{align}
	This defines a function $\sigma_F \in \mathcal{C}^{\infty}_e(P,G)$ and the induced map
	\DefMap{\sigma:(\mathcal{G}(P),\circ)}{(\mathcal{C}^{\infty}_e(P,G), \cdot)}{F}{\sigma_F,}
	is a group isomorphism.
\end{Lem}

\begin{Prf} $ $
	\begin{steplist}
		\step[existence] 
			The existence of $\sigma_F(p)$ follows from the fact that $F$ preserves the fibres and that the $G$-action on $P$ is simply-transitive on the fibres. 
		\step[equivariance]
			We simply have to verify
			\begin{align*}
				p.(g \sigma_{F}(p.g))
				=(p.g).\sigma_F(p.g)
				= F(p.g)
				= F(p).g
				= p.\sigma_F(p)g,
			\end{align*}
			thus $g \sigma_F(p.g) = \sigma_F(p)g$ as claimed.
		\step[smoothness]
			Let $s \in \Gamma_U(P)$ be a local section and $\psi_s:P_U \to U \times G$ be the associated local trivialization, i.e.
			\DefMap{\psi_s:P_U}{U \times G}{s(x).g}{(x, g).}
			We obtain for any $x \in U$,
			\begin{align*}
				\sigma_F(s(x)) = \pr_2(\psi_s(s(x).\sigma_F(s(x))))= \pr_2(\psi_s(F(s(x)))),
			\end{align*}
			thus $\sigma_F \circ s = \pr_2 \circ \psi_s \circ F \circ s \in \mathcal{C}^\infty(U,G)$. This implies
			\begin{align}
				\label{EqGaugeToFunctionsLocal}
				\forall (x, g) \in U \times G: \sigma_F(\psi_s^{-1}(x, g)) = \sigma_F(s(x).g) = g^{-1} \sigma_F(s(x)) g,
			\end{align}
			thus $\sigma_F \circ \psi_s^{-1} \in \mathcal{C}^{\infty}(U \times G, G)$. Consequently, $\sigma_F:P_U \to G$ is smooth.
		\step[homomorphism]
			Let $F, F' \in \mathcal{G}(P)$, $p \in P$, $g \in G$. We calculate
			\begin{align*}
				p.\sigma_{F' \circ F}(p)
				=(F' \circ F)(p)
				=F'(p.\sigma_F(p))
				=F'(p).\sigma_F(p)
				=p.\sigma_{F'}(p)\sigma_F(p),
			\end{align*}
			thus $\sigma_{F' \circ F} = \sigma_{F'} \cdot \sigma_F$. Clearly, $\sigma_{\id_P} \equiv e \in G$.
		\step[injectivity]
			If $\sigma_F(g)=e$, then $F=\id$ by \cref{EqDefGaugeSigma}.
		\step[surjectivity]
			Let $\sigma \in \mathcal{C}^{\infty}_e(P,G)$ be arbitrary. Define $F$ by \cref{EqDefGaugeSigma}. By construction, for any $p \in P$, $g \in G$
			\begin{align*}
				F(p.g)
				= (p.g).\sigma(pg)
				= p.(gg^{-1}\sigma(p)g)
				=p.\sigma(p)g
				=F(p).g,
			\end{align*}
			thus $F \in \G(P)$ and $\sigma_F = \sigma$.
	\end{steplist}
\end{Prf}

\begin{Lem}
	\label{CorWeakTopGBundle}
	Let $P \to M$ be a principal $G$-bundle. With respect to any $\mathcal{C}^k_w$-topology, the following hold.
	\begin{enumerate}
		\item
			The left translation $L_g$, the right translation $R_g$ and the conjugation $C_g$ define continuous maps $L, R, C: G \to \Diff(G)$.
		\item 
			The gauge group $\G(P) \subset \Diff(P)$ is a topological group.
		\item
			The group action of $G$ on $P$ defines a continuous map $G \to \Diff(P)$.
		\item
			\label{ItWeakTopSigmaHoem}
			The isomorphism $\sigma:(\mathcal{G}(P),\circ) \to (\mathcal{C}^{\infty}_e(P,G), \cdot)$ from \cref{LemGaugeToFunctions} is a homeomorphism.		
	\end{enumerate}
\end{Lem}

\begin{Thm}[association is continuous]
	\label{ThmContGaugeAssocIsos}
	Let $P \to M$ be a principal $G$-bundle, $\rho:G \to \Aut(V)$ be any representation of $G$ and $E := P \times_{\rho} V \to M$ be the associated vector bundle. For any $F \in \G(P)$, we define
		\DefMap{\bar F: E}{E}{{[p,v]}}{[F(p), v].}	
	The map
		\DefMap{(\G(P), \mathcal{C}^k_w)}{(\Gamma(\Iso E), \mathcal{C}^k)}{F}{\bar F}
	is continuous for any $k \in \N \cup \{\infty\}$.
\end{Thm}

\begin{Prf}
	By \cref{CorWeakTopGBundle}, the map $\sigma$ is a homeomorphism. Consequently, the composition
		\DefMap{(\G(P), \mathcal{C}^k_w)}{(\mathcal{C}^{\infty}(P, \Aut(V)), \mathcal{C}^k_w)}{F}{\rho^{-1} \circ \sigma_F}
	is continuous.
	We calculate for any $[p,v] \in E$,
	\begin{align*}
		\bar F([p,v])
		=[F(p),v]
		=[p.\sigma(p), v]
		=[p, \rho^{-1}(\sigma_{F}(p))v].
	\end{align*}
	Let $b_1, \ldots, b_n$ be a basis of $V$, $b^1, \ldots, b^n$ be its dual basis in $V^*$, and $\psi_s$ be a trivialization of $P$ associated to a local section $s \in \Gamma_U(P)$. We obtain that the coordinate matrix of $\bar F$ is locally given by the continuous map
		\DefMap{(\G(P), \mathcal{C}^k_w)}{(\mathcal{C}^{\infty}(U,\GL_n), \mathcal{C}^k)}{F}{(x \mapsto \bar F^i_j(x) := b^j(\rho^{-1}(\sigma_{F}(s(x)))(b_i))).}		
	This implies the claim.
\end{Prf}

\section{Banach and Hilbert Bundles}
\label{SecHilbertBundles}
The notion of a Banach respectively a Hilbert bundle is a straightforward generalization of a vector bundle. But since this is less standard, we give a formal definition.

\begin{Def}[vector bundle]
	\index{vector bundle}\index{Hilbert bundle}\index{Banach bundle}
	\label{DefVectBanHilBundle}
	A continuous map of topological spaces $\pi_{\mathcal{V}}:\mathcal{V} \to X$ is a \emph{vector bundle (over $\K$)}, if
	\begin{enumerate}
		\item 
			For any $x \in X$, the fibre $\mathcal{V}_x := \pi^{-1}_{\mathcal{V}}(x)$ has the structure of a vector space over $\K$ (not necessarily finite dimensional).
		\item
			There exists a normed vector space $(V, \| \_ \|)$ over $\K$ and for each $x \in X$, there exists an open neighborhood $Y \subset X$ of $x$ and a homeomorphism $\phi:\pi^{-1}_{\mathcal{V}}(Y) \to Y \times V$ such that for each $y \in W$, the restriction $\phi:\pi^{-1}_{\mathcal{V}}(y) \to \{y\} \times V$ is a linear isomorphism.
	\end{enumerate}
	We say $\mathcal{V}$ is \emph{modelled} on $V$.	We say $\pi_{\mathcal{V}}$ is a \emph{Banach}, \emph{pre-Hilbert} respectively \emph{Hilbert bundle}, if in addition $V$ and each $\mathcal{V}_x$ have the structure of a Banach, pre-Hilbert, respectively Hilbert space.
\end{Def}

\begin{Rem}
	A Hilbert space structure $\langle \_, \_ \rangle_{\mathcal{H}}$ on a Hilbert bundle $\mathcal{H}$ is part of the data and in no way unique. Therefore, we also denote a Hilbert bundle by $(\mathcal{H}, \langle \_, \_ \rangle_{\mathcal{H}}) \to X$, if we want to stress its dependence on the Hilbert structure.
\end{Rem}

\begin{Def}[morphism]
	Let $\pi_i:\mathcal{V}_i \to X_i$, $i=1,2$, be vector bundles. A \emph{morphism of vector bundles} is a tuple $(f, T)$ of continuous maps between topological spaces such that
	\begin{align*}
		\xymatrix{
			\mathcal{V}_1
				\ar[r]^-{T}
				\ar[d]^{\pi_1}
			& \mathcal{V}_2
				\ar[d]^-{\pi_2}
			\\
			X_1
				\ar[r]^-{f}
			&X_2
		}
	\end{align*}
	commutes and such that for each $x_1 \in X_1$, the restriction $T|_{\pi_1^{-1}(x_1)}$ is a linear map $\mathcal{V}_1|_{x_1} \to \mathcal{V}_2|_{f(x_1)}$. We denote by $\B(\mathcal{V}_1, \mathcal{V}_2)$ the space of such maps. 	
	The notion of a morphism of Banach, pre-Hilbert respectively Hilbert bundles is defined identically (i.e. the morphisms are not required to preserve any additional structure).
\end{Def}

Of course, one can form categories out of these notions.

\nomenclature[HilB]{$\cat{HilB}_{\K}$}{category of Hilbert bundles}
\nomenclature[BanB]{$\cat{BanB}_{\K}$}{category of Banach bundles}
\begin{Def}
	The \emph{category $\cat{BanB}_{\K}$ ($\cat{HilB}_{\K}$) of Banach (Hilbert) Bundles over $\K$} is given by:
	\begin{category}
		\item[Objects]
			Banach (Hilbert) bundles $\pi_{\mathcal{V}}:\mathcal{V} \to X$ in the sense of \cref{DefVectBanHilBundle}.
			
		\item[Morphisms]
			Morphisms of Banach (Hilbert) bundles in the sense of \cref{DefVectBanHilBundle}.
	\end{category}	
\end{Def}

We would like to generalize the notion of an unbounded operator on a Hilbert space to an unbounded morphism between Hilbert bundles.

\begin{Lem}
	\label{LemMorBanachBundlesBanach}
	Let $\mathcal{V}, \mathcal{V}' \to X$ be two Banach bundles over $X$. Then $\B(\mathcal{V}, \mathcal{V}')$ can also be regarded as a map
		\DefMap{\B(\mathcal{V}, \mathcal{V}') = \coprod\limits_{x \in X}{\B(\mathcal{V}_x, \mathcal{V}'_x)}}{X}{T_x}{x,}
	which is a Banach bundle as well.
\end{Lem}

\begin{Prf}
	Local trivializations can be constructed as follows: Let $Y \subseto X$ be an open neighborhood sufficiently small such that there exist trivializations $\psi: \mathcal{V}_Y \to Y \times V$ and $\psi': \mathcal{V}'_Y \to Y \times V'$. For any $T \in \B(\mathcal{V}, \mathcal{V}')$ we obtain a commutative diagram
	\begin{align*}
		\xymatrixcolsep{3.5em}
		\xymatrix{
			\mathcal{V}_Y
				\ar[d]^-{\psi}
				\ar[r]^-{T}
			& \mathcal{V}_Y'
				\ar[d]^-{\psi'}
			\\
			Y \times V
				\ar[r]^-{\id \times \tau}
			&Y \times V'.
		}
	\end{align*}
	where $\tau:Y \to \B(V, V')$ is a continuous map. This defines a local trivialization of $\B(\mathcal{V}, \mathcal{V}')$ by
		\DefMap{\B(\mathcal{V}, \mathcal{V}')_Y}{Y \times \B(V, V')}{T_y}{(y, \tau_y).}
	Since $V$ and $V'$ are Banach spaces, $\B(V, V')$ and all the $\B(\mathcal{V}_x, \mathcal{V}'_x)$ are Banach spaces as well.
\end{Prf}
	
\begin{Def}[domain subspace]
	Let $(H, \langle \_, \_ \rangle_H)$ be a Hilbert space. A \emph{domain subspace} of $H$ is an injection
	\begin{align*}
		\iota: (U, \langle \_, \_ \rangle_U) \hookrightarrow (H, \langle \_, \_ \rangle_H),
	\end{align*}
	where $U \subset H$ is a linear subspace, $\langle \_, \_ \rangle_U$ is a (possibly different) scalar product on $U$ such that $(U, \langle \_, \_ \rangle_U)$ is a Hilbert space in its own right and such that $\iota$ is continuous with dense image.
\end{Def}

\begin{Exm}
	The first order Sobolev space $H^1$ together with the Sobolev norm is a domain subspace for $L^2$ with the $L^2$-norm. Typically, this space serves as the domain of definition for a first order differential operator that is viewed as an unbounded operator $L^2 \to L^2$.
\end{Exm}

\begin{Def}[domain subbundle]
	\label{DefDomainSubbundle}
	Let $(U, \langle \_, \_ \rangle_U) \hookrightarrow (H, \langle \_, \_ \rangle_H)$ be a domain subspace and $\pi_{\mathcal{H}}:(\mathcal{H}, \langle \_, \_ \rangle_{\mathcal{H}}) \to X$ be a Hilbert bundle modelled on $(H, \langle \_, \_ \rangle_H)$. A \emph{domain subbundle} is a subset $\mathcal{U} \subset \mathcal{H}$ together with a Hilbert space structure $\langle \_, \_ \rangle_{\mathcal{U}}$ such that $\pi_{\mathcal{U}} := \pi_{\mathcal{H}}|_{\mathcal{U}}:(\mathcal{U}, \langle \_, \_ \rangle_{\mathcal{U}}) \to X$ is a Hilbert bundle modelled on $(U, \langle \_, \_ \rangle_U)$ and such that near each point $x \in X$ there exists a trivialization $\phi:\mathcal{H}_Y \to Y \times H$ satisfying $\phi(\mathcal{U}_Y) = Y \times U$. Such a trivialization $\phi$ is called \emph{adapted to $U$}.
\end{Def}

\begin{Rem}
	If $\mathcal{U} \subset \mathcal{H}$ is a domain subbundle, then 
	\begin{align*}
		(\mathcal{U}_x, \langle \_, \_ \rangle_{\mathcal{U}_x}) \hookrightarrow (\mathcal{H}_x, \langle \_, \_ \rangle_{\mathcal{H}_x})		
	\end{align*}
	is a domain subspace for any $x \in X$.
\end{Rem}

\begin{Cor}
	Let $(\mathcal{U}, \langle \_, \_ \rangle_{\mathcal{U}}) \to X$ be a domain subbundle of $(\mathcal{H}, \langle \_, \_ \rangle_{\mathcal{H}}) \to X$. Then $\B(\mathcal{U}, \mathcal{H}) \to X$ is a Banach bundle.
\end{Cor}

\begin{Prf}
	This follows immediately from \cref{LemMorBanachBundlesBanach}.
\end{Prf}

\begin{Def}
	\label{DefHilbertBundleUnboundedMorphism}
	Let $(\mathcal{U}, \langle \_, \_ \rangle_{\mathcal{U}})$ be a domain subbundle of $(\mathcal{H}, \langle \_, \_ \rangle_{\mathcal{H}})$. An \emph{unbounded morphism $T:\mathcal{H} \to \mathcal{H}$} is a section $T \in \Gamma(\B(\mathcal{U}, \mathcal{H}))$.
\end{Def}

\section{Spectral Separation of Unbounded Operators}
In this section, we discuss a spectral separation theorem by Kato, see \cref{ThmSpectralSeparation}. We also provide some important definitions to formulate this theorem and draw some useful conclusions. For the entire section, let $X$ be a Banach space.

\begin{Def}[spectrum]
	\nomenclature[DT]{$\dom(T)$}{domain of $T$}
	Let $T:\dom(T) \subset X \to X$ be an unbounded operator with \emph{domain} $T$. The \emphi{resolvent set} $\rho(T) \subset \C$ is the set of all $\lambda \in \C$ such that the operator
	\begin{align*}
		\lambda-T:\dom(T) \to X
	\end{align*}
	is bijective and the inverse
	\begin{align*}
		\xymatrixcolsep{4em}
		\xymatrix{ 
			X \ar[r]^-{(\lambda-T)^{-1}} & \dom(T) \ar@{^(->}[r] & X
			}
	\end{align*}
	belongs to $\B(X)$, the bounded operators on $X$. The complement $\sigma(T) := \C \setminus \rho(T) $ is the \emph{spectrum of $T$}. The subsets
	\begin{align*}
		\sigma_p(T) &:= \{ \lambda \in \C \mid (\lambda - T) \text{ is not injective} \}, \\
		\sigma_c(T) &:= \{ \lambda \in \C \mid (\lambda - T) \text{ is injective, not surjective but has dense image} \}, \\
		\sigma_r(T) &:= \{ \lambda \in \C \mid (\lambda - T) \text{ is injective and does not have dense image} \},
	\end{align*}
	are the \emph{point spectrum}, the \emph{continuous spectrum} and the \emph{rest spectrum}.\footnote{For Dirac operators, $\sigma(\Dirac^g) = \sigma_p(\Dirac^g)$, which we usually denote by $\spec \Dirac^g$. However, the notation $\sigma$ is also very common in functional analysis.}
\end{Def}

\begin{Def}[commutativity]
	\index{commutativity!of unbounded operators}
	Let $T: \dom(T) \subset X \to X$ be an unbounded operator and let $S \in \B(X)$. Then \emph{$T$ commutes with $S$}, if
	\begin{align*}
		\forall x \in \dom(T): Sx \in \dom(T) \text{ and } TSx = STx.
	\end{align*}
\end{Def}

\begin{Def}[decomposition] 
	\label{DefDecompositionUnboundedOperator} \index{decomposition!of unbounded operators}
	Let $T:\dom(T) \subset X \to X$ be an unbounded operator. For any subspace $U \subset X$, we define
	\begin{align*}
		T(U) := T(U \cap \dom(T)) = \{ Tx \mid x \in U \cap \dom(T) \}.
	\end{align*} 
	A \emph{decomposition of $T$} is a pair of subspaces $X_1,X_2 \subset X$ such that $X = X_1 \oplus X_2$ (as Banach spaces, i.e. $X_1$, $X_2$ are closed) and
	\begin{align*}
		T(X_1) \subset X_1, && T(X_2) \subset X_2, && P(\dom(T)) \subset \dom(T),
	\end{align*}
	where $P:X \to X_1$ is the projection onto $X_1$ along $X_2$. For any such decomposition, we obtain the \emph{parts}
	\begin{align*}
		T_i: \dom(T_i) = \dom(T) \cap X_i \to X_i, \; \; T_i := T|_{\dom(T_i)}, \; \;  i=1,2.
	\end{align*}
	Notice that if $T$ is closed, its parts are also closed.
\end{Def}

We will also need the following standard result from functional analysis.

\begin{Thm}
	\label{ThmSpectralTheoremUnboundedOp}
	Let $T:\dom(T) \subset H \to H$ be a closed, densely defined operator on a Hilbert space $H$ with compact resolvent. If $T$ is self-adjoint, then
	\begin{enumerate}
		\item $\sigma(T)=\sigma_p(T) \subset \R$ is discrete and consists entirely of isolated eigenvalues of finite multiplicity.
		\item There exists a Hilbert ONB of eigenvectors of $T$.
	\end{enumerate} 
\end{Thm}

\begin{Thm}[\protect{\cite[Theorem 6.17, p. 178]{kato}}]
	\label{ThmSpectralSeparation}
	Let $T:\dom(T) \subset X \to X$ be a closed operator on a Banach space $X$ and assume that the spectrum $\sigma(T)$ can be written as a disjoint union $\sigma(T) = \sigma_1 \dot \cup \sigma_2$. Assume, there even exists a simple closed curve $\Gamma \subset \C$ such that $\sigma_1$ lies in the interior of $\Gamma$ and $\sigma_2$ lies in the exterior of $\Gamma$. Then there exists a decomposition $X = X_1 \oplus X_2$ (see \cref{DefDecompositionUnboundedOperator}) such that the parts $T_i := T|_{X_i}: \dom(T_i) \subset X_i \to X_i$, $i=1,2$, satisfy $\sigma(T_i) = \sigma_i$. In addition $T_1 \in B(X_1)$. The projection $P \in B(X)$ onto $X_1$ along $X_2$ commutes with $T$ and satisfies
	\begin{align}
		\label{EqSpecralProj}
		\forall x \in X: Px = -\frac{1}{2 \pi i} \oint_{\Gamma}{(z-T)^{-1}xdz} .
	\end{align}
\end{Thm}
 
\begin{Cor}
	\label{CorSpectralEigenProjection}
	Under the hypothesis of \cref{ThmSpectralSeparation} it holds in addition that 
	\begin{enumerate}
		\item $\sigma_p(T_i) = \sigma_p(T) \cap \sigma_i$, $i=1,2$.
		\item 
			If $E_\lambda$ is the eigenspace of $T$ with respect to a $\lambda \in \sigma_p(T)$, then 
			\begin{align*}
				E^1 := \bigoplus_{\lambda \in \sigma_p(T_1)}{E_\lambda} \subset X_1, 
				&&E^2 := \mathop{\overline{\bigoplus}}\limits_{\lambda \in \sigma_p(T_2)}{E_\lambda} \subset X_2, 
			\end{align*}
		\item 
			If in addition $T$ is self-adjoint and has compact resolvent, $X_1 = E^1$ and $X_2 = E^2$.
	\end{enumerate}
\end{Cor}
 
\begin{Prf} $ $ 
	\begin{steplist}
	\step[$(\sigma_p(T) \cap \sigma_1) \subset \sigma_p(T_1)$]
		Let $\lambda \in \sigma_p(T)$ be arbitrary. By definition
		\begin{align}
			\label{EqSpectralEigenProjectionDefEV}
			\exists 0 \neq x \in \dom(T): Tx = \lambda x.
		\end{align}
		Since $P$ commutes with $T$ (by \cref{ThmSpectralSeparation}), so does $I-P$, which is the projection onto $X_2$ along $X_1$. Decompose $x = Px + (I-P)x =:x_1 + x_2$ and calculate
		 \begin{align} 
				\label{EqSpectralEigenProjectionx1}
				 \lambda x_1 &= \lambda Px = P Tx = T Px = T_1 x_1, \\ \label{EqSpectralEigenProjectionx2}
				 \lambda x_2 &= \lambda (I-P) x = (I-P)Tx = T(I-P)x = T_2 x_2.
		 \end{align}
		We claim that
		\begin{align}
			\label{EqSpectralEigenProjectionZerox2}
				\forall \lambda \in \sigma_p(T) \cap \sigma_1: \forall x \in E_\lambda: x_2 &= (I-P)x = 0.\\ \label{EqSpectralEigenProjectionZerox1}
				\forall \lambda \in \sigma_p(T) \cap \sigma_2: \forall x \in E_\lambda: x_1 &= Px = 0.
		\end{align}
		To see that \cref{EqSpectralEigenProjectionZerox2} holds, we argue by contradiction: If it does not hold, then $x_2$ would be an eigenvector of $T_2$ by \cref{EqSpectralEigenProjectionx2}. This implies $\lambda \in \sigma_2$, which cannot hold by assumption. By the same reasoning, we obtain \cref{EqSpectralEigenProjectionZerox1}.
				 
		This implies the claim: By \cref{EqSpectralEigenProjectionZerox2}, we obtain $x_2=0$. Therefore,
		\begin{align}
			\label{EqSpecProjFinal}
			x_1 = x_1 + x_2 = x \neq 0, 
		\end{align}
		which implies $\lambda \in \sigma_p(T_1)$ by \cref{EqSpectralEigenProjectionx1}. This also implies $E_\lambda \subset P(X)=X_1$. Analogously, if $\lambda \in \sigma_p(T) \cap \sigma_2$, we obtain $\lambda \in \sigma_p(T_2)$ and $E_{\lambda} \subset X_2$ by the same reasoning.
		 
	 \step[$\sigma_p(T_1) \subset (\sigma_p(T) \cap \sigma_1) $]
			Conversely, assume that $\lambda \in \sigma_p(T_1)$. Then there exists $x_1 \in \dom(T_1) \subset \dom(T)$ such that $x_1 \neq 0$ and $\lambda x_1 = T_1 x_1 = T x_1$. Thus $\lambda \in \sigma_p(T)$. Using \cref{ThmSpectralSeparation}, we obtain $\lambda \in \sigma_p(T) \cap \sigma_1$. Analogously, if $\lambda \in \sigma_p(T_2)$ it follows that $\lambda \in \sigma_p(T) \cap \sigma_2$. 
			
	\step 
		We have already proven along the way that
		\begin{align*}
			\bigoplus_{\lambda \in \sigma_p(T_i)}{E_\lambda} \subset X_i, \qquad i=1,2.
		\end{align*}
		By \cref{DefDecompositionUnboundedOperator} of a decomposition, $X_2$ is closed, thus $E^2 \subset X_2$.

	\step[self-adjoint case]
		By \cref{ThmSpectralTheoremUnboundedOp}, the hypothesis implies $\sigma(T)=\sigma_p(T)$, all eigenspaces are finite-dimensional and there exists an orthonormal system of eigenvectors. In particular, any vector $x \in X_1$ can be decomposed into
		\begin{align*}
			x = u + v \in E^1 \oplus E^2.
		\end{align*}
		Since $E^1 \subset X_1$ and $E^2 \subset X_2$, this implies $x, u \in X_1$ and $v \in X_2$. Therefore, 
		\begin{align*}
			v = (I-P)v = (I-P)(x-u) = 0,
		\end{align*}
		thus $x=u \in E^1$. This proves $X_1 = E^1$. The equality $X_2 = E^2$ is proven analogously.
	\end{steplist}
\end{Prf}
 
\begin{Thm}
	\label{ThmImageBundle}
	Let $X$ be a topological space, $H$ be a Hilbert space and $T:X \to B(H)$, $x \mapsto T_x$, be a continuous map. Assume there exists $k \in \N$ such that for all $x \in X$, the operator $T_x$ has rank $k$. Denote by $H_x := T_x(H)$ the image of $T_x$. Then the union
	\begin{align*}
		E:= \coprod_{x \in X}{H_x} \to X
	\end{align*}
	is a continuous vector subbundle of rank $k$ of the trivial bundle $X \times H \to X$.
\end{Thm}

\begin{Prf}
	Local trivializations can be constructed as follows: Fix any $x \in X$. By construction, there exist $v_1, \ldots, v_k \in H$ such that $(w_1,\ldots,w_k):=(T_x(v_1),\ldots,T_x(v_k))$ is a basis of $H_x$. In particular, $w_1, \ldots, w_k$ are linearly independent. This can also be expressed by saying that the continuous map
		\DefMap{X}{\Lambda^k H}{y}{T_y(v_1) \wedge  \ldots \wedge T_y(v_k)}
	is nonzero at $x$. Consequently, there exists a neighborhood $U$ of $x$ where this map is non-zero. Consequently, for all $y \in U$, the vectors $(T_y(v_1), \ldots, T_y(v_k))$ are linearly independent. Therefore, they must be a basis of $H_y$. Consequently,
		\DefMap{\Psi: U \times \R^k}{E_U}{(y,c_1, \ldots, c_k)}{\sum_{i=1}^{k}{c_i T_y(v_i)} }
	defines a local trivialization of $E_U$.
\end{Prf}

%% file: higher.appendix.tex
\section{Equivariant Covering Spaces}
In this section, we fix our notion of covering spaces and extend some lifting resuts to the equivariant setting.

\begin{DefI}[covering space]
	A continuous surjection $p:\hat X \to X$ between topological spaces is a \emph{covering}, if for each $x \in X$, there exists an open neighborhood $U$ of $x$, an index set $A$ and a family of disjoint open subsets $\{\hat U_{\alpha} \subset \hat X \}_{\alpha \in A}$, such that $p^{-1}(U) = \dot \cup_{\alpha \in A}{\hat U_{\alpha}}$ and $p|_{\hat U_{\alpha}} \to U$ is a homeomorphism.
\end{DefI}

\begin{Rem}
	The definition of a \emph{covering space} varies slightly in the literature, in particular with respect to the connectedness assumptions on the spaces and with respect to the surjectivity of $p$. We do not assume the spaces to be connected, but the map to be surjective. It is also common to choose basepoints $\hat x_0 \in \hat X$, $x_0 \in X$ and view $p$ as a map of pointed spaces $p:(\hat X, \hat x_0) \to (X, x_0)$.
\end{Rem}

\begin{Thm}[Homotopy Lifting Property, {\cite[Prop. 1.30]{HatcherAlgTop}}]
	\label{ThmHomotopyLifting}
	Let $p : \hat X \to X$ be a covering space. Let $f:Y \to X$ be a map and assume there exists a lift $\hat f:Y \to \hat X$, i.e. $p \circ \hat f = f$. Let $H:Y \times I \to X$ be a homotopy starting at $f$, i.e. $H(\_,0)=f$. Then $H$ lifts to a unique homotopy $\hat H: Y \times I \to \hat X$ satisfying $\hat H(\_,0)=\hat f$, i.e. the following diagram commutes
	\begin{align*}
		\xymatrixcolsep{3.5em}
		\xymatrix{
			 &
			 &
			 \hat X
				\ar[d]^-{p}
			\\
			Y
				\ar@/^2pc/[rru]_-{\hat f}
				\ar@/_2pc/[rr]^-{f}
				\ar@{^(->}[r]^-{\id \times \{0\}}
			& Y \times I
				\ar@{..>}[ru]^-{\hat H}
				\ar[r]^-{H}
			& X.
		}
	\end{align*}
\end{Thm}

\begin{Thm}[Lifting Theorem, {\cite[Prop. 1.33]{HatcherAlgTop}}]
	\label{ThmLifting}
	Let $p:(\hat X, x_0) \to (X,x_0)$ be a covering and $f:(Y,y_0) \to (X,x_0)$ be any continuous map. We assume that $Y$ is path-connected and locally path-connected. There exists a lift $\hat f:(Y,y_0) \to (\hat X,\hat x_0)$, i.e. a map such that
	\begin{align*}
		\xymatrix{
			& (\hat X, \hat x_0)
				\ar[d]^-{p}
			\\
			(Y,y_0)
				\ar@{..>}[ur]^-{\hat f}
				\ar[r]^-{f}
			&(X,x_0)
		}
	\end{align*}
	commutes, if and only if $f_{\sharp}(\pi_1(Y,y_0)) \subset p_{\sharp}(p_1(\hat X, \hat x_0))$.
\end{Thm}

\begin{Lem}[equivariant lifts]
	\label{LemLiftEquiv}
	Let $p:\hat X \to X$ be a covering map and $f:Y \to X$. Assume that there are topological groups $G$, $H$, $K$ acting from the right on $Y$, $X$, $\hat X$. Assume that there exists a lift $\hat f:Y \to \hat X$ of $f$ against $p$ and assume that there are group homomorphisms $\alpha$, $\beta$, $\gamma$ such that $\gamma=\beta \circ \alpha$. If $f$ and $p$ are equivariant with respect to the group actions and these homomorphisms, i.e.
	\begin{align*}
		\forall y \in Y: \forall g \in G: f(y.g) &= f(y).\gamma(g) \\
		\forall \hat x \in \hat X: \forall k \in K: p(\hat x.k) &= p(\hat x).\beta(k),
	\end{align*}
	and if $G$ is connected, then $\hat f$ is equivariant as well, i.e.
	\begin{align*}
	\forall y \in Y: \forall g \in G: \hat f(y.g) =\hat f(y).\alpha(g).
	\end{align*}
	The situation can be visualized in the following diagram:
	\begin{align*}
		\xymatrixrowsep{1.5em}
		\xymatrixcolsep{1.5em}
		\xymatrix{ 
			&&& K
				\ar@/^2pc/[ddddr]^-{\beta}
				\ar@/_/@{..>}[dd]
			\\
			\\
			&&& \hat X
				\ar[d]^-{p}
				\ar@/_/@{..}[uu]
			\\
			&&Y
				\ar[r]_-{f}
				\ar[ur]^-{\hat f}
				\ar@/_/@{..}[dll]
			& X
				\ar@/_/@{..}[dr]
			\\
			G
				\ar@/_2pc/[rrrr]^-{\gamma}
				\ar@/^2pc/[uuuurrr]^-{\alpha}
				\ar@/_/@{..>}[urr]
			&&&&H.
				\ar@/_/@{..>}[ul]
		}
	\end{align*}	
\end{Lem}

\begin{Prf}
	Let $y \in Y$ and $g \in G$ be arbitrary. Since $G$ is connected, there exists a continuous path $\sigma:[0,1] \to G$ joining $e_G \in G$ and $g \in G$. This implies for any $t \in [0,1]$,
	\begin{align*}
		p(\hat f(y.\sigma(t)).\alpha(\sigma(t)^{-1}))
		&=p(\hat f(y.\sigma(t))).\beta(\alpha(\sigma(t)^{-1}))
		=f(y.\sigma(t)).\gamma(\sigma(t)^{-1}) \\
		&=f(y).\gamma(\sigma(t)) \gamma(\sigma(t)^{-1})
		=f(y).\gamma(\sigma(t) \sigma(t)^{-1}) 
		=f(y).
	\end{align*}
	Consequently, the continuous path
	\DefMap{[0,1]}{\hat X}{t}{\hat f(y . \sigma(t)) . \alpha(\sigma(t)^{-1})}
	maps into the discrete fibre over $f(y)$ and is therefore constant. In particular, for $t=0$ and $t=1$, we obtain
	\begin{align*}
	& \hat f(y.\sigma(1)). \alpha(\sigma(1)^{-1})=\hat f(y.\sigma(0)).\alpha(\sigma(0)^{-1}) \\
	\Longrightarrow & \hat f(y. g).\alpha(g^{-1})=\hat f(y . e_G) . \alpha(e_G^{-1}) \\
	\Longrightarrow & \hat f(y. g)=\hat f(y).\alpha(g) .
	\end{align*}
\end{Prf}

\section{The Round Sphere as a Spin Manifold}
In this section, we give an explicit description of the spin structure on the round sphere, see \cref{ThmSpinStructureRoundSphere}. This will be useful for some explicit calculations. We denote by $\bar g$ the Euclidean metric on $\R^{m+1}$ and by $g\degree$ the round metric induced by $\bar g$ on $\S^m$.

\begin{Rem}[tangent bundle and orientation]
	Recall that the map
		\DefMap{\{ (p,v) \in \S^m \times \R^{m+1} \mid \bar g( p, v)  = 0 \}}{T \S^m}{(p,v)}{\sum_{i=1}^{m+1}{v^i \partial_i|_p}}
	is an isomorphism of vector bundles. With this identification, the vector field $N \in \mathcal{N}(\S^m)$, defined by $N(p):=p$, is the outward pointing unit normal field. We enumerate the canonical basis of $\R^{m+1}$ by $(e_0, \ldots, e_m)$ and define this basis to be positively oriented. A basis $(b_1, \ldots, b_m)$ of $T_p\S^m$ is \emph{positively oriented} if and only if $(p,b_1,\ldots, b_m) \in \R^{m+1}$ is positive. This defines an orientation on $\S^m$, which is sometimes called the \emph{Stokes orientation}, see \cite[Prop. 13.17]{LSM}.
\end{Rem}

\begin{Rem}[frame bundle]
	\label{RemSOSmChar}
	Let $\SO(\S^m)$ be the bundle of positively oriented orthonormal frames. Then the map
		\DefMap{\Psi:\SO_{m+1}}{\SO(\S^m)}{A}{(A_1, \ldots, A_m) \in \SO_{A_0}(\S^m),}
	where $A_i$ is the $i$-th column of $A$, is a diffeomorphism. Clearly, its inverse is given by
	\DefMap{\Psi^{-1}:\SO(\S^m)}{\SO_{m+1}}{(b_1, \ldots, b_m) \in T_p\S^m}{(p,b_1, \ldots, b_m).}
	Therefore, we will identify $\SO(\S^m)$ with $\SO_{m+1}$. The group action is given by $\SO_{m+1} \times \SO_m \to \SO_{m+1}$, $(b,A) \mapsto b j(A)$, where 
	$j$ is the canonical inclusion
		\DefMap{j:\SO_m}{\SO_{m+1}}{A}{\begin{pmatrix}1 & 0 \\ 0 & A\end{pmatrix}.}	
\end{Rem}

To describe $\Spin(\S^m)$ it will be useful to have the following auxiliary map available.

\begin{Lem} \label{ThmSOmSpinmp}
	Let $\tilde e_1, \ldots, \tilde e_m \in \R^{m}$ and $e_0, \ldots, e_m \in \R^{m+1}$ be the canonical bases of $\R^m$ respectively $\R^{m+1}$. The canonical inclusion
		\DefMap{\iota:\R^{m}}{\R^{m+1}}{x = \sum_{i=1}^{m}{x^i \tilde e_i}}{\sum_{i=1}^{m}{x^i e_i}}
	is isometric and therefore lifts to a map $\Cl_{m} \to \Cl_{m+1}$, which we also denote by $ \iota$. This map is an injective homomorphism of algebras and makes the following diagram commutative
	\begin{align}
		\label{EqIotaJSOSpin}
		\begin{split}
			\xymatrixcolsep{3.5em}
			\xymatrix{
				\Spin_m
					\ar@{->>}[d]_-{\vartheta_m}
					\ar@{^(->}[r]^-{ \iota} 
				& \Spin_{m+1}
					\ar@{->>}[d]^-{\vartheta_{m+1}} \\
				\SO_{m}
					\ar@{^(->}[r]_-{j} 
				& \SO_{m+1}. 
				}
		\end{split}
	\end{align}

\end{Lem}

\begin{Prf}
	It is clear that $\iota:\R^m \to \R^{m+1}$ is an isometric embedding. This implies the existence of the lift $\iota: \Cl_m \to \Cl_{m+1}$, see \cref{ThmCliffFunctor}. Since $\iota:\R^{m} \to \R^{m+1}$ is an injective homomorphism of vector spaces, its lift is an injective homomorphism of algebras, see \cite[Prop. 10.7.2]{greub2}.
	By construction $\iota(\Spin_{m}) \subset \Spin_{m+1}$. To see that \cref{EqIotaJSOSpin} commutes, recall that for any $v \neq 0$ the reflection $\rho_v \in \SO_m$ along the hyperplane orthogonal to $v$ is given by
		\DefMap{\rho_v:\R^m}{\R^m}{w}{w - 2\frac{\bar g(v,w)}{\bar g(v,v)}v}
	and that $\vartheta_m(v) = (\rho_v(\tilde e_1), \ldots, \rho_v(\tilde e_m))$. Consequently, for any $1 \leq i \leq m$,
	\begin{align*}
		j(\vartheta_m(\tilde e_i))
		&=j \left(  \tilde e_1, \ldots, -\tilde e_i, \ldots, \tilde e_m \right)
		=( e_0,  e_1, \ldots, - e_i, \ldots, e_m ) \\
		&=\vartheta_{m+1}(e_i)
		=\vartheta_{m+1}( \iota (\tilde e_i)).
	\end{align*}
\end{Prf}

\begin{Rem}
	\label{RemSpinmActsSpinmp1}
	\cref{ThmSOmSpinmp} above allows us to identify $\Spin_m$ with its image under $\iota$ in $\Spin_{m+1}$. Consequently, $\Spin_m$ acts on $\Spin_{m+1}$ by right multiplication.
\end{Rem}

\begin{Rem}
	Recall that $\Spin_m = \Spin(\R^{m},\bar g)$, $m \geq 2$, is the connected non-trivial double cover of $\SO_m$ and that the covering map $\vartheta_{m}:\Spin_{m} \to \SO_{m}$ satisfies the following property: Whenever $v_1 \cdot \ldots \cdot v_r \in \Spin_{m}$ then $\vartheta_m(v_1 \cdot \ldots \cdot v_r)$ is the coordinate matrix (with respect to the canonical basis) of the map $\rho_{v_1} \circ \ldots \circ \rho_{v_r}$, where each $\rho_{v_i}$ is the reflection along the hyperplane $v_i^\perp$. 
\end{Rem}

\begin{Thm}[spin structure on spheres] 
	\label{ThmSpinStructureRoundSphere}
	Let $m \geq 2$ and define $\Spin(\S^m) := \Spin_{m+1}$. The map			
		\DefMap{\pi:\Spin(\S^m)}{\S^m}{w}{\vartheta_{m+1}(w) e_0,}
	is a principal $\Spin_{m}$-bundle over $\S^m$ and $\Theta_{m} := \Psi \circ \vartheta_{m+1}:\Spin(\S^m) \to \SO(\S^m)$ is a spin structure for $\S^m$. Here, $\Psi$ is the isomorphism from \cref{RemSOSmChar}.
\end{Thm}

\begin{Prf} $ $
	\begin{steplist}
	\step[group action]
		By \cref{RemSpinmActsSpinmp1}, we already have a group action of $\Spin_m$ on $\Spin_{m+1}$. We have to show that the orbits of this action are precisely the fibres of the map $\pi$. This requires a simple lemma from Linear Algebra, see \cref{LemABCSOn} below.
		So take $s \in \Spin_{m}$, $w \in \Spin_{m+1}$ and calculate (using the terminology of \cref{ThmSOmSpinmp})
		\begin{align*}
			\pi(w.s)
			&=\vartheta_{m+1}(w.s)e_0
			=\vartheta_{m+1}(w \cdot \iota(s))e_0
			=\vartheta_{m+1}(w) \vartheta_{m+1}(\iota(s))e_0\\
			&=\vartheta_{m+1}(w) j(\vartheta_m(s))e_0
			\jeq{\ref{LemABCSOn}}{=}\vartheta_{m+1}(w)e_0
			=\pi(w).
		\end{align*}
		Thus, the action preserves the fibres of $\pi$. To see that the action is transitive on the fibres, let $w' \in \Spin_{m+1}$ such that $\pi(w) = \pi(w')$. This implies
		\begin{align*}
			e_0 = \underbrace{\vartheta_{m+1}(w)^{-1}}_{=:B^{-1}} \underbrace{\vartheta_{m+1}(w')}_{=:A}e_0 
		\end{align*}
		and by \cref{LemABCSOn}, there exists $C \in \SO_{m}$ such that $B^{-1} A = j(C)$. There exists $s \in \Spin_{m}$ such that $\vartheta_m(\pm s)=C$. This implies
			\begin{align*}
			\vartheta_{m+1}(w.s)
			=\vartheta_{m+1}(w \cdot \iota(s))
			=\vartheta_{m+1}(w) j(\vartheta_m(s))
			=B \; j(C)
			=A
			=\vartheta_{m+1}(w'),
		\end{align*}
		thus $\pm w' = w.s$. In case $-w' = w.s$ replace $s$ with $- s$. Finally, to see that the action is simply-transitive on the fibres, let $w \in \Spin_{m+1}$ and $s \in \Spin_m$ such that
		\begin{align*}
			w = w.s = w \cdot \iota(s).
		\end{align*}
		Since $w$ is invertible, $\iota(s) = 1$ and since $\iota$ is injective, $s=1$.
	\step[trivialization]
		To see that $\pi$ is locally trivial, we first construct a trivialization for $\SO(\S^m)$ as follows: For a matrix $A \in \R^{(m+1) \times (m+1)}$ denote by $\varphi(A) \in \R^{m \times m}$ the projection onto the lower right square. Define the map
		\DefMap{\Phi: \SO_{m+1}}{\S^m \times \SO_m}{A}{(Ae_0,\varphi(A))}
		For any matrix $C \in \SO_m$, let $\psi(C) \in \R^{1 \times m}$ be defined by
		\begin{align*}
			\psi(C)_j^2 := 1 - \sum_{i=1}^{m}{c_{ij}^2}.
		\end{align*}
		Since any matrix $A \in \SO_{m+1}$ satisfies
		\begin{align*}
			\forall 1 \leq j \leq m: 1 = a_{0j}^2 + \sum_{i=1}^{m}{a_{ij}^2},
		\end{align*}
		the map
		\DefMap{\Psi:\S^m \times \SO_m}{\SO_{m+1}}{(p,C)}{\begin{pmatrix}p & \psi(C) \\
		& C \end{pmatrix}}
		is an inverse to $\Phi$. Thus, $\Phi$ is a global trivialization for $\SO(\S^m)$. Since $\Spin_{m+1}$ is simply connected, we obtain a lift $\tilde \Phi$ such that
		\begin{align*}
			\xymatrix{
				\Spin_{m+1}
					\ar@{..>}[r]^-{\tilde \Phi}
					\ar[d]^-{\vartheta_{m+1}}
				&\S^m \times \Spin_m
					\ar[d]^-{\id \times \vartheta_m}_{2:1}
				\\
				\SO_{m+1}
					\ar[r]^-{\cong}_-{\Phi}
				&\S^m \times \SO_m
			}
		\end{align*}
		commutes.  This implies that $\tilde \Phi$ has full rank, is $\Spin_m$-equivariant by \cref{LemLiftEquiv} and satisfies for any $w \in \Spin_{m+1}$,
		\begin{align*}
			(\pr_1 \circ \tilde \Phi)(w)
			=\pr_1 \circ \id \times \vartheta_{m} \circ \tilde \Phi(w)
			=\vartheta_{m+1}(w)e_0
			=\pi(w).
		\end{align*}
		Now, if $\tilde \Phi(w_1) = \tilde \Phi(w_2)$, this implies that $w_1$ and $w_2$ are in the same fibre, hence there exists a unique $s \in \Spin_m$ such that $w_2 = w_1.s$. Since $\tilde \Phi$ is equivariant, we obtain
		\begin{align*}
			\tilde \Phi(w_1) = \tilde \Phi(w_2) = \tilde \Phi(w_1.s) = \tilde \Phi(w_1).s,
		\end{align*}
		thus $s= 1$ and $w_2 = w_1$. Therefore, $\tilde \Phi$ is injective. Since $\vartheta_{m+1}$, $\Phi$ and $\id \times \vartheta_m$ are all surjective, $\tilde \Phi$ is also surjective. We have shown that $\pi$ is a principal $\Spin_m$-bundle.
		
	\step[spin structure] 
		By construction $\Theta_m = \Psi \circ \vartheta_{m+1}$ is a two-fold covering map. It follows from \cref{EqIotaJSOSpin} that
		\begin{align*}
			\xymatrix{
				\Spin_{m+1} \times \Spin_m
					\ar[r]
					\ar[d]^-{\vartheta_{m+1} \times \vartheta_m}
				& \Spin_{m+1}
					\ar[d]^-{\vartheta_{m+1}}
				\\
				\SO_{m+1} \times \SO_m
					\ar[r]
				&
				\SO_{m+1}
			}
		\end{align*}
		commutes. This implies that $\Theta_m$ is a spin structure.
	\end{steplist}
\end{Prf} 

\begin{Lem}
	\label{LemABCSOn}
	For any $A,B \in \SO_{m+1}$
	\begin{align*}
		Ae_0 = B e_0 \Longleftrightarrow \exists!  C \in \SO_m: A^{-1} B = j(C),
	\end{align*}
	where $j$ is the map from \cref{ThmSOmSpinmp}.
\end{Lem}

\begin{Prf} $ $ \\
	'''$\Longrightarrow$'': By hypothesis $M:=A^{-1}B \in \SO_{m+1}$ satisfies $Me_0 = e_0$. Consequently, 
	\begin{align*}
		\exists ! C \in \R^{m \times m}: 
		M =
		\begin{pmatrix}
			1 & * \\
			0 & C
		\end{pmatrix}.
	\end{align*}
	By hypothesis, $M \in \SO_{m+1}$ and therefore $1 = \det(M) = \det(C)$, so $C \in \SL_m$. Since the columns of $M$ form an orthonormal system, we automatically obtain $*=0$. Consequently, the columns of $C$ are an orthonormal system as well, so $C \in \SO_m$ and $M=j(C)$. \\
	''$\Longleftarrow$'': This is obvious.
\end{Prf}

\begin{Lem}
	\label{LemRoundSphereLift}
	Define $f:\R \to \SO_{m+1}$, $\alpha \mapsto R_{\alpha}|_{\S^m}$, where $R_{\alpha}$ are the rotations from \cref{DefRotations}. Let
		\DefMap{\hat f :\R}{\Spin_{m+1}}{\alpha}{\cos(\alpha) + \sin(\alpha)e_{m-1} \cdot e_m.}
	and denote by $\cdot 2:\R \to \R, \alpha \mapsto 2\alpha$, the multiplication by $2$.
	Then the diagram
	\begin{align}
		\label{EqRoundSphereLift}
		\begin{split}
			\xymatrix{
				\R
					\ar[d]^-{\cdot 2}
					\ar[r]^-{\hat f}
				& \Spin_{m+1}
					\ar[d]^-{\vartheta_{m+1}}
				\\
				\R
					\ar[r]^-{f}
				& \SO_{m+1}
			}
		\end{split}
	\end{align}
	commutes.
\end{Lem}

\begin{Prf} 
	We split the calculation into several steps. In \cref{StpTransPath}, we show that $\hat f$ can be factored into $\hat f(\alpha) = v_{\hat \alpha} w_{\hat \alpha}$. In \cref{StpCalcvw}, we calculate $\vartheta_{m+1}(v_{\hat \alpha})$ and $\vartheta_{m+1}(w_{\hat \alpha})$ separately. Finally, in \cref{StpFinal}, we put everything together and obtain the result.
	\begin{steplist}
	\step[transform the path]
		\label{StpTransPath}
		Let $\alpha \in \R$. We define
		\begin{align*}
			v_\alpha := \cos(\alpha)e_{m-1} + \sin(\alpha)e_m, &&
			w_\alpha := -\cos(\alpha)e_{m-1} + \sin(\alpha)e_m.
		\end{align*}
		Using the double angle formulas
		\begin{align}
			\label{EqDblAngle}
			\cos(2 \alpha) = \cos(\alpha)^2 - \sin(\alpha)^2, && 
			\sin(2\alpha) = 2 \cos(\alpha) \sin( \alpha),
		\end{align}
		we obtain
		\begin{align*}
			v_\alpha w_\alpha
			&= (\cos(\alpha)e_{m-1} + \sin(\alpha)e_m)(-\cos(\alpha)e_{m-1} + \sin(\alpha)e_m) \\
			&=-\cos(\alpha)^2e_{m-1}^2 
			+ \sin(\alpha)^2e_m^2
			+\cos(\alpha)\sin(\alpha)e_{m-1} \cdot e_m
			-\sin(\alpha) \cos(\alpha) e_m \cdot e_{m-1} \\
			&=\cos(\alpha)^2
			- \sin(\alpha)^2
			+2\cos(\alpha)\sin(\alpha)e_{m-1} \cdot e_m \\
			&\jeq{\cref{EqDblAngle}}{=}\cos(2 \alpha) + \sin(2 \alpha)e_{m-1} \cdot e_m
			=\hat f(2 \alpha).
		\end{align*}
		Consequently, by setting $\hat \alpha := \alpha / 2$, we obtain
		\begin{align*}
			\hat f(\alpha) = v_{\hat \alpha} w_{\hat \alpha}.
		\end{align*}

	\step[calculate $\vartheta_{m+1}(v_{\hat \alpha})$, $\vartheta_{m+1}(w_{\hat \alpha})$]
		\label{StpCalcvw}
		By definition, for any $v \neq 0$, $\vartheta_{m+1}(v)=\rho_v$, where $\rho_v$ is the reflection along the hyperplane $v^\perp$. This map is explicitly given by
		\begin{align*}
			\forall w \in \R^{m+1}: \rho_v(w) = w - 2 \frac{\bar g(v,w)}{\bar g(v,v)}v.
		\end{align*}
		We calculate
		\begin{align*}
			\bar g(v_{\hat \alpha},v_{\hat \alpha}) = 1, &&
			\bar g(w_{\hat \alpha},w_{\hat \alpha}) = 1, \\
			\bar g(v_{\hat \alpha},e_{m-1}) = \cos(\hat \alpha), &&
			\bar g(w_{\hat \alpha},e_{m-1}) = - \cos(\hat \alpha) \\
			\bar g(v_{\hat \alpha},e_{m}) = \sin(\hat \alpha), &&
			\bar g(w_{\hat \alpha},e_{m}) = \sin(\hat \alpha).
		\end{align*}
		Combining this with
		\begin{align}
			\label{EqDoubleAngleCorcos}
			1 - 2 \cos(\hat \alpha)^2
			& = \sin(\hat \alpha)^2 + \cos(\hat \alpha)^2 - 2 \cos(\hat \alpha)^2
			=\sin(\hat \alpha)^2 - \cos(\hat \alpha)^2
			\jeq{\cref{EqDblAngle}}{=} - \cos(\alpha) \\ \label{EqDoubleAngleCorsin}
			1 - 2 \sin(\hat \alpha)^2
			&= \sin(\hat \alpha)^2 + \cos(\hat \alpha)^2 - 2 \sin(\hat \alpha)^2 
			= \cos(\hat \alpha)^2 - \sin(\hat \alpha)^2 
			\jeq{\cref{EqDblAngle}}{=} \cos(\alpha),
		\end{align}
		we obtain
		\begin{align*}
			\vartheta_{m+1}(v_{\hat \alpha})(e_{m-1})
			&=e_{m-1} - 2 \cos(\hat \alpha)(\cos(\hat \alpha)e_{m-1} + \sin(\hat \alpha)e_m), \\
			&=(1 - 2 \cos(\hat \alpha)^2)e_{m-1} - 2 \cos(\hat \alpha) \sin(\hat \alpha)e_m \\
			&=-\cos(\alpha) e_{m-1} - \sin(\alpha) e_m,
		\end{align*}
		\begin{align*}
			\vartheta_{m+1}(v_{\hat \alpha})(e_{m})
			&=e_m - 2 \sin(\hat \alpha)(\cos(\hat \alpha) e_{m-1} + \sin(\hat \alpha) e_m) \\
			&=(1 - 2 \sin(\hat \alpha)^2)e_m - 2 \sin(\hat \alpha) \cos(\hat \alpha) e_{m-1}\\
			&=\cos( \alpha) e_m - \sin(\alpha) e_{m-1},
		\end{align*}
		\begin{align*}
			\vartheta_{m+1}(w_{\hat \alpha})(e_{m-1})
			&=e_{m-1} + 2 \cos(\hat \alpha)(- \cos(\hat \alpha)e_{m-1} + \sin(\hat \alpha) e_m)\\
			&=(1-2 \cos(\hat \alpha)^2)e_{m-1} + 2 \cos(\hat \alpha) \sin(\hat \alpha) e_m \\
			&\jeq{}{=}-\cos(\alpha) e_{m-1} + \sin( \alpha) e_m,
		\end{align*}
		\begin{align*}
			\vartheta_{m+1}(w_{\hat \alpha})(e_{m})
			&=e_m - 2 \sin(\hat \alpha)(-\cos(\hat \alpha)e_{m-1} + \sin(\hat \alpha) e_m)\\
			&=(1 - 2 \sin(\hat \alpha)^2)e_m + 2 \sin(\hat \alpha)\cos(\hat \alpha) e_{m-1}\\
			&=\cos(\alpha) e_m + \sin( \alpha) e_{m-1}.
		\end{align*}
		We obtain the matrices
		\begin{align*}
			\vartheta_{m+1}(v_{\hat \alpha}) 
			=\begin{pmatrix}
				I_{m-1} & 0 & 0 \\
				0 & -\cos(\alpha) & -\sin(\alpha) \\
				0 & -\sin(\alpha) & \cos(\alpha)
			\end{pmatrix},
		\end{align*}
		\begin{align*}
			\vartheta_{m+1}(w_{\hat \alpha}) 
			=\begin{pmatrix}
				I_{m-1} & 0 & 0 \\
				0 & -\cos(\alpha) & \sin(\alpha) \\
				0 & \sin(\alpha) & \cos(\alpha)
			\end{pmatrix}
		\end{align*}

	\step[final calculation]
		\label{StpFinal}
		Since
		\begin{align*}
			\begin{pmatrix}
				-\cos(\alpha) & -\sin(\alpha) \\
				- \sin(\alpha) & \cos(\alpha) \end{pmatrix}
			\begin{pmatrix}
				-\cos(\alpha) & \sin(\alpha) \\
				\sin(\alpha) & \cos(\alpha)
			\end{pmatrix}
			&=
			\begin{pmatrix}
				\cos(\alpha)^2 - \sin(\alpha)^2 & - 2 \cos(\alpha)\sin(\alpha) \\
				2 \cos(\alpha) \sin(\alpha) & \cos(\alpha)^2 - \sin(\alpha)^2
			\end{pmatrix} \\
			&=
			\begin{pmatrix}
				\cos(2 \alpha) & - \sin(2 \alpha) \\
				\sin(2 \alpha) & \cos(2 \alpha)
			\end{pmatrix},
		\end{align*}
		we obtain
		\begin{align*}
			\vartheta_{m+1}(\hat f(\alpha))
			&=\vartheta_{m+1}(v_{\hat \alpha} w_{\hat \alpha})
			=\vartheta_{m+1}(v_{\hat \alpha}) \vartheta_{m+1}( w_{\hat \alpha}) \\
			&=
			\begin{pmatrix}
				I_{m-1} & 0 & 0 \\
				0 & \cos(2 \alpha) & - \sin(2 \alpha) \\
				0 & \sin(2 \alpha) & \cos(2 \alpha)
			\end{pmatrix}
			=f(2 \alpha).
		\end{align*}
		\end{steplist}
\end{Prf}

\section{Rayleigh Quotients and the Min-Max Principle}
\label{SectRayleighRitz}

In this section, we collect some well known facts about Rayleigh quotients and the Min-Max principle. Our treatment of the subject will be similar to \cite[XIII.1]{Reed4}.

\begin{Thm}[Min-Max principle]
	\index{min-max principle}
	\label{ThmMinMaxPrinciple}
	Let $H$ be a Hilbert space and $T: \dom(T) \subset H \to H$ be a densely defined self-adjoint operator with compact resolvent. Assume that the spectrum of $T$ satisfies $b \leq \lambda_1 \leq \lambda_2 \leq \ldots$ for some lower bound $b \in \R$. Then for any $k \in \N$
	\begin{align}
		\label{EqMinMaxEV}
		\lambda_k = \min_{\substack{U \subset \dom(T), \\ \dim(U) = k}}{\max_{\substack{x \in U, \\ \|x\|=1}}{\langle Tx, x \rangle }},
	\end{align}
	where the $\min$ is taken over all linear subspaces $U \subset \dom(T)$ of dimension $k$.
\end{Thm}

\begin{Prf} 
	By \cref{ThmSpectralTheoremUnboundedOp}, there exists a complete orthonormal system $\{u_k\}_{k \in \N}$ such that $Tu_k = \lambda_k u_k$ and
	\begin{align*}
		H = \overline{\Lin \{ u_k \mid k \in \N \} }.
	\end{align*}
	Let $U \subset \dom(T)$, $\dim(U)=k$, be arbitrary. The subspace
	\begin{align*}
		V_k := \Lin\{ u_i \mid i \geq k \}
	\end{align*}
	satisfies $\codim(V_k) = k-1$. Therefore, by \cref{LemDimCodimk},
	\begin{align*}
		\exists \hat x \in V_k \cap U, && \|\hat x\|=1.
	\end{align*}
	By construction, this $\hat x$ has a representation
	\begin{align*}
		\hat x = \sum_{i=k}^{\infty}{\alpha_i u_i}.
	\end{align*}
	This implies
	\begin{align*}
		\langle T\hat x, \hat x \rangle
		&= \sum_{i=k}^{\infty}{\alpha_i \langle Tu_i, \hat x \rangle }
		= \sum_{i=k}^{\infty}{\alpha_i \lambda_i \langle u_i, \hat x \rangle }
		= \sum_{i=k}^{\infty}{\alpha_i \lambda_i \sum_{j=k}^{\infty}{\alpha_j \langle u_i, u_j \rangle} } \\
		&= \sum_{i=k}^{\infty}{\alpha_i^2 \lambda_i }
		\geq \lambda_k \sum_{i=k}^{\infty}{\alpha_i^2 }
		=\lambda_k \|\hat x\|
		=\lambda_k.
	\end{align*}
	Therefore,
	\begin{align*}
		\sup_{\substack{x \in U,\\ \|x\|=1}}{\langle Tx, x \rangle } \geq \langle T \hat x, \hat x \rangle  \geq \lambda_k.
	\end{align*}
	Since $U$ was arbitrary,
	\begin{align*}
		\inf_{\substack{U \subset \dom(T), \\ \dim(U) = k}}{\sup_{\substack{x \in U, \\ \|x\|=1}}{\langle Tx, x \rangle }} \geq  \lambda_k.
	\end{align*}
	On the other hand, the subspace $U_k:=\Lin(u_1, \ldots, u_k) \subset \dom(T)$ is $k$-dimensional and satisfies
	\begin{align*}
		\forall x = \sum_{i=1}^{k}{\alpha_i u_i} \in U_k: \langle Tx, x \rangle = \sum_{i=1}^{k}{ \alpha_i^2 \lambda_i } \leq \lambda_k \|x\|.
	\end{align*}
	Since $\langle Tu_k, u_k \rangle = \lambda_k$, this implies
	\begin{align} \label{EqSupBound}
		\sup_{\substack{x \in U_k, \\ \|x\| = 1}}{\langle Tx, x \rangle } = \langle Tu_k,u_k \rangle  =\lambda_k.
	\end{align}
	Consequently, the infimum is attained at $U_k$ and equal to $\lambda_k$. Since any unit sphere of a $k$-dimensional subspace $U$ is compact and any linear map on a finite dimensional space is continuous, we may also replace the ''$\sup$'' by a ''$\max$''. 
\end{Prf}

\begin{Lem}
	\label{LemDimCodimk}
	Let $X$ be any vector space and let $U,V \subset X$ be two subspaces such that $\dim U = k$ and $\codim V = k-1$. Then
	\begin{align*}
		U \cap V \neq \{0\}.
	\end{align*}
\end{Lem}

\begin{Prf}
	Define the map
	\DefMap{\Psi:U+V}{U/(U \cap V)}{u+v}{[u].}
	This map is well-defined: If $u,u' \in U$ and $v,v' \in V$ such that $u+v=u'+v'$, this implies 
	\begin{align*}
		u-u' = v' - v \in U \cap V \Longrightarrow [u-u'] = [0] \Longrightarrow [u] = [u'].
	\end{align*}
	Clearly, $\Psi$ is surjective. We claim that $\ker \Psi = V$. On the one hand $\Psi(0+v) = [0]$, so $V \subset \ker \Psi$ and on the other hand
	\begin{align*}
		u+v \in \ker \Psi \Longrightarrow [0] = [u] = \Psi(u+v) \Longrightarrow u \in U \cap V \Longrightarrow u + v \in V.
	\end{align*}
	Consequently, $\Psi$ descends to an isomorphism
	\begin{align*}
		(U+V) / V \longrightarrow U/(U \cap V).
	\end{align*}
	Now assume by contradiction that $U \cap V = {0}$. Then this is an isomorphism
	\begin{align*}
		(U+V) / V \longrightarrow U.
	\end{align*}
	But this implies
	\begin{align*}
		k = \dim(U) = \dim((U+V)/V) \leq \dim(X/V) = \codim(V) = k-1,
	\end{align*}
	which is a contradiction.
\end{Prf}

\begin{Def}[Rayleigh quotient]
	In the situation of \cref{ThmMinMaxPrinciple}, for any $x \in \dom(T)$, the expression
	\begin{align*}
		R_T(x) := \frac{\langle Tx, x \rangle }{\langle x, x\rangle }
	\end{align*}
	is called \emphi{Rayleigh quotient}.
\end{Def}

\begin{Cor}
	\label{CorMinMaxEVEstimate}
	In the situation of \cref{ThmMinMaxPrinciple}, assume there exists a subspace $W \subset \dom(T)$ of dimension $k$ such that
	\begin{align}
		\label{EqMinMaxEstimate}
		\forall x \in W: R_T(x) \leq \Lambda,
	\end{align}
	for some $\Lambda > 0$. Then $T$ has $k$ eigenvalues in $[b,\Lambda]$.
\end{Cor}

\begin{Prf}
	By \cref{EqMinMaxEV}, we obtain
	\begin{align*}
		\lambda_k = \min_{\substack{U \subset \dom(T), \\ \dim(U) = k}}{\max_{\substack{x \in U, \\ \|x\|=1}}{\langle Tx, x \rangle }} 
		\leq {\max_{\substack{x \in W, \\ \|x\|=1}}{\langle Tx, x \rangle }}
		\leq \Lambda.
	\end{align*}
\end{Prf}

\begin{Rem}
	In spin geometry, \cref{CorMinMaxEVEstimate} is usually applied to $T = (\Dirac^g_{\C})^2$. An equation of the form \cref{EqMinMaxEstimate} implies that $\Dirac^g_{\C}$ has $k$ eigenvalues in $[-\Lambda, \Lambda]$.
\end{Rem}

\begin{Thm}
	\label{ThmRayleighRitzDistance}
	Let $H$ be a Hilbert space, $T:H \to H$ be a densely defined operator, self-adjoint with compact resolvent. Let $\Lambda> 0$, $k \in \N$ and assume that the first $k+1$ distinct eigenvalues of $T$ satisfy
	\begin{align*}
		0 \leq \lambda_1 < \ldots < \lambda_k < \Lambda < \lambda_{k+1}.
	\end{align*}
	Define
	\begin{align*}
		E^{(\nu)} := \ker(T-\lambda_{\nu}), && V := \bigoplus_{\nu=1}^{k}{E^{(\nu)}},
	\end{align*}
	and let $x \in \dom(T)$, $\|x\|=1$, such that
	\begin{align*}
		R_T(x) = \langle Tx, x \rangle  \leq \Lambda + \varepsilon.
	\end{align*}
	Then the distance between $x$ and $V$ satisfies
	\begin{align*}
		d(V,x)^2 \leq \frac{\Lambda + \varepsilon}{\lambda_{k+1}}.
	\end{align*}
\end{Thm}

\begin{Prf}
	Consider the orthogonal decomposition
	\begin{align*}
		x = \sum_{\nu=1}^{\infty}{x^{(\nu)}}, && x^{(\nu)} \in E^{(\nu)}.
	\end{align*}
	By hypothesis
	\begin{align*}
		\varepsilon + \Lambda \geq 
		R_T(x) &= \langle Tx, x \rangle \\
		&=\langle \sum_{\nu=1}^{k}{Tx^{(\nu)}} + \sum_{\nu=k+1}^{\infty}{Tx^{(\nu)}}, \sum_{\nu=1}^{k}{x^{(\nu)}} + \sum_{\nu=k+1}^{\infty}{x^{(\nu)}} \rangle \\
		&=\langle \sum_{\nu=1}^{k}{\lambda_{\nu} x^{(\nu)}} + \sum_{\nu=k+1}^{\infty}{\lambda_{\nu} x^{(\nu)}}, \sum_{\nu=1}^{k}{x^{(\nu)}} + \sum_{\nu=k+1}^{\infty}{x^{(\nu)}} \rangle \\	
		&=\sum_{\nu=1}^{k}{\lambda_{\nu} \|x^{(\nu)}\|^2} + \sum_{\nu=k+1}^{\infty}{\lambda_{\nu} \|x^{(\nu)}\|^2} \\
		& \geq \lambda_{k+1} \sum_{\nu=k+1}^{\infty}{\|x^{(\nu)}\|^2}.
	\end{align*}
	Let $P_V:H \to H$ be the orthogonal projection onto $V$. We obtain
	\begin{align*}
		d(V, x)^2
		= \|P_V(x) - x\|^2 
		=\sum_{\nu=k+1}^{\infty}{\|x^{(\nu)}\|^2}
		\leq \frac{\Lambda + \varepsilon}{\lambda_{k+1}}.
	\end{align*}
\end{Prf}

\section{Extension of the Surgery Theorem}
\label{SecSurgeryExtension}

In this section, we give some more details on how to modify the proof of \cref{ThmBaerDahlSurgeryPimped} to get \cref{ThmSurgeryFamilyPimped}. The two major aspects are the treatment of unsymmetric intervals, see \cref{RemLambdaLambda}, and the case of compact families of metrics, see \cref{RemCompactFamilies}. The reader is assumed to be very familiar with \cite{BaerDahlSurgery}. We will use the notation from \cref{ThmSurgeryFamilyPimped}, \cref{RemCuttingOffSpinorFields} and \cref{DefSurgeryCuttOffFunctions}.

\begin{Rem}[$\Lambda_1 < \Lambda_2$]
	\label{RemLambdaLambda}
	The first part of the proof does not have to be modified at all, because it is shown that if $\lambda_i$ is an eigenvalue of $\Dirac^g$, then $\Dirac^{\tilde g}$ has an eigenvalue $\tilde \lambda_i$ in $]\lambda_i-\varepsilon, \lambda_i + \varepsilon[$. Of course this is still true, if $\lambda_i \in [\Lambda_1, \Lambda_2]$. 
	
	In the second part, one has to apply the Min-Max principle, see \cref{ThmMinMaxPrinciple}. Now, one has to bound the Rayleigh quotient of $(\Dirac^{g}-c)^2$ in terms of $l^2$ instead of bounding the Rayleigh quotient of $(\Dirac^g)^2$ in terms of $\Lambda^2$. Again, $c:= \tfrac{1}{2}(\Lambda_1 + \Lambda_2)$, $l := \tfrac{1}{2}|\Lambda_2 - \Lambda_1|$. This is due to the fact that if $\Dirac^g$ has eigenvalues $\lambda_1, \ldots, \lambda_n$ in $[\Lambda_1, \Lambda_2]$, then $(\Dirac-c)^2$ has eigenvalues $(\lambda_1 -c)^2, \ldots, (\lambda_n-c)^2$ in $[0,l^2]$. Consequently, one has to prove
	\begin{align}
		\label{EqL1L2AppendixClaim}
		\frac{\| (\Dirac^g_{\K}-c) \psi^{\varepsilon}\|_{L^2(\Sigma_{\K}^g M)}^2}{\|\psi^{\varepsilon}\|_{L^2(\Sigma^g_{\K} M)}^2} < (l + \varepsilon)^2.
	\end{align}	
	To bound the denominator, one first has to generalize \cite[Lemma 2.2]{BaerDahlSurgery}, see \cref{LemEigenSpinorsScalarCurvatureAvoidMod}. Due to the slightly different constant in \cref{EqScalCurvAvoidMod}, one has to choose $S_1$ sufficiently large such that
	\begin{align*}
		1 - \frac{4 l^2 + 8c\Lambda_2 - S_0}{S_1 - S_0} = \frac{S_1 - 4 l^2 - 8c\Lambda_2}{S_1 - S_0} \geq \left( \frac{l + \tfrac{\varepsilon}{2}}{l + \varepsilon} \right)^2.
	\end{align*}
	Now, the calculation \cite[p. 67]{BaerDahlSurgery} gives
	\begin{align*}
		\|\psi^{\varepsilon}\|_{L^2(\Sigma^g_{\K} M)} \geq \left( \frac{l + \tfrac{\varepsilon}{2}}{l + \varepsilon} \right)^2 \|\tilde \psi^{\varepsilon}\|_{L^2(\Sigma^{\tilde g^{\varepsilon}}_{\K} \tilde M^{\varepsilon})},
	\end{align*}
	which is the desired bound on the denominator. 
	
	For the nominator, we have to modify \cite[Eq. (7)]{BaerDahlSurgery} to
	\begin{align*}
		\| (\Dirac^{g}_{\K} - c ) \psi^{\varepsilon}\|_{L^2(\Sigma^{\tilde g^{\varepsilon}}_{\K} \tilde M)}
		&= \| \nabla \chi^{\varepsilon} \tilde  \psi^{\varepsilon} +   \chi^{\varepsilon} (\Dirac^{g}_{\K} - c ) \tilde \psi^{\varepsilon}\|_{L^2(\Sigma^{\tilde g^{\varepsilon}}_{\K} \tilde M^{\varepsilon})} \\
		\leq & \tfrac{c}{r_{\varepsilon}}  \|\tilde \psi^{\varepsilon}\|_{L^2(\Sigma^g_{\K} A_{S}(r_{\varepsilon}, 2 r_{\varepsilon}))}  + l \|  \tilde \psi^{\varepsilon}\|_{L^2(\Sigma^{\tilde g^{\varepsilon}}_{\K} \tilde M^{\varepsilon})} .
	\end{align*}
	Now, the first term can be estimated just as in \cite{BaerDahlSurgery}, which gives
	\begin{align*}
		\|(\Dirac^g_{\K} -c) \psi^{\varepsilon}\|_{L^2(\Sigma^g_{\K} M)} < (l + \tfrac{\varepsilon}{2}) \|\tilde \psi^{\varepsilon} \|_{L^2(\Sigma^{\tilde g^{\varepsilon}}_{\K} \tilde M^{\varepsilon})}.
	\end{align*}
	This implies \cref{EqL1L2AppendixClaim}. 
\end{Rem}

\begin{Lem}
	\label{LemEigenSpinorsScalarCurvatureAvoidMod}
	Let $(M,g,\Theta^g)$ be a closed Riemannian spin manifold, $\Lambda_1 < \Lambda_2$, $S_0 < S_1$, $c:= \tfrac{1}{2}(\Lambda_1 + \Lambda_2)$, $l := \tfrac{1}{2}|\Lambda_2 - \Lambda_1|$. Assume the scalar curvature of $M$ satisfies $\scal^g \geq S_0$. Define
	\begin{align*}
		M_+ := \{x \in M \mid \scal^g \geq S_1\}.
	\end{align*}
	Then for any spinor field $\psi \in L^2_{[\Lambda_1, \Lambda_2]}(\Sigma_{\K}^g M)$, the following inequality holds:
	\begin{align}
		\label{EqScalCurvAvoidMod}
		\int_{M_+}{|\psi|^2 \dvol_{g}} \leq \frac{4 l^2 + 8c\Lambda_2 - S_0}{S_1 - S_0} \int_{M}{|\psi|^2 \dvol_{g}}.
	\end{align}
\end{Lem}

\begin{Prf}
	Let $\lambda_1 \leq \ldots \leq \lambda_k$ be the eigenvalues of $\Dirac^g_{\K}$ in $[\Lambda_1, \Lambda_2]$. Then $\psi$ can be decomposed into
	\begin{align*}
		\psi = \sum_{i=1}^{k}{\psi_i}, && \Dirac^g_{\K} \psi_{i} = \lambda_i \psi_{i}. 
	\end{align*}
	This implies
	\begin{align*}
		\|(\Dirac^g_{\K} - c)\psi\|^2_{L^2(\Sigma_{\K}^gM)}
		& =\sum_{i,j=1}^k{\langle (\Dirac^g_{\K} - c) \psi_i, (\Dirac^g_{\K} - c)\psi_j \rangle } \\
		& =\sum_{i,j=1}^k{\langle (\lambda_i - c) \psi_i, (\lambda_j - c)\psi_j \rangle } \\
		& =\sum_{i=1}^k{(\lambda_i - c)^2 \|\psi_i\|^2_{L^2(\Sigma_{\K}^g M)}} \\
		& \leq l^2 \|\psi\|_{L^2(\Sigma_{\K}^g M)}
	\end{align*}	
	and
	\begin{align*}
		\langle \Dirac^g_{\K} \psi, \psi \rangle
		= \sum_{i=1}^k{\lambda_i \|\psi_i\|^2_{L^2(\Sigma_{\K}^g M)}}
		\leq \Lambda_2 \|\psi\|_{L^2(\Sigma_{\K}^g M)}^2.
	\end{align*}
	Using the Schrödinger--Lichnerowicz formula, we obtain
	\begin{align*}
		|(\Dirac^g_{\K}-c)\psi|^2
		&=\langle \Dirac^g_{\K} \psi, \Dirac^g_{\K}\psi \rangle  - 2 c  \langle \Dirac^g_{\K} \psi, \psi \rangle  + c^2 |\psi|^2 \\
		&=|\nabla \psi|^2 + \tfrac{1}{4} \scal^g |\psi|^2  - 2 c \langle \Dirac^g_{\K} \psi, \psi \rangle  + c^2 |\psi|^2 .
	\end{align*}
	Altogether, this implies
	\begin{align*}
		0 &\leq \int_{M}{|\nabla \psi|^2 + c^2 |\psi|^2 \dvol_{g}} \\
		&=\int_{M}{|(\Dirac^g_{\K} - c)^2\psi|^2 + 2c \langle \Dirac^g_{\K} \psi, \psi \rangle - \frac{1}{4} \scal^g |\psi|^2} \dvol_{g} \\
		& \leq l^2 \int_{M}{|\psi|^2 \dvol_{g}} + 2c \Lambda_2 \int_{M}{|\psi|^2 \dvol_{g}} 
		- \frac{S_0}{4} \int_{M \setminus M_+}{|\psi|^2\dvol_{g}} - \frac{S_1}{4} \int_{M_+}{|\psi|^2 \dvol_{g}}.
	\end{align*}
	Therefore,
	\begin{align*}
		\left( \frac{S_1 - S_0}{4} \right) \int_{M_+}{|\psi|^2 \dvol_{g}}
		\leq \left( l^2 + 2c \Lambda_2 - \frac{S_0}{4} \right) \int_{M}{|\psi|^2 \dvol_{g}},
	\end{align*}
	which gives the result.
\end{Prf}

\begin{Rem}[compact families]
	\label{RemCompactFamilies}
	When generalizing \cref{ThmBaerDahlSurgeryPimped} to compact $\mathcal{C}^2$-continuous families $\mathbf{g}:A \to \Rm(M)$, the neighborhood $U_S(r_{\varepsilon})$ in the proof will depend on the family parameter $\alpha \in A$, but this does not matter as long as one chooses $r_{\varepsilon}$ so small such that $U_{S}(r_{\varepsilon}) \subset V$, where $V$ is the open neighborhood of the surgery sphere from the statement of \cref{ThmSurgeryFamilyPimped}. The cut-off functions from \cref{DefSurgeryCuttOffFunctions} will then depend on the family parameter and yield families of functions $\chi^{\varepsilon}_{\alpha}$. Nevertheless, if one chooses the constants $C_1$, $C_2$ on \cite[p. 65]{BaerDahlSurgery} uniformly, the first part of the proof goes through. 
	
	For the second part of the proof, one first has to generalize \cite[Thm. 2.1]{BaerDahlSurgery} to a family version, see \cref{ThmProp21DahlExtended}. This version is then used to argue that one can assume that $\mathbf{g}$ satisfies
	\begin{itemize}
		\item 
			$\forall \alpha \in A: \scal^{g_{\alpha}} \geq S_0$ on $M$.
		\item
			$\forall \alpha \in A: \scal^{g_{\alpha}} \geq 2 S_1$ on an open neighborhood $U_0 \subset V$ of $S$ in $M$ with $S_1$ arbitrarily large.
	\end{itemize}
	Now, for each $\alpha \in A$, one has to choose $r_{\varepsilon, \alpha} > 0$ sufficiently small such that
	\begin{itemize}
		\item 
			$r_{\varepsilon,\alpha} < \frac{\varepsilon^4}{2^8 \cdot 100 \cdot (m+1)^4}$.
		\item
			$U_S(r_{\varepsilon,\alpha})((2r_{\varepsilon,\alpha})^{1/11}) \subset U_0 \subset V$.
		\item
			$(2r_{\varepsilon,\alpha})^{1/11})  \leq R$, where $R$ is from \cite[Lem. 2.4]{BaerDahlSurgery}. One observes from the proof of that lemma that this $R$ can be chosen uniformly for the family.
	\end{itemize}	
	Now, we perform the surgery on the neighborhood $U_{S}(r_{\varepsilon,\alpha})$. This yields manifolds of the form $\tilde M^{\varepsilon,\alpha} = (M \setminus U_S(r_{\varepsilon,\alpha})) \cup \tilde U_{\varepsilon,\alpha}$. The family of metrics $\mathbf{\tilde g^{\varepsilon}}$ is obtained from the Gromov-Lawson construction, see \cite{GromowLawson,RosenbergSurgery}. This construction actually does hold in the family case as well, but a detailed discussion of this fact is beyond the scope of this text. Readers interested in family versions of this construction are referred to \cite[Sect. 3.7]{Walsh}. Strictly speaking, the resulting manifolds $\tilde M^{\varepsilon,\alpha}$ depend on $\alpha$. But for different $\alpha$, these manifolds are mutually diffeomorphic. So, by pulling back the resulting metrics to one reference manifold $\tilde M^{\varepsilon,\alpha_0}=: \tilde M^{\varepsilon}$ one obtains a manifold that does not depend on the family parameter and a family of metrics $\mathbf{\tilde g^{\varepsilon}}$ as required. Now, the rest of the proof can be applied to all $(M,g_{\alpha})$, $\alpha \in A$, which gives \cref{ThmBaerDahlSurgeryPimped}.
\end{Rem}

\begin{Thm}
	\label{ThmProp21DahlExtended}
	Let $M$ be a closed manifold, $A$ be a compact space, $\mathbf{g}:A \to (\Rm(M), \mathcal{C}^2)$ be a continuous family of Riemannian metrics. Let $N \subset M$ be a compact submanifold of positive codimension and $U$ be a compact neighborhood of $N$ in $M$. There exist continuous families of metrics $\mathbf{g}^{(j)}:A \to (\Rm(M), \mathcal{C}^2)$, $j \in \N$, such that
	\begin{enumerate}
		\item
			$\forall \alpha \in A: \forall j \in \N: g^{(j)}_{\alpha}$ is conformally equivalent to $g_{\alpha}$.
		\item
			$\forall \alpha \in A: \forall j \in \N: g^{(j)}_{\alpha} = g_{\alpha}$ on $M \setminus U$.
		\item
			$\sup_{\alpha \in A}{\|g^{(j)}_{\alpha} - g_{\alpha}\|_{\mathcal{C}^1}} \to 0$ as $j \to \infty$.
		\item
			$\min_{\alpha \in A}{\min_{N} \scal^{g_{\alpha}^{(j)}}} \to \infty$ as $j \to \infty$.
	\end{enumerate}
\end{Thm}

\begin{Prf}
	In this proof, one only has to choose the constants $C_1, C_2$ in \cite[p. 59]{BaerDahlSurgery} uniformly, which is possible, since $A$ is compact.
\end{Prf}

%% file: evpaper.appendix.tex
\section{Some Fundamental Results}
\label{SectEvpaperAppendixFundRes}
Here we collect some results for convenient reference.

\begin{Thm}[\protect{identification of spinor bundles, cf. \cite{MaierGen,BourgGaud}}]
	\label{ThmSpinorIdentification}
	\nomenclature[Dgh]{$\Dirac^h_g$}{Dirac operator for $h$ pulled back to $\Sigma^g M$}
	\nomenclature[Df]{$\Dirac_g$}{family of Dirac operators $\Dirac^h_g$}
	Let $g \in \Rm(M)$ be a fixed metric. For every $h \in \Rm(M)$, there exists an isometry of Hilbert spaces $\bar \beta_{g,h}: L^2(\Sigma_{\C}^gM) \to L^2(\Sigma_{\C}^hM)$, such that the operator
	\begin{align*}
	\Dirac^h_g := \bar \beta_{g,h} \circ \Dirac_{\C}^h \circ \bar \beta_{g,h}: L^2(\Sigma_{\C}^gM) \to L^2(\Sigma_{\C}^gM)
	\end{align*}
	is closed, densely defined on $H^1(\Sigma_{\C}^g M)$, isospectral to $\Dirac_{\C}^h$, and such that the map
	\begin{align*}
	\Dirac_g:\Rm(M) \to B(H^1(\Sigma_{\C}^gM),L^2(\Sigma_{\C}^gM)), && h \mapsto \Dirac_g^h,
	\end{align*}
	is continuous. (Here $B(\_)$ denotes the space of bounded linear operators endowed with the operator norm.)
\end{Thm}

\begin{Thm}[\protect{\cite[Thm. 1.3]{Hanke}}]
	\label{ThmHanke}
	Given $k,l \geq 0$ there is an $N=N(k,l) \in \N_{\geq 0}$ with the following property: For all $n \geq N$, there is a $4n$-dimensional smooth closed spin manifold $P$ with non-vanishing $\widehat{A}$-genus and which fits into a smooth fibre bundle
	\begin{align*}
		X \to P \to S^k.
	\end{align*}
	In addition, we can assume that the following conditions are satisfied:
	\begin{enumerate}
		\item The fibre $X$ is $l$-connected.
		\item The bundle $P \to S^k$ has a smooth section $s:S^k \to P$ with trivial normal bundle.
	\end{enumerate}
\end{Thm}

\begin{Thm}[\protect{\cite[Thm 1.1]{AmmannSurgery}}]
	\label{ThmAmmannSurgery}
	Let $M$ be a closed spin manifold of dimension $m$. Assume $m \not \equiv 0 \mod 4$ and $m \not \equiv 1,2 \mod 8$. Then there exists a metric $g \in \Rm(M)$ such that $\Dirac_{\C}^g$ is invertible. 
\end{Thm}

\begin{Lem}[classification of fibre bundles over $S^1$]
	\label{LemClassBundlesS1}
	Let $M$ be any smooth manifold.
	\begin{enumerate}
		\item For any $f \in \Diff(M)$, the space obtained by setting
		\begin{align*}
			P_f := [0,1] \times M / \sim, &&
			\forall (t,x) \in [0,1] \times M: (1,x) \sim (0,f(x))
		\end{align*}
		is a smooth $M$-bundle over $S^1 = [0,1] / (0 \sim 1)$. 
		\item Let $\mathcal{I}$ denote the isotopy classes of $\Diff(M)$, and $\mathcal{P}$ the isomorphism classes of $M$-bundles over $S^1$. The map
			\DefMap{\mathcal{I}}{\mathcal{P}}{{[f]}}{{[P_f]}}
		is well-defined, surjective, and if $[P_{f}] = [P_{f'}]$, then $f$ is isotopic to a conjugate of~$f'$. 

		\item Let $M$ be oriented. Then $P_f$ is orientable if and only if $f \in \Diff^+(M)$.
		\item Let $M$ be spin and simply connected. Then $P_f$ is spin if and only if $f \in \Diff^{\spin}(M)$. 
	\end{enumerate}
\end{Lem} 

The proof of this is elementary. 

\begin{Thm}[\protect{\cite[Thm. 8.17, Rem. 8.18a)]{BaerBallmann}}] 
	\label{ThmBB}
	Let $(M,g)$ be an even-dimensional smooth complete spin manifold with volume element $\mu$ and a fixed spin structure. Let $N$ be a closed, two-sided hypersurface in $M$. Cut $M$ open along $N$ to obtain a manifold $M'$ with two isometric boundary components $N_1$ and $N_2$. Consider the pullbacks $\mu'$, $\Sigma_{\C} M'$, $\Dirac_{\C}'$ of $\mu$, $\Sigma_{\C} M$ and $\Dirac_{\C}$. Let $\Dirac_{\C}$ be \emph{coercive at infinity}, i.e. assume there exists a compact subset $K \subset M$ and a $C>0$ such that 
	\begin{align*}
		\|\psi\|_{L^2(\Sigma_{\C} M)} \leq C \| \Dirac_{\C} \psi \|_{L^2(\Sigma_{\C} M)}
	\end{align*}
	for all $\psi \in \Gamma(\Sigma_{\C} M)$, which are compactly supported in $M \setminus K$. Then $\Dirac_{\C}'$ is Fredholm and
	\begin{align*}
		\ind \Dirac'_+ = \ind \Dirac_+.
	\end{align*}
	Here $\Dirac_{\C}'$ is to be understood as the Dirac operator with APS-boundary conditions, i.e.
	\begin{align*}
		\dom(\Dirac_{\C}') = \{ \psi \in H^1(\Sigma_{\C} M ') \mid \psi|_{N_1} \in H^{1/2}_{\mathopen{]}-\infty,0\mathclose{]}}(\tilde \Dirac_{\C}), \psi|_{N_2} \in H^{1/2}_{\mathopen{[}0, \infty \mathclose{[}}(\tilde \Dirac_{\C})\},
	\end{align*}
	where $\tilde \Dirac_{\C} $ is the Dirac operator on the boundary (resp. its two-fold copy as in \eqref{EqDiracCopy}).
\end{Thm}

\begin{Thm}[\protect{\cite[Thm A]{Salamon}}]
	\label{ThmSalamon}
	Assume we are given the following data.
	\begin{enumerate}
		\item A complex separable Hilbert space $(H,\|\_\|_H)$.
		\item A dense subspace $W \subset H$ and a norm $\| \_ \|_W$ on $W$ such that $W$ is also a Hilbert space and such that the injection $W \hookrightarrow H$ is compact.
		\item A family of unbounded self-adjoint operators $\{A(t)_{t\in \R}\}$ on $H$ with time independent domain $W$, such that for each $t \in \R$ the graph norm of $A(t)$ is equivalent to $\| \_ \|_W$.
		\item A map $\R \to L(W,H)$, $t \mapsto A(t)$, which is continuously differentiable with respect to the weak operator topology.
		\item Invertible operators $A^\pm \in L(W,H)$ such that $\lim_{t \to \pm \infty}{A(t)} = A^\pm$ in norm topology.
	\end{enumerate}
	 Then the operator
	 \begin{align*}
		D_A := \tfrac{d}{dt} - A(t): W^{1,2}(\R,H) \cap L^2(\R,W) \to L^2(\R,H)
	 \end{align*}
	 is Fredholm, and its Fredholm index is equal to the spectral flow $\specfl(A)$ of the operator family $A=(A(t))_{t \in \R}$. 
\end{Thm}